\documentclass[preprint,3p,times]{elsarticle}

\usepackage[section]{placeins}
\usepackage{multirow}
\usepackage{makecell}
\usepackage{longtable}
\usepackage{booktabs}
\usepackage{supertabular}
\usepackage[table]{xcolor}
\usepackage{setspace}
\usepackage{bm}
\usepackage[bf]{caption}

\usepackage[ruled,linesnumbered,lined,commentsnumbered]{algorithm2e}

\newcommand{\tabincell}[2]{\begin{tabular}{@{}#1@{}}#2\end{tabular}} 

\usepackage{amssymb}
\usepackage{amsmath}
\newtheorem{theorem}{Theorem}
\newtheorem{lemma}{Lemma}

\newtheorem{definition}{Definition}
\newtheorem{corollary}{Corollary}

\newtheorem{remark}{Remark}
\newtheorem{example}{Example}

\journal{Elsevier}

\begin{document}

\begin{frontmatter}

\title{An extension of the order-preserving mapping to the 
WENO-Z-type schemes}


\author[a]{Ruo Li}
\ead{rli@math.pku.edu.cn}

\author[b,c]{Wei Zhong\corref{cor1}}
\ead{zhongwei2016@pku.edu.cn}

\cortext[cor1]{Corresponding author}

\address[a]{CAPT, LMAM and School of Mathematical Sciences, Peking
University, Beijing 100871, China}

\address[b]{School of Mathematical Sciences, Peking University,
Beijing 100871, China}

\address[c]{Northwest Institute of Nuclear Technology, Xi'an 
710024, China}

\begin{abstract}

  In our latest studies, by introducing the novel order-preserving 
  (OP) criterion, we have successfully addressed the widely 
  concerned issue of the previously published mapped weighted essentially non-oscillatory (WENO) schemes that it is rather 
  difficult to achieve high resolutions on the premise of removing
  spurious oscillations for long-run simulations of the hyperbolic systems. In the present study, we extend 
  the OP criterion to the WENO-Z-type schemes as the forementioned 
  issue has also been extensively observed numerically for these 
  schemes. Firstly, we innovatively present the concept of the 
  generalized mapped WENO schemes by rewriting the Z-type weights in a uniform formula 
  from the perspective of the mapping relation. Then, we naturally 
  introduce the OP criterion to improve the WENO-Z-type schemes, and 
  the resultant schemes are denoted as MOP-GMWENO-X. Finally, 
  extensive numerical experiments have been conducted to demonstrate 
  the benefits of these new schemes. We draw the conclusion that, 
  the convergence propoties of the proposed schemes are 
  equivalent to the corresponding WENO-X schemes. The major 
  benefit of the new schemes is that they have the 
  capacity to achieve high resolutions and simultaneously remove 
  spurious oscillations for long simulations. The new schemes have the additional 
  benefit that they can greatly decrease the post-shock 
  oscillations on solving 2D Euler problems with strong shock waves.

\end{abstract}



\end{frontmatter}

\section{Introduction}
\label{secIntroduction}
Over the past several decades, the WENO methods 
\cite{ENO-Shu1988,ENO-Shu1989,WENO-LiuXD,WENO-JS,WENOoverview}
have received considerable scholarly attention. The first WENO scheme that can obtain the designed convergence order 
of accuracy was proposed by Jiang and Shu \cite{WENO-JS}, dubbed 
WENO-JS. By using the information of all $r$-point substencils of the ENO scheme \cite{ENO1987JCP71,ENO1987V24, ENO1986, 
ENO1987JCP83}, WENO-JS maintains the ENO property near the region with discontinuities or large gradients and in the 
meantime achieves the designed convergence rates of accuracy. It was pointed out by Henrick et al. in \cite{WENO-M} that the fifth-order WENO-JS scheme can not recover the designed accuracy at critical points of order $n_{\mathrm{cp}} = 1$. For function $g$, one has $g'=0, g'' \neq 0$ for $n_{\mathrm{cp}} = 1$. Similarly, $g'=0, g'' = 0, g''' \neq 0$ for $n_{\mathrm{cp}} = 2$, etc. Then, the sufficient condition for optimality of the convergence rates of accuracy were derived in \cite{WENO-Z}, and this sufficient condition can be extended to higher order cases trivially \cite{WENO-IM}. 
In the work of Henrick et al. \cite{WENO-M}, a mapping function, namely $(g^{\mathrm{M}})_{s}(\omega^{\mathrm{JS}})$, was designed and the resultant mapped WENO scheme, dubbed WENO-M, can achieve the designed convergence properties even in the presence of critical points. It is 
since the work of Henrick et al. \cite{WENO-M} that the study of 
different mapped WENO methods has gained momentum, and a series of new 
mapping functions \cite{WENO-PM,WENO-IM,WENO-PPM5,WENO-RM260,
WENO-MAIMi,WENO-ACM} have been proposed by obeying the similar 
principles proposed by Henrick et al. \cite{WENO-M}. 

Later, the work of Henrick et al. \cite{WENO-M} inspired the 
development of a new family of nonlinear weights, dubbed Z-type 
weights. From a different perspective, Borges et al. \cite{WENO-Z} 
proposed another version of nonlinear weights by using available and 
previously unused information of the WENO-JS scheme. In other words, 
a global smoothness indicator (GSI) of higher order, obtained via a 
linear combination of the original smoothness indicators of the 
WENO-JS scheme, was proposed and employed to devise the new 
nonlinear weights. The resultant scheme was denoted as WENO-Z. 
Because of the success of the WENO-Z scheme that its nonlinear 
weights can satisfy the sufficient conditions for optimality of the 
convergence order without any costly mapping processes, leading to 
superior results with almost the same computational effort of the 
WENO-JS method, different researchers have developed a multitude of 
techniques to design their Z-type weights \cite{WENO-NS,WENO-eta,
WENO-eta-02,WENO-P,WENO-Zplus,WENO-ZplusI,WENO-ZA,P-WENO,MWENO-P,
WENO-D_WENO-A,WENO-NIP} by obeying the similar priciples proposed by 
Borges et al. \cite{WENO-Z}. In this paper, all the WENO 
schemes using Z-type weights is collectively called
WENO-Z-type schemes. We will give a brief review of several 
WENO-Z-type schemes in subsection \ref{subsec:WENO-Z-type}.

Despite the success mainly for short-output-time simulations, the 
family of mapped WENO schemes has a serious and ubiquitous problem 
in calculations with long output times, that is, they can hardly 
avoid spurious oscillations and meanwhile preserve high resolutions 
for long simulations. This disadvantage in long simulations of the 
mapped WENO methods was firstly noticed and carefully studied by 
Feng et al. \cite{WENO-PM} and it has attracted considerable 
attention over the past decade \cite{WENO-IM,WENO-RM260,
WENO-RM-Vevek2018,WENO-AIM,WENO-MAIMi,MOP-WENO-ACMk,MOP-WENO-X,
PoAOP-WENO-X}. However, up to now, far too little attention has been 
paid to the long-output-time simulations of the WENO-Z-type schemes. 
Indeed, our extensive calculations (see subsection 
\ref{subsec:Long-run} below) show that the WENO-Z-type schemes also 
terribly suffer from either losing high resolutions or generating 
numerical oscillations on long-run calculations. Nevertheless, nodoubtly, this issue is worthy of scholarly attention.

Accordingly, in this article, we would like to focus on the theme of 
addressing the aforementioned drawback of the WENO-Z-type schemes. 
First of all, we give the important observation and analysis of the 
implicit relationship between the nonlinear weights of the WENO-JS 
scheme and the Z-type weights, denoted as IMR (standing for 
\textit{implicit mapping relation}) for simplicity. It can be found 
that the profiles of IMRs for various WENO-Z-type schemes are very 
similar to the traditional designed mapping curves (of the mapped 
WENO schemes) that embrace evident optimal weight intervals where 
the nonlinear weights are replaced by the ideal weights. It is 
well known that this kind of replacements appears favorable for 
reducing the dissipation and improving the resolution at least for 
short-output-time simulations. However, these replacements also 
generate the non-order-preserving (non-OP) points where the order of 
the nonlinear weights is disrupted. It has been demonstrated 
\cite{MOP-WENO-ACMk,MOP-WENO-X,PoAOP-WENO-X} that the non-OP points 
are extremely harmful for the WENO schemes to preserving high 
resolutions and meanwhile avoiding spurious oscillations for long-run calculations. Therefore, after advancing the new 
concept of the generalized mapped WENO schemes by reformulating the 
nonnormalized nonlinear Z-type weights in a uniform form from the 
perspective of the mapping relation, the OP property is introduced 
to modify the previously published WENO-Z-type schemes for the purpose of improving their performance in long simulations. Necessary theoretical analysis will be provided and extensive numerical experiments of 1D linear advection equation with various intial conditions for long output times will be conducted and carefully discussed to examine the enhancement of the new proposed schemes. In addition, we also carry out several problems modeled via 1D and 2D Euler euqations to show the good performance of these schemes.

We organize the remainder of this article as follows. The preliminaries are reviewed in Section \ref{secWENO-Z-type}, where we briefly recall the procedures of several typical WENO schemes, including WENO-JS \cite{WENO-JS}, WENO-Z
\cite{WENO-Z} and some other WENO methods with Z-type weights. In Section \ref{sec:Z-mapping}, the analysis of the Z-type nonlinear weights from the perspective of the mapping relation is given, and the concept of the generalized mapped WENO schemes is proposed. The method to generally improve the previously published WENO-Z-type schemes by requiring the use of the OP property is devised in Section \ref{sec:MOP-GMWENO-X}, and some examples of 1D linear advection equation are performed to manifest the major benefits of the proposed schemes in this Section. In Section 
\ref{NumericalExperiments}, more numerical 
results of 1D and 2D Euler equations are provided for illustration. 
The conclusions are given in Section \ref{secConclusions}.


\section{Brief review of the WENO-Z-type schemes}
\label{secWENO-Z-type}
It was reported \cite{WENO-Z01} that there are three versions of the 
odd-order WENO methods with $r \geq 3$ when applied to 
hyperbolic systems. They are the classical WENO schemes 
(e.g., WENO-JS \cite{WENO-JS}, WENO-ZS \cite{WENO-ZS}, WENO-NIS 
\cite{WENO-NIS}), the mapped WENO schemes (e.g., WENO-M
\cite{WENO-M}, WENO-PM$k$ \cite{WENO-PM}, WENO-MAIMi
\cite{WENO-MAIMi}, WENO-ACM \cite{WENO-ACM}) and the WENO-Z-type 
schemes (e.g., WENO-Z \cite{WENO-Z}, WENO-Z+ \cite{WENO-Zplus}, 
WENO-NIP \cite{WENO-NIP}). As the version of the mapped WENO schemes 
has been discussed carefully in our previous work \cite{MOP-WENO-X}, 
we mainly describe the other two versions in this section. In other 
words, we first introduce the classical WENO-JS 
scheme \cite{WENO-JS}, along with the original WENO-Z scheme as 
designed in \cite{WENO-Z} and then several WENO schemes with improved Z-type weights.  

For brevity's sake, we denote our discussion to the 
following 1D scalar equation
\begin{equation}
u_{t} + f(u)_{x} = 0.
\label{eq:1D-scalar}
\end{equation}
Let $\bar{u}(x_{j}, t)= \dfrac{1}{\Delta x}\displaystyle\int_{x - \Delta x/2}^{x + \Delta x/2}u(\eta,t)\mathrm{d}\eta$, within the Finite Volume Method (FVM) framework, we can transform Eq. \eqref{eq:1D-scalar} into the semi-discretized form
\begin{equation}\left\{
\begin{array}{l}
\dfrac{\mathrm{d}\bar{u}_{j}(t)}{\mathrm{d}t}
\approx \mathcal{L}(u_{j}), \\
\mathcal{L}(u_{j}) = -\dfrac{1}{\Delta x}\Big( \hat{f}_{j+1/2} - 
\hat{f}_{j-1/2} \Big),
\end{array}
\right.
\label{eq:1D-scalar:semiFVM}
\end{equation}
where $\bar{u}_{j}(t)$ is the approximation to $\bar{u}(x_{j}, t)$, and $\hat{f}_{j \pm 1/2} = \hat{f}(u_{j \pm 1/2}^{-}, u_{j \pm 1/2}^{+})$ is the numerical flux used to replace $f(u)$. The WENO 
reconstruction techniques are used to calculate the
$u_{j \pm 1/2}^{\pm}$. We will choose the global LF flux in this article.

\begin{remark}
In the discussion below, we only show the reconstruction of 
$u_{j+1/2}^{-}$ as that of $u_{j+1/2}^{+}$ is symmetric to it with respect to $x_{j+1/2}$. And for brevity, the ``-'' sign in the superscript will be dropped without causing any confusion.
\end{remark}

\subsection{WENO-JS}\label{subsec:WENO-JS}
In the 5th-order WENO-JS, $u_{j+1/2}$ is approximated by 
\begin{equation}
u_{j + 1/2} = \omega_{0}u_{j + 1/2}^{0} + \omega_{1}u_{j + 1/2}^{1} + \omega_{2}u_{j + 1/2}^{2},
\label{eq:approx_WENO}
\end{equation}
with
\begin{equation}
\begin{array}{l}
\begin{aligned}
&u_{j+1/2}^{0} = \dfrac{1}{6}(2\bar{u}_{j-2} - 7\bar{u}_{j-1}
+ 11\bar{u}_{j}), \\
&u_{j+1/2}^{1} = \dfrac{1}{6}(-\bar{u}_{j-1} + 5\bar{u}_{j}
+ 2\bar{u}_{j+1}), \\
&u_{j+1/2}^{2} = \dfrac{1}{6}(2\bar{u}_{j} + 5\bar{u}_{j+1}
- \bar{u}_{j+2}),
\end{aligned}
\end{array}
\label{eq:approx_ENO}
\end{equation}
and (for $s = 0,1,2$)
\begin{equation} 
\omega_{s}^{\mathrm{JS}} = \dfrac{\alpha_{s}^{\mathrm{JS}}}{\sum_{l =
 0}^{2} \alpha_{l}^{\mathrm{JS}}}, \alpha_{s}^{\mathrm{JS}} = \dfrac{
 d_{s}}{(\epsilon + \beta_{s})^{2}},
\label{eq:weights:WENO-JS}
\end{equation} 
where $\epsilon$ is a small positive number, $d_{0} = \frac{1}{10}, 
d_{1} = \frac{3}{5}, d_{2} = \frac{3}{10}$. In \cite{WENO-JS}, $\beta_{s}$, dubbed LSI, 
are explicitly given as follows
\begin{equation*}
\begin{array}{l}
\begin{aligned}
\beta_{0} &= \dfrac{13}{12}\big(\bar{u}_{j - 2} - 2\bar{u}_{j - 1} + 
\bar{u}_{j} \big)^{2} + \dfrac{1}{4}\big( \bar{u}_{j - 2} - 4\bar{u}_
{j - 1} + 3\bar{u}_{j} \big)^{2}, \\
\beta_{1} &= \dfrac{13}{12}\big(\bar{u}_{j - 1} - 2\bar{u}_{j} + \bar
{u}_{j + 1} \big)^{2} + \dfrac{1}{4}\big( \bar{u}_{j - 1} - \bar{u}_{
j + 1} \big)^{2}, \\
\beta_{2} &= \dfrac{13}{12}\big(\bar{u}_{j} - 2\bar{u}_{j + 1} + \bar
{u}_{j + 2} \big)^{2} + \dfrac{1}{4}\big( 3\bar{u}_{j} - 4\bar{u}_{j 
+ 1} + \bar{u}_{j + 2} \big)^{2}.
\end{aligned}
\end{array}
\end{equation*}

It was indicated by Henrick et al. \cite{WENO-M} that the convergence rates of accuracy of WENO-JS may drop to 3rd-order at critical points. In the same article, they indicated that to ensure 
the overall scheme retaining fifth-order accuracy, it is sufficient 
to require
\begin{equation}
\omega_{s}^{\pm} - d_{s} = O(\Delta x^{3}).
\label{eq:sufficient:optimal}
\end{equation} 

\subsection{WENO-Z}\label{subsec:WENO-Z}
Borges et al. \cite{WENO-Z} devised the first WENO-Z-type scheme, 
dubbed WENO-Z. By introducing the GSI, 
say, $\tau_{5} = \lvert \beta_{0} - \beta_{2} \rvert$, they compute 
the the nonlinear weights by
\begin{equation}
\omega_{s}^{\mathrm{Z}}=\dfrac{\alpha_{s}^{\mathrm{Z}}}{\sum_{j=0}^{2
}\alpha_{j}^{\mathrm{Z}}}, \quad \alpha_{s}^{\mathrm{Z}} = d_{s}
\bigg(1 + \Big(\dfrac{\tau_{5}}{\beta_{s} + \epsilon} \Big)^{p} 
\bigg), \quad s = 0,1,2,
\label{eq:weights:Z}
\end{equation}
where the parameters $\beta_{s}$ and $\epsilon$ are the same as in 
the classical WENO-JS scheme, and $p$ is a tunable parameter. Borges 
et al. demonstrated that \cite{WENO-Z}, if $p = 1$, the WENO-Z 
scheme only achieves fourth-order convergence orders in the presence of critical points; if $p = 2$, it can recover the designed 5th-order convergence rates. Unless indicated otherwise, we choose $p = 2$ in the present study.

\subsection{Other WENO-Z-type schemes}
\label{subsec:WENO-Z-type}
The significant contribution of WENO-Z is the inclusion 
of higher-order information by introducing the GSI, $\tau_{5}$, in the definition of the nonlinear 
weights. This helps the WENO-Z scheme get closer to the central 
scheme than the WENO-JS scheme and hence generate less dissipation.
The formation of the nonlinear weights of WENO-Z gives 
an effective way for enhancing the capacity of WENO-JS. Thus, by obeying the similar principles of WENO-Z, a series 
of new Z-type nonlinear weights with various kinds of GSI are proposed. In this subsection, we briefly 
review other WENO schemes with different Z-type weights that will be further 
studied in the rest of this paper. It should be noted that there are 
so many different types of WENO-Z-type schemes \cite{WENO-Z,
WENO-NS,WENO-eta,WENO-eta-02,WENO-P,WENO-Zplus,WENO-ZplusI,WENO-ZA,
P-WENO,MWENO-P,WENO-D_WENO-A,WENO-NIP}, and we only considered a 
limited number of them just for illustrative purposes in the present 
study.

WENO-Z$\eta(\tau_{l})$: Fan et al. \cite{WENO-eta,WENO-eta-02} 
proposed the WENO-Z$\eta(\tau_{l})$ scheme with the following Z-type 
nonlinear weights
\begin{equation}
\omega_{s}^{\mathrm{Z}\eta} = \dfrac{\alpha_{s}^{\mathrm{Z}\eta}}{
\sum_{l = 0}^{2}\alpha_{l}^{\mathrm{Z}\eta}},\alpha_{s}^{\mathrm{Z}
\eta} = d_{s}\bigg( 1 + \bigg(\dfrac{\tau_{l}}{\eta_{s} + \epsilon}
\bigg)^{2} \bigg), \quad s = 0,1,2,
\label{eq:weights:Z-eta}
\end{equation}
where the LSI is calculated using the flexible and simple formula 
suggested by Shen and Zha \cite{WENO-eta-LSI}, and we state them 
explicitly here
\begin{equation}\left\{
\begin{array}{l}
\begin{aligned}
\eta_{0} &= \dfrac{1}{4}\big( \bar{u}_{j - 2} - 4\bar{u}_{j - 1} + 
3\bar{u}_{j} \big)^{2} + \big(\bar{u}_{j - 2} - 2\bar{u}_{j - 1} + 
\bar{u}_{j} \big)^{2}, \\
\eta_{1} &= \dfrac{1}{4}\big( \bar{u}_{j - 1} - \bar{u}_{j + 1} 
\big)^{2} + \big(\bar{u}_{j - 1} - 2\bar{u}_{j} + \bar{u}_{j + 1} 
\big)^{2}, \\
\eta_{2} &= \dfrac{1}{4}\big( 3\bar{u}_{j} - 4\bar{u}_{j + 1} + 
\bar{u}_{j + 2} \big)^{2} + \big(\bar{u}_{j} - 2\bar{u}_{j + 1} + 
\bar{u}_{j + 2} \big)^{2}.
\end{aligned}
\end{array}
\right.
\label{eq:ISk:Z-eta}
\end{equation}
Fan et al. \cite{WENO-eta} proposed a fifth-, a sixth- and two 
eighth- order GSIs as
\begin{equation*}\left\{
\begin{array}{l}
\tau_{5} = \lvert \eta_{0} - \eta_{2} \rvert, \quad
\tau_{6} = \bigg\lvert \eta_{5} - \dfrac{\eta_{0} + 4\eta_{1} +
\eta_{2}}{6} \bigg\rvert, \\
\tau_{81} = \bigg\lvert \Big( \big\lvert P_{0}^{(1)} 
\big\rvert - \big\lvert P_{2}^{(1)} \big\rvert \Big)\Big( P_{0}^{(2)}
-2P_{1}^{2} + P_{2}^{(2)} \Big)\bigg\rvert, \\
\tau_{82} = \Big( \big\lvert P_{0}^{(1)} \big\rvert - 
\big\lvert P_{2}^{(1)} \big\rvert \Big)^{2} + \Big( P_{0}^{(2)} -
2P_{1}^{2} + P_{2}^{(2)} \Big)^{2}, \\
\end{array}
\right.
\label{eq:tau:Z-eta}
\end{equation*}
where
\begin{equation*}
\eta_{5} = \dfrac{1}{144}\bigg( \big( \bar{u}_{j-2} -8\bar{u}_{j-1} + 8b\bar{u}_{j+1} - \bar{u}_{j+2} \big)^{2} + \big( \bar{u}_{j-2} - 16\bar{u}_{j-1} + 30\bar{u}_{j} - 16\bar{u}_{j+1} + \bar{u}_{j+2} \big)^{2} \bigg), \\
\end{equation*}
and
\begin{equation*}
\centering
\begin{array}{l}
\left\{
\begin{array}{l}
\begin{aligned}
&P_{0}^{(1)} = \dfrac{1}{2}\bar{u}_{j-2} - 2\bar{u}_{j-1} + \dfrac{3}{2}\bar{u}_{j}, \\
&P_{1}^{(1)} = -\dfrac{1}{2}\bar{u}_{j-1} + \dfrac{1}{2}\bar{u}_{j+1}, \\
&P_{2}^{(1)} = -\dfrac{3}{2}\bar{u}_{j} + 2\bar{u}_{j+1} - \dfrac{1}{2}\bar{u}_{j+2}, \\
\end{aligned}
\end{array}
\right.  
\quad \quad \quad
\left\{
\begin{array}{l}
\begin{aligned}
&P_{0}^{(2)} = \bar{u}_{j-2} - 2\bar{u}_{j-1} + \bar{u}_{j}, \\
&P_{1}^{(2)} = \bar{u}_{j-1} - 2\bar{u}_{j} + \bar{u}_{j+1}, \\
&P_{2}^{(2)} = \bar{u}_{j} - 2\bar{u}_{j+1} + \bar{u}_{j+2}.
\end{aligned}
\end{array}
\right.
\end{array}
\label{eq:eta_5:Z-eta}
\end{equation*}

WENO-Z+: Acker et al. \cite{WENO-Zplus} proposed the WENO-Z+ scheme 
with the following Z-type nonlinear weights
\begin{equation}
\omega_{s}^{\mathrm{Z+}} = \dfrac{\alpha_{s}^{\mathrm{Z+}}}{
\sum_{l = 0}^{2}\alpha_{l}^{\mathrm{Z+}}},\alpha_{s}^{\mathrm{Z+}} 
= d_{s}\bigg( 1 + \bigg( \dfrac{\tau_{5} + \epsilon}{\beta_{s} + 
\epsilon} \bigg)^{2} + \lambda\bigg( \dfrac{\beta_{s} + \epsilon}{
\tau_{5} + \epsilon} \bigg) \bigg), \quad s = 0,1,2,
\label{eq:weights:Z+}
\end{equation}
where $\lambda$ is recommended to be $\Delta x^{2/3}$ in their work.

WENO-ZA: Shen et al. \cite{WENO-ZA} proposed the WENO-ZA scheme with 
the following Z-type nonlinear weights
\begin{equation}
\omega_{s}^{\mathrm{ZA}} = \dfrac{\alpha_{s}^{\mathrm{ZA}}}{
\sum_{l = 0}^{2}\alpha_{l}^{\mathrm{ZA}}},\alpha_{s}^{\mathrm{ZA}} = 
d_{s}\bigg( 1 + \dfrac{A\tau_{6}}{\beta_{s} + \epsilon} \bigg), 
\quad s = 0,1,2,
\label{eq:weights:ZA}
\end{equation}
where the LSI is defined by
\begin{equation}
\beta_{s} = \gamma_{1}\big(\bar{u}_{s}^{(1)}\big)^{2} + 
\gamma_{2}\big(\bar{u}_{s}^{(2)}\big)^{2}, \quad s = 0,1,2,
\label{eq:beta_s:ZA}
\end{equation}
with
\begin{equation}
\begin{array}{ll}
\left\{
\begin{array}{l}
\bar{u}_{0}^{(1)}=(\bar{u}_{j-2}-4\bar{u}_{j-1} + 3\bar{u}_{j})/2, \\
\bar{u}_{1}^{(1)}=(- \bar{u}_{j-1}+\bar{u}_{j+1})/2, \\
\bar{u}_{2}^{(1)}=(-3\bar{u}_{j}+4\bar{u}_{j+1} - \bar{u}_{j+2})/2, 
\end{array}\right.
& \quad
\left\{
\begin{array}{l}
\bar{u}_{0}^{(2)} = \bar{u}_{j-2} - 2\bar{u}_{j-1} + \bar{u}_{j}, \\
\bar{u}_{1}^{(2)} = \bar{u}_{j-1} - 2\bar{u}_{j} + \bar{u}_{j+1}, \\
\bar{u}_{2}^{(2)} = \bar{u}_{j} - 2\bar{u}_{j+1} + \bar{u}_{j+2}. 
\end{array}\right.
\end{array}
\label{eq:bar_u:ZA}
\end{equation}
In Eq. \eqref{eq:weights:ZA}, the parameter $A$ is recommended as
\begin{equation}
A = \dfrac{\tau_{6}}{\beta_{0} + \beta_{2} - \tau_{6} + \epsilon},
\label{eq:functionA:ZA}
\end{equation}
with
\begin{equation}
\tau_{6} = \gamma_{1}\big( \lvert \bar{u}_{0}^{(1)} \rvert - \lvert 
\bar{u}_{2}^{(1)} \rvert \big)^{2} + \gamma_{2}\big( \lvert 
\bar{u}_{0}^{(2)} \rvert - \lvert \bar{u}_{2}^{(2)} \rvert \big)^{2}.
\label{eq:tau6:ZA}
\end{equation}
According to the recommendation of Shen et al. \cite{WENO-ZA}, 
$\gamma_{1} = 1$, $\gamma_{2} = 13/12$ are used.

WENO-D/WENO-A: Wang et al. \cite{WENO-D_WENO-A} proposed  
WENO-D and WENO-A with the following Z-type nonlinear weights
\begin{equation}\left\{
\begin{array}{ll}
\text{WENO-D:} &
\omega_{s}^{\mathrm{D}} = \dfrac{\alpha_{s}^{\mathrm{D}}}{\sum_{l = 0}^{2}\alpha_{l}^{\mathrm{D}}},\alpha_{s}^{\mathrm{D}} = d_{s}\bigg( 1 + \Phi\bigg(\dfrac{\tau_{5}}{\beta_{s} + \epsilon}\bigg)^{p} \bigg), \\
\text{WENO-A:} &
\omega_{s}^{\mathrm{A}} = \dfrac{\alpha_{s}^{\mathrm{A}}}{\sum_{l = 0}^{2}\alpha_{l}^{\mathrm{A}}},\alpha_{s}^{\mathrm{A}} = d_{s}\bigg( \max\bigg(1, \Phi\bigg(\dfrac{\tau_{5}}{\beta_{s} + \epsilon}\bigg)^{p}\bigg) \bigg),
\end{array}
\right. 
\quad s = 0,1,2,
\label{eq:weights:D}
\end{equation}
where $p=2$ is used in their numerical examples, and $\tau_{5} = 
\lvert \beta_{0} - \beta_{2} \rvert$, $\Phi$ is computed by
\begin{equation}
\Phi = \min\{ 1,\phi \}, \quad \phi = \sqrt{\lvert \beta_{0} - 
2\beta_{1} + \beta_{2} \rvert}.
\label{eq:D:phi}
\end{equation}

WENO-NIP: Yuan \cite{WENO-NIP} proposed the WENO-NIP scheme with the 
following Z-type nonlinear weights
\begin{equation}
\omega_{s}^{\mathrm{NIP}} = \dfrac{\alpha_{s}^{\mathrm{NIP}}}{
\sum_{l = 0}^{2}\alpha_{l}^{\mathrm{NIP}}},\alpha_{s}^{\mathrm{NIP}} 
= d_{s}\bigg( 1 + \dfrac{\tau}{(\epsilon + \chi_{s})^{2}}\bigg), 
\quad s = 0,1,2,
\label{eq:weights:NIP}
\end{equation}
where 
\begin{equation}
\begin{array}{l}
\chi_{0} = \theta \big\lvert \bar{u}_{j-2} - 3\bar{u}_{j-1} + 
2\bar{u}_{j} \big\rvert + \big\lvert \bar{u}_{j-2} -2\bar{u}_{j-1} + 
\bar{u}_{j} \big\rvert, \\
\chi_{1} = \theta \big\lvert \bar{u}_{j+1} - \bar{u}_{j} \big\rvert +
\big\lvert \bar{u}_{j-1} -2\bar{u}_{j} + \bar{u}_{j+1} \big\rvert, \\
\chi_{2} = \theta \big\lvert \bar{u}_{j+1} - \bar{u}_{j} \big\rvert +
\big\lvert \bar{u}_{j} -2\bar{u}_{j+1} + \bar{u}_{j+2} \big\rvert,
\end{array}
\label{eq:ISk:NIP}
\end{equation}
and $\tau = \lvert \bar{u}_{j-2} - 4\bar{u}_{j-1} + 6\bar{u}_{j} 
- 4\bar{u}_{j+1} + \bar{u}_{j+2} \rvert^{J}$. As mentioned in \cite{WENO-NIP}, the requirement 
$J \geq \frac{2n_{\mathrm{cp}}+4}{4}$ needs to be satisfied. 
According to Yuan's solutions \cite{WENO-NIP}, $\theta = 0.1$ seems 
to be a better choice.



\section{Study on the Z-type weights from the perspective of 
the mapping relation}\label{sec:Z-mapping}

It is well known
\cite{WENO-M,WENO-Z,WENO-IM,WENO-MAIMi,WENO-ACM} that the WENO-JS 
scheme can not achieve the designed convergence order at critical 
points, and there are usually two ways to address this issue: one is 
to introduce a mapping function \cite{WENO-M,WENO-PM,WENO-IM,
WENO-RM260,WENO-PPM5,WENO-MAIMi,WENO-ACM}, the other is to introduce 
a GSI \cite{WENO-Z,WENO-eta,WENO-Zplus,WENO-ZA,WENO-D_WENO-A,
WENO-NIP}. Recently, the long-run simulations of the mapped WENO 
methods have been widely concerned \cite{WENO-PM,WENO-IM,WENO-AIM,
WENO-RM260,WENO-MAIMi,MOP-WENO-ACMk,MOP-WENO-X}. The mapped WENO 
methods have one great drawback in common that it is very difficult 
to get high resolutions and remove spurious oscillations 
simultaneously for long-run simulations and many researches 
were devoted to amending it. However, any researches hardly focused on the performances of the WENO-Z-type schemes for this topic. In this study, this will be studied and discussed carefully. Actually, our numerical calculations indicate that the drawback of long-run 
simulations mentioned above also exist for the WENO-Z-type schemes.
From Figs. \ref{fig:Z:N1600:01}-\ref{fig:SLP:N3200:02} below, 
we can see that most of the existing WENO-Z-type schemes suffer from 
either losing high resolutions or generating spurious oscillations. 
This issue will be addressed in this paper.

In practical calculations, for any output time $t$, we can calculate 
uniquely the $\omega_{s}^{\mathrm{JS}}$ of each location $x_{j}$ 
according to Eq. \eqref{eq:weights:WENO-JS}. And also, for any 
WENO-X scheme, we can uniquely get its nonlinear weights, say 
$\omega_{s}^{\mathrm{X}}$, according to its nonlinear weights 
formula, like the formulas of the various Z-type weights as shown in 
subsection \ref{subsec:WENO-Z-type}. Trivially, the relation of 
$\omega_{s}^{\mathrm{JS}} \sim \omega_{s}^{\mathrm{X}}$ is uniquely 
determined. Thus, we call it {\it{Implicit Mapping Relation (IMR)}}
of the WENO-X scheme as an explicit mapping function is not 
necessary for this case.

\begin{example}\rm{
(SLP) The following \textit{\textbf{L}inear \textbf{P}roblem} 
proposed by \textit{\textbf{S}hu} et al. in \cite{WENO-JS}, denoted 
as SLP, is widely used to examine the performance of high resolution 
schemes. It is defined by}
\label{ex:SLP}
\end{example}
\begin{equation}\left\{
\begin{array}{l}
u_{t} + u_{x} = 0, \quad -1 < x < 1,\\
u(x, 0) = u_{0}(x).
\end{array}
\right.
\end{equation}
Let $G(x, \beta, z) = \exp\Big(-\beta \big(x - z\big)^{2}\Big),
F(x, \alpha, a) = \Big(\max \big(1 - \alpha ^{2}(x - a)^{2}, 0 \big)\Big)^{1/2}$, and the initial condition $u_{0}(x)$ is given as
\begin{equation}
\begin{array}{l}
u_{0}(x) = \left\{
\begin{array}{ll}
\dfrac{1}{6}\big[ G(x, \beta, z - \hat{\delta}) + 4G(x, \beta, z) + G
(x, \beta, z + \hat{\delta}) \big], & -0.8 \leq x \leq -0.6, \\
1, & -0.4 \leq x \leq -0.2, \\
1 - \big\lvert 10(x - 0.1) \big\rvert, & 0.0 \leq x \leq 0.2, \\
\dfrac{1}{6}\big[ F(x, \alpha, a - \hat{\delta}) + 4F(x, \alpha, a) +
 F(x, \alpha, a + \hat{\delta}) \big], & 0.4 \leq x \leq 0.6, \\
0, & \mathrm{otherwise},
\end{array}\right. 
\end{array}
\label{eq:LAE:SLP}
\end{equation}
with $a = \frac{1}{2}, \alpha = 10, \hat{\delta} = \frac{1}{200}, 
z = -\frac{7}{10}$ and $\beta = \log (2 / 36\hat{\delta} ^{2})$.

\begin{figure}[ht]
\flushleft
  \includegraphics[height=0.28\textwidth]
  {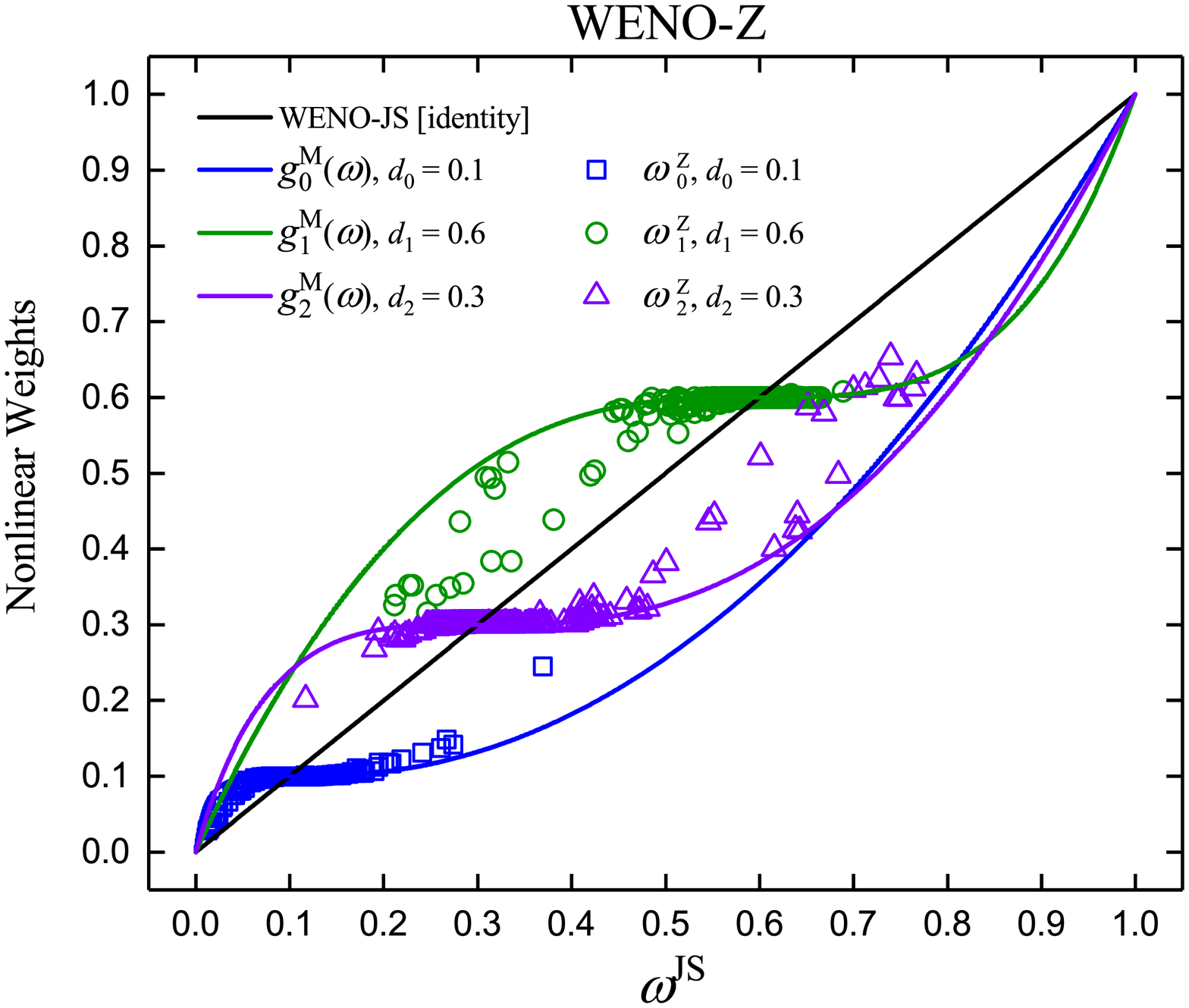}
  \includegraphics[height=0.28\textwidth]
  {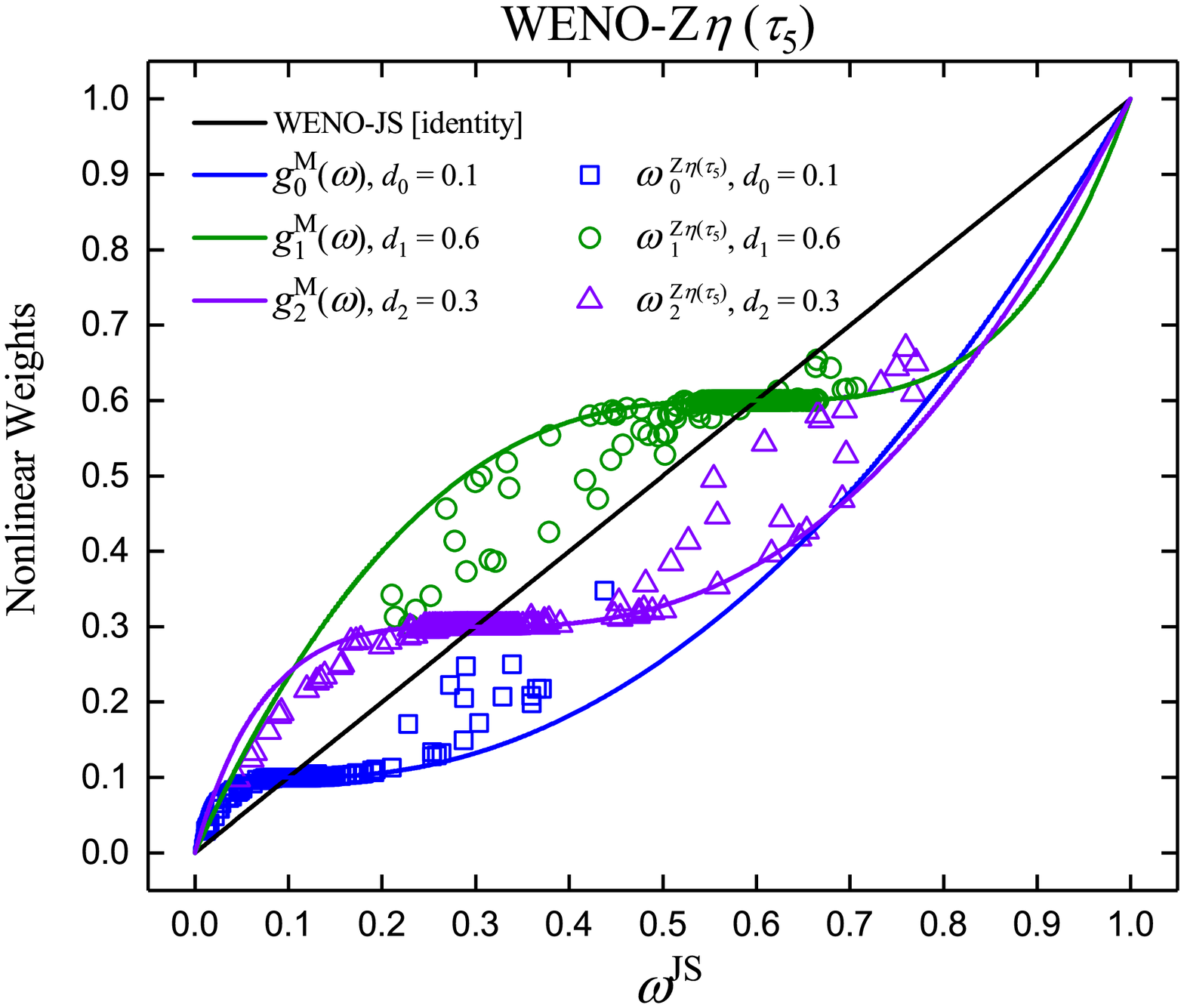}
  \includegraphics[height=0.28\textwidth]
  {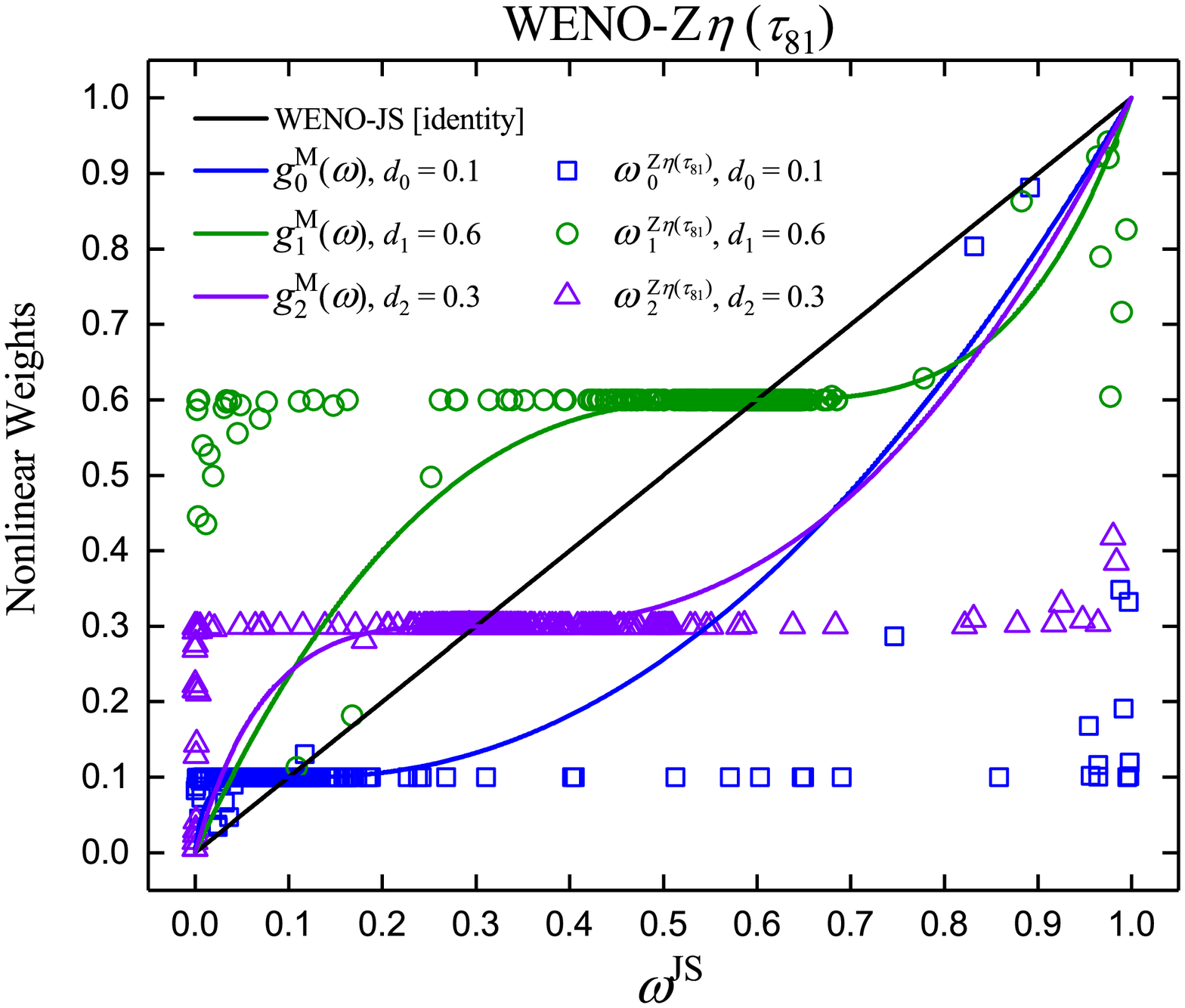}\\
  \includegraphics[height=0.28\textwidth]
  {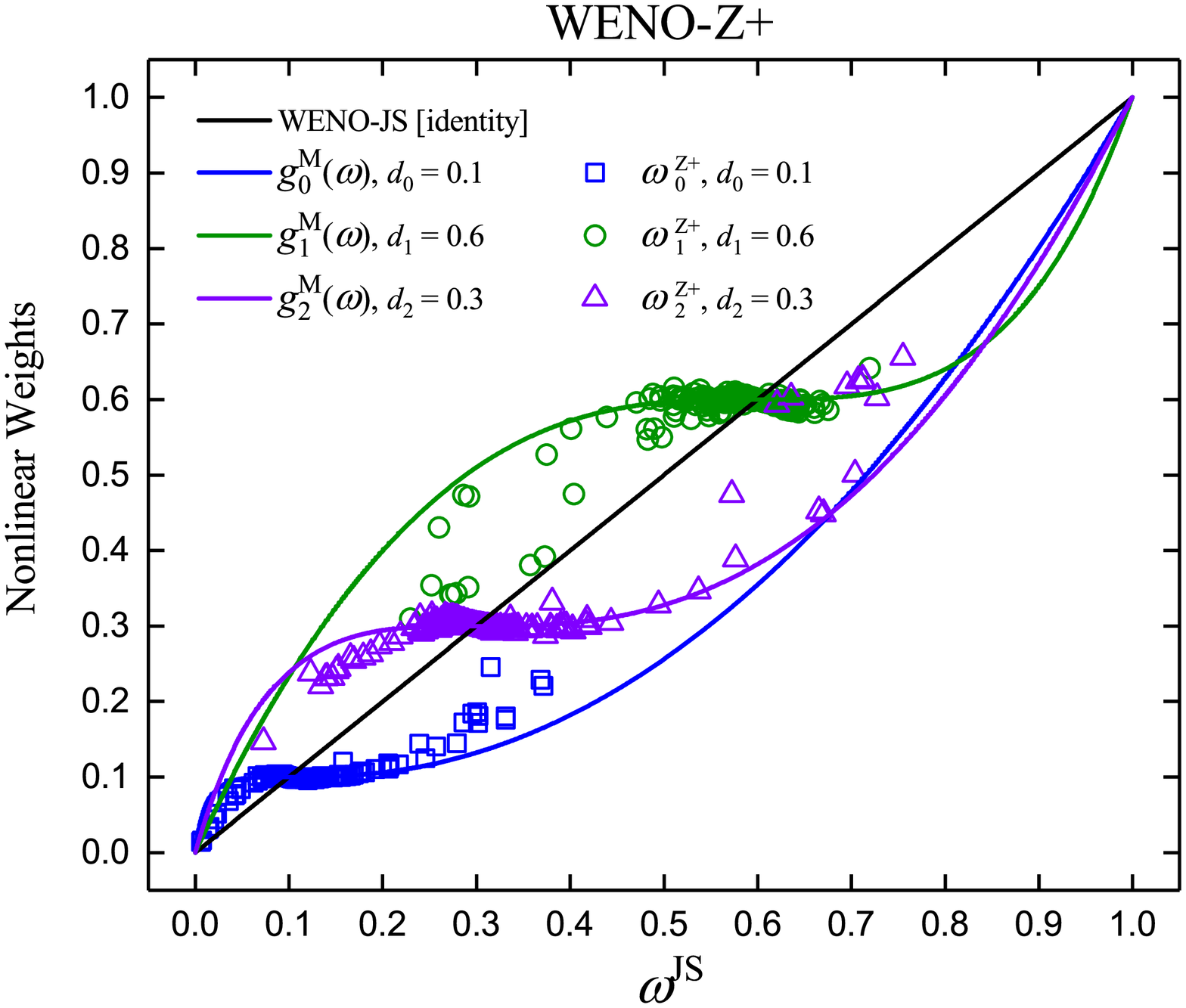}
  \includegraphics[height=0.28\textwidth]
  {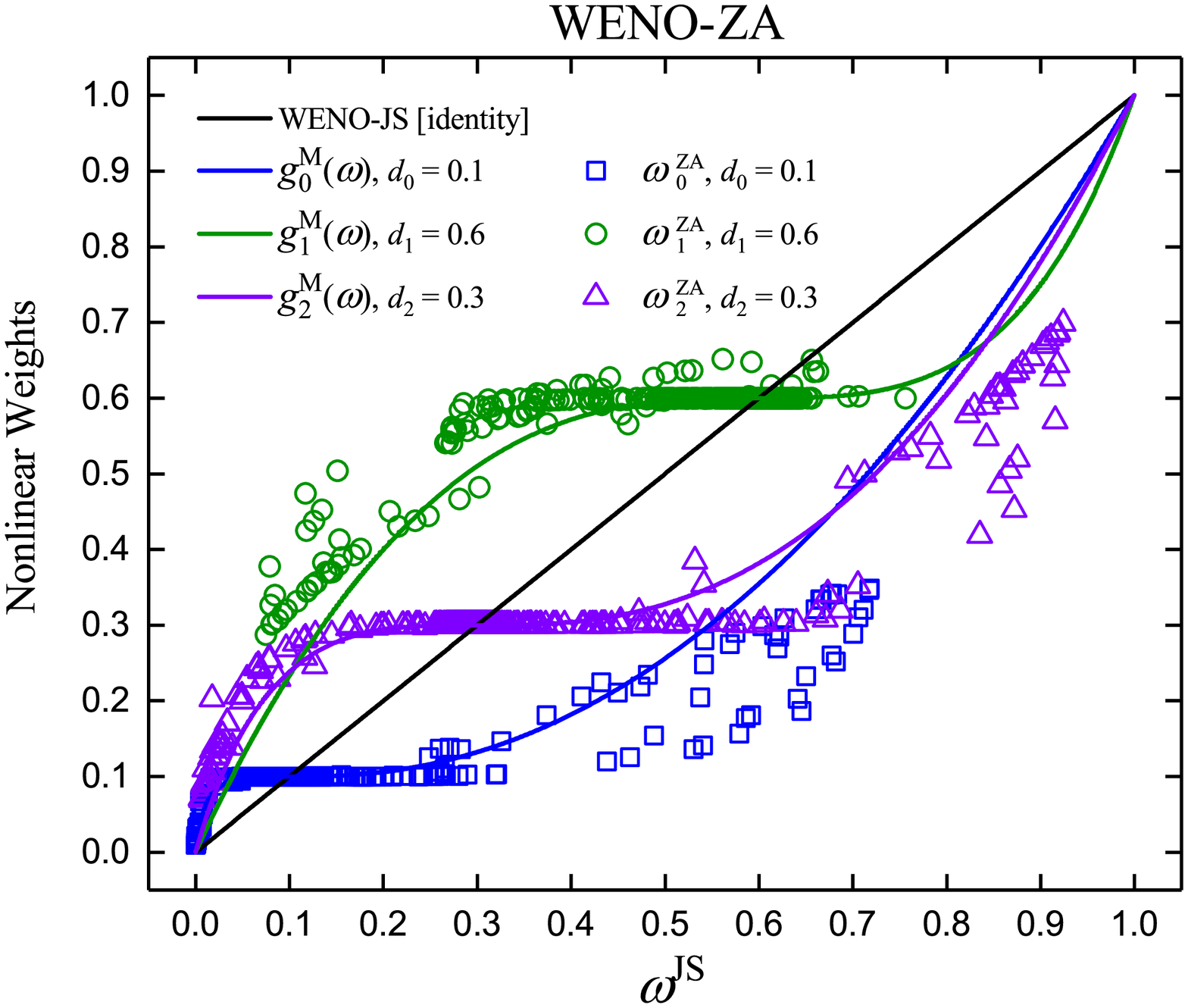}
  \includegraphics[height=0.28\textwidth]
  {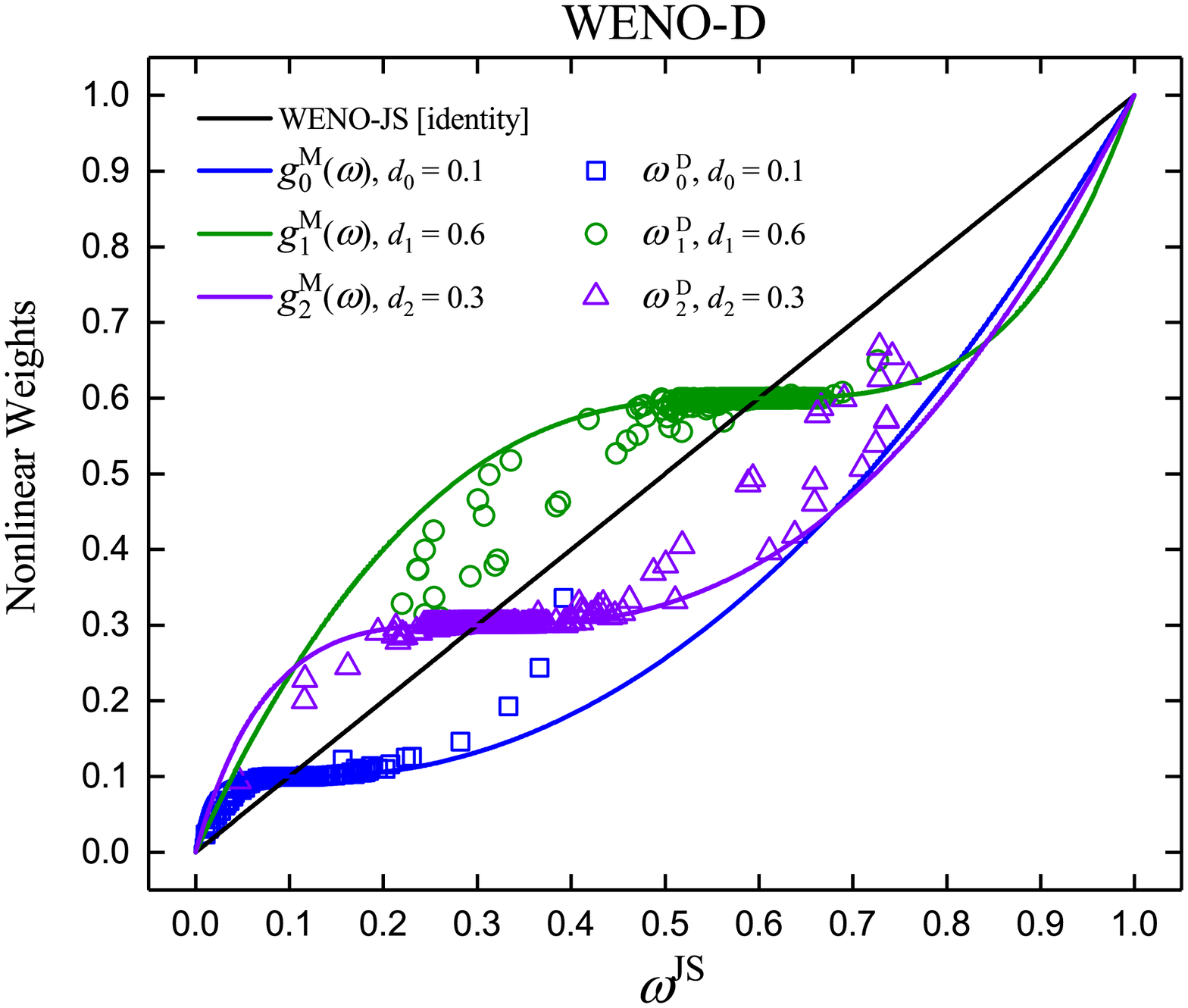}\\
  \includegraphics[height=0.28\textwidth]
  {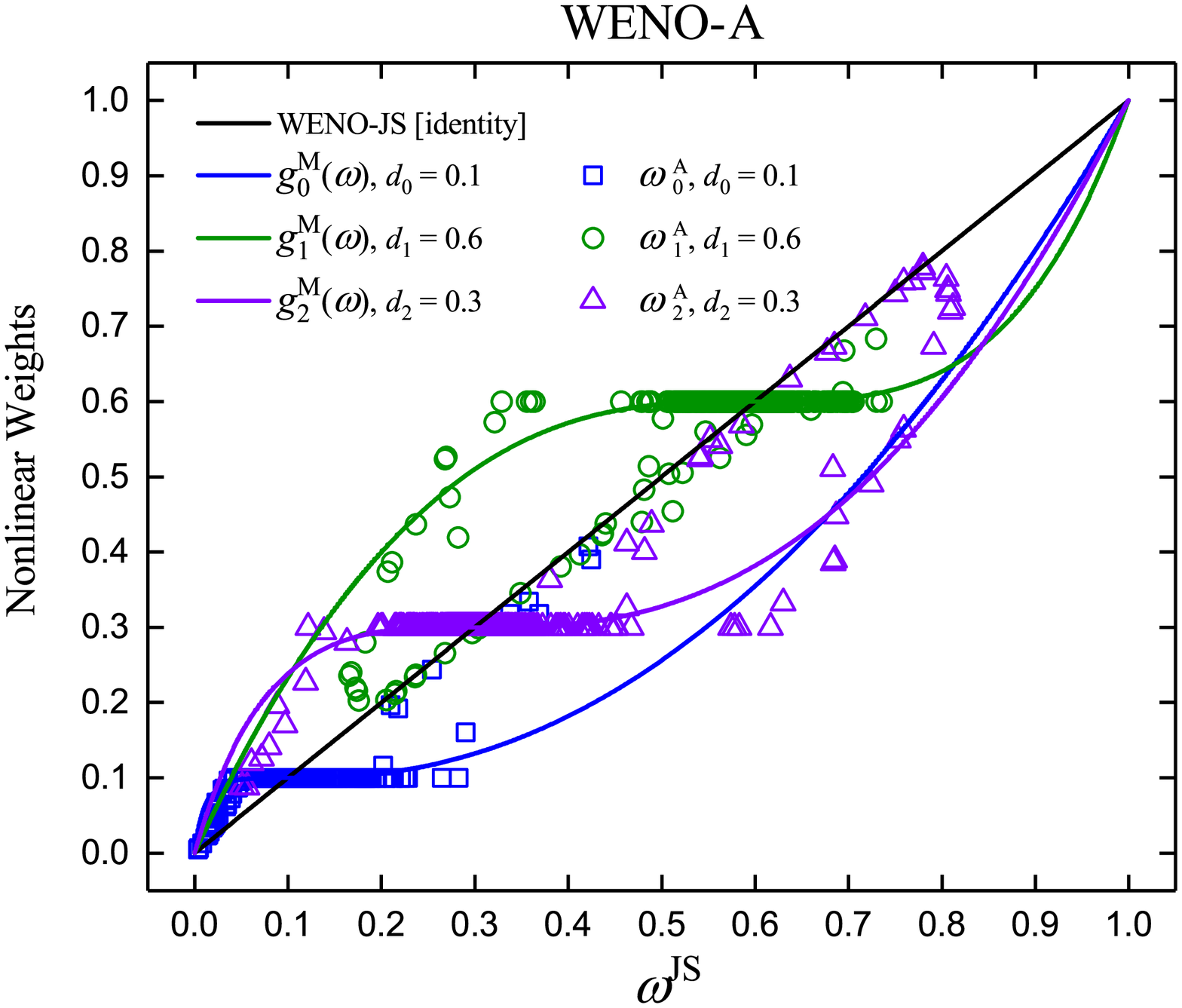}
  \includegraphics[height=0.28\textwidth]
  {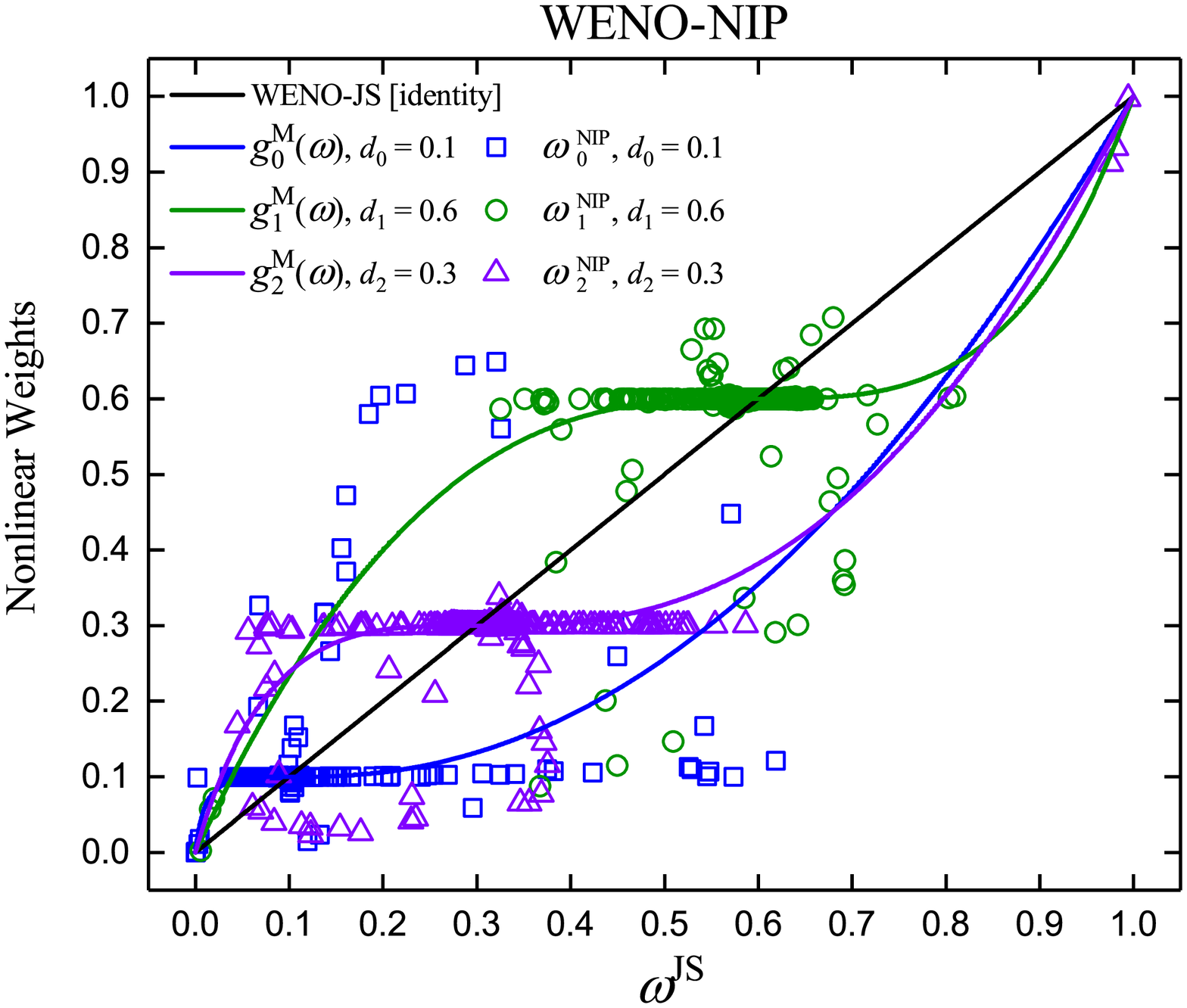}  
\caption{The IMRs for various WENO-Z-type schemes.}
\label{fig:IMR:Z-type}
\end{figure}

Taking SLP as an example, by choosing $t = 2$ and $N = 800$, we 
plot the IMRs of the WENO-Z, WENO-Z$\eta(\tau_{5})$, 
WENO-Z$\eta(\tau_{81})$, WENO-Z+, WENO-ZA, WENO-D, WENO-A and 
WENO-NIP schemes in Fig. \ref{fig:IMR:Z-type}. For comparison 
purpose, we also plot the designed mapping of WENO-JS, or 
say identity mapping, and that of WENO-M. It is very 
interesting that the IMRs of all the considered WENO-Z-type schemes 
are similar to the designed mapping of the WENO-M scheme in general. 
In other words, it is quite obvious that there are ``optimal weight 
intervals'' where the mappings, as well as the Z-type weight formulas, replace the nonlinear weights using the associated ideal weights. We refer to subsection 3.4.2 of \cite{WENO-MAIMi} for more details about the ``optimal weight intervals''. As far as we know, these features are not accounted for in presiously published articles. However, it is indicated that \cite{WENO-IM,WENO-AIM,WENO-MAIMi} the {\it{optimal weight interval}} plays a critical part in recovering the designed convergence orders in the presence of critical points and improving 
the resolution of the corresponding WENO scheme. Therefore, based on 
this meaningful observation, we innovatively build the concept of 
{\it{Generalized Mapped WENO schemes}}.

\begin{definition}
({\rm{Generalized Mapped WENO scheme}}) If the IMR of WENO-X 
exhibits the key characteristic of the traditional mapped 
WENO scheme, that is, existing apparent {\it{optimal weight 
intervals}},we say WENO-X is a Generalized Mapped WENO scheme.
\end{definition}

In fact, the various Z-type weights can be written in the following 
generalized form
\begin{equation}
\begin{array}{ll}
\left\{
\begin{array}{l}
\begin{aligned}
&\omega_{s}^{\mathrm{X}} = \dfrac{\alpha_{s}^{\mathrm{X}}}{\sum_{l=0}^{2}\alpha_{l}^{\mathrm{X}}}, \\
&\alpha_{s}^{\mathrm{X}} = \psi_{s,1}^{\mathrm{X}}(d_{s}) + \dfrac{d_{s}/(\beta_{s}^{\mathrm{JS}} + \epsilon)^{2}}{\sum_{l=0}^{2}\alpha_{l}^{\mathrm{JS}}}\cdot\psi_{s,2}^{\mathrm{X}}(d_{s}) + 
\psi_{s,3}^{\mathrm{X}}(d_{s}) = \psi_{s,1}^{\mathrm{X}}(d_{s}) + 
\omega_{s}^{\mathrm{JS}}\cdot\psi_{s,2}^{\mathrm{X}}(d_{s}) + 
\psi_{s,3}^{\mathrm{X}}(d_{s}), 
\end{aligned}
\end{array}
\right. & s = 0,1,2,
\end{array}
\label{eq:omegaWENO-Z_general}
\end{equation}
where the $\psi_{s,1}^{\mathrm{X}}(d_{s})$, 
$\psi_{s,2}^{\mathrm{X}}(d_{s})$, and $\psi_{s,3}^{\mathrm{X}}(d_{s})$ 
can be obtained trivially by some simple mathematical manipulations. 
As examples, we provide the $\psi_{s,1}^{\mathrm{X}}(d_{s})$, 
$\psi_{s,2}^{\mathrm{X}}(d_{s})$, and $\psi_{s,3}^{\mathrm{X}}(d_{s})$ 
of the WENO-Z, WENO-Z$\eta(\tau_{5})$, WENO-Z$\eta(\tau_{81})$, 
WENO-Z+, WENO-ZA, WENO-D, WENO-A and WENO-NIP schemes in Table 
\ref{tab:psi:Z-type}.

Furthermore, we give a uniform formula of the nonlinear weights of 
all three versions of the odd-order WENO schemes, taking 
the form
\begin{equation}
\omega_{s}^{\mathrm{X}} = \dfrac{\alpha_{s}^{\mathrm{X}}}{
\sum_{l=0}^{2}\alpha_{l}^{\mathrm{X}}},\alpha_{s}^{\mathrm{X}} = 
\Big(g^{\mathrm{GMX}}\Big)_{s}(\omega_{s}^{\mathrm{Y}})=
\psi_{s,1}^{\mathrm{X}}(d_{s}) + \mathcal{H}^{\mathrm{X}}(\omega_{s}^{
\mathrm{Y}}, d_{s})\cdot\psi_{s,2}^{\mathrm{X}}(d_{s}) + \psi_{s,3}^{
\mathrm{X}}(d_{s}), \quad s = 0,1,2,
\label{eq:omega_general:final}
\end{equation}
where $\omega_{s}^{\mathrm{Y}}$ is the nonlinear weights of some 
classical WENO scheme, and $\omega_{s}^{\mathrm{Y}} = 
\omega_{s}^{\mathrm{JS}}$ is used usually. Also, 
$\psi_{s,i}^{\mathrm{X}}(d_{s}), i=1,2,3$ and 
$\mathcal{H}^{\mathrm{X}}(\omega_{s}^{\mathrm{Y}})$ can easily be 
obtained. As examples, we provide the $\psi_{s,i}^{\mathrm{X}}(d_{s})
, i=1,2,3$ and $\mathcal{H}^{\mathrm{X}}(\omega_{s}^{\mathrm{Y}})$ 
of the WENO-JS, WENO-M, WENO-Z and WENO-Z+ schemes in Table 
\ref{tab:psi:Z-type:final}.

\begin{definition}
({\rm{Generalized Mapping Function}}) We call 
$\Big(g^{\mathrm{GMX}}\Big)_{s}(\omega_{s}^{\mathrm{Y}})$ in Eq.
\eqref{eq:omega_general:final} the Generalized Mapping Function of WENO-X, and for simplicity, $\Big(g^{\mathrm{GMX}}\Big)_{s}(\omega)$ is used without causing any confusion.
\end{definition}

\begin{corollary}
If WENO-X is a traditional mapped WENO scheme, then its generalized 
mapping function is exactly its designed mapping function; if WENO-X
is a WENO-Z-type WENO scheme, then its generalized mapping function 
is the nonnormalized form of the associated IMR.
\label{corollary:GMF}
\end{corollary}

The proof of Corollary \ref{corollary:GMF} is very trivial and we 
omit it here just for the sake of brevity.

\begin{table}[ht]
\renewcommand\arraystretch{1.5}
\footnotesize
\centering
\caption{$\psi_{s,1}^{\mathrm{X}}(d_{s})$, 
$\psi_{s,2}^{\mathrm{X}}(d_{s})$, and $\psi_{s,3}^{\mathrm{X}}(d_{s})$ 
of various Z-type weights.}
\label{tab:psi:Z-type}
\begin{tabular*}{\hsize}
{@{}@{\extracolsep{\fill}}llll@{}}
\hline
Scheme, WENO-X & $\psi_{s,1}^{\mathrm{X}}(d_{s})$ 
               & $\psi_{s,2}^{\mathrm{X}}(d_{s})$ 
               & $\psi_{s,3}^{\mathrm{X}}(d_{s})$\\
\hline
WENO-Z    & $d_{s}$ 
          & $\sum_{l=0}^{2}\alpha_{l}^{\mathrm{JS}}\cdot\tau_{5}^{2}$
          & $0$\\
WENO-Z$\eta(\tau_{5})$   
          & $d_{s}$ 
          & $\sum_{l=0}^{2}\alpha_{l}^{\mathrm{JS}}\cdot\dfrac{(
            \beta_{s}^{\mathrm{JS}} + \epsilon)^{2}}{(\eta_{s} + 
            \epsilon)^{2}}\cdot\tau_{5}^{2}$  
          & $0$\\
WENO-Z$\eta(\tau_{81})$   
          & $d_{s}$ 
          & $\sum_{l=0}^{2}\alpha_{l}^{\mathrm{JS}}\cdot\dfrac{(
            \beta_{s}^{\mathrm{JS}} + \epsilon)^{2}}{(\eta_{s} + 
            \epsilon)^{2}}\cdot\tau_{81}^{2}$  
          & $0$\\
WENO-Z$+$ & $d_{s}$ 
          & $\sum_{l=0}^{2}\alpha_{l}^{\mathrm{JS}}\cdot(\tau_{5} + 
            \epsilon)^{2}$     
          & $d_{s}\cdot\lambda\bigg( \dfrac{\beta_{s} + 
            \epsilon}{\tau_{5} + \epsilon} \bigg) $\\
WENO-ZA   & $d_{s}$ 
          & $\sum_{l=0}^{2}\alpha_{l}^{\mathrm{JS}}\cdot\dfrac{(
            \beta_{s}^{\mathrm{JS}} + \epsilon)^{2}}{(\beta_{s}^{
            \mathrm{ZA}} + \epsilon)}\cdot A\tau_{6}$  
          & $0$\\
WENO-D    & $d_{s}$ 
          & $\sum_{l=0}^{2}\alpha_{l}^{\mathrm{JS}}\cdot 
          \Phi\tau_{5}^{p}\cdot(\beta_{s}^{\mathrm{JS}} + 
          \epsilon)^{2-p}$     
          & $0$ \\
WENO-A    & $d_{s}\cdot B$ 
          & $\sum_{l=0}^{2}\alpha_{l}^{\mathrm{JS}}\cdot 
          \Phi\tau_{5}^{p}\cdot(\beta_{s}^{\mathrm{JS}} + 
          \epsilon)^{2-p}(1 - B)$
          & $0$\\
WENO-NIP  & $d_{s}$ 
          & $\sum_{l=0}^{2}\alpha_{l}^{\mathrm{JS}}\cdot\tau\cdot
            \dfrac{(\beta_{s}^{\mathrm{JS}} + 
            \epsilon)^{2}}{(\chi_{s} + \epsilon)^{2}}$  
          & $0$\\
\hline
\multicolumn{4}{l}{\tabincell{l}{(1) $B = \mathrm{BOOL}\Bigg(\max
\bigg(1, \Phi\bigg(\dfrac{\tau_{5}}{\beta_{s} + \epsilon}\bigg)^{p}
\bigg) = 1 \Bigg)$, $\Phi$ is computed by Eq. \eqref{eq:D:phi}; \\ 
(2) the other parameters are the same as in subsection 
\ref{subsec:WENO-Z} and subsection \ref{subsec:WENO-Z-type}.}} \\ 
\hline
\end{tabular*}
\end{table}

\begin{table}[ht]
\renewcommand\arraystretch{1.5}
\footnotesize
\centering
\caption{$\psi_{s,i}^{\mathrm{X}}(d_{s}), i=1,2,3$ and 
$\mathcal{H}^{\mathrm{X}}(\omega_{s}^{\mathrm{Y}})$ 
for WENO-JS, WENO-M, WENO-Z and WENO-Z+.}
\label{tab:psi:Z-type:final}
\begin{tabular*}{\hsize}
{@{}@{\extracolsep{\fill}}llllll@{}}
\toprule
Scheme, WENO-X & $\omega_{s}^{\mathrm{Y}}$ 
               & $\psi_{s,1}^{\mathrm{X}}(d_{s})$ 
               & $\psi_{s,2}^{\mathrm{X}}(d_{s})$ 
               & $\psi_{s,3}^{\mathrm{X}}(d_{s})$ 
               & $\mathcal{H}^{\mathrm{X}}(\omega_{s}^{\mathrm{Y}}, 
                 d_{s})$\\
\hline
WENO-JS        & $\omega_{s}^{\mathrm{JS}}$     
               & $0$        
               & $1$  
               & $0$  
               & $\omega_{s}^{\mathrm{Y}}$   \\
WENO-M         & $\omega_{s}^{\mathrm{JS}}$     
               & $d_{s}$    
               & $1$  
               & $0$  
               & $\Big(g^{\mathrm{M}}\Big)_{s}(\omega_{s}^{
               \mathrm{Y}})=\dfrac{(\omega_{s}^{\mathrm{Y}}-d_{s}
               )^{3}}{(\omega_{s}^{\mathrm{Y}}-d_{s})^{2}+\omega_{s
               }^{\mathrm{Y}} (1-\omega_{s}^{\mathrm{Y}})}$\\
WENO-Z         & $\omega_{s}^{\mathrm{JS}}$     
               & $d_{s}$    
               & $\sum_{l=0}^{2}\alpha_{l}^{\mathrm{JS}}\cdot
                 \tau_{5}^{2}$  
               & $0$  
               & $\omega_{s}^{\mathrm{Y}}$   \\
WENO-Z+        & $\omega_{s}^{\mathrm{JS}}$     
               & $d_{s}$    
               & $\sum_{l=0}^{2}\alpha_{l}^{\mathrm{JS}}\cdot
                 (\tau_{5}+\epsilon)^{2}$   
               & $d_{s}\cdot\lambda\bigg( \dfrac{\beta_{s} + 
               \epsilon}{\tau_{5} + \epsilon} \bigg) $ 
               & $\omega_{s}^{\mathrm{Y}}$ \\
\bottomrule
\end{tabular*}
\end{table}


\section{Design and properties of WENO-Z-type schemes with  
order-preserving generalized mappings}
\label{sec:MOP-GMWENO-X}
\subsection{The new WENO-Z-type schemes}\label{subsec:MOP-GMWENO-X}
In order to provide convenience to the reader and clarify our major 
concern, we firstly restate the following notations and 
preliminaries, originally proposed in \cite{MOP-WENO-X}, that will 
be needed for the present study.

\begin{definition}
Construct the {\rm{\textbf{Min\_DIST}}} function taking the form
\begin{equation}\left\{
\begin{array}{ll}
\bm{\mathrm{Min\_DIST}}\big(v_{0},\cdots,v_{r-1};\bar{\omega}_{0},\cdots,\bar{\omega}_{r-1};\omega\big) = v_{l^{*}}, \\
l^{*} = \min\Bigg(\mathrm{INDEX}\bigg(\min\Big\{\lvert \omega-\bar{\omega}_{0}\rvert, \cdots, \lvert \omega-\bar{\omega}_{r-1}
\rvert\Big\}\bigg)\Bigg),
\label{eq:minDist}
\end{array}\right.
\end{equation}
where $\bar{\omega}_{s}$ stands for the ideal weight and 
\begin{equation}
\mathrm{INDEX}\bigg(\min\Big\{\lvert \omega-\bar{\omega}_{0}
\rvert, \cdots, \lvert \omega-\bar{\omega}_{r-1}
\rvert\Big\}\bigg) = \Big\{q_{1},\cdots,q_{M}\Big\},
\label{eq:IndexOf}
\end{equation}
if $\min\Big\{\lvert \omega-\bar{\omega}_{0}\rvert, \cdots, \lvert \omega - \bar{\omega}_{r-1}\rvert\Big\} = \lvert \omega - 
\bar{\omega}_{q_{1}} \rvert = \cdots = \lvert \omega - \bar{\omega}_{q_{M}} \rvert$.
\label{definition:minDist}
\end{definition}

\begin{definition}
Denote $\mathcal{S} = \Big\{d_{0},\cdots,d_{r-1}\Big\}$ and
$\mathcal{\widehat{S}} = \Big\{\widehat{d}_{0}, \cdots,\widehat{d}_{r-1}\Big\}$, and they satisfy: (1) the elements of $\mathcal{S}$ and $\mathcal{\widehat{S}}$ are only different in turn; (2) $ 0 < \widehat{d}_{0} < \cdots < \widehat{d}_{r-1} < 1$. Let $\mathcal{G} = \Big\{ \big( g^{\mathrm{GMX}} \big)_{0}(\omega), \cdots, \big( 
g^{\mathrm{GMX}} \big)_{r-1}(\omega) \Big\}$, we define $\mathcal{\widehat{G}}=\Big\{\widehat{\big(g^{\mathrm{GMX}}\big)}_{0}(\omega),\widehat{\big(g^{\mathrm{GMX}}\big)}_{1}(\omega), \cdots, \widehat{\big(g^{\mathrm{GMX}}\big)}_{r-1}(\omega) \Big\}$, where $\widehat{\big(g^{\mathrm{GMX}}\big)}_{s}(\omega)$ is computed by substituting $\widehat{d}_{s}$ into the corresponding $\big(g^{\mathrm{GMX}}\big)_{s}(\omega)$.
\label{def:newMappingSet}
\end{definition}

\begin{lemma}
Assume $\widehat{d}_{-1} = 0, \widehat{d}_{r} = 1$, and denote 
$\ddot{d}_{-1} = \widehat{d}_{-1}, \ddot{d}_{0} = 
(\widehat{d}_{0} + \widehat{d}_{1})/2, \cdots,
\ddot{d}_{r-2}=(\widehat{d}_{r-2}+\widehat{d}_{r-1})/2
, \ddot{d}_{r-1}$ $= \widehat{d}_{r}$, then we have
\begin{equation*}
\min\Bigg(\mathrm{INDEX}\bigg(\min\Big\{\lvert \omega-
\widehat{d}_{0}\rvert, \cdots, \lvert \omega-\widehat{d}_{r-1}\rvert\Big\}\bigg)\Bigg)= j, \quad \mathrm{for} \quad \forall\omega \in (\ddot{d}_{j-1}, \ddot{d}_{j}].
\end{equation*}
\label{lemma:Omega_i:minDist}
\end{lemma}
\textbf{Proof.} We refer to the proof of the Lemma 3 in 
\cite{MOP-WENO-X}.
$\hfill\square$ \\

\begin{lemma}
Let $\Omega_{i}=\Big\{\omega \lvert \bm{\mathrm{Min\_DIST}}(
\widehat{d}_{0}, \cdots, \widehat{d}_{r-1}; 
\widehat{d}_{0}, \cdots, \widehat{d}_{r-1}; 
\omega) = \widehat{d}_{i} \Big\}$ and $a, b\in\{0,1,\cdots,r-1\}$, 
if $a > b$, then for $\forall \omega_{\alpha} \in \Omega_{a}$ and $\forall \omega_{\beta} \in \Omega_{b}$, one can get

C1. $\Omega_{i} = (\ddot{d}_{i-1}, \ddot{d}_{i}]$; 

C2. $\omega_{\alpha}>\omega_{\beta}, \widehat{d}_{\alpha}>\widehat{d}_{\beta}$; 

C3. $\widehat{\big(g^{\mathrm{GMX}}\big)}_{a}(\omega_{\alpha}) > \widehat{\big(g^{\mathrm{GMX}}\big)}_{b}(\omega_{\beta})$ .
\label{lemma:Omega_i}
\end{lemma}
\textbf{Proof.}
According to Definition \ref{definition:minDist} and Lemma \ref{lemma:Omega_i:minDist}, we can prove \textit{C1} trivially. As $a > b$, then $\Omega_{b}$ is on the left of $\Omega_{a}$. It is already given that $\omega_{\alpha} \in \Omega_{a}, \omega_{\beta} 
\in \Omega_{b}$, so 
\begin{equation}
\omega_{\alpha} > \omega_{\beta}, \quad
\widehat{d}_{\alpha}>\widehat{d}_{\beta}. 
\label{eq:proof:lemma:Omega_i:01}
\end{equation}
Thus, the \textit{C2} is proved.

The proof of \textit{C3} is case by case. For simplicity but without loss of 
generality, we present the proof for the case of ``WENO-X = 
WENO-Z'' as an example.
According to Eq. \eqref{eq:omega_general:final}, we can get
\begin{equation}
\begin{array}{l}
\begin{aligned}
&\widehat{\big(g^{\mathrm{GMZ}}\big)}_{a}(\omega_{\alpha}) = 
\psi_{1}^{\mathrm{Z}}(\widehat{d}_{\alpha}) + 
\mathcal{H}^{\mathrm{Z}}(\omega_{\alpha}, \widehat{d}_{\alpha})
\cdot\psi_{2}^{\mathrm{Z}}(\widehat{d}_{\alpha}) + 
\psi_{3}^{\mathrm{Z}}(\widehat{d}_{\alpha}), \\
&\widehat{\big(g^{\mathrm{GMZ}}\big)}_{b}(\omega_{\beta}) = 
\psi_{1}^{\mathrm{Z}}(\widehat{d}_{\beta}) + 
\mathcal{H}^{\mathrm{Z}}(\omega_{\beta}, \widehat{d}_{\beta})
\cdot\psi_{2}^{\mathrm{Z}}(\widehat{d}_{\beta}) + 
\psi_{3}^{\mathrm{Z}}(\widehat{d}_{\beta}).
\end{aligned}
\end{array}
\label{eq:proof:lemma:Omega_i:02}
\end{equation}
From Table \ref{tab:psi:Z-type:final}, we can obtain
\begin{equation}
\begin{array}{l}
\begin{aligned}
&\mathcal{H}^{\mathrm{Z}}(\omega_{\alpha},
\widehat{d}_{\alpha}) = \omega_{\alpha}, \quad
\mathcal{H}^{\mathrm{Z}}(\omega_{\beta},
\widehat{d}_{\beta}) = \omega_{\beta}, \\
&\psi_{1}^{\mathrm{Z}}(\widehat{d}_{\alpha})=
\widehat{d}_{\alpha},\quad 
\psi_{1}^{\mathrm{Z}}(\widehat{d}_{\beta})=\widehat{d}_{\beta},\\
&\psi_{2}^{\mathrm{Z}}(\widehat{d}_{\alpha})=
\psi_{2}^{\mathrm{Z}}(\widehat{d}_{\beta})=
\sum\nolimits_{l=0}^{2}\alpha_{l}^{\mathrm{JS}}\cdot\tau_{5}^{2}>0,\\
&\psi_{3}^{\mathrm{Z}}(\widehat{d}_{\alpha})=
\psi_{3}^{\mathrm{Z}}(\widehat{d}_{\beta})= 0.
\end{aligned}
\end{array}
\label{eq:proof:lemma:Omega_i:03}
\end{equation}
Then, from Eqs. \eqref{eq:proof:lemma:Omega_i:01}\eqref{eq:proof:lemma:Omega_i:02}\eqref{eq:proof:lemma:Omega_i:03}, 
we get 
\begin{equation*}
\widehat{\big(g^{\mathrm{GMZ}}\big)}_{a}(\omega_{\alpha}) >
\widehat{\big(g^{\mathrm{GMZ}}\big)}_{b}(\omega_{\beta}).
\end{equation*}
Thus, \textit{C3} for the case of ``WENO-X = WENO-Z'' is proved. The other cases can 
also be proved analogously and they are omitted here to save 
space. Actually, from Fig. \ref{fig:IMR:GMZ} below, it can be intuitively observed from the IMR curve of MOP-GMWENO-X, where MOP-GMWENO-X stands for the new WENO-Z-type scheme 
derived from the WENO-X scheme taking the OP property and we will 
present it in the following. $\hfill\square$ \\

It has been reported \cite{MOP-WENO-ACMk,MOP-WENO-X} that the OP 
property of the mapping plays an essential 
role in preserving high resolutions and meanwhile avoiding spurious 
oscillations for long-run calculations. Now, we extend the concept of the 
order-preserving mapping (see Definition 2 of \cite{MOP-WENO-ACMk} 
or Definition 1 of \cite{MOP-WENO-X}) to the generalized mapping 
function as follows.

\begin{definition} For $\forall m > n$, and $\omega_{l_{m}} \in \Omega_{m}$, $\omega_{l_{n}} \in \Omega_{n}$, if
\begin{equation}
\big(g^{\mathrm{GMX}}\big)_{l_{m}}(\omega_{l_{m}})
>\big(g^{\mathrm{GMX}}\big)_{l_{n}}(\omega_{l_{n}}),
\label{def:order-preserving_mapping}
\end{equation} 
then the set of generalized mappings \Big\{$\big( g^{\mathrm{GMX}}\big)_{s}(\omega), s = 0,\cdots,r-1$\Big\} is called to be \textbf{order-preserving (OP) mapping}. Otherwise, it is called to be \textbf{non-order-preserving (non-OP) mapping}.
\label{def:OPGM}
\end{definition}


In Algorithm \ref{alg:minDist}, we introduce the \textit{OP} 
property into the previously published WENO-Z-type schemes whose 
generalized mappings are \textit{non-OP}. Before giving Algorithm 
\ref{alg:minDist}, we rewrite the generalized mapping 
$\big( g^{\mathrm{GMX}} \big)_{s}(\omega)$ defined in Eq. 
\eqref{eq:omega_general:final} in a more meaningful form
\begin{equation}
\big( g^{\mathrm{GMX}} \big)_{s}(\omega_{s}^{\mathrm{JS}}) = 
g^{\mathrm{GMX}}\Big(\omega^{\mathrm{JS}}_{s}; 
  \psi_{s,1},\psi_{s,2},\psi_{s,3}\Big).
\label{eq:rewrite:omega_general:final}
\end{equation}

\begin{algorithm}[htb]
\caption{The method to get new Z-type weights satisfying 
the \textit{OP} property.}
\label{alg:minDist}
\SetKwInOut{Input}{input}\SetKwInOut{Output}{output}
\Input{$s$, $d_{s}$, $\omega^{\mathrm{JS}}_{s}$, $\psi_{s,m}$ with
$m = 1,2,3$}
\Output{the new \textit{OP} generalized mappings, namely $\Big\{\big(g^{\mathrm{MOP-GMX}}\big)_{s}(\omega^{\mathrm{JS}}_{s}), s=0,1,\cdots,r-1\Big\}$}
\BlankLine
\emph{$\big( g^{\mathrm{GMX}} \big)_{s}(\omega)$
is a generalized mapping function, $\Big\{\big(g^{\mathrm{GMX}} \big)_{s}(\omega), s = 0, 1,\cdots,r-1 \Big\}$ is \textit{non-OP}}\;
set $s = 0$\;
\While{$s \leq r - 1$}{  
  compute $\omega_{s}^{\mathrm{JS}}$\;
  set $j = 1, d_{\mathrm{MIN}}=\lvert\omega_{s}^{\mathrm{JS}}-d_{0}\rvert, l_{s}^{*} = 0$\;  
  \While{$j \leq r - 1$}{
    \While{$-d_{\mathrm{MIN}} < \omega_{s}^{\mathrm{JS}}-d_{j}<
    d_{\mathrm{MIN}}$}{
      $d_{\mathrm{MIN}}=\lvert\omega_{s}^{\mathrm{JS}}-d_{j}\rvert$\;
      $l_{s}^{*} = j$\;
    } 
    $j ++$\;
  }
  set $m = 1$\;
  \While{$m \leq 3$}{
    $\overline{\psi}_{s,m} = \psi_{l_{s}^{*},m}$\;
    $m ++$\;
  }
  $s ++$\;
}
set $s=0$\;
\While{$s \leq r - 1$}{
  $\big(g^{\mathrm{MOP-GMX}}\big)_{s}(\omega^{\mathrm{JS}}_{s}) = 
  g^{\mathrm{GMX}}\Big(\omega^{\mathrm{JS}}_{s}; 
  \overline{\psi}_{s,1},\overline{\psi}_{s,2},
  \overline{\psi}_{s,3}\Big)$\;
  $s ++$\;
}
\end{algorithm}

\begin{theorem}
$\Big\{\big(g^{\mathrm{MOP-GMX}}\big)_{s}(\omega^{\mathrm{JS}}_{s}), 
s=0,1,\cdots,r-1\Big\}$ computed by Algorithm \ref{alg:minDist} 
is \textit{OP}.
\label{theorem:g_MOP}
\end{theorem}
\textbf{Proof.} 
Without loss of generality, assume $\omega_{m}^{\mathrm{JS}}, \omega_{n}^{\mathrm{JS}} \in [0,1]$ and $\omega_{m}^{\mathrm{JS}} \in \Omega_{l_{m}^{*}}, \omega_{n}^{\mathrm{JS}} \in \Omega_{l_{n}^{*}}$, then
\begin{equation*}\left\{
\begin{array}{l}
\big(g^{\mathrm{MOP-GMX}}\big)_{m}(\omega^{\mathrm{JS}}_{m}) = 
g^{\mathrm{GMX}}\Big(\omega^{\mathrm{JS}}_{m}; \psi_{l_{m}^{*},1}, 
\psi_{l_{m}^{*},2}, \psi_{l_{m}^{*},3}\Big), \\
\big(g^{\mathrm{MOP-GMX}}\big)_{n}(\omega^{\mathrm{JS}}_{n}) = 
g^{\mathrm{GMX}}\Big(\omega^{\mathrm{JS}}_{n}; \psi_{l_{n}^{*},1}, 
\psi_{l_{n}^{*},2}, \psi_{l_{n}^{*},3}\Big).
\end{array}
\right.
\label{eq:gMOP_m-n}
\end{equation*}
Trivially, we have
\begin{equation*}\left\{
\begin{array}{l}
g^{\mathrm{GMX}}\Big(\omega^{\mathrm{JS}}_{m}; \psi_{l_{m}^{*},1}, 
\psi_{l_{m}^{*},2}, \psi_{l_{m}^{*},3}\Big) = 
\widehat{\big(g^{\mathrm{GMX}}\big)}_{l_{m}^{*}}
(\omega_{m}^{\mathrm{JS}}), \\
g^{\mathrm{GMX}}\Big(\omega^{\mathrm{JS}}_{n}; \psi_{l_{n}^{*},1}, 
\psi_{l_{n}^{*},2}, \psi_{l_{n}^{*},3}\Big) = 
\widehat{\big(g^{\mathrm{GMX}}\big)}_{l_{n}^{*}}
(\omega_{n}^{\mathrm{JS}}).
\end{array}
\right.
\end{equation*}
Up to now, the proof can be finished easily as Lemma \ref{lemma:Omega_i} holds true. 
$\hfill\square$ \\

Now, we give the improved Z-type weights satisfying the \textit{OP} property by
\begin{equation}
\omega_{s}^{\mathrm{MOP-GMX}} = 
\dfrac{\alpha_{s}^{\mathrm{MOP-GMX}}}{\sum_{l=0}^{r-1}\alpha_{l}^{
\mathrm{MOP-GMX}}}, \quad 
\alpha_{s}^{\mathrm{MOP-GMX}}=\big(g^{\mathrm{MOP-GMX}}\big)_{s}(
\omega^{\mathrm{JS}}_{s}), \quad s = 0,\cdots,r-1,
\label{eq:MOP-mapping}
\end{equation}
where $\big(g^{\mathrm{MOP-GMX}}\big)_{s}(\omega^{\mathrm{JS}}_{s})$ 
is determined by Algorithm \ref{alg:minDist} and the resultant 
scheme is named MOP-GMWENO-X.

In Fig. \ref{fig:IMR:GMZ}, we plot the profiles of the generalized 
mappings. It can be seen that, for the improved WENO-Z-type schemes: 
(1) the \textit{OP} property is obtained; (2) the widths of the 
``optimal weight intervals'' are decreased compared to the WENO-X schemes. 

Next, we conduct numerical examples to show that the 
convergence properties of the MOP-GMWENO-X schemes.

\begin{figure}[ht]
\flushleft
  \includegraphics[height=0.28\textwidth]
  {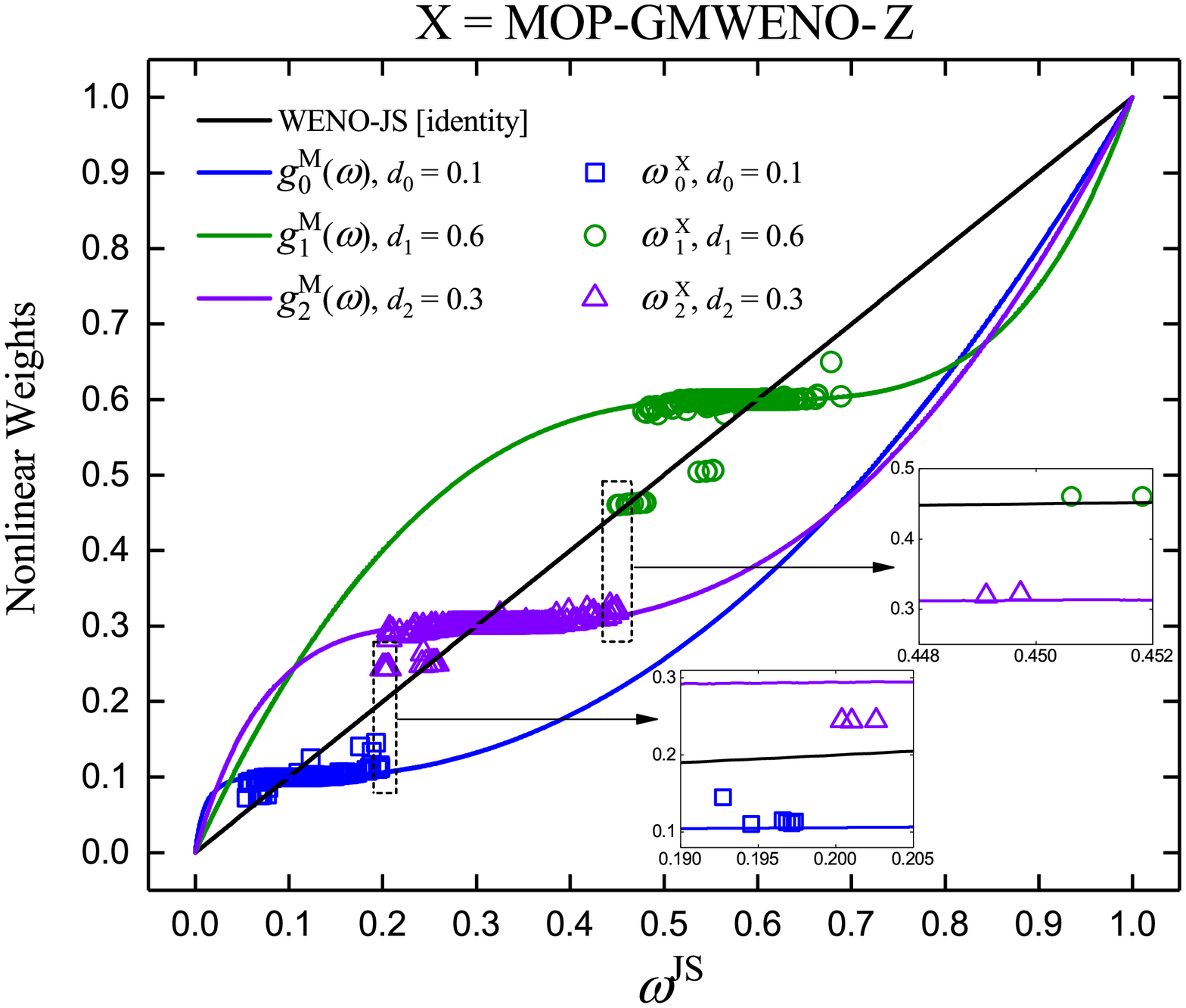}
  \includegraphics[height=0.28\textwidth]
  {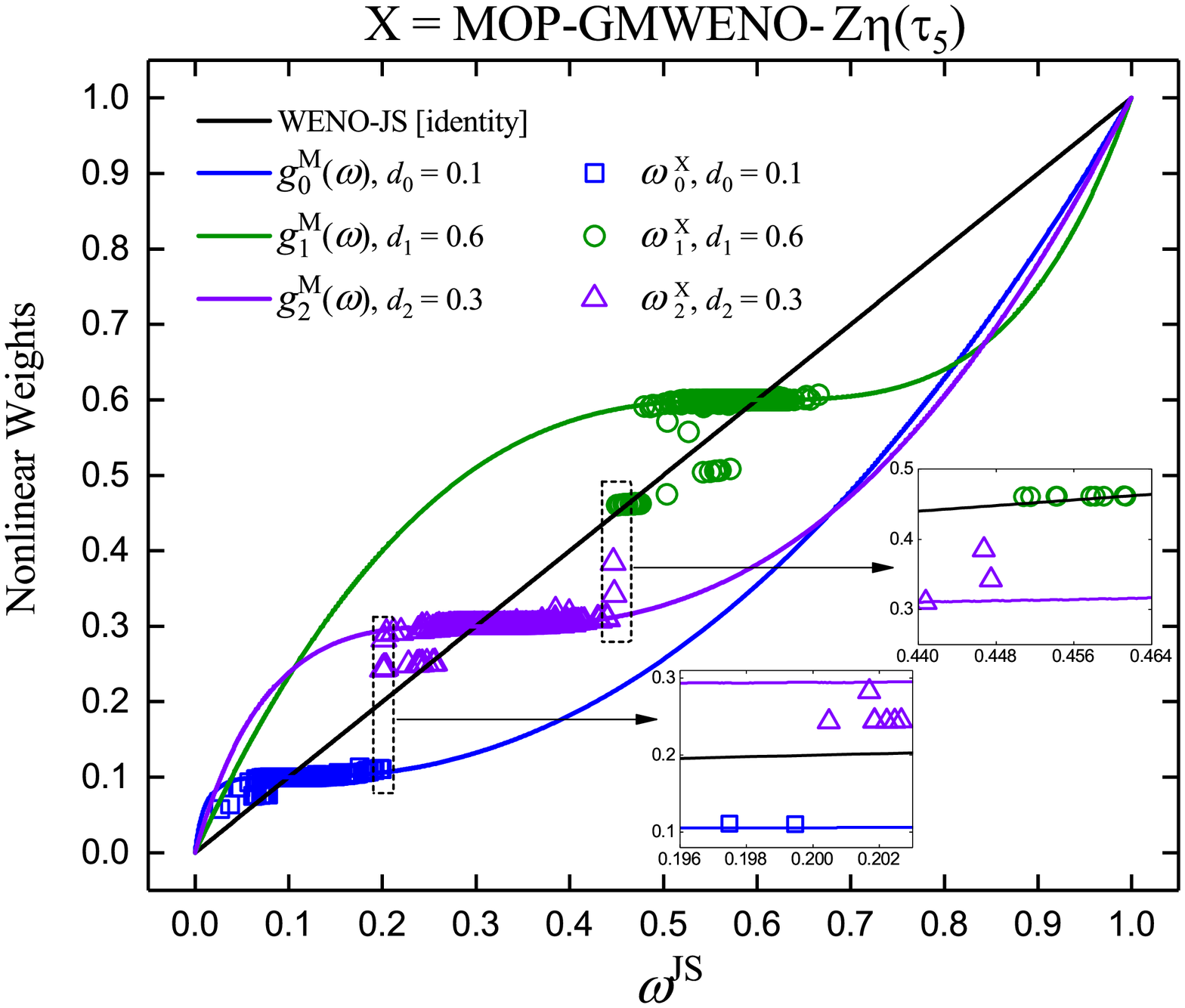}
  \includegraphics[height=0.28\textwidth]
  {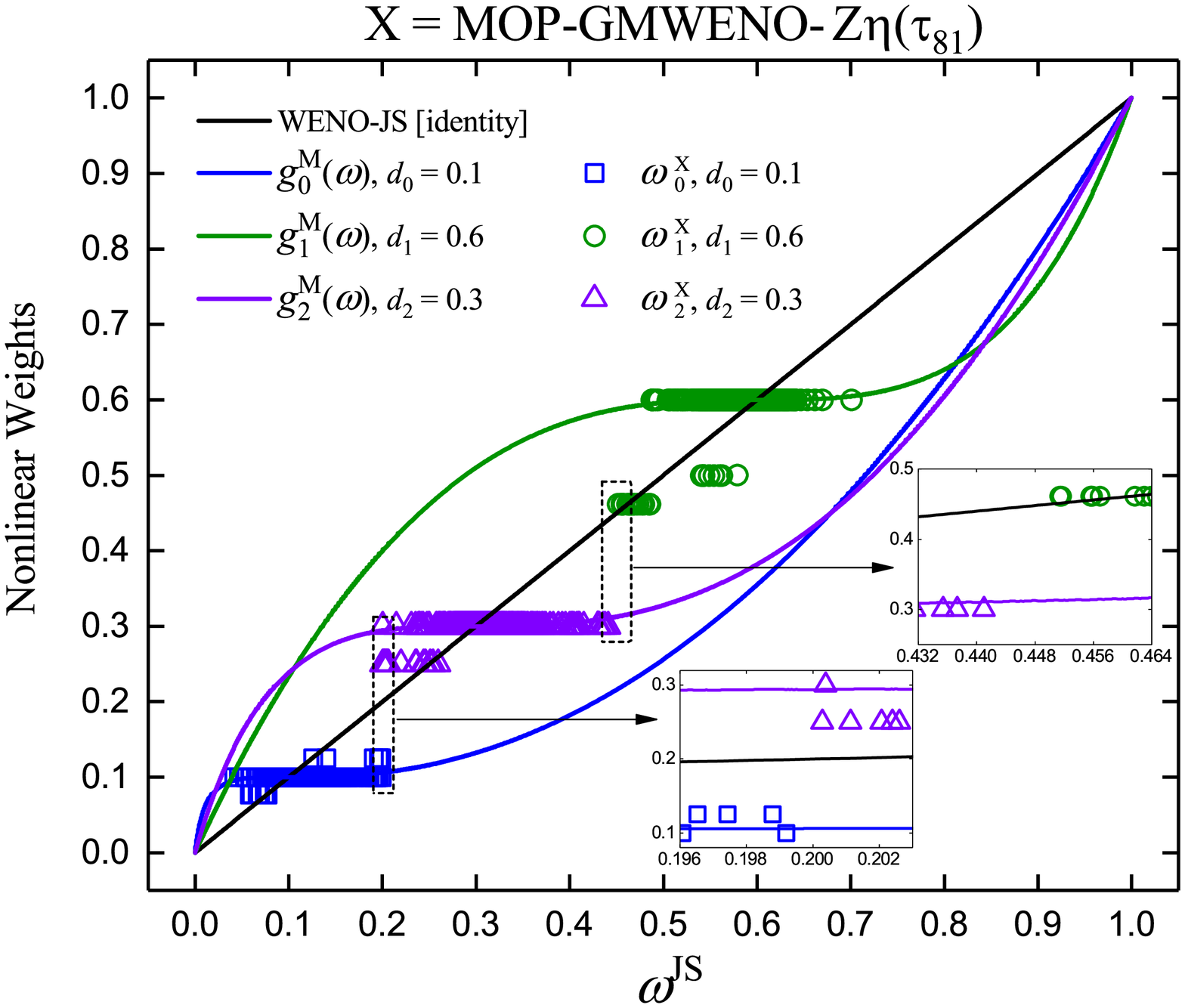}\\
  \includegraphics[height=0.28\textwidth]
  {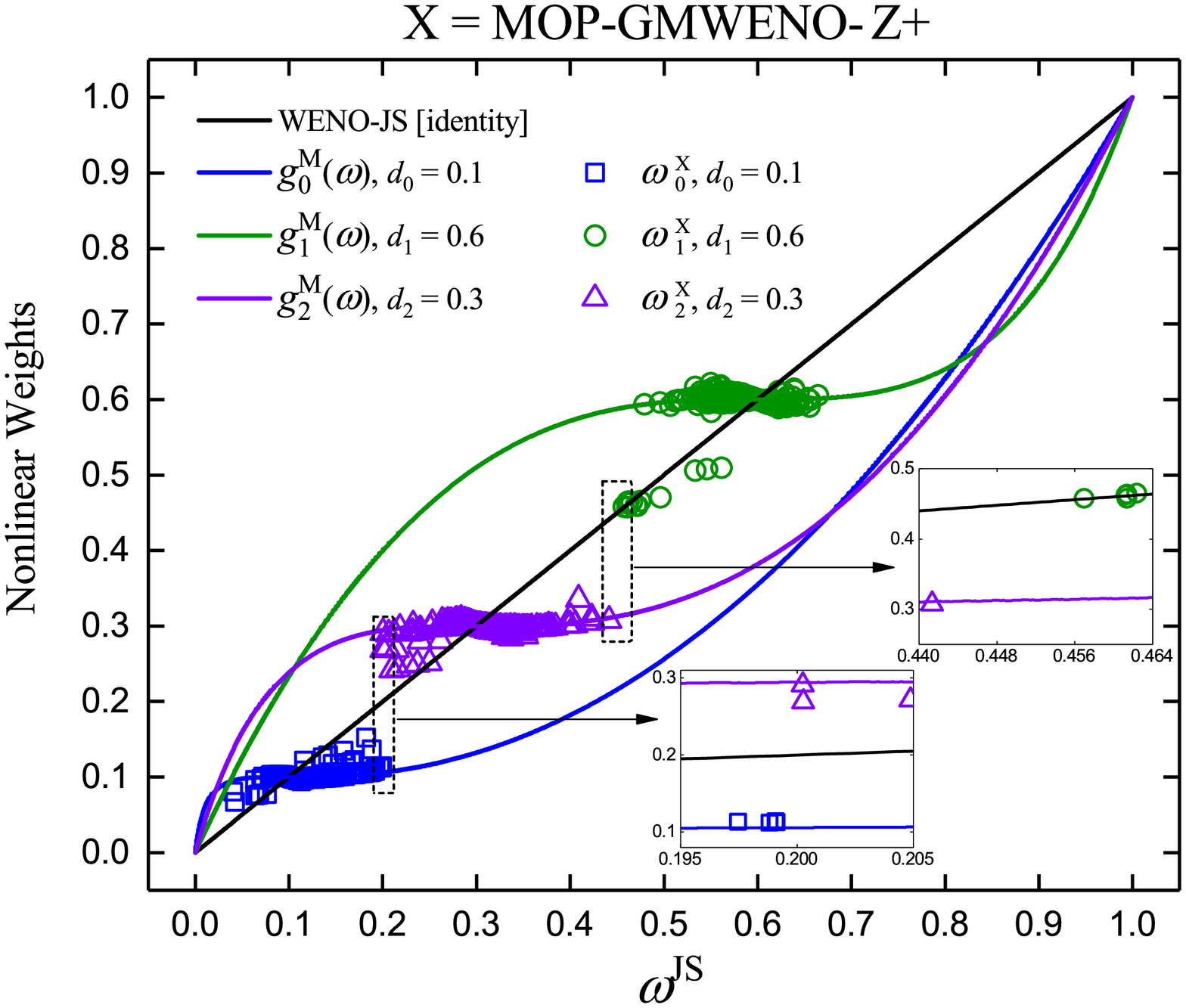}
  \includegraphics[height=0.28\textwidth]
  {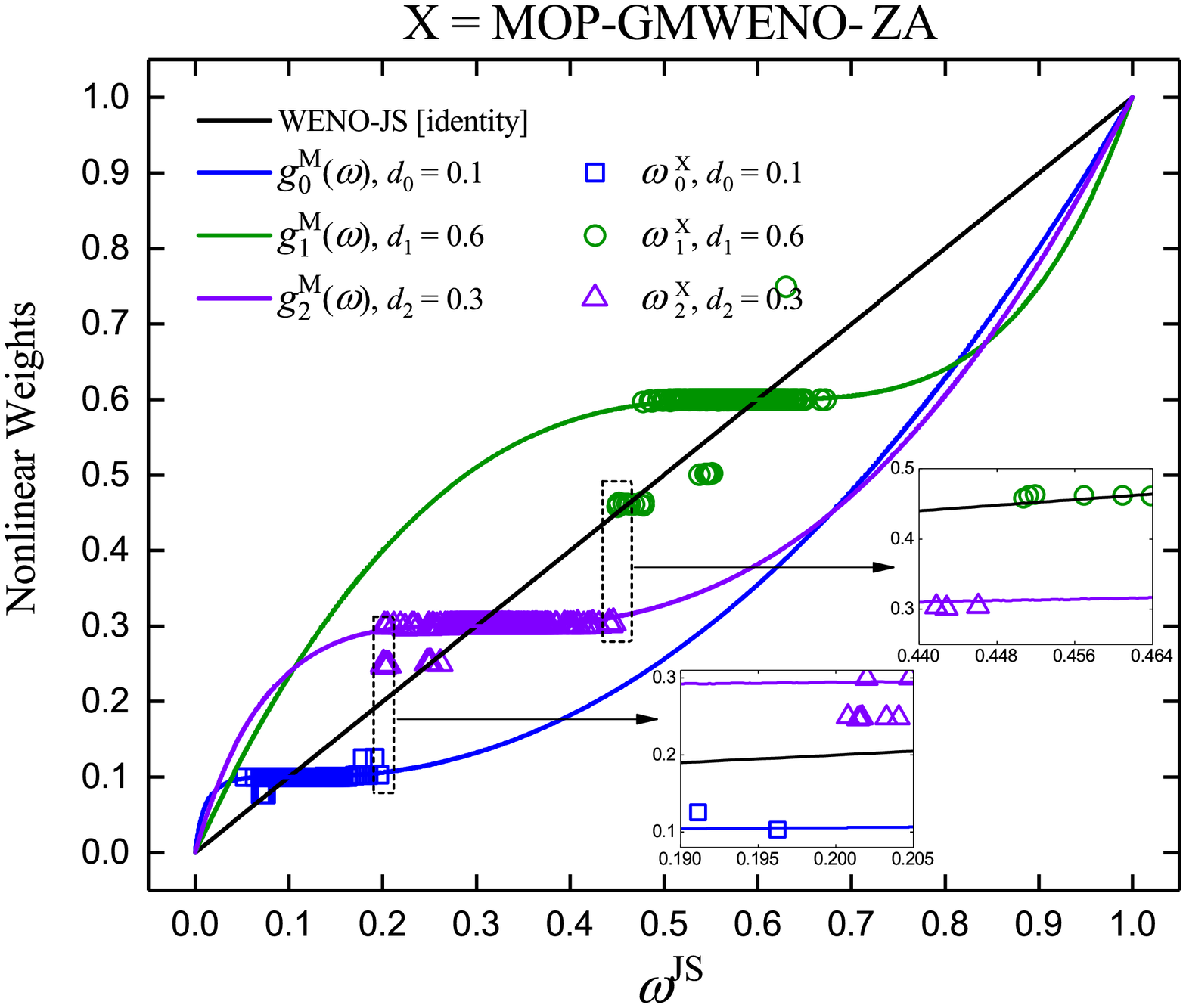}
  \includegraphics[height=0.28\textwidth]
  {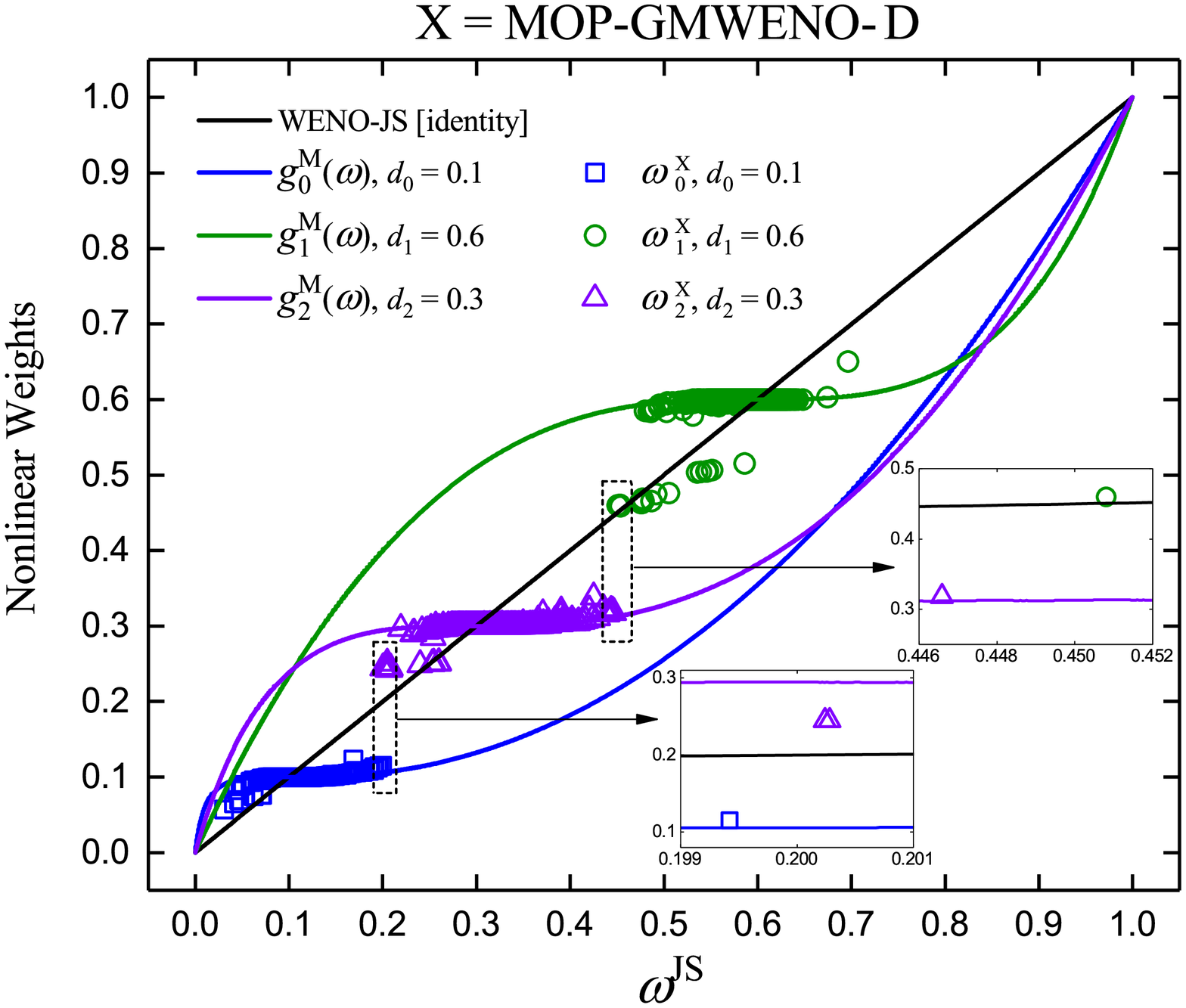}\\
  \includegraphics[height=0.28\textwidth]
  {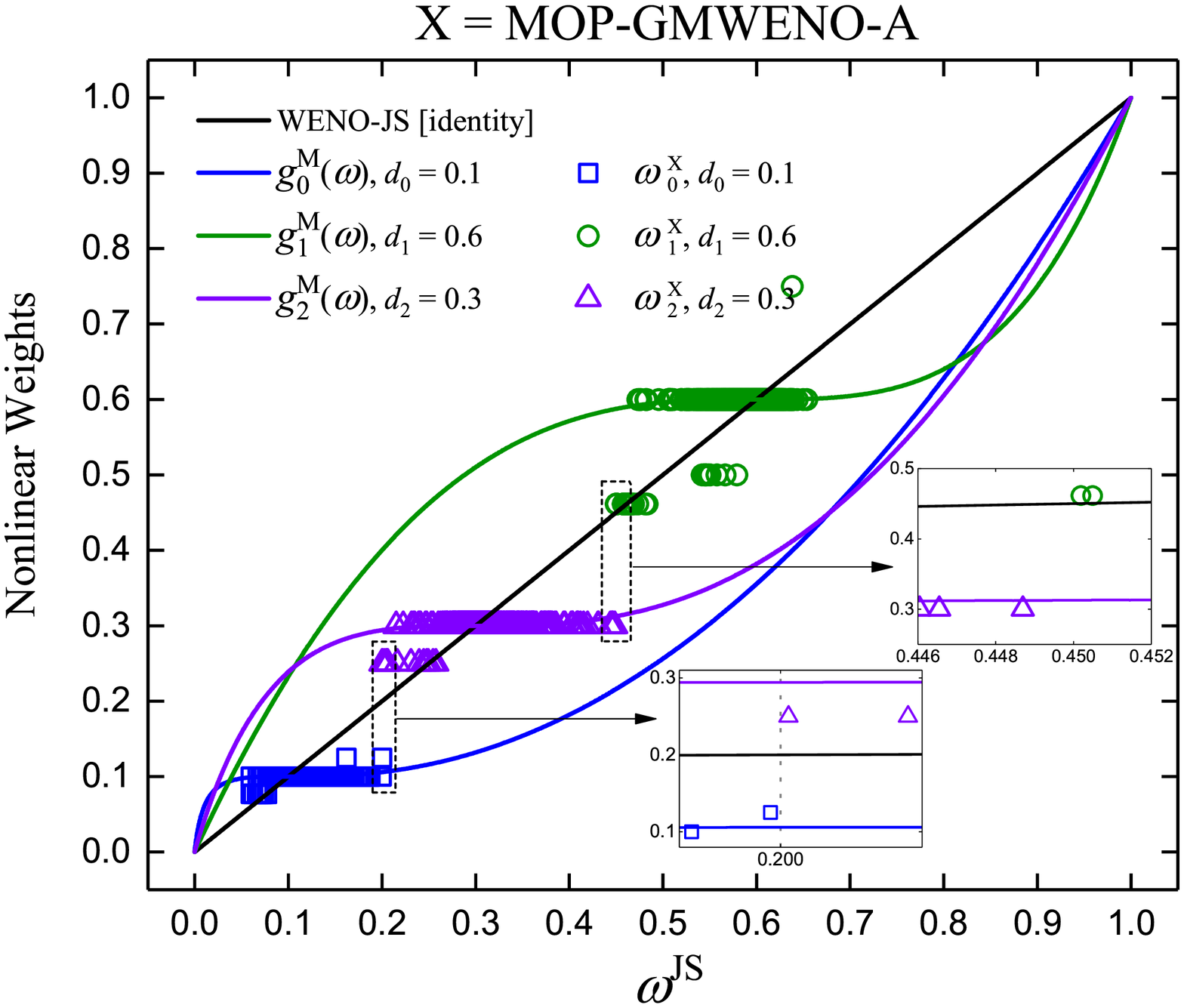}
  \includegraphics[height=0.28\textwidth]
  {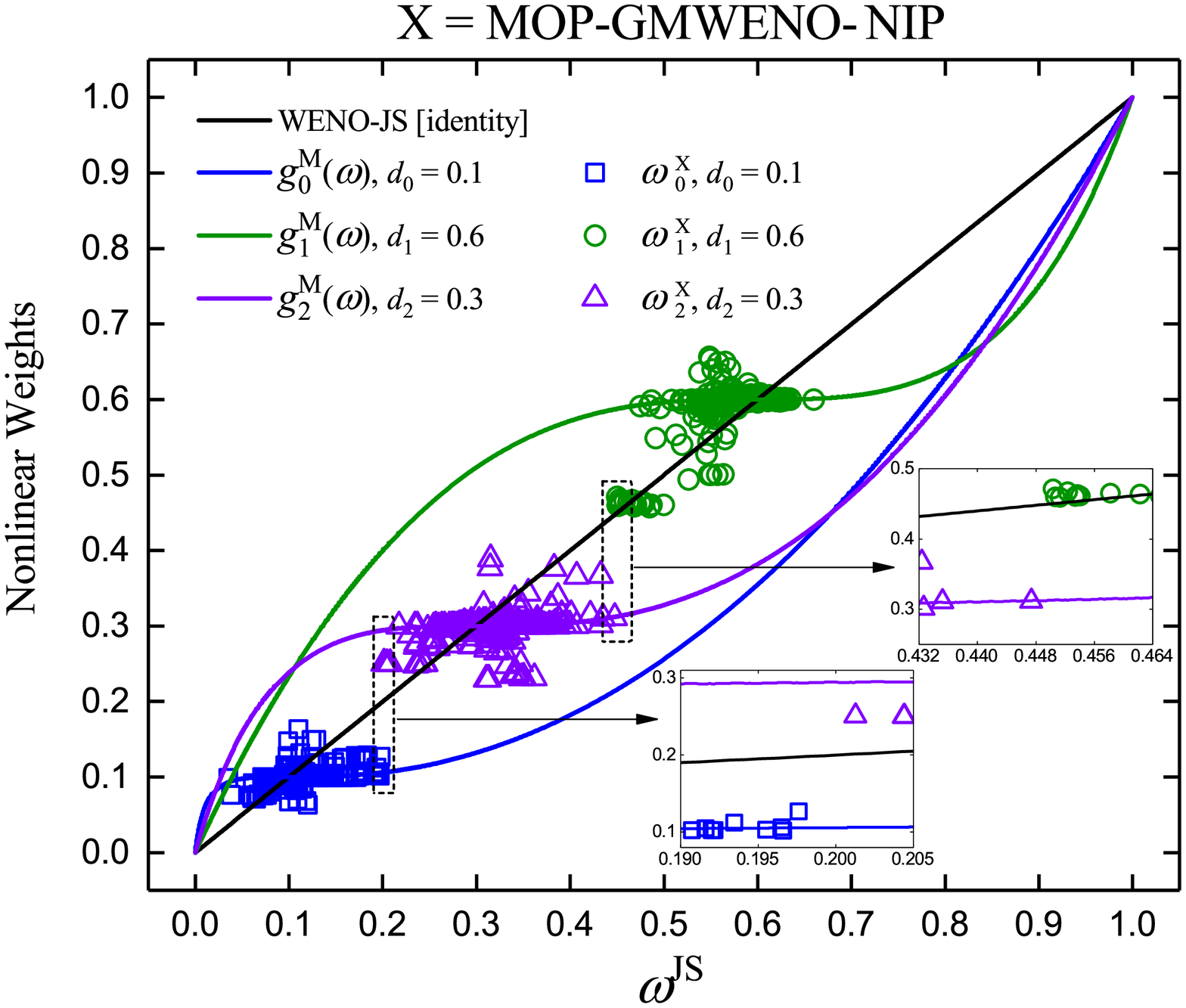}  
\caption{The IMRs for various MOP-GMWENO-X schemes.}
\label{fig:IMR:GMZ}
\end{figure}

\subsection{Convergence at critical points}\label{subsec:Critical}
Since Henrick et al. \cite{WENO-M} pointed out that the WENO-JS 
scheme achieves only third-order convergence rate of accuracy at 
critical points of smooth solutions, it has become a focus of 
discussion \cite{WENO-IM,WENO-PM,WENO-RM260,WENO-PPM5,WENO-MAIMi,
WENO-ACM,MOP-WENO-ACMk,MOP-WENO-X,WENO-Z}. 

We compute the function $f(x) = x^{3} + \cos(x)$ to test the 
convergence property of the MOP-GMWENO-X schemes at critical points. 
It is trivial to know that $f'_{x = 0}= 0$ and $f'''_{x = 0}\neq 0$. 
For comparison purpose, the WENO-X schemes, as well as 
WENO-JS and WENO5 using ideal linear weights (denoted as WENO5-ILW 
in this paper for simplicity), are also conducted. 

Table \ref{tab:Convergence:Critical} shows the $L_{\infty}$ 
convergence behaviors for the considered schemes at the critical 
point $x = 0$. As expected, the WENO-JS scheme only gets third-order 
accuracy. It should be noted that, in this test, the WENO-Z+ scheme 
can also only obtain third-order accuracy (in consistency with the 
results reported in \cite{WENO-Zplus}), leading to the fact that the MOP-GMWENO-Z+ scheme only gets third-order accuracy. 
Moreover, the WENO-JS, WENO-Z+ and MOP-GMWENO-Z+ schemes obtain the 
errors of at least 4 orders of magnitude larger than those of the 
the WENO-ILW scheme. The other WENO-X schemes, say, WENO-Z, 
WENO-Z$\eta(\tau_{5})$, WENO-Z$\eta(\tau_{81})$, WENO-ZA, WENO-D, 
WENO-A, WENO-NIP, and the MOP-GMWENO-X schemes can 
get fifth-order accuracy, but only the WENO-Z$\eta(\tau_{81})$, 
WENO-ZA, WENO-A, WENO-NIP schemes, and the
MOP-GMWENO-X schemes, can obtain the errors with the same order of 
magnitude as those of the WENO-ILW scheme. The errors generated 
by the WENO-Z, WENO-Z$\eta(\tau_{5})$, WENO-D schemes, and the
MOP-GMWENO-X schemes, are about 3 orders of 
magnitude larger than those of the WENO-ILW scheme. It is worthy to 
note that all the new WENO schemes perform similarly as the 
WENO-X schemes for this test.

\begin{table}[ht]
\begin{myFontSize}
\centering
\caption{Convergence rate of accuracy at the critical point.}
\label{tab:Convergence:Critical}
\begin{tabular*}{\hsize}
{@{}@{\extracolsep{\fill}}lllllllll@{}}
\hline
\space    &\multicolumn{2}{l}{\cellcolor{gray!35}{WENO5-ILW}}  
          &\multicolumn{2}{l}{\cellcolor{gray!35}{WENO-JS}}  
          &\multicolumn{2}{l}{\cellcolor{gray!35}{WENO-Z}}  
          &\multicolumn{2}{l}{\cellcolor{gray!35}{MOP-GMWENO-Z}}\\
\cline{2-3}  \cline{4-5}   \cline{6-7}  \cline{8-9}
$\Delta x$            & Error              & Order 
                      & Error              & Order 
		  			  & Error              & Order 
      				  & Error              & Order \\
\Xhline{0.65pt}
0.0125                & 2.72851e-11        & -
                      & 3.57689e-06        & -
                      & 5.36240e-08        & -
                      & 5.36245e-08        & -\\
0.00625               & 1.11260e-12        & 4.6161
                      & 3.64804e-07        & 3.2935
                      & 1.14381e-09        & 5.5510
                      & 1.14381e-09        & 5.5510\\
0.003125              & 3.47692e-14        & 5.0000
                      & 3.81177e-08        & 3.2586
                      & 2.50723e-11        & 5.5116
                      & 2.50723e-11        & 5.5116\\
0.0015625             & 1.08654e-15        & 5.0000
                      & 3.99276e-09        & 3.2550
                      & 5.58613e-13        & 5.4881
                      & 5.58613e-13        & 5.4881\\
\hline
\space    &\multicolumn{2}{l}{\cellcolor{gray!35}{WENO-Z$\eta(\tau_{5})$}}  
          &\multicolumn{2}{l}{\cellcolor{gray!35}{MOP-GMWENO-Z$\eta(\tau_{5})$}}  
          &\multicolumn{2}{l}{\cellcolor{gray!35}{WENO-Z$\eta(\tau_{81})$}}  
          &\multicolumn{2}{l}{\cellcolor{gray!35}{MOP-GMWENO-Z$\eta(\tau_{81})$}}\\
\cline{2-3}  \cline{4-5}   \cline{6-7}  \cline{8-9}
$\Delta x$            & Error      & Order 
                      & Error      & Order 
		  			  & Error      & Order 
      				  & Error      & Order \\
\Xhline{0.65pt}
0.0125                & 5.23932e-08        & -
                      & 5.23937e-08        & -
                      & 3.50808e-11        & -
                      & 3.65801e-11        & -\\
0.00625               & 1.11430e-09        & 5.5552
                      & 1.11430e-09        & 5.5552
                      & 1.11260e-12        & 4.9787
                      & 1.11260e-12        & 5.0391\\
0.003125              & 2.43877e-11        & 5.5138
                      & 2.43877e-11        & 5.5138
                      & 3.47692e-14        & 5.0000
                      & 3.47692e-14        & 5.0000\\
0.0015625             & 5.42433e-13        & 5.4906
                      & 5.42433e-13        & 5.4906
                      & 1.08654e-15        & 5.0000
                      & 1.08654e-15        & 5.0000\\
\hline
\space    &\multicolumn{2}{l}{\cellcolor{gray!35}{WENO-Z+}}  
          &\multicolumn{2}{l}{\cellcolor{gray!35}{MOP-GMWENO-Z+}}  
          &\multicolumn{2}{l}{\cellcolor{gray!35}{WENO-ZA}}  
          &\multicolumn{2}{l}{\cellcolor{gray!35}{MOP-GMWENO-ZA}}\\
\cline{2-3}  \cline{4-5}   \cline{6-7}  \cline{8-9}
$\Delta x$            & Error      & Order 
                      & Error      & Order 
		  			  & Error      & Order 
      				  & Error      & Order \\
\Xhline{0.65pt}
0.0125                & 6.57316e-07        & -
                      & 6.57315e-07        & -
                      & 1.07615e-10        & -
                      & 1.08743e-10        & -\\
0.00625               & 8.77288e-08        & 2.9055
                      & 8.77288e-08        & 2.9055
                      & 1.50193e-12        & 6.1629
                      & 1.50193e-12        & 6.1780\\
0.003125              & 1.08696e-08        & 3.0128
                      & 1.08696e-08        & 3.0128
                      & 3.69247e-14        & 5.3461
                      & 3.69247e-14        & 5.3461\\
0.0015625             & 1.29563e-09        & 3.0686
                      & 1.29563e-09        & 3.0686
                      & 1.09865e-15        & 5.0708
                      & 1.09865e-15        & 5.0708\\
\hline
\space    &\multicolumn{2}{l}{\cellcolor{gray!35}{WENO-D}}  
          &\multicolumn{2}{l}{\cellcolor{gray!35}{MOP-GMWENO-D}}  
          &\multicolumn{2}{l}{\cellcolor{gray!35}{WENO-A}}  
          &\multicolumn{2}{l}{\cellcolor{gray!35}{MOP-GMWENO-A}}\\
\cline{2-3}  \cline{4-5}   \cline{6-7}  \cline{8-9}
$\Delta x$            & Error      & Order 
                      & Error      & Order 
		  			  & Error      & Order 
      				  & Error      & Order \\
\Xhline{0.65pt}
0.0125                & 5.36240e-08        & -
                      & 5.36245e-08        & -
                      & 3.62665e-11        & -
                      & 3.67745e-11        & -\\
0.00625               & 1.14328e-09        & 5.5516
                      & 1.14328e-09        & 5.5516
                      & 1.11260e-12        & 5.0266
                      & 1.11260e-12        & 5.0467\\
0.003125              & 2.50723e-11        & 5.5109
                      & 2.50723e-11        & 5.5109
                      & 3.47692e-14        & 5.0000
                      & 3.47692e-14        & 5.0000\\
0.0015625             & 5.58613e-13        & 5.4881
                      & 5.58613e-13        & 5.4881
                      & 1.08654e-15        & 5.0000
                      & 1.08654e-15        & 5.0000\\
\hline
\space    &\multicolumn{2}{l}{\cellcolor{gray!35}{WENO-NIP}}  
          &\multicolumn{2}{l}{\cellcolor{gray!35}{MOP-GMWENO-NIP}}  
          &\multicolumn{2}{l}{}  
          &\multicolumn{2}{l}{}\\
\cline{2-3}  \cline{4-5}
$\Delta x$            & Error      & Order 
                      & Error      & Order 
		  			  & {}         & {} 
      				  & {}         & {} \\
\cline{1-5}
0.0125                & 3.56050e-11        & -
                      & 3.56511e-11        & -
                      & {}                 & {}
                      & {}                 & {}\\
0.00625               & 1.11264e-12        & 5.0000
                      & 1.11264e-12        & 5.0019
                      & {}                 & {}
                      & {}                 & {}\\
0.003125              & 3.47695e-14        & 5.0000
                      & 3.47695e-14        & 5.0000
                      & {}                 & {}
                      & {}                 & {}\\
0.0015625             & 1.08654e-15        & 5.0000
                      & 1.08654e-15        & 5.0000
                      & {}                 & {}
                      & {}                 & {}\\
\cline{1-5}
\end{tabular*}
\end{myFontSize}
\end{table}

\subsection{Long-run simulations of 1D linear advection equation for 
comparison}\label{subsec:Long-run}
\subsubsection{With high-order critical points}
\label{subsubsec:LongCritical}
For the purpose of demonstrating the improvement of the new WENO 
schemes that they can preserve high resolutions for the problem with 
high-order critical points at long output times, we perform the 
following test.

\begin{example}
\rm{We compute}
\label{ex:Long-run:HighOrderCritical}
\end{example}
\begin{equation*}
u_{t} + u_{x} = 0, \quad x \in (7.5, 10.5),
\end{equation*} 
with the following initial condition
\begin{equation}
u(x, 0) = \exp\Big(-\big(x - 9.0\big)^{10}\Big)\cos^{9}\Big(\pi\big(x - 9.0\big)\Big) .
\label{eq:Long-run:IC:HighOrderCritical}
\end{equation}
Here, we set CFL = $(\Delta x)^{2/3}$. 

To compare the dissipations, the following 
$L_{1}$ and $L_{\infty}$ errors are computed
\begin{equation}
\displaystyle
\begin{aligned}
&L_{1} = h \cdot \displaystyle\sum\nolimits_{i=1}^{N}\Big\lvert 
u_{i}^{\mathrm{exact}} - (u_{h})_{i} \Big\rvert, \quad \\
&L_{\infty} = \displaystyle\max_{1\leq i\leq N} \Big\lvert 
u_{i}^{\mathrm{exact}} -(u_{h})_{i}
\Big\rvert,
\end{aligned}
\label{normsDefinition}
\end{equation}
where $N$ is the number of the cells, $h$ the mesh size,
$(u_{h})_{i}$ the computing result and $u_{i}^{\mathrm{exact}}$ the 
exact solution that can be easily computed by $u(x,t) = \exp\Big(-\big((x - t) - 9.0\big)^{5}\Big)\cos^{9}\Big(\pi\big((x - t) - 9.0\big)\Big) $. 

For the scheme ``Y'', its increased errors compared to WENO5-ILW are 
also been computed by
\begin{equation*}
\begin{aligned}
\chi_{1}=\frac{L_{1}^{\mathrm{Y}}(t)-L_{1}^{\mathrm{ILW}}(t)}{L_{1}^{\mathrm{ILW}}(t)}\times100\%, \quad 
\chi_{\infty}=\frac{L_{\infty}^{\mathrm{Y}}(t)-L_{\infty}^{\mathrm{ILW}}(t)}{L_{\infty}^{\mathrm{ILW}}(t)}\times100\%,
\end{aligned}
\end{equation*}
where $L_{1}^{\mathrm{ILW}}(t)$ and $L_{\infty}^{\mathrm{ILW}}(t)$ 
stand for the $L_{1}$ and $L_{\infty}$ errors of WENO5-ILW 
respectively. 

In Table \ref{tab::long-run:N300}, we show the results computed by 
considered schemes with $N=300$ and $t = 300, 600, 900, 1200$. It 
can be seen that: (1) WENO-JS produces the largest $L_{1}$ and $L_{\infty}$
errors, leading to the largest increased errors, among all 
considered schemes for all output times; (2) the WENO-Z, 
WENO-Z$\eta(\tau_{5})$, WENO-Z+, WENO-D and WENO-A 
schemes also generate very large $L_{1}$ and $L_{\infty}$ errors, 
although slightly smaller than the WENO-JS scheme, 
and their associated increased errors are extremely large naturally; 
(3) however, the MOP-GMWENO-Z, 
MOP-GMWENO-Z$\eta(\tau_{5})$, MOP-GMWENO-Z+, MOP-GMWENO-D and 
MOP-GMWENO-A schemes can significantly decrease the increased errors 
to a tolerable level, as they can achieve much smaller $L_{1}$ and 
$L_{\infty}$ errors that actually get very close to that of 
WENO-ILW; (4) moreover, the 
MOP-GMWENO-Z$\eta(\tau_{81})$, MOP-GMWENO-ZA and MOP-GMWENO-NIP 
schemes can also ensure that their $L_{1}$ and $L_{\infty}$
errors get close to that of WENO-ILW, 
and thus their increased errors are also at a tolerable level.
It seems that the WENO-Z$\eta(\tau_{81})$, WENO-ZA and WENO-NIP 
schemes can get solutions almost as accurate as, or even more 
accurate than, that of WENO-ILW, and this is good for 
this test. However, it should be pointed out that these schemes may 
suffer from lack of robustness as their dissipations are too small for long-run calculations. Indeed, we will demonstrate this in 
detail by numerical examples below.

\begin{figure}[!ht]
\centering
  \includegraphics[height=0.33\textwidth]
  {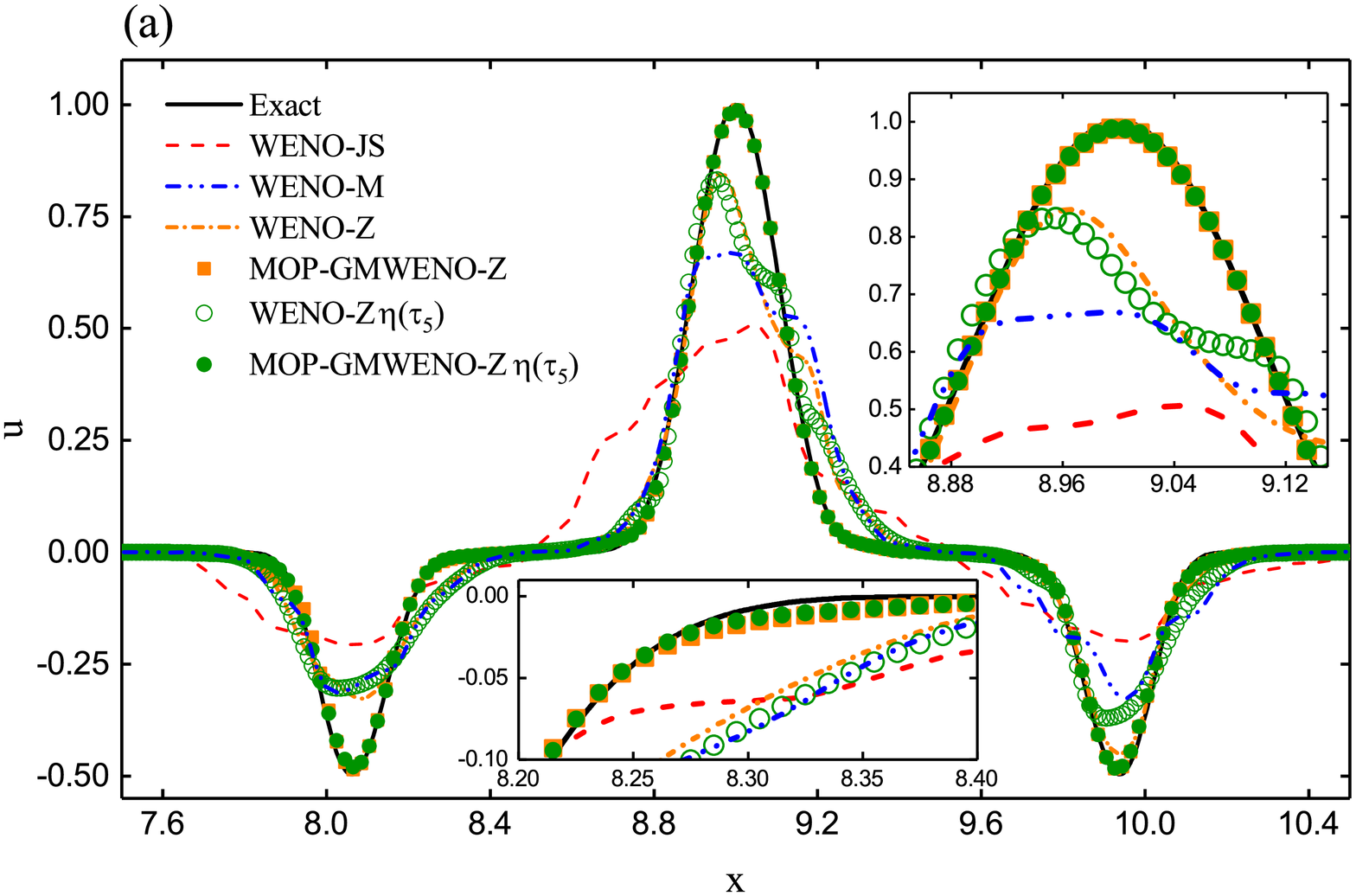}
  \includegraphics[height=0.33\textwidth]
  {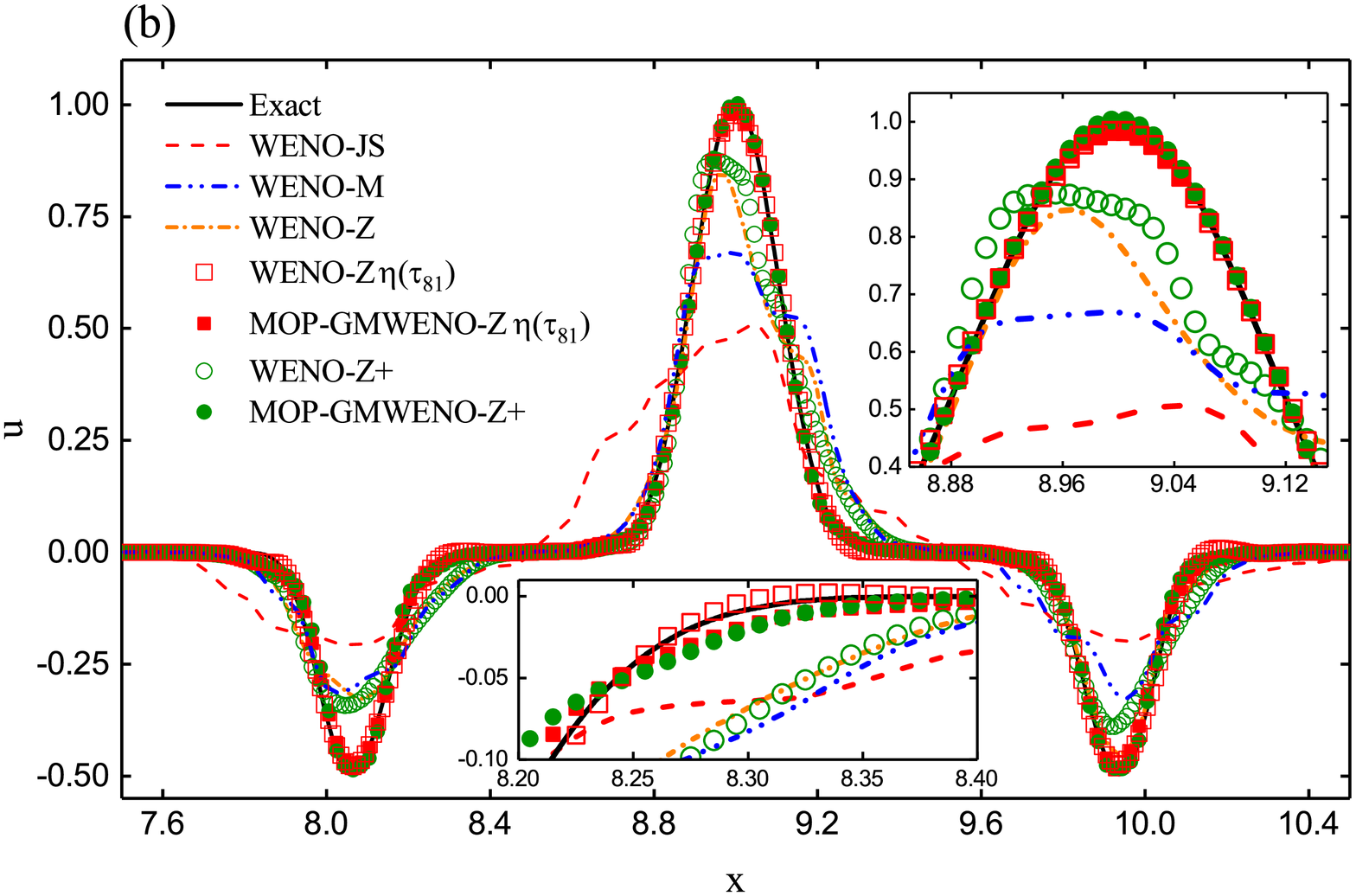}\\
  \includegraphics[height=0.33\textwidth]
  {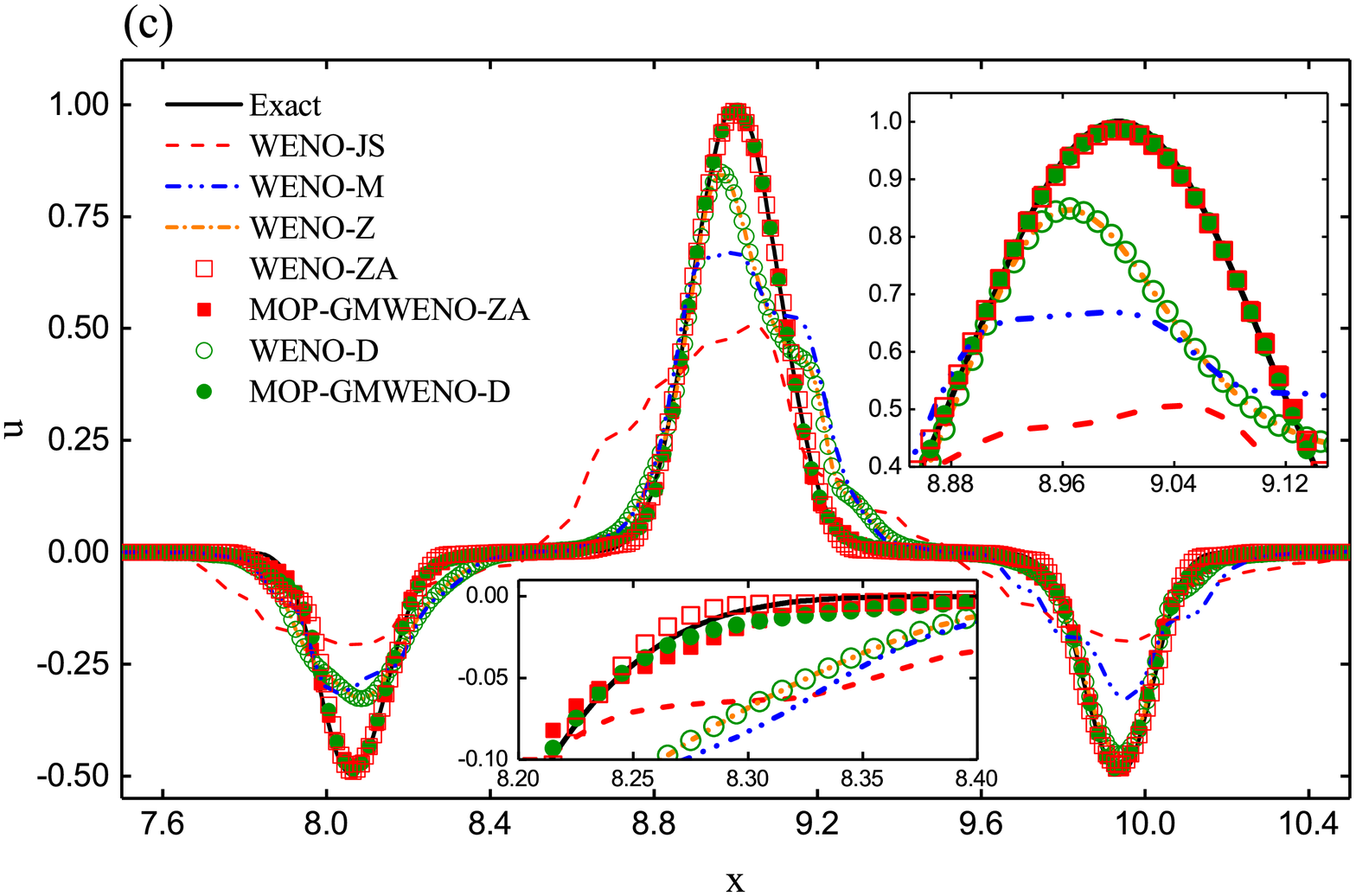}  
  \includegraphics[height=0.33\textwidth]
  {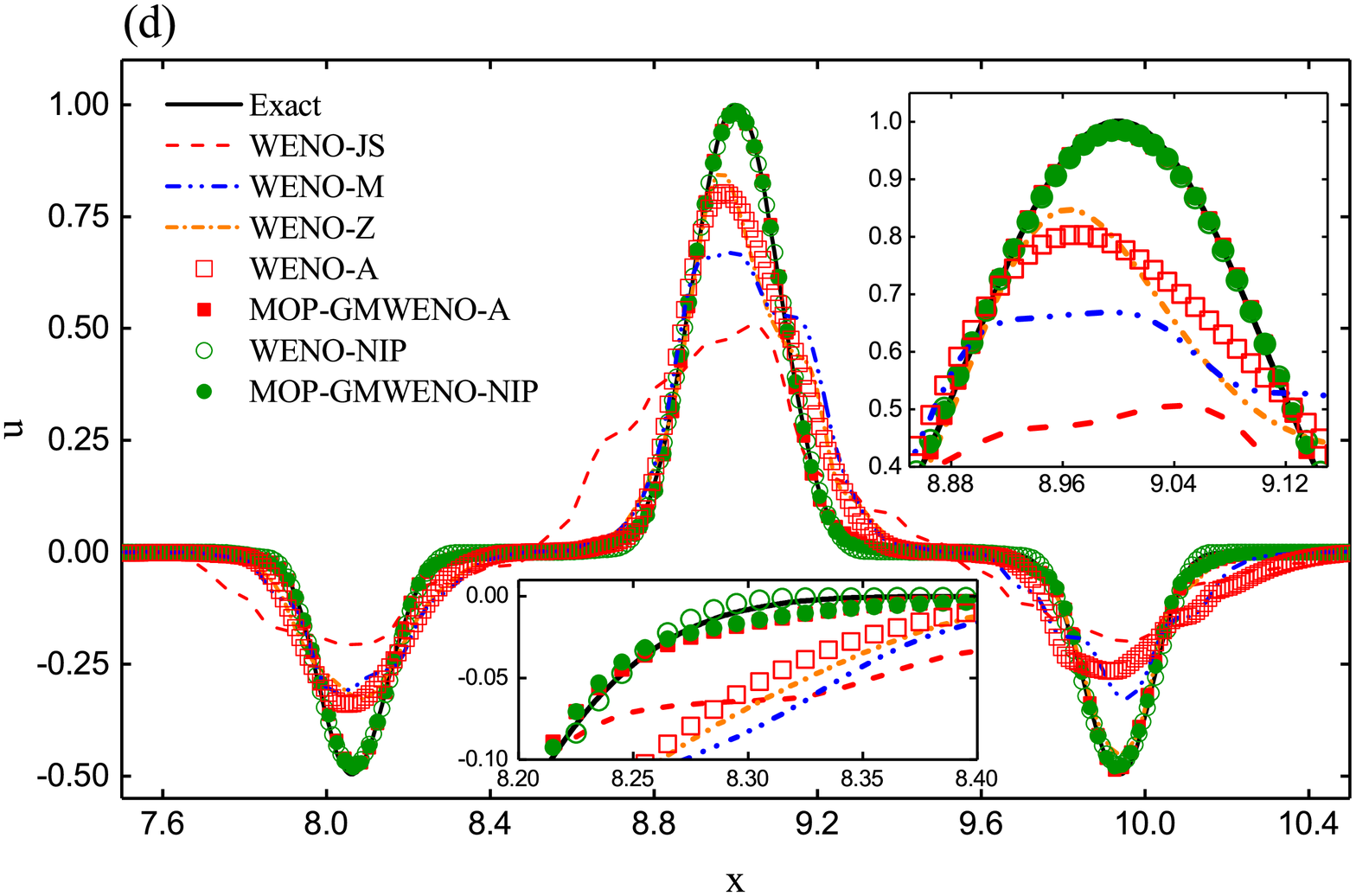}  
  \caption{The solutions on solving Example \ref{ex:Long-run:HighOrderCritical}.}
\label{fig:ex:long-run:N300}
\end{figure}

Fig. \ref{fig:ex:long-run:N300} shows the solutions at $t = 1200$.
For comparison purpose, we also plot the results computed by the 
WENO-JS \cite{WENO-JS} and WENO-M \cite{WENO-M} schemes. From Fig. 
\ref{fig:ex:long-run:N300}, it can be observed that: (1) WENO-JS
shows the lowest resolution, followed by WENO-M; (3) the 
WENO-Z, WENO-Z$\eta(\tau_{5})$, WENO-Z+, WENO-D and WENO-A schemes 
also show very low resolutions but the improved 
MOP-GMWENO-X schemes can significantly improve the resolutions; (4) 
the resolutions of MOP-GMWENO-Z$\eta(\tau_{81})$, MOP-GMWENO-ZA and 
MOP-GMWENO-NIP are slightly smaller than that of the
WENO-Z$\eta(\tau_{81})$, WENO-ZA and 
WENO-NIP schemes but they are still far better than the WENO-JS, 
WENO-M, and the other WENO-Z-type schemes.

\begin{table}[!ht]
\begin{myFontSize}
\centering
\caption{$L_{1}$, $L_{\infty}$ errors and the increased errors (in percentage) on solving Example \ref{ex:Long-run:HighOrderCritical}.}
\label{tab::long-run:N300}
\begin{tabular*}{\hsize}
{@{}@{\extracolsep{\fill}}rlrlrlrlr@{}}
\hline
\space    &\multicolumn{4}{l}{\cellcolor{gray!35}{WENO5-ILW}}  
          &\multicolumn{4}{l}{\cellcolor{gray!35}{WENO-JS}}  \\
\cline{2-5}  \cline{6-9}
Time, $t$             & $L_{1}$ error      & $\chi_{1}$
                      & $L_{\infty}$ error & $\chi_{\infty}$
		  			  & $L_{1}$ error      & $\chi_{1}$
      				  & $L_{\infty}$ error & $\chi_{\infty}$ \\
\Xhline{0.65pt}
300                   & 5.39974E-03        & -
                      & 8.81363E-03        & -
                      & 7.93589E-02        & 1370\%
                      & 1.34321E-01        & 1424\% \\
600                   & 9.94133E-03        & -
                      & 1.50917E-02        & -
                      & 2.10016E-01        & 2013\%
                      & 3.04860E-01        & 1920\% \\
900                   & 1.38061E-02        & -
                      & 1.96281E-02        & -
                      & 2.84632E-01        & 1962\%
                      & 4.20080E-01        & 2040\% \\
1200                  & 1.74067E-02        & -
                      & 2.39652E-02        & -
                      & 3.26687E-01        & 1777\%
                      & 5.14072E-01        & 2045\% \\
\hline
\space    &\multicolumn{4}{l}{\cellcolor{gray!35}{WENO-Z}}
          &\multicolumn{4}{l}{\cellcolor{gray!35}{MOP-GMWENO-Z}}\\
\cline{2-5}  \cline{6-9}
Time, $t$             & $L_{1}$ error      & $\chi_{1}$         
                      & $L_{\infty}$ error & $\chi_{\infty}$     
		  			  & $L_{1}$ error      & $\chi_{1}$         
      				  & $L_{\infty}$ error & $\chi_{\infty}$      \\
\Xhline{0.65pt}
300                   & 3.56420e-02        & 560\%
                      & 1.23424e-01        & 1300\%
                      & 1.18556e-02        & 120\%
                      & 2.21312e-02        & 151\% \\
600                   & 8.58998e-02        & 777\%
                      & 2.36754e-01        & 1530\%
                      & 1.56335e-02        & 60\%
                      & 2.43793e-02        & 68\% \\
900                   & 1.13762e-01        & 726\%
                      & 2.77206e-01        & 1324\%
                      & 2.57924e-02        & 87\%
                      & 3.33097e-02        & 71\% \\
1200                  & 1.27002e-01        & 631\%
                      & 2.78133e-01        & 1070\%
                      & 2.51627e-02        & 45\%
                      & 4.06869e-02        & 71\% \\
\hline
\space    &\multicolumn{4}{l}{\cellcolor{gray!35}{WENO-Z$\eta(\tau_{5})$}}
          &\multicolumn{4}{l}{\cellcolor{gray!35}{MOP-GMWENO-Z$\eta(\tau_{5})$}}\\
\cline{2-5}  \cline{6-9}
Time, $t$             & $L_{1}$ error      & $\chi_{1}$         
                      & $L_{\infty}$ error & $\chi_{\infty}$     
		  			  & $L_{1}$ error      & $\chi_{1}$         
      				  & $L_{\infty}$ error & $\chi_{\infty}$      \\
\Xhline{0.65pt}
300                   & 4.52391e-02        & 738\%
                      & 1.14337e-01        & 1197\%
                      & 1.16734e-02        & 116\%
                      & 2.42074e-02        & 175\% \\
600                   & 9.44265e-02        & 864\%
                      & 2.36673e-01        & 1530\%
                      & 1.56751e-02        & 60\%
                      & 2.51991e-02        & 74\% \\
900                   & 1.09605e-01        & 695\%
                      & 2.12852e-01        & 993\%
                      & 2.52547e-02        & 83\%
                      & 3.20434e-02        & 65\% \\
1200                  & 1.38434e-01        & 697\%
                      & 3.05553e-01        & 1185\%
                      & 2.15191e-02        & 24\%
                      & 3.08777e-02        & 30\% \\
\hline
\space    &\multicolumn{4}{l}{\cellcolor{gray!35}{WENO-Z$\eta(\tau_{81})$}}
          &\multicolumn{4}{l}{\cellcolor{gray!35}{MOP-GMWENO-Z$\eta(\tau_{81})$}}\\
\cline{2-5}  \cline{6-9}
Time, $t$             & $L_{1}$ error      & $\chi_{1}$         
                      & $L_{\infty}$ error & $\chi_{\infty}$     
		  			  & $L_{1}$ error      & $\chi_{1}$         
      				  & $L_{\infty}$ error & $\chi_{\infty}$      \\
\Xhline{0.65pt}
300                   & 5.57811e-03        & 3\%
                      & 9.60994e-03        & 9\%
                      & 9.38301e-03        & 74\%
                      & 1.47213e-02        & 67\% \\
600                   & 1.04455e-02        & 7\%
                      & 1.50933e-02        & 4\%
                      & 1.70600e-02        & 74\%
                      & 2.32755e-02        & 60\% \\
900                   & 1.47709e-02        & 7\%
                      & 2.00720e-02        & 3\%
                      & 2.09101e-02        & 52\%
                      & 3.42640e-02        & 76\% \\
1200                  & 1.83129e-02        & 5\%
                      & 2.40373e-02        & 1\%
                      & 2.66533e-02        & 54\%
                      & 3.22354e-02        & 36\% \\
\hline
\space    &\multicolumn{4}{l}{\cellcolor{gray!35}{WENO-Z+}}
          &\multicolumn{4}{l}{\cellcolor{gray!35}{MOP-GMWENO-Z+}}\\
\cline{2-5}  \cline{6-9}
Time, $t$             & $L_{1}$ error      & $\chi_{1}$         
                      & $L_{\infty}$ error & $\chi_{\infty}$     
		  			  & $L_{1}$ error      & $\chi_{1}$         
      				  & $L_{\infty}$ error & $\chi_{\infty}$      \\
\Xhline{0.65pt}
300                   & 5.09235e-02        & 843\%
                      & 1.59585e-01        & 1710\%
                      & 9.45604e-03        & 75\%
                      & 2.50726e-02        & 184\% \\
600                   & 8.20317e-02        & 737\%
                      & 1.84579e-01        & 1171\%
                      & 1.50895e-02        & 54\%
                      & 3.40525e-02        & 135\% \\
900                   & 1.03748e-01        & 653\%
                      & 1.76101e-01        & 804\%
                      & 2.18063e-02        & 58\%
                      & 6.96188e-02        & 258\% \\
1200                  & 1.15534e-01        & 565\%
                      & 2.22997e-01        & 838\%
                      & 3.00894e-02        & 73\%
                      & 5.66315e-02        & 138\% \\
\hline
\space    &\multicolumn{4}{l}{\cellcolor{gray!35}{WENO-ZA}}
          &\multicolumn{4}{l}{\cellcolor{gray!35}{MOP-GMWENO-ZA}}\\
\cline{2-5}  \cline{6-9}
Time, $t$             & $L_{1}$ error      & $\chi_{1}$         
                      & $L_{\infty}$ error & $\chi_{\infty}$     
		  			  & $L_{1}$ error      & $\chi_{1}$         
      				  & $L_{\infty}$ error & $\chi_{\infty}$      \\
\Xhline{0.65pt}
300                   & 5.10578e-03        & -5\%
                      & 8.86265e-03        & 1\%
                      & 8.87687e-03        & 64\%
                      & 1.54090e-02        & 75\% \\
600                   & 9.05646e-03        & -8\%
                      & 1.48176e-02        & 2\%
                      & 1.82371e-02        & 86\%
                      & 2.52104e-02        & 74\% \\
900                   & 1.41249e-02        & 3\%
                      & 2.36763e-02        & 22\%
                      & 1.99799e-02        & 45\%
                      & 3.82169e-02        & 96\% \\
1200                  & 2.32315e-02        & 34\%
                      & 6.80578e-02        & 186\%
                      & 2.81372e-02        & 62\%
                      & 3.76059e-02        & 58\% \\
\hline
\space    &\multicolumn{4}{l}{\cellcolor{gray!35}{WENO-D}}
          &\multicolumn{4}{l}{\cellcolor{gray!35}{MOP-GMWENO-D}}\\
\cline{2-5}  \cline{6-9}
Time, $t$             & $L_{1}$ error      & $\chi_{1}$         
                      & $L_{\infty}$ error & $\chi_{\infty}$     
		  			  & $L_{1}$ error      & $\chi_{1}$         
      				  & $L_{\infty}$ error & $\chi_{\infty}$      \\
\Xhline{0.65pt}
300                   & 3.56420e-02        & 560\%
                      & 1.23424e-01        & 1300\%
                      & 1.18556e-02        & 120\%
                      & 2.21312e-02        & 151\% \\
600                   & 8.58998e-02        & 777\%
                      & 2.36754e-01        & 1530\%
                      & 1.56335e-02        & 60\%
                      & 2.43793e-02        & 68\% \\
900                   & 1.13762e-01        & 726\%
                      & 2.77206e-01        & 1324\%
                      & 2.57924e-02        & 87\%
                      & 3.33097e-02        & 71\% \\
1200                  & 1.27002e-01        & 631\%
                      & 2.78133e-01        & 1070\%
                      & 2.51627e-02        & 45\%
                      & 4.06869e-02        & 71\% \\
\hline
\space    &\multicolumn{4}{l}{\cellcolor{gray!35}{WENO-A}}
          &\multicolumn{4}{l}{\cellcolor{gray!35}{MOP-GMWENO-A}}\\
\cline{2-5}  \cline{6-9}
Time, $t$             & $L_{1}$ error      & $\chi_{1}$         
                      & $L_{\infty}$ error & $\chi_{\infty}$     
		  			  & $L_{1}$ error      & $\chi_{1}$         
      				  & $L_{\infty}$ error & $\chi_{\infty}$      \\
\Xhline{0.65pt}
300                   & 1.18154e-01        & 2088\%
                      & 2.14177e-01        & 2330\%
                      & 9.18713e-03        & 70\%
                      & 1.47167e-02        & 67\% \\
600                   & 1.33707e-01        & 1265\%
                      & 2.38106e-01        & 1540\%
                      & 1.68488e-02        & 72\%
                      & 2.48995e-02        & 71\% \\
900                   & 1.50686e-01        & 993\%
                      & 2.33834e-01        & 1101\%
                      & 2.07575e-02        & 51\%
                      & 3.50088e-02        & 80\% \\
1200                  & 1.56130e-01        & 799\%
                      & 2.35669e-01        & 891\%
                      & 2.41207e-02        & 39\%
                      & 3.23873e-02        & 36\% \\
\hline
\space    &\multicolumn{4}{l}{\cellcolor{gray!35}{WENO-NIP}}
          &\multicolumn{4}{l}{\cellcolor{gray!35}{MOP-GMWENO-NIP}}\\
\cline{2-5}  \cline{6-9}
Time, $t$             & $L_{1}$ error      & $\chi_{1}$         
                      & $L_{\infty}$ error & $\chi_{\infty}$     
		  			  & $L_{1}$ error      & $\chi_{1}$         
      				  & $L_{\infty}$ error & $\chi_{\infty}$      \\
\Xhline{0.65pt}
300                   & 5.06954e-03        & -6\%
                      & 8.83927e-03        & 0\%
                      & 9.18661e-03        & 70\%
                      & 1.14320e-02        & 30\% \\
600                   & 8.65059e-03        & -12\%
                      & 1.48766e-02        & 2\%
                      & 1.51265e-02        & 54\%
                      & 2.28329e-02        & 57\% \\
900                   & 1.13469e-02        & -18\%
                      & 1.91798e-02        & -2\%
                      & 1.72611e-02        & 25\%
                      & 2.19778e-02        & 13\% \\
1200                  & 1.36442e-02        & -21\%
                      & 2.29721e-02        & -3\%
                      & 2.13328e-02        & 23\%
                      & 2.68344e-02        & 13\% \\
\hline
\end{tabular*}
\end{myFontSize}
\end{table}

\subsubsection{With discontinuities}
To examine the major benefit of the new schemes that they are able to achieve high resolutions and meanwhile avoid spurious oscillations for long output time simulations, we solve 
\begin{equation*}\left\{
\begin{array}{l}
\dfrac{\partial u}{\partial t} + \dfrac{\partial u}{\partial x} = 0,  \quad -1 < x < 1,\\
u(x, 0) = u_{0}(x).
\end{array}
\right.
\end{equation*}
with two different $u_{0}(x)$.

\begin{example}\rm{
(Case 1) The initial condition $u_{0}(x)$ is computed by}
\label{ex:AccuracyTest:Z}
\end{example}
\begin{equation}
\begin{array}{l}
u_{0}(x) = \left\{
\begin{array}{ll}
1,   & -1 < x \leq 0, \\
0,   & 0 < x < 1.
\end{array}\right. 
\end{array}
\label{eq:LAE:Z}
\end{equation}
This problem simply consists of two constant states separated by 
sharp discontinuities at $x = 0.0, \pm 1.0$.

\begin{example}\rm{
(Case 2) The SLP defined by Eq. \eqref{eq:LAE:SLP}.}
\label{ex:AccuracyTest:SLP}
\end{example}

In order to test the convergence properties, we compute both Case 1 
and Case 2 to the final time $t = 2000$ with the CFL = 0.1. For the 
purpose of comparison, we also solve Case 1 and Case 2 by WENO-JS 
and WENO5-ILW.

In Tables \ref{table:AccuracyTest:Z:t2000}, 
\ref{table:AccuracyTest:SLP:t2000}, the $L_{1}$, $L_{\infty}$ errors 
and the convergence orders are presented. It can be observed that: 
(1) because of its highest 
dissipation, the WENO-JS scheme generates significantly larger 
numerical errors than all other schemes, leading to the smallest 
$L_{1}$ convergence orders; (2) the $L_{1}$ convergence orders of 
the MOP-GMWENO-X schemes are 
distinctly higher than the WENO-X schemes, and 
the $L_{1}$ errors produced by the MOP-GMWENO-Z, 
MOP-GMWENO-Z$\eta(\tau_{5})$, MOP-GMWENO-Z+, MOP-GMWENO-ZA, 
MOP-GMWENO-D and MOP-GMWENO-A schemes are slightly smaller than 
the WENO-X schemes in general; (3) the $L_{1}$ errors 
produced by the MOP-GMWENO-NIP scheme on solving 
Case 1 are also slightly smaller than the
WENO-NIP scheme, while this holds true only for the computing case 
of $N = 800$ on solving Case 2 as the MOP-GMWENO-NIP scheme 
generates slightly and evidently larger $L_{1}$
errors for the computing cases of $N = 400$ and $N = 200$ 
respectively; (4) the $L_{1}$ errors produced by the 
MOP-GMWENO-Z$\eta(\tau_{81})$ scheme are slightly larger on solving 
both Case 1 and Case 2; (5) the $L_{\infty}$ errors of the 
MOP-GMWENO-X schemes on solving Case 1 are smaller than
the WENO-X schemes; (6) for Case 2, the 
$L_{\infty}$ errors of the MOP-GMWENO-X schemes are very close 
to, or even smaller for many cases than, the WENO-X 
schemes. In addition, by taking a view of the $x - u$ profiles, we 
find that the resolutions of the WENO-Z, 
WENO-Z$\eta(\tau_{5})$, WENO-Z+, WENO-ZA, WENO-D and WENO-A schemes 
are significantly lower than the MOP-GMWENO-Z, 
MOP-GMWENO-Z$\eta(\tau_{5})$, MOP-GMWENO-Z+, MOP-GMWENO-ZA, 
MOP-GMWENO-D and MOP-GMWENO-A schemes, and the 
WENO-Z$\eta(\tau_{81})$ and WENO-NIP schemes produce spurious 
oscillations but the MOP-GMWENO-Z$\eta(\tau_{81})$ and 
MOP-GMWENO-NIP schemes, as well as all other MOP-GMWENO-X schemes, 
can remove these oscillations properly. To demonstrate this, 
we next perform detailed tests and show the solutions carefully.

We re-run both Case 1 and Case 2 by considered WENO schemes. Just 
for the purpose of providing better illustrations but without loss 
of generality, we set the output time to be $t = 200$ and use the 
uniform meshes of $N = 1600$ and $N = 3200$ this time. 

\begin{figure}[ht]
\centering
\includegraphics[height=0.33\textwidth]
{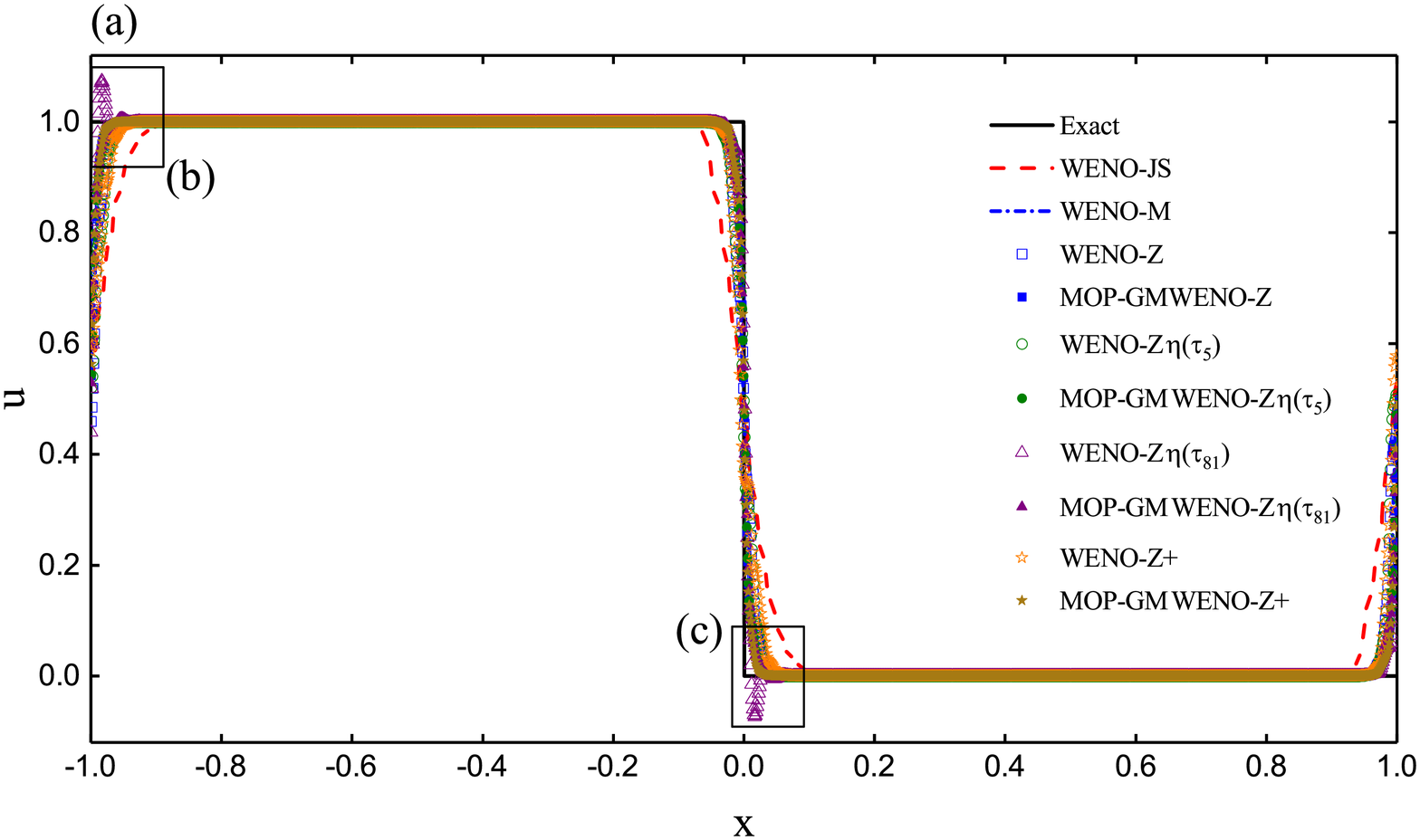}
\includegraphics[height=0.33\textwidth]
{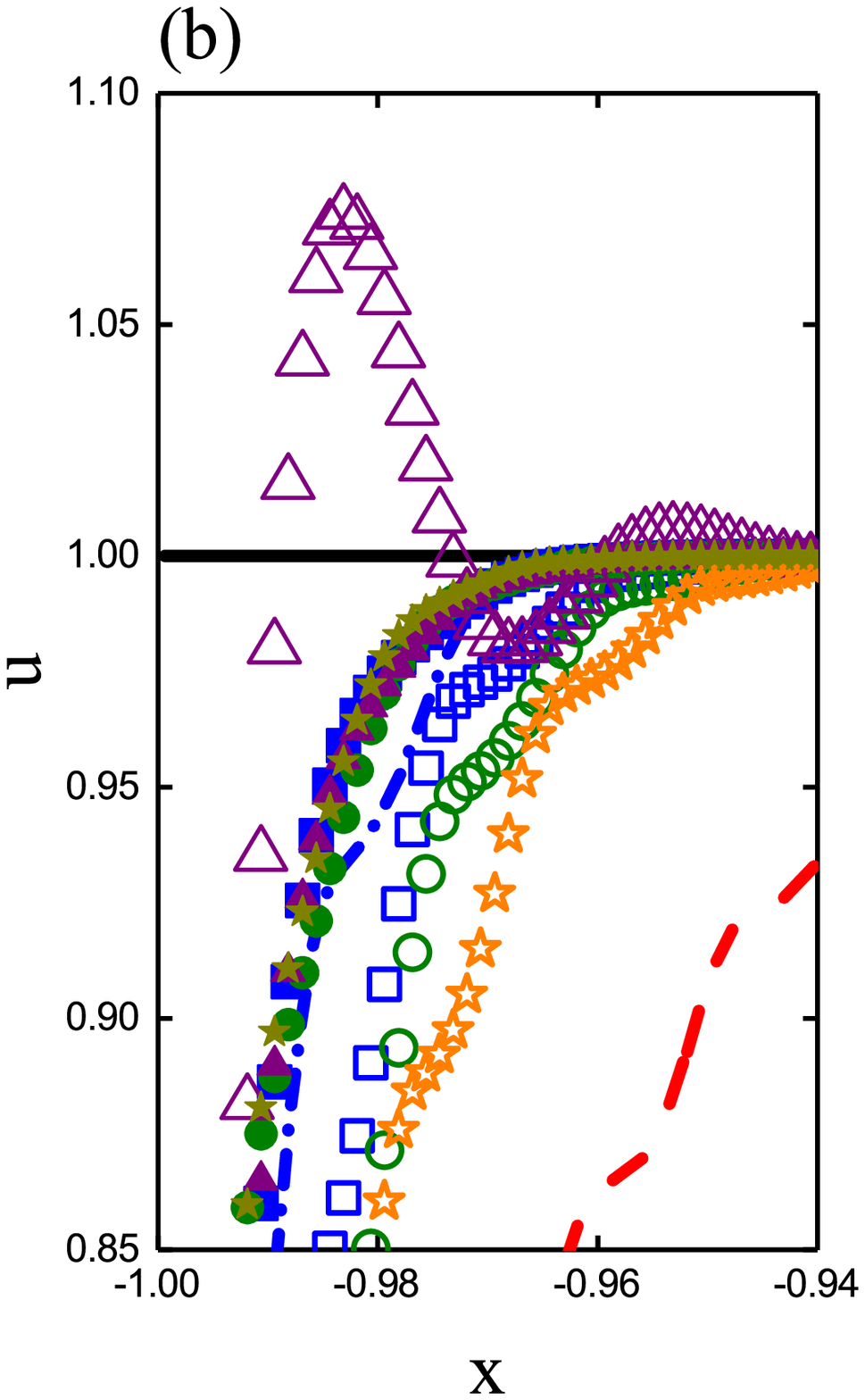}
\includegraphics[height=0.33\textwidth]
{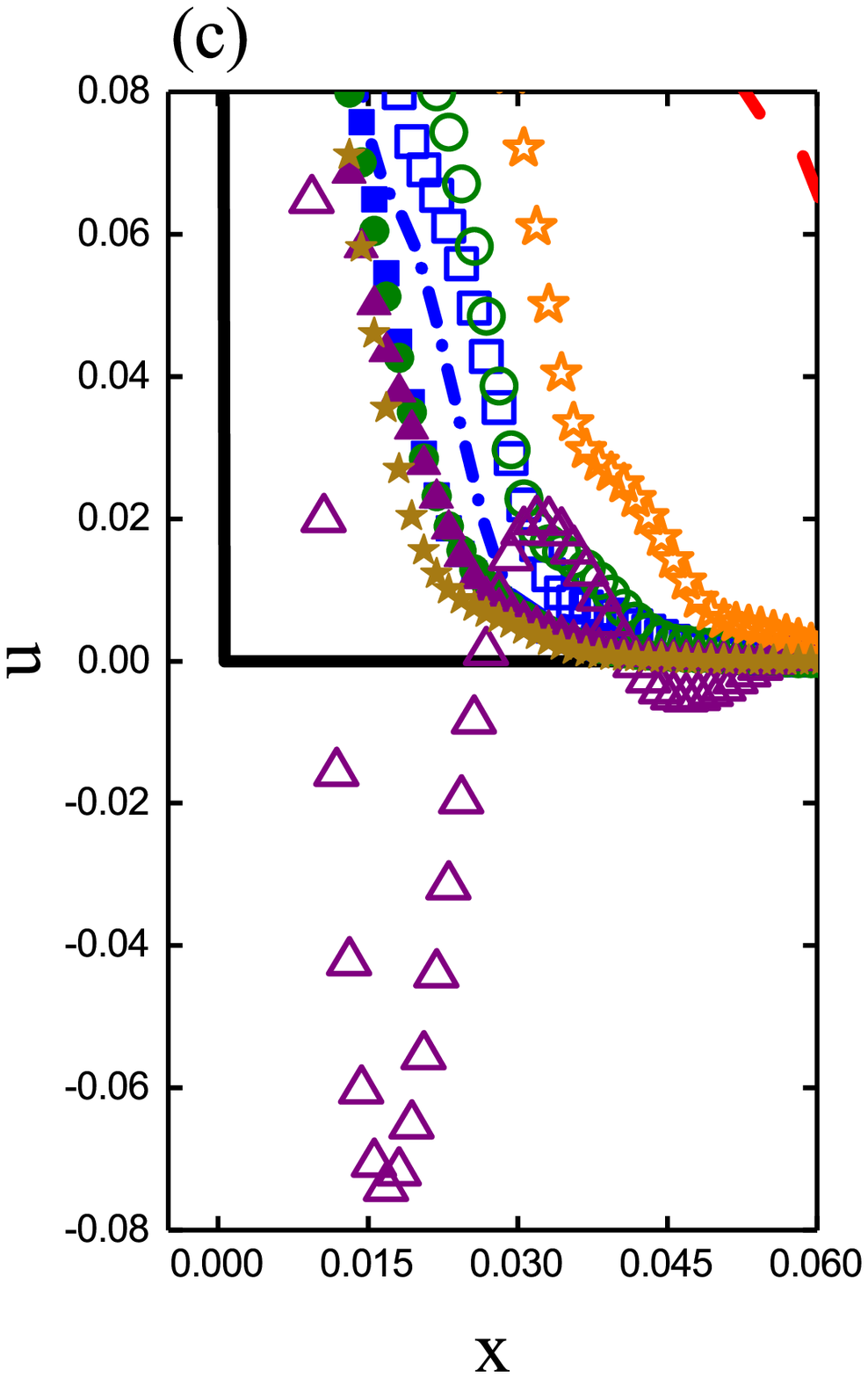}
\caption{Results of different schemes on solving Example 
\ref{ex:AccuracyTest:Z} (Case 1) with $t = 200, N = 1600$.}
\label{fig:Z:N1600:01}
\end{figure}

\begin{figure}[ht]
\centering
\includegraphics[height=0.33\textwidth]
{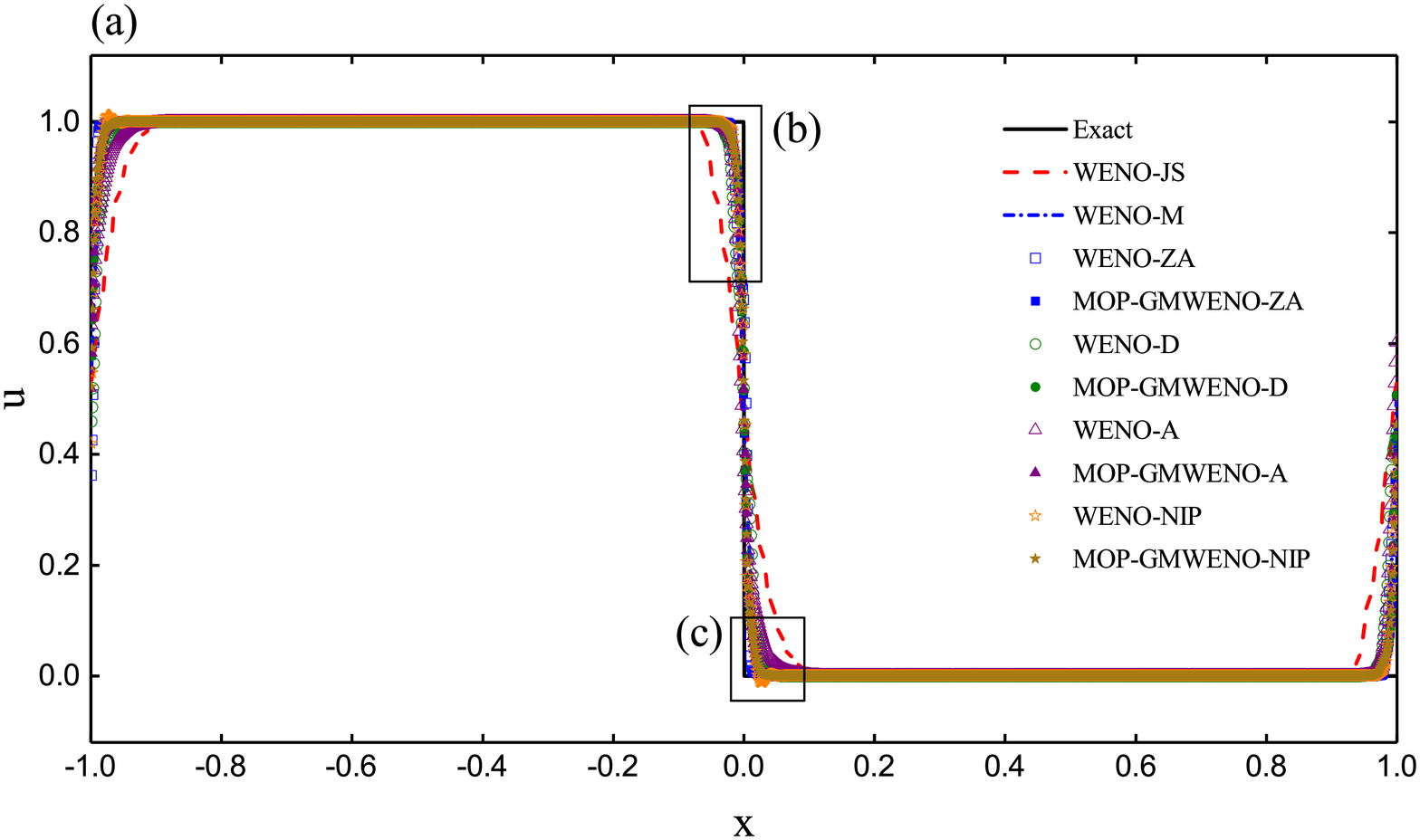}
\includegraphics[height=0.33\textwidth]
{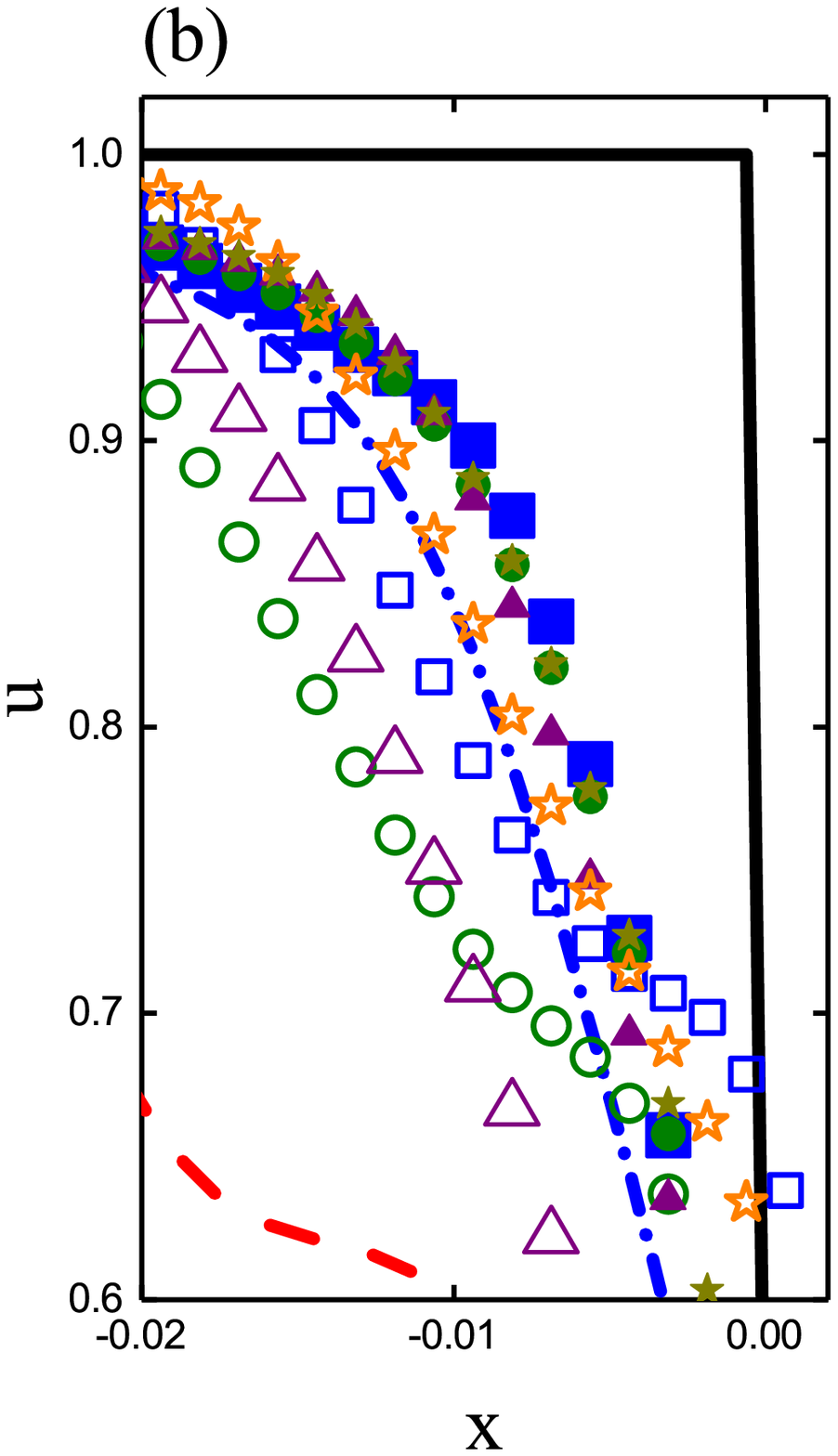}
\includegraphics[height=0.33\textwidth]
{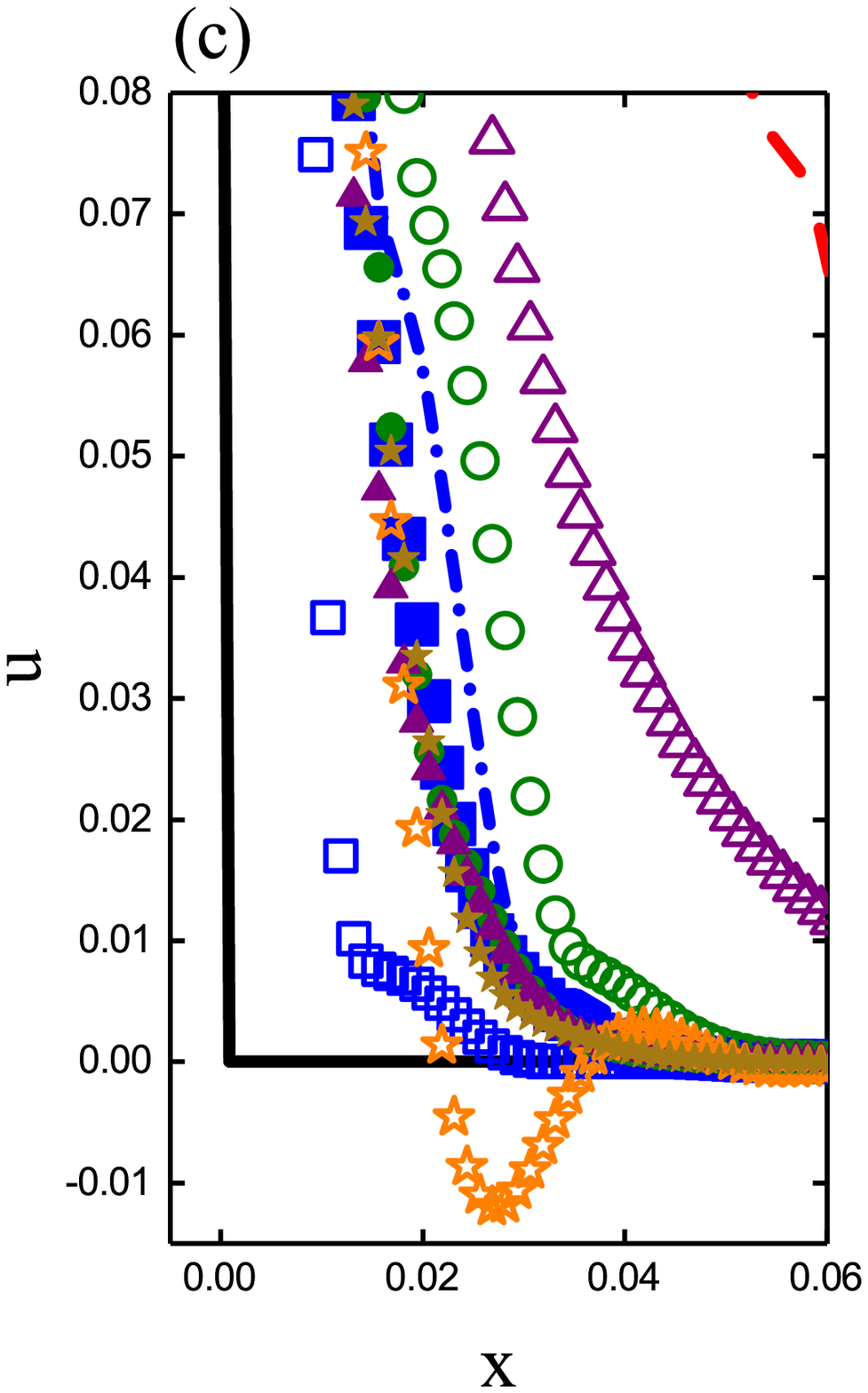}
\caption{Results of different schemes on solving Example 
\ref{ex:AccuracyTest:Z} (Case 1) with $t = 200, N = 1600$.}
\label{fig:Z:N1600:02}
\end{figure}

\begin{figure}[ht]
\centering
\includegraphics[height=0.33\textwidth]
{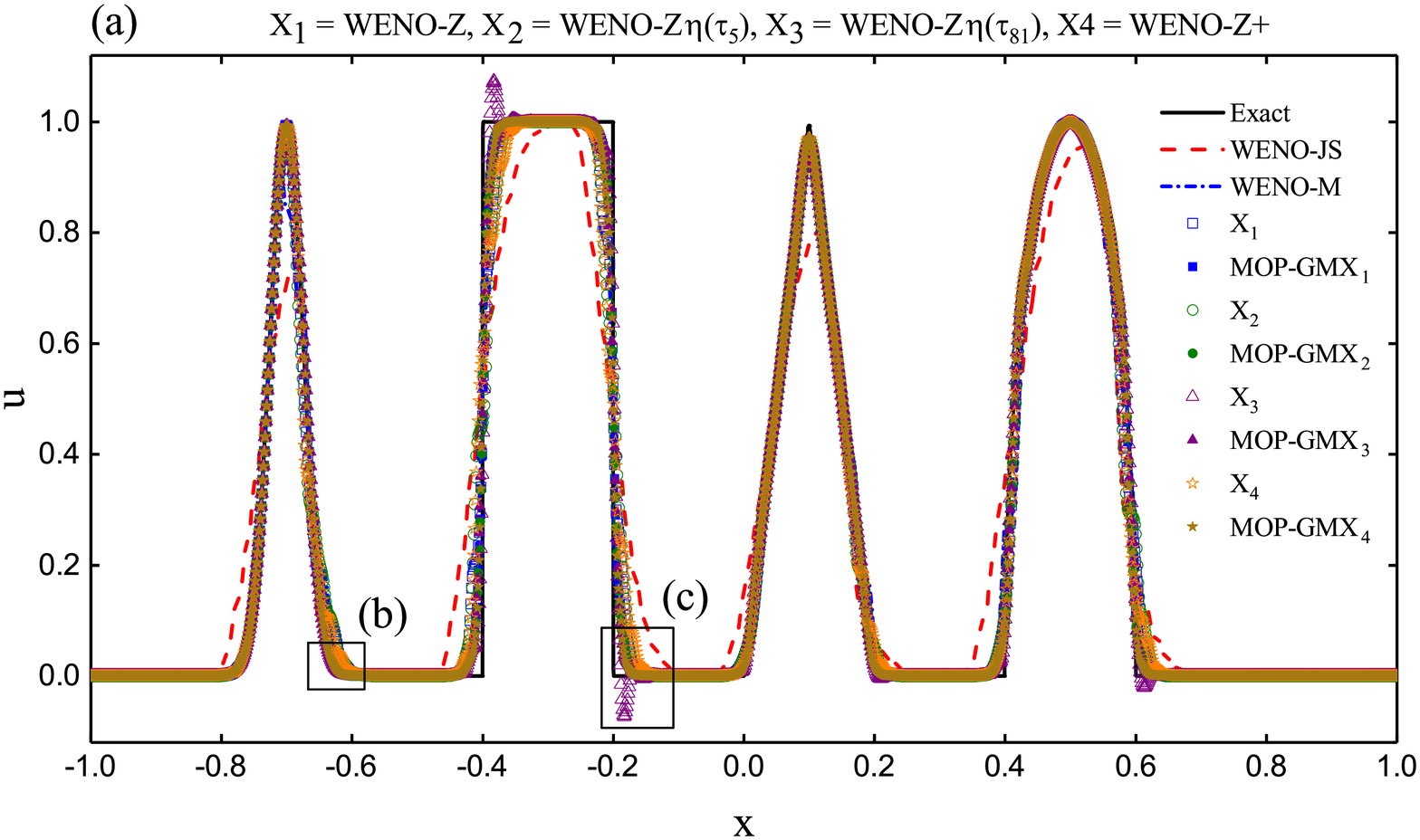}
\includegraphics[height=0.33\textwidth]
{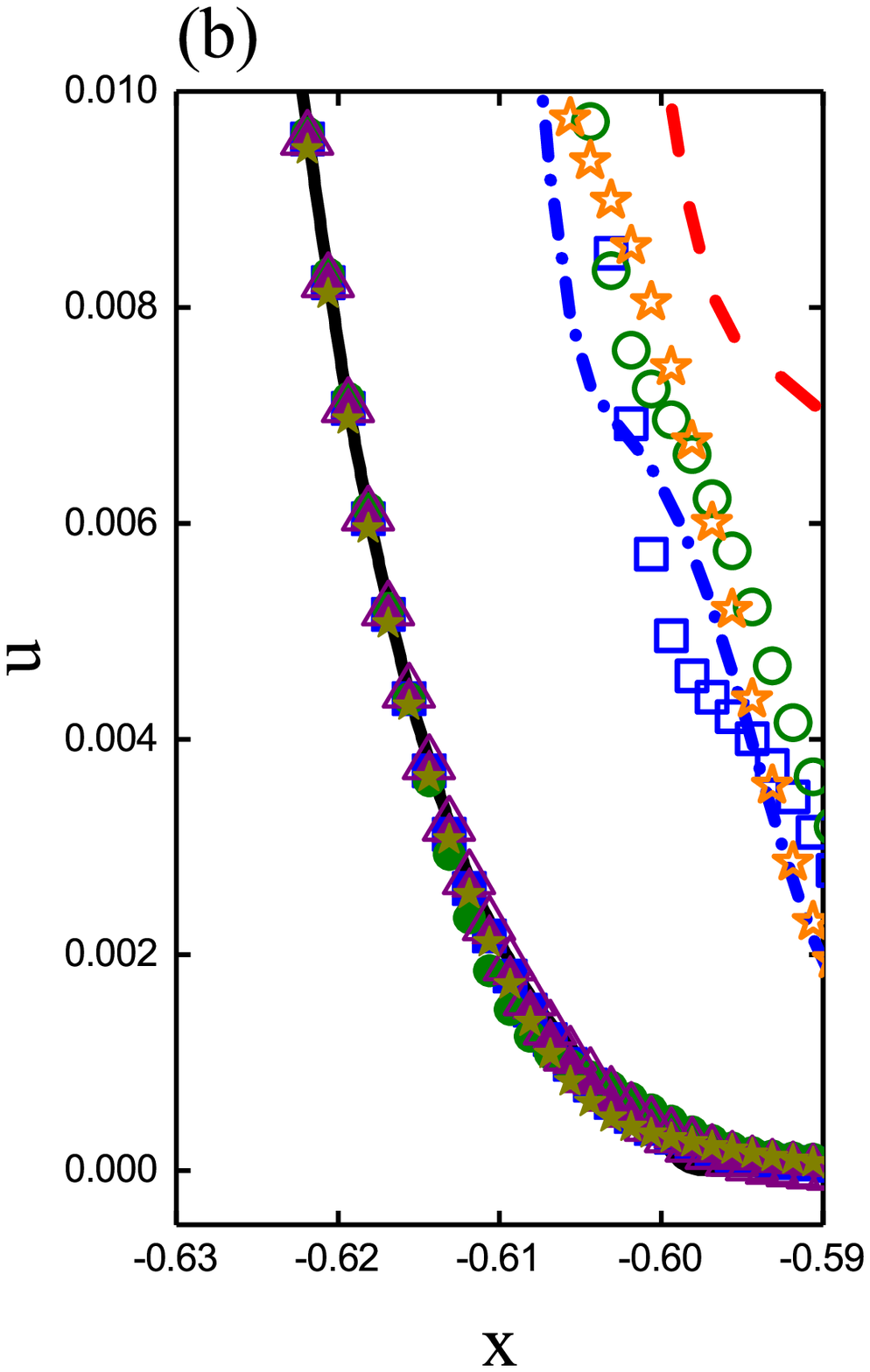}
\includegraphics[height=0.33\textwidth]
{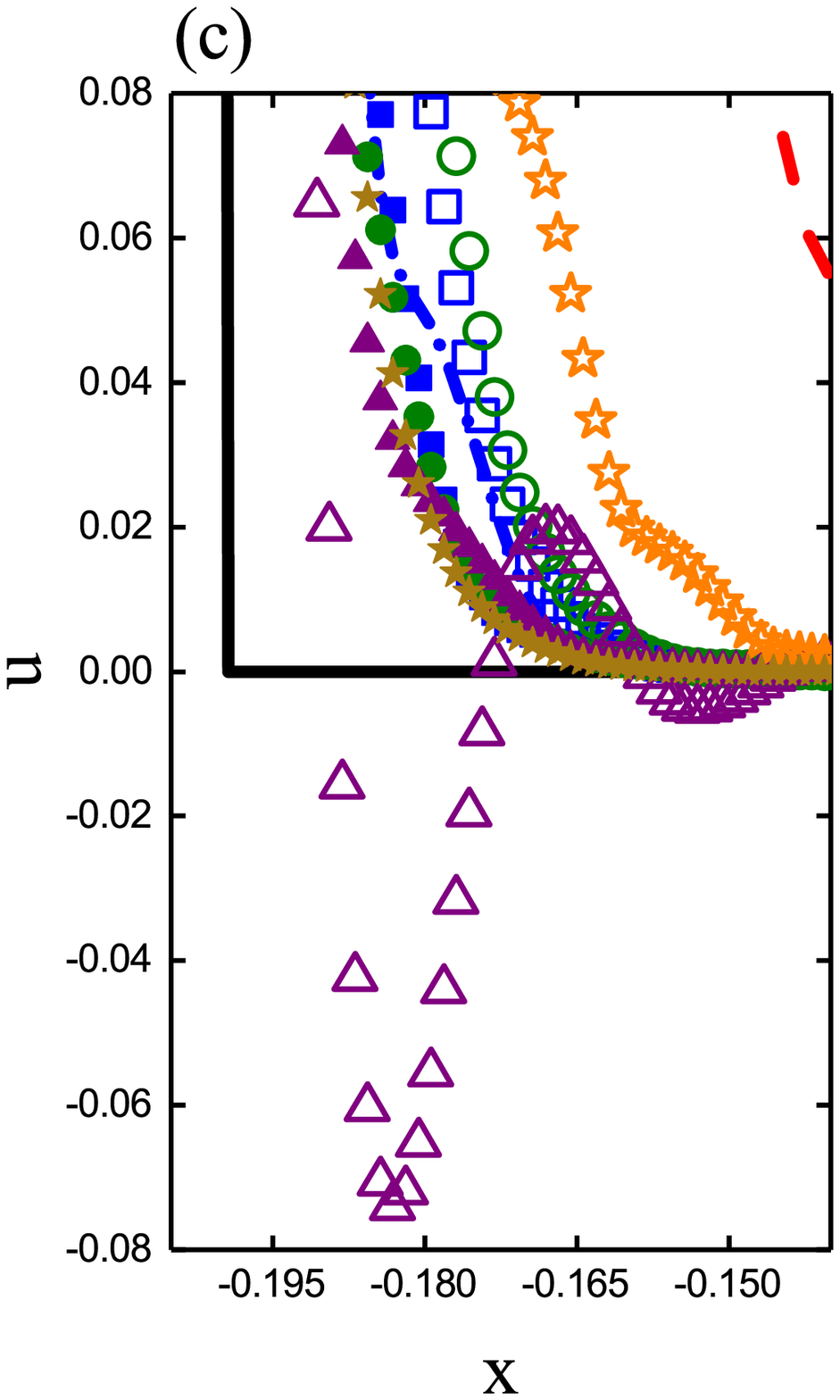}
\caption{Results of different schemes on solving Example 
\ref{ex:AccuracyTest:Z} (Case 2) with $t = 200, N = 1600$.}
\label{fig:SLP:N1600:01}
\end{figure}

\begin{figure}[ht]
\centering
\includegraphics[height=0.33\textwidth]
{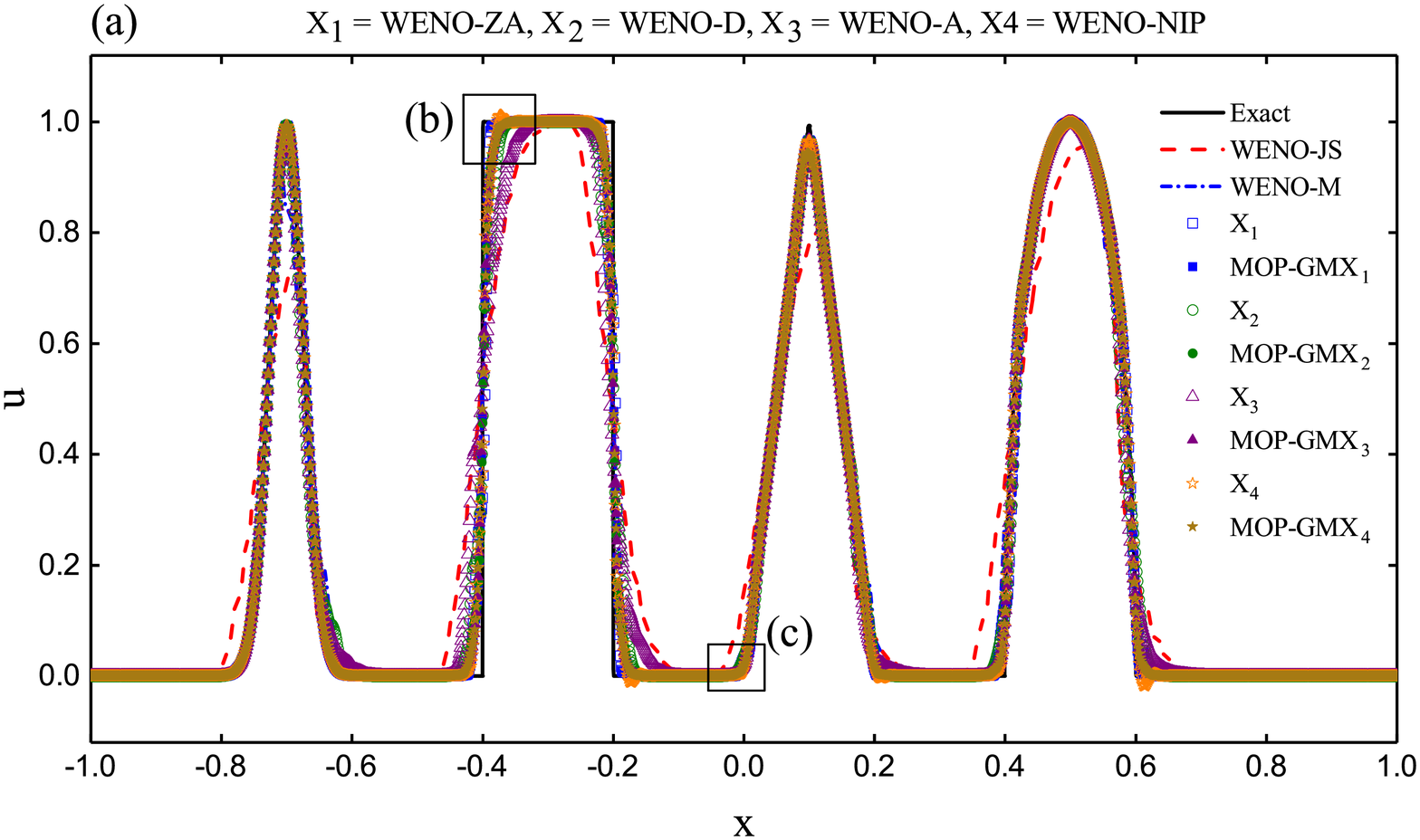}
\includegraphics[height=0.33\textwidth]
{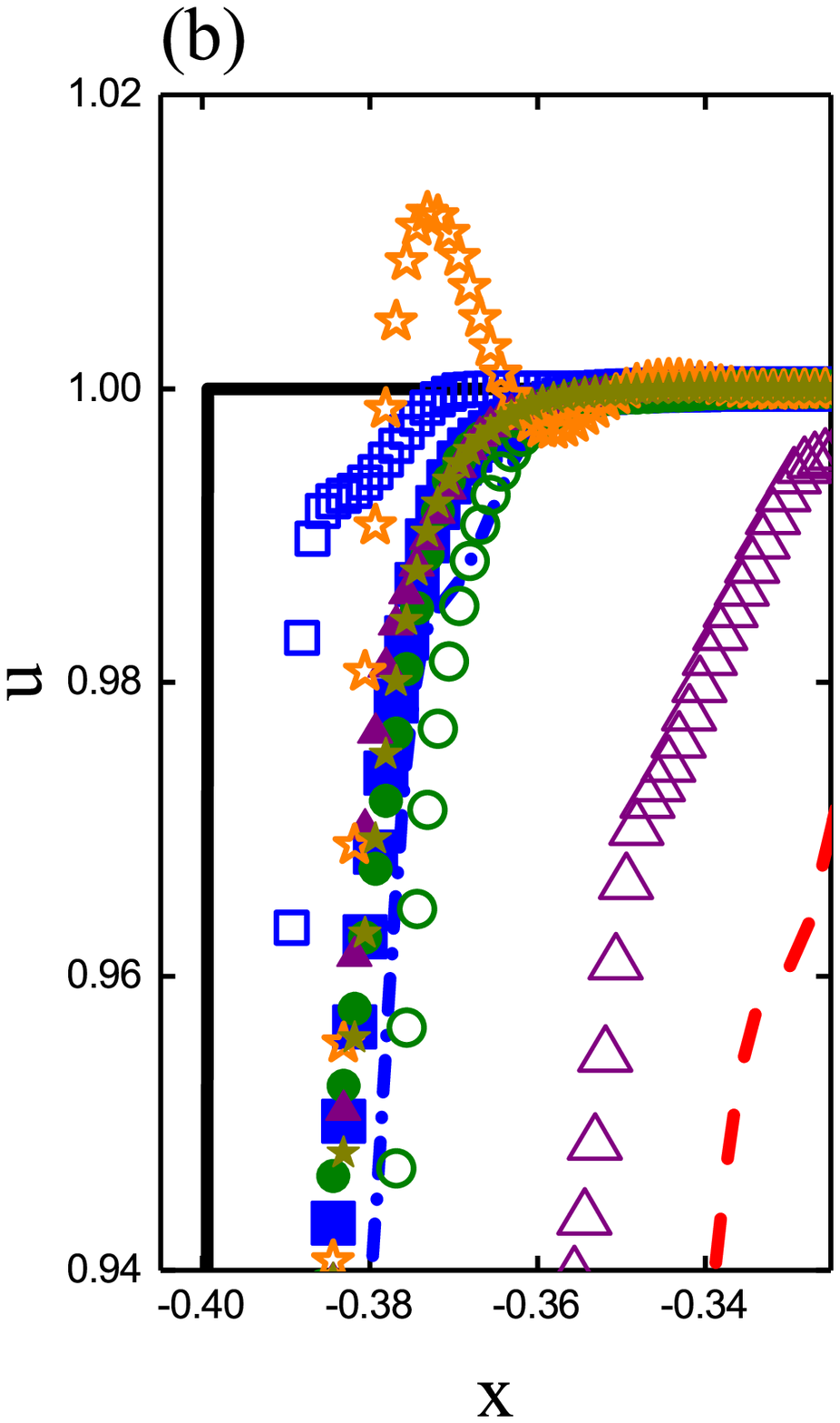}
\includegraphics[height=0.33\textwidth]
{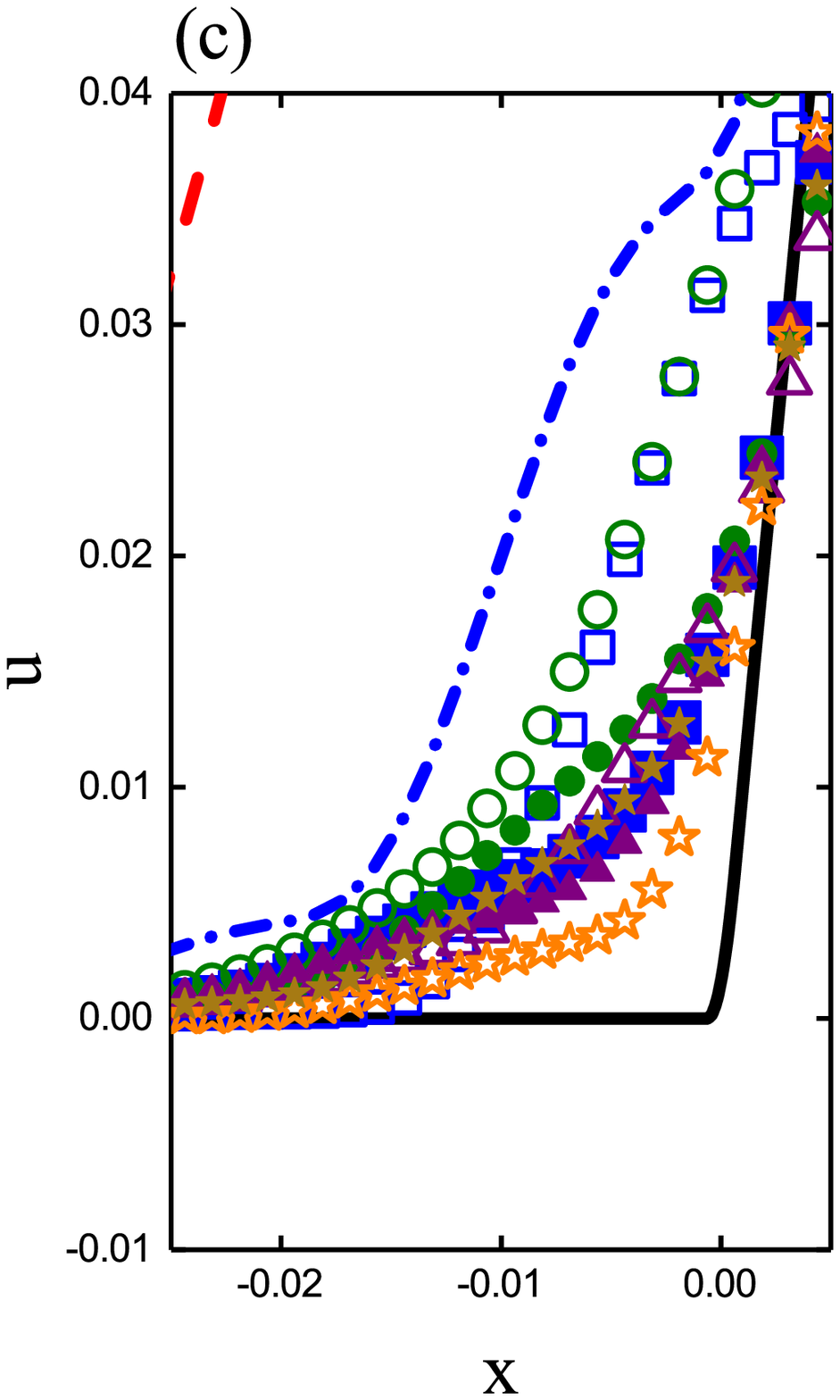}
\caption{Results of different schemes on solving Example 
\ref{ex:AccuracyTest:Z} (Case 2) with $t = 200, N = 1600$.}
\label{fig:SLP:N1600:02}
\end{figure}

Figures \ref{fig:Z:N1600:01}-\ref{fig:Z:N1600:02} and Figures 
\ref{fig:SLP:N1600:01}-\ref{fig:SLP:N1600:02} give the solutions of different schemes with $t = 200, N = 1600$ for 
Case 1 and Case 2, respectively. It can be observed that: (1) the 
MOP-GMWENO-Z, MOP-GMWENO-Z$\eta(\tau_{5})$, MOP-GMWENO-Z+, 
MOP-GMWENO-ZA, MOP-GMWENO-D and MOP-GMWENO-A schemes generate no any 
spurious oscillations and obtain numerical results with 
significantly better resolutions than the
WENO-Z, WENO-Z$\eta(\tau_{5})$, WENO-Z+, WENO-ZA, WENO-D, WENO-A 
schemes, and the WENO-JS and WENO-M schemes; (2) the 
WENO-Z$\eta(\tau_{81})$ and WENO-NIP schemes inevitably generate 
severe spurious oscillations while the
MOP-GMWENO-Z$\eta(\tau_{81})$ and MOP-GMWENO-NIP schemes can avoid 
these spurious oscillations successfully and at the same time they 
can achieve considerablely high resolutions, as accurate as the 
other MOP-GMWENO-X schemes.

\begin{figure}[ht]
\centering
\includegraphics[height=0.33\textwidth]
{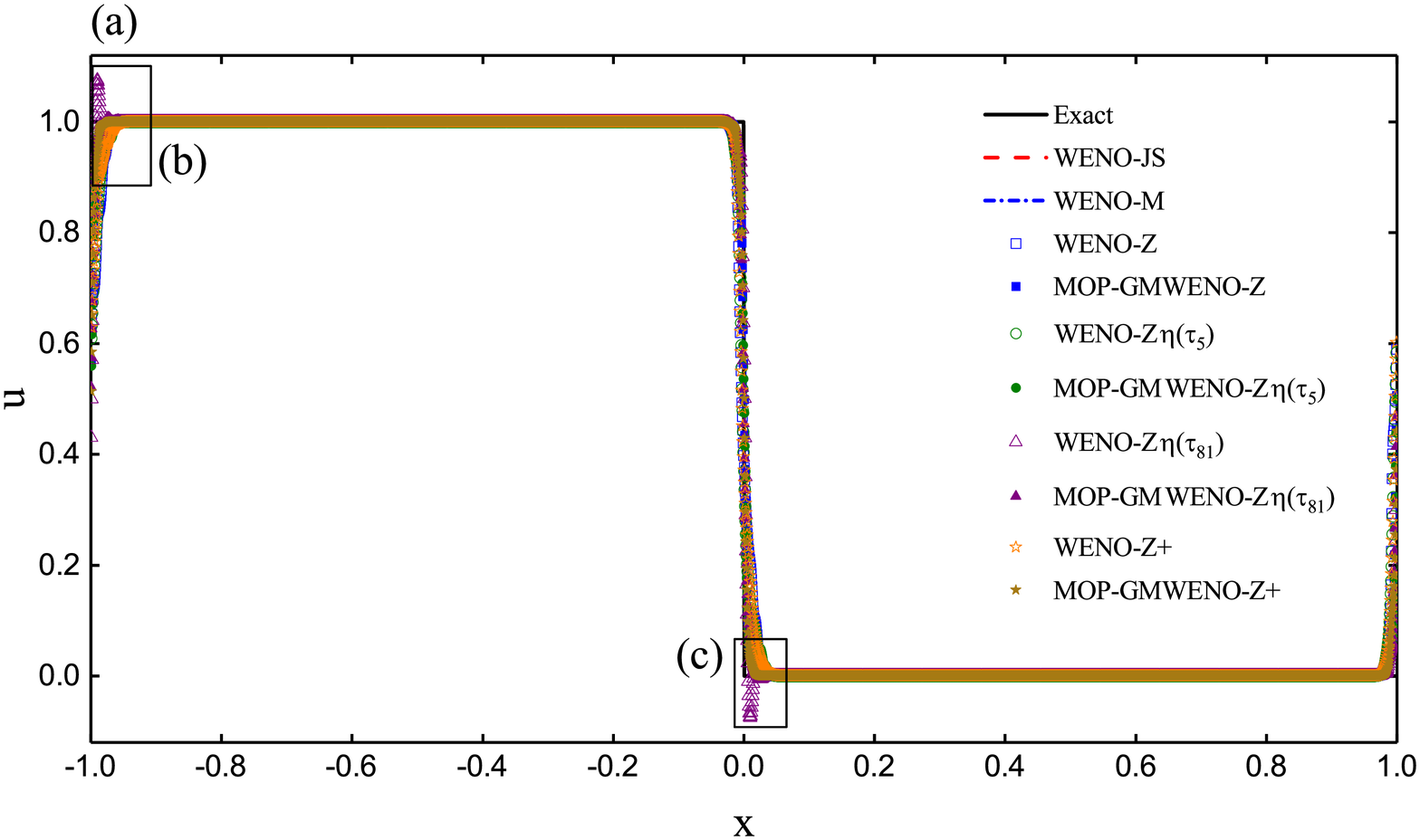}
\includegraphics[height=0.33\textwidth]
{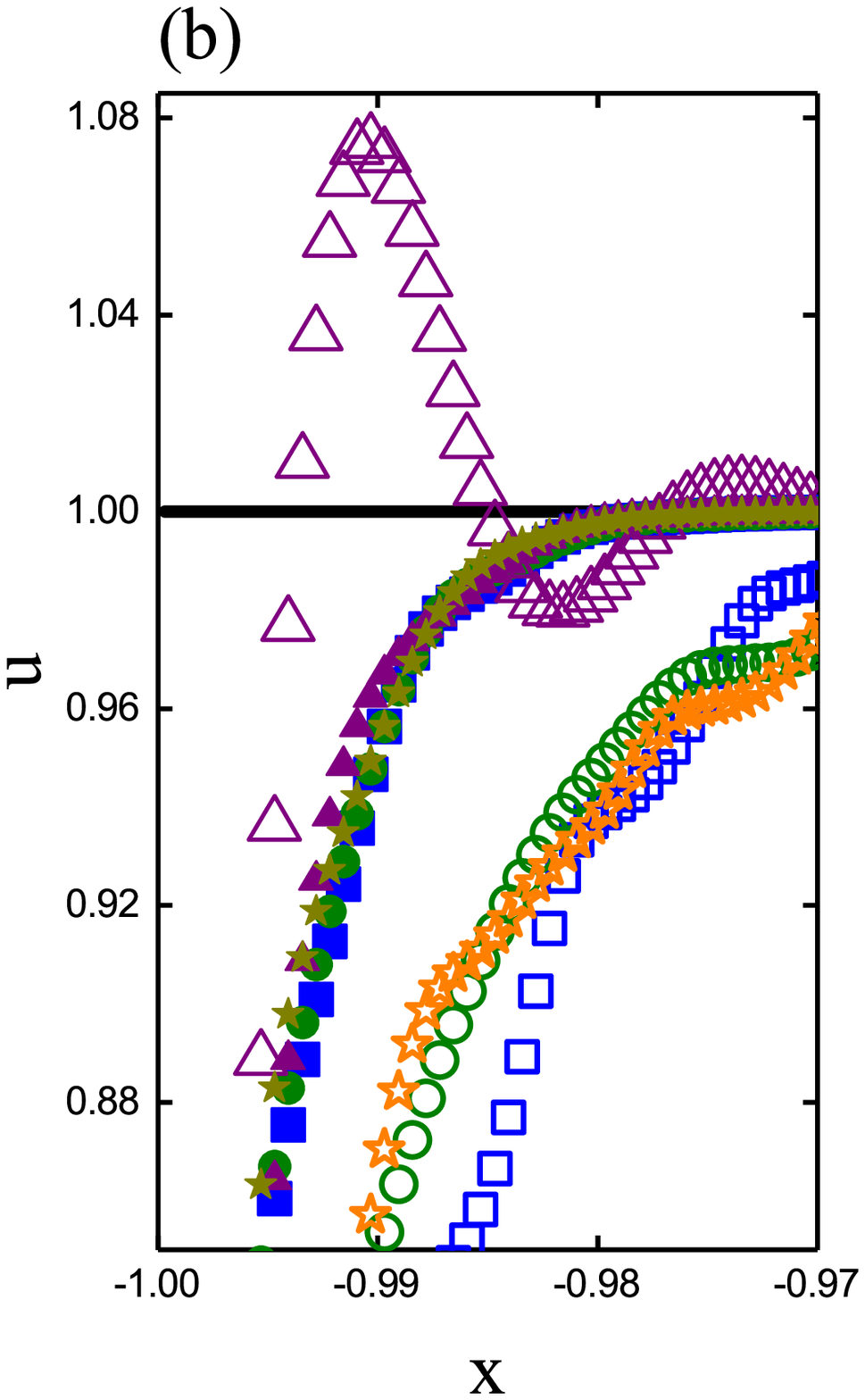}
\includegraphics[height=0.33\textwidth]
{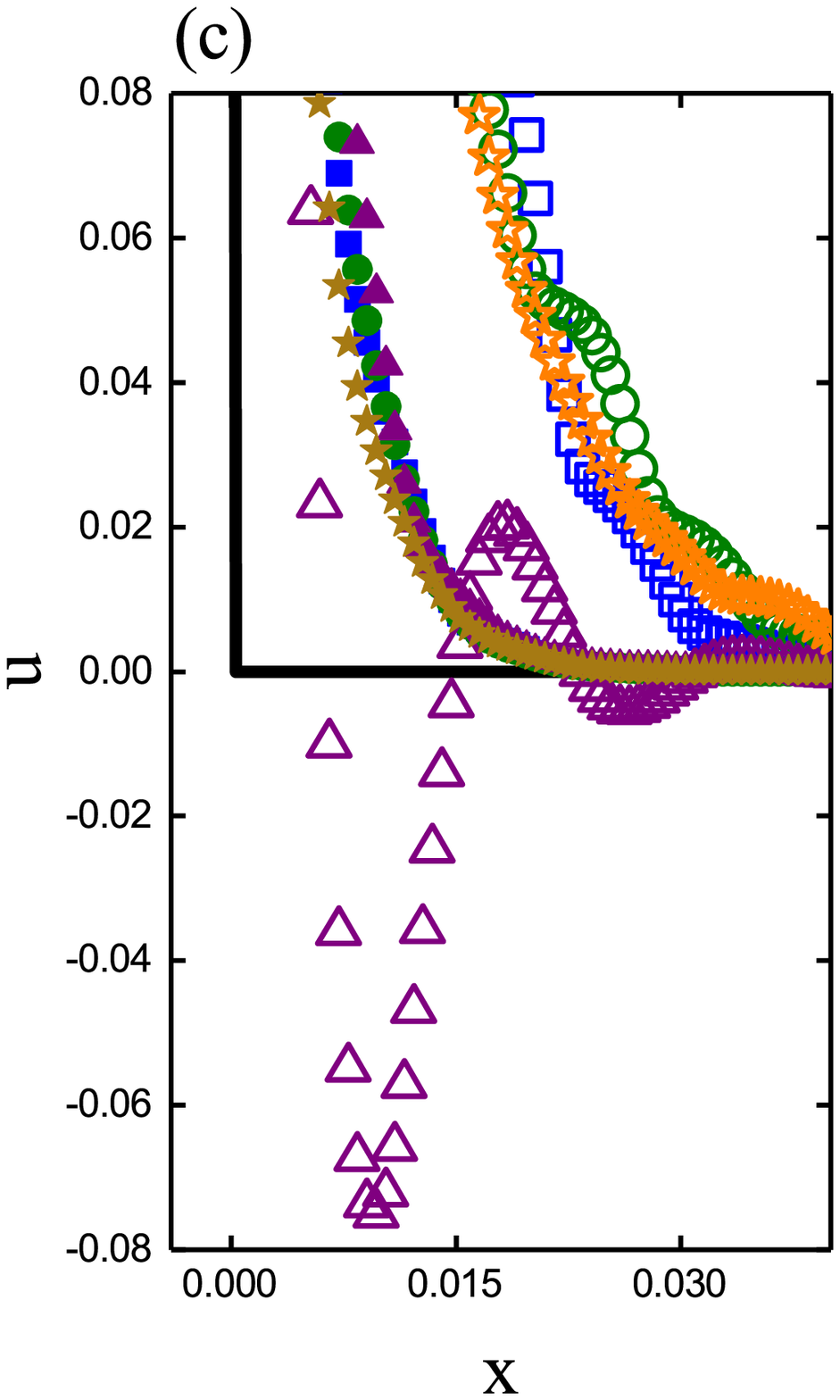}
\caption{Results of different schemes on solving Example 
\ref{ex:AccuracyTest:Z} (Case 1) with $t = 200, N = 3200$.}
\label{fig:Z:N3200:01}
\end{figure}

\begin{figure}[ht]
\centering
\includegraphics[height=0.33\textwidth]
{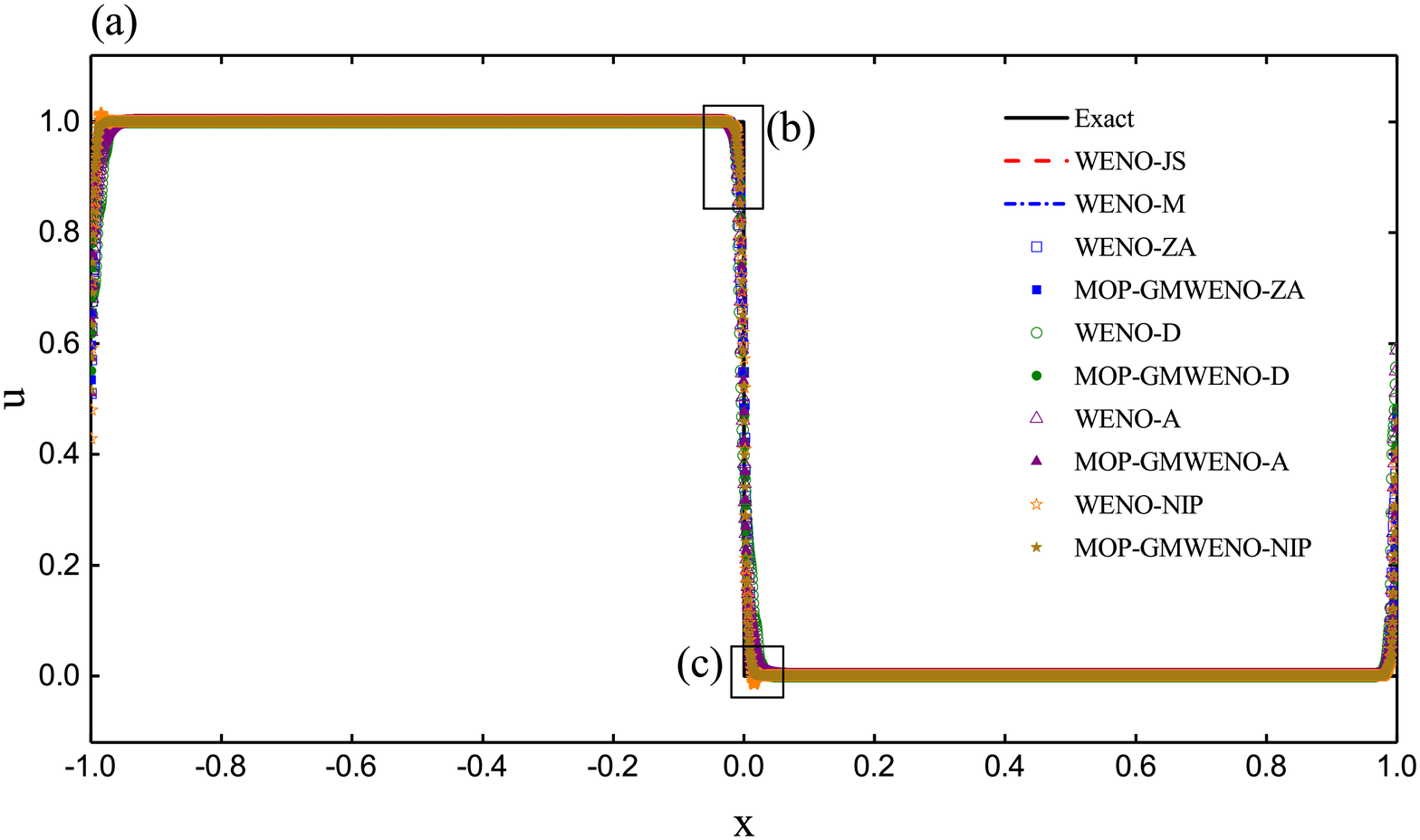}
\includegraphics[height=0.33\textwidth]
{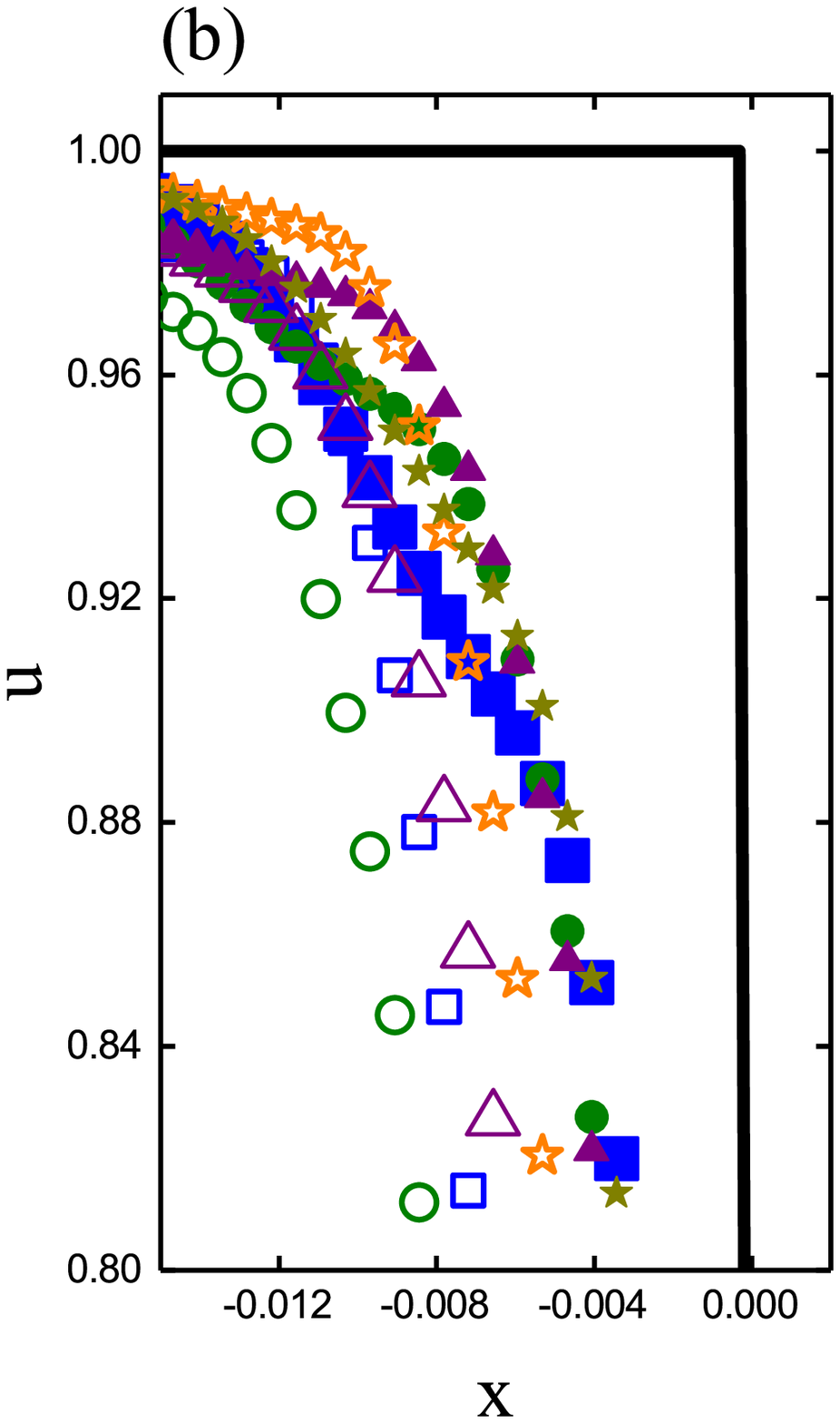}
\includegraphics[height=0.33\textwidth]
{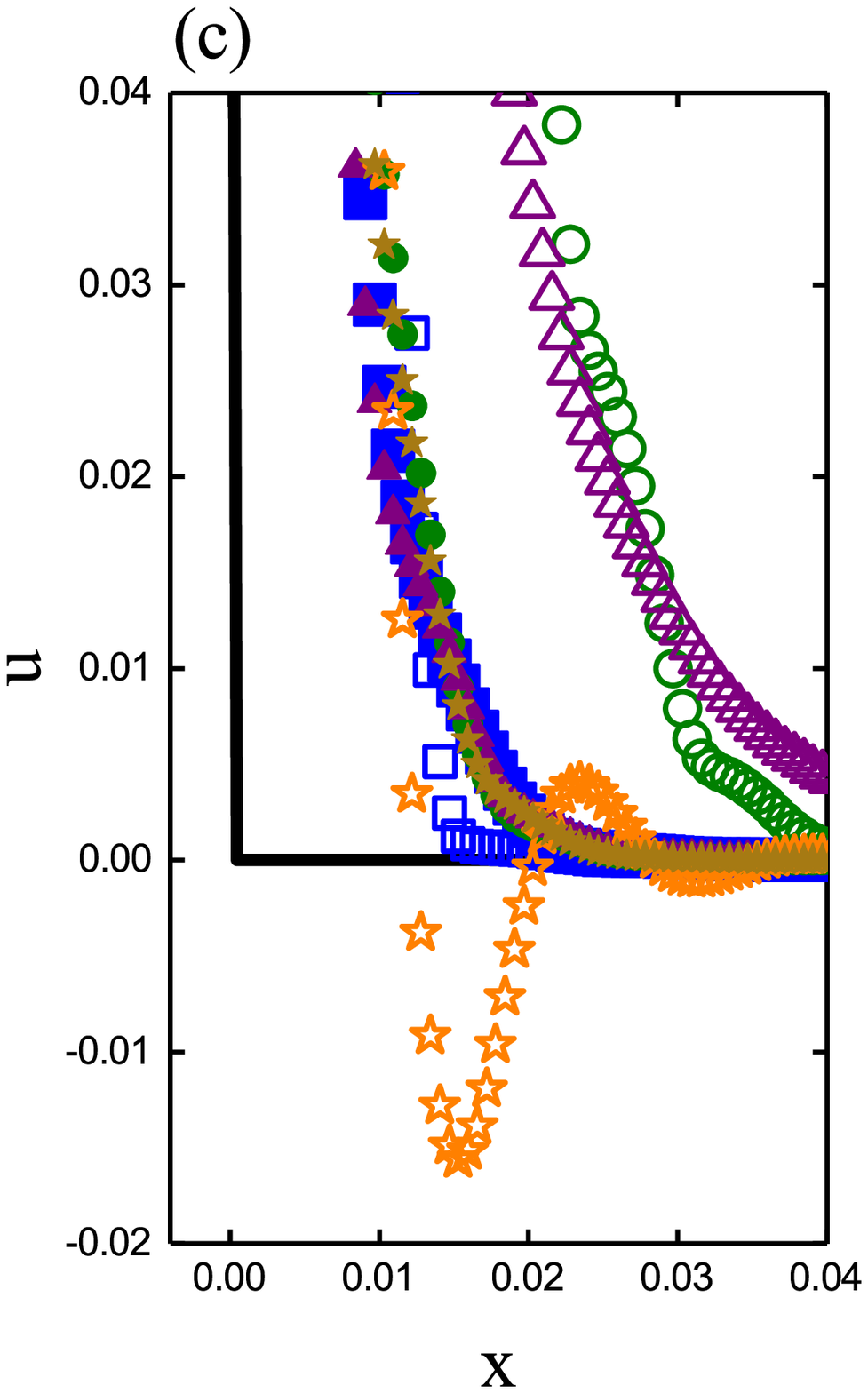}
\caption{Results of different schemes on solving Example 
\ref{ex:AccuracyTest:Z} (Case 1) with $t = 200, N = 3200$.}
\label{fig:Z:N3200:02}
\end{figure}

\begin{figure}[ht]
\centering
\includegraphics[height=0.33\textwidth]
{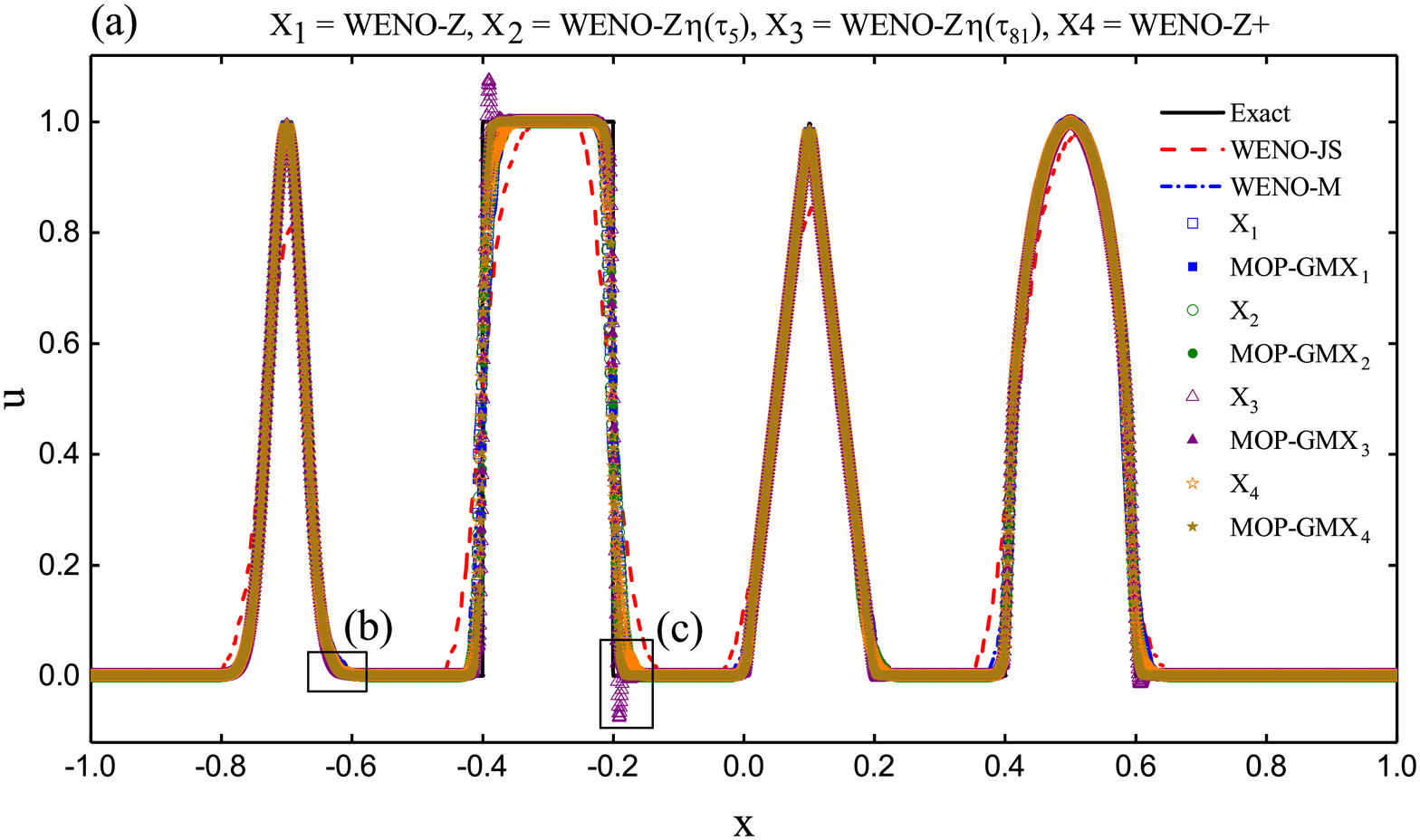}
\includegraphics[height=0.33\textwidth]
{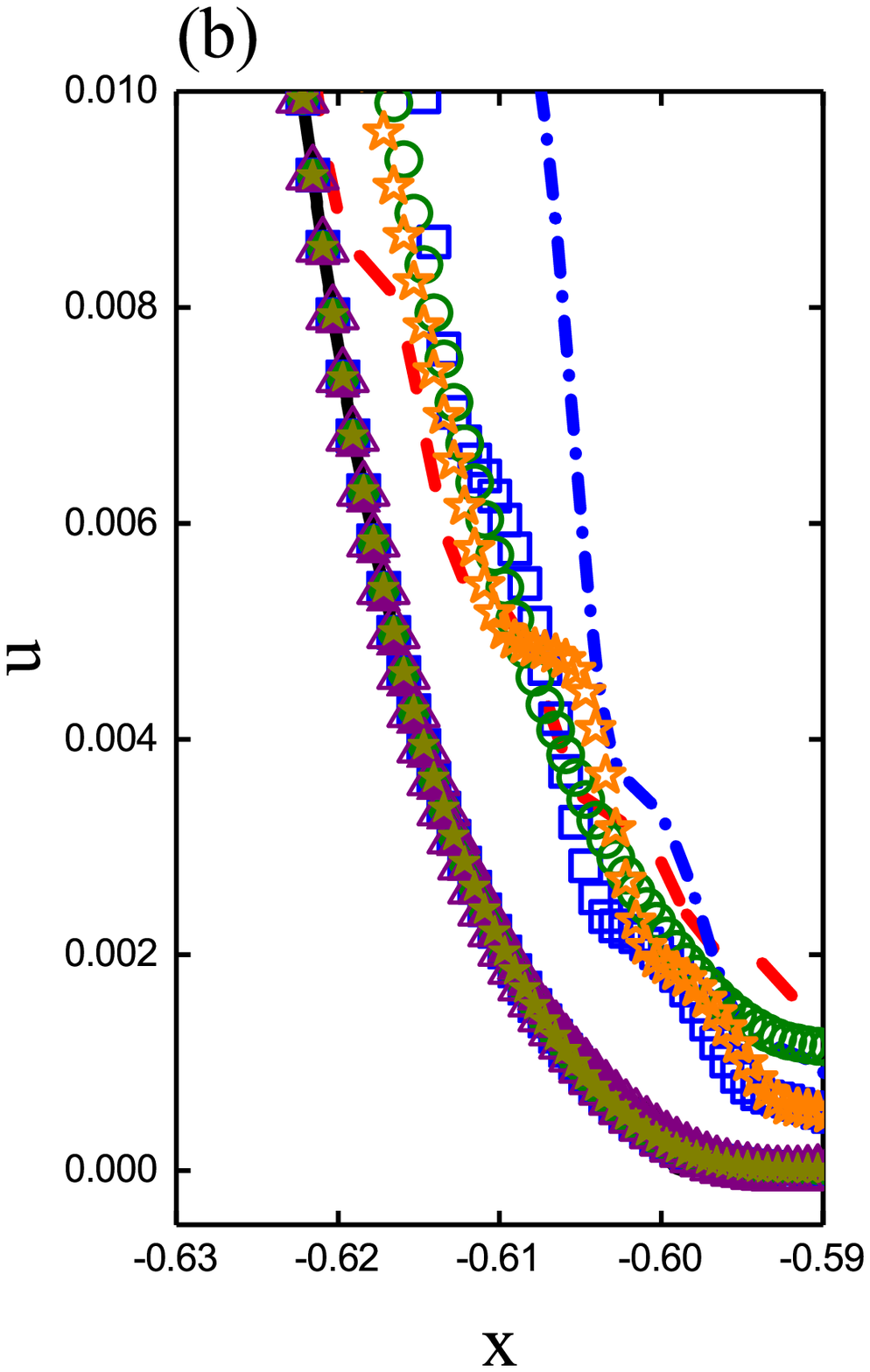}
\includegraphics[height=0.33\textwidth]
{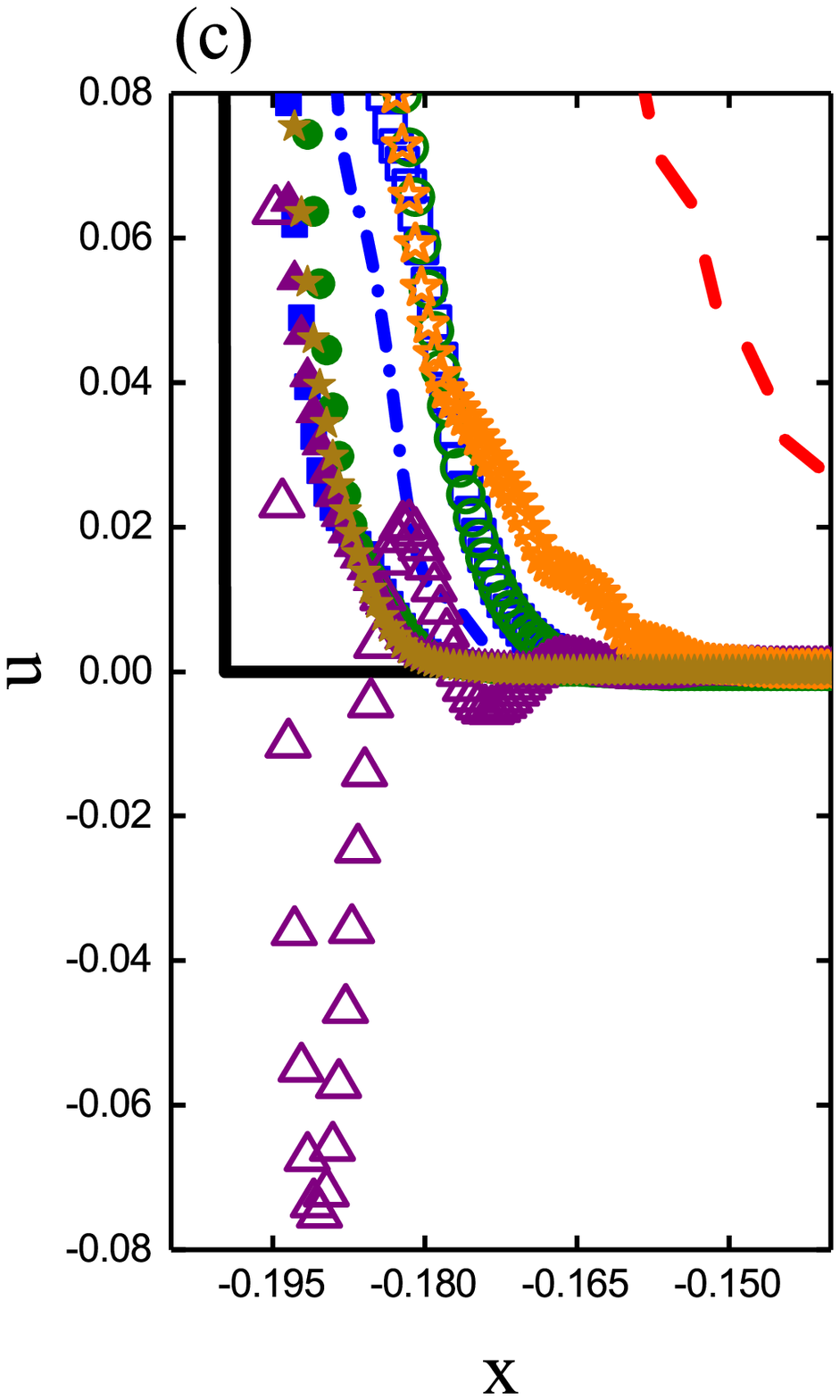}
\caption{Results of different schemes on solving Example 
\ref{ex:AccuracyTest:Z} (Case 2) with $t = 200, N = 3200$.}
\label{fig:SLP:N3200:01}
\end{figure}

\begin{figure}[ht]
\centering
\includegraphics[height=0.33\textwidth]
{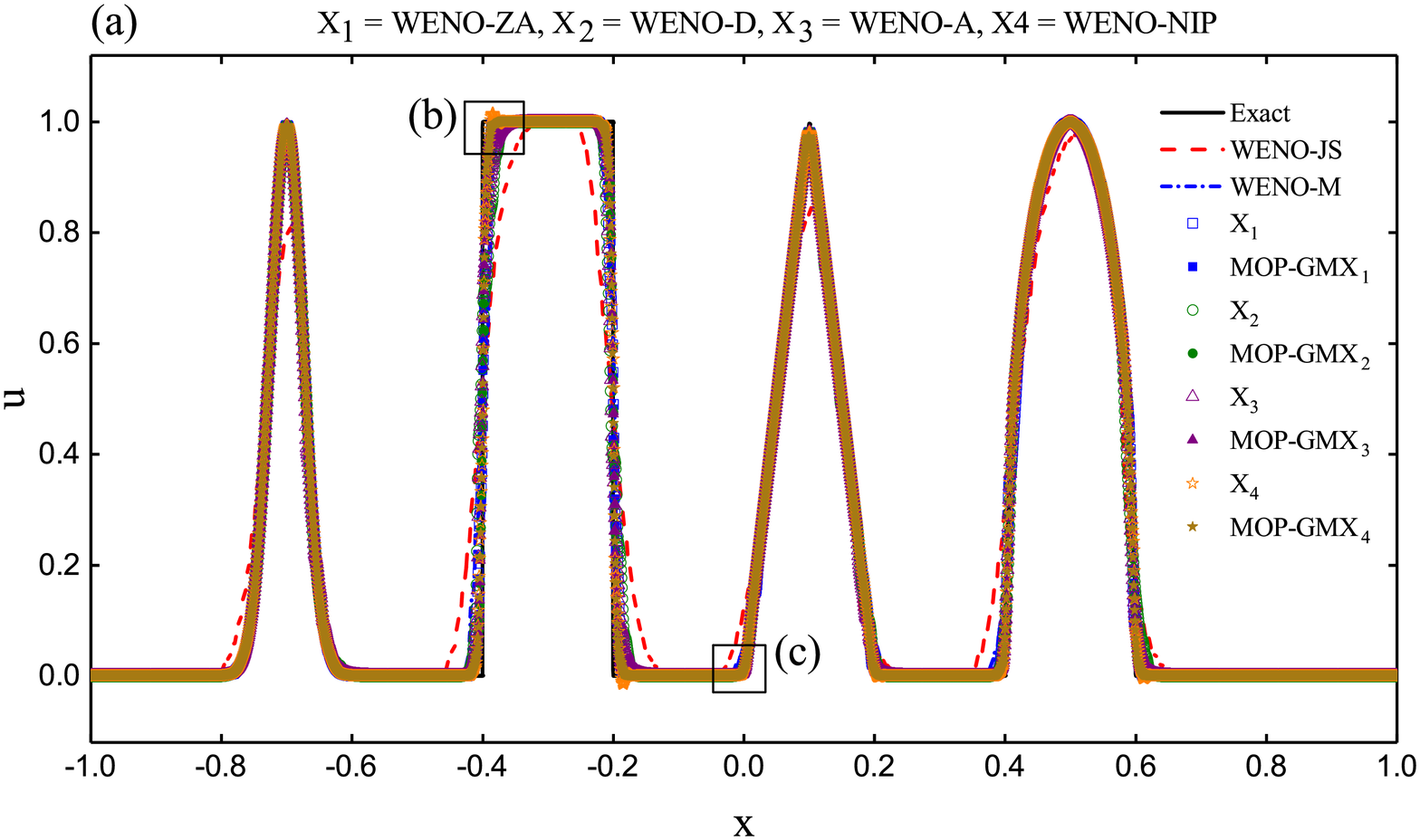}
\includegraphics[height=0.33\textwidth]
{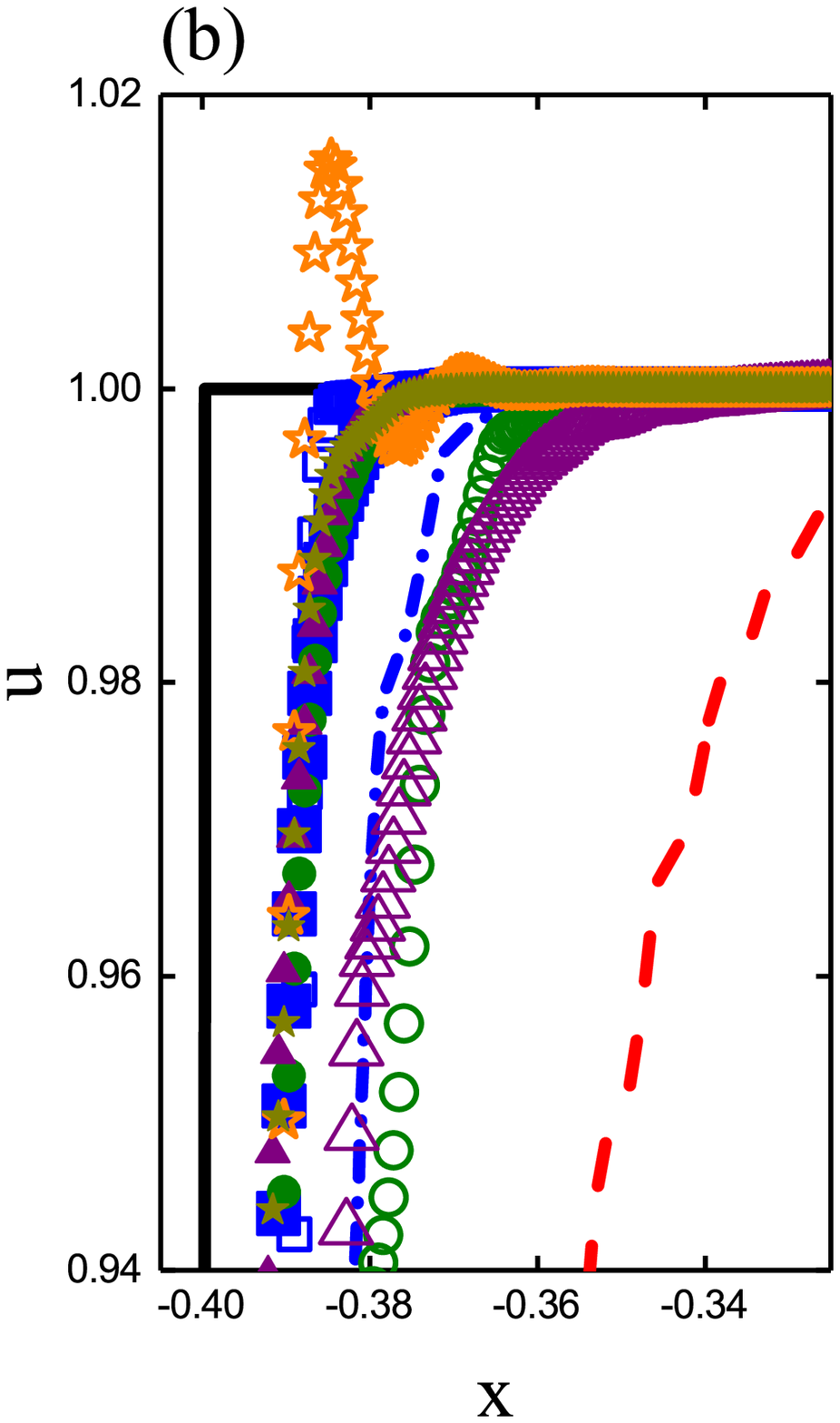}
\includegraphics[height=0.33\textwidth]
{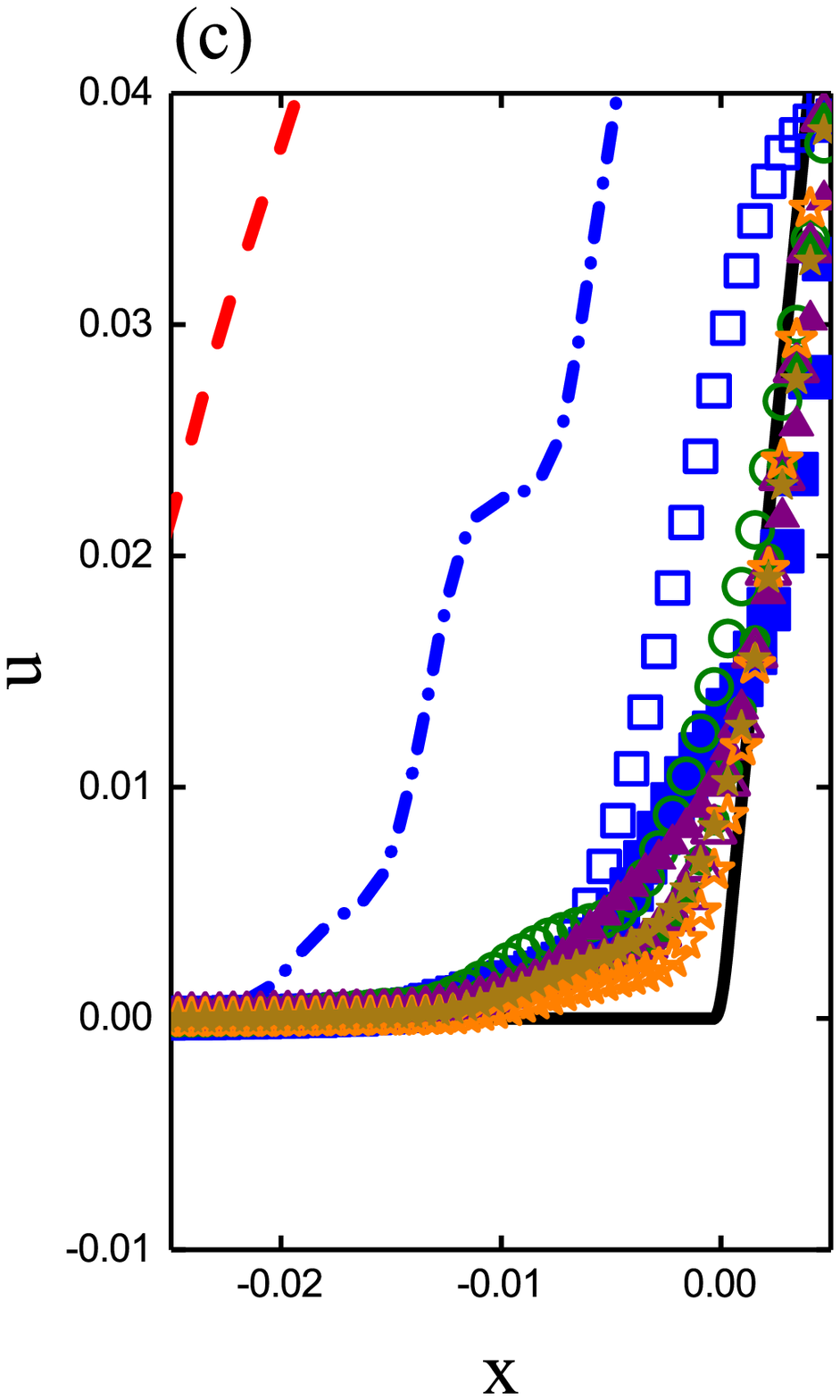}
\caption{Results of different schemes on solving Example 
\ref{ex:AccuracyTest:Z} (Case 2) with $t = 200, N = 3200$.}
\label{fig:SLP:N3200:02}
\end{figure}

In Figures \ref{fig:Z:N3200:01}-\ref{fig:Z:N3200:02} and 
Figures \ref{fig:SLP:N3200:01}-\ref{fig:SLP:N3200:02}, we give the 
solutions of different schemes respectively for the case of $t = 200,
N = 3200$. We can see 
that: (1) with the increase of grid number, the spurious oscillations 
produced by WENO-Z$\eta(\tau_{81})$ and WENO-NIP appear 
to be closer to the discontinuities, and the amplitudes of these 
spurious oscillations seem to become larger, however, 
MOP-GMWENO-Z$\eta(\tau_{81})$ and MOP-GMWENO-NIP are still able to
prevent the spurious oscillations while provide greatly improved 
resolutions; (2) all the other MOP-GMWENO-X schemes still generate 
no any spurious oscillations and evidently provide much better 
resolutions than the WENO-X schemes and the 
WENO-JS and WENO-M schemes.

To sumerize, we may conclude that the \textit{order-preserving} 
property introduced in the present study, can help the WENO-Z-type 
schemes to get high resolutions and remove spurious 
oscillations at the same time for long-run simulations.

\begin{table}[!ht]
\begin{myFontSize}
\centering
\caption{Numerical errors and convergence orders of accuracy on 
Example \ref{ex:AccuracyTest:Z} at $t = 2000.0$.}
\label{table:AccuracyTest:Z:t2000}
\begin{tabular*}{\hsize}
{@{}@{\extracolsep{\fill}}cllllllll@{}}
\hline
\space    &\multicolumn{4}{l}{\cellcolor{gray!35}{WENO5-ILW}}  
          &\multicolumn{4}{l}{\cellcolor{gray!35}{WENO-JS}}  \\
\cline{2-5}  \cline{6-9}
$N$                   & $L_{1}$ error      & $L_{1}$ order 
                      & $L_{\infty}$ error & $L_{\infty}$ order
		  			  & $L_{1}$ error      & $L_{1}$ order 
      				  & $L_{\infty}$ error & $L_{\infty}$ order \\
\Xhline{0.65pt}
200                   & 1.03240E-01        & -
                      & 4.67252E-01        & -
                      & 4.48148E-01        & -
                      & 5.55748E-01        & - \\
400                   & 5.79848E-02        & 0.8323
                      & 4.70837E-01        & -0.0110
                      & 3.37220E-01        & 0.4103
                      & 5.77105E-01        & -0.0544 \\
800                   & 3.25843E-02        & 0.8315
                      & 4.74042E-01        & -0.0098
                      & 2.93752E-01        & 0.1991
                      & 5.17829E-01        & 0.1564 \\
\hline
\space    &\multicolumn{4}{l}{\cellcolor{gray!35}{WENO-Z}}
          &\multicolumn{4}{l}{\cellcolor{gray!35}{MOP-GMWENO-Z}}\\
\cline{2-5}  \cline{6-9}
$N$                   & $L_{1}$ error      & $L_{1}$ order 
                      & $L_{\infty}$ error & $L_{\infty}$ order
		  			  & $L_{1}$ error      & $L_{1}$ order 
      				  & $L_{\infty}$ error & $L_{\infty}$ order \\
\Xhline{0.65pt}
200                   & 2.08722e-01        & -
                      & 4.98364e-01        & -
                      & 1.18245e-01        & -
                      & 4.88876e-01        & - \\
400                   & 1.25878e-01        & 0.7296
                      & 5.94047e-01        & -0.2534
                      & 6.67048e-02        & 0.8259
                      & 5.32787e-01        & -0.1241\\
800                   & 8.56252e-02        & 0.5559
                      & 6.02088e-01        & -0.0194
                      & 3.79314e-02        & 0.8144
                      & 5.30365e-01        & 0.0066\\
\hline
\space    &\multicolumn{4}{l}{\cellcolor{gray!35}{WENO-Z$\eta(\tau_{5})$}}
          &\multicolumn{4}{l}{\cellcolor{gray!35}{MOP-GMWENO-Z$\eta(\tau_{5})$}}\\
\cline{2-5}  \cline{6-9}
$N$                   & $L_{1}$ error      & $L_{1}$ order 
                      & $L_{\infty}$ error & $L_{\infty}$ order
              & $L_{1}$ error      & $L_{1}$ order 
                & $L_{\infty}$ error & $L_{\infty}$ order \\
\Xhline{0.65pt}
200                   & 2.35848e-01        & -
                      & 5.46963e-01        & -
                      & 1.18115e-01        & -
                      & 4.82687e-01        & - \\
400                   & 1.34410e-01        & 0.8112
                      & 5.77563e-01        & -0.0785
                      & 6.57450e-02        & 0.8452
                      & 4.99738e-01        & -0.0501\\
800                   & 8.88830e-02        & 0.5967
                      & 6.10969e-01        & -0.0811
                      & 3.74280e-02        & 0.8128
                      & 5.14095e-01        & -0.0409\\
\hline
\space    &\multicolumn{4}{l}{\cellcolor{gray!35}{WENO-Z$\eta(\tau_{81})$}}
          &\multicolumn{4}{l}{\cellcolor{gray!35}{MOP-GMWENO-Z$\eta(\tau_{81})$}}\\
\cline{2-5}  \cline{6-9}
$N$                   & $L_{1}$ error      & $L_{1}$ order 
                      & $L_{\infty}$ error & $L_{\infty}$ order
              & $L_{1}$ error      & $L_{1}$ order 
                & $L_{\infty}$ error & $L_{\infty}$ order \\
\Xhline{0.65pt}
200                   & 1.06673e-01        & -
                      & 5.57758e-01        & -
                      & 1.16619e-01        & -
                      & 4.77852e-01        & - \\
400                   & 6.02857e-02        & 0.8233
                      & 5.67245e-01        & -0.0243
                      & 6.51575e-02        & 0.8398
                      & 4.95518e-01        & -0.0524\\
800                   & 3.50816e-02        & 0.7811
                      & 5.87196e-01        & -0.0499
                      & 3.70925e-02        & 0.8128
                      & 5.15105e-01        & -0.0559\\
\hline
\space    &\multicolumn{4}{l}{\cellcolor{gray!35}{WENO-Z+}}
          &\multicolumn{4}{l}{\cellcolor{gray!35}{MOP-GMWENO-Z+}}\\
\cline{2-5}  \cline{6-9}
$N$                   & $L_{1}$ error      & $L_{1}$ order 
                      & $L_{\infty}$ error & $L_{\infty}$ order
              & $L_{1}$ error      & $L_{1}$ order 
                & $L_{\infty}$ error & $L_{\infty}$ order \\
\Xhline{0.65pt}
200                   & 2.35835e-01        & -
                      & 5.53876e-01        & -
                      & 9.90921e-02        & -
                      & 4.99725e-01        & - \\
400                   & 1.46265e-01        & 0.6892
                      & 6.02883e-01        & -0.1223
                      & 5.76023e-02        & 0.7826
                      & 4.90011e-01        & 0.0283\\
800                   & 6.51259e-02        & 1.1673
                      & 5.83954e-01        & 0.0460
                      & 3.24177e-02        & 0.8293
                      & 4.81302e-01        & 0.0259\\
\hline
\space    &\multicolumn{4}{l}{\cellcolor{gray!35}{WENO-ZA}}
          &\multicolumn{4}{l}{\cellcolor{gray!35}{MOP-GMWENO-ZA}}\\
\cline{2-5}  \cline{6-9}
$N$                   & $L_{1}$ error      & $L_{1}$ order 
                      & $L_{\infty}$ error & $L_{\infty}$ order
              & $L_{1}$ error      & $L_{1}$ order 
                & $L_{\infty}$ error & $L_{\infty}$ order \\
\Xhline{0.65pt}
200                   & 1.31936e-01        & -
                      & 5.50979e-01        & -
                      & 1.17961e-01        & -
                      & 5.10189e-01        & - \\
400                   & 1.13400e-01        & 0.2184
                      & 6.01658e-01        & -0.1269
                      & 6.61903e-02        & 0.8336
                      & 5.02192e-01        & 0.0228\\
800                   & 6.70946e-02        & 0.7572
                      & 5.18541e-01        & 0.2145
                      & 3.65564e-02        & 0.8565
                      & 4.80328e-01        & 0.0642\\
\hline
\space    &\multicolumn{4}{l}{\cellcolor{gray!35}{WENO-D}}
          &\multicolumn{4}{l}{\cellcolor{gray!35}{MOP-GMWENO-D}}\\
\cline{2-5}  \cline{6-9}
$N$                   & $L_{1}$ error      & $L_{1}$ order 
                      & $L_{\infty}$ error & $L_{\infty}$ order
              & $L_{1}$ error      & $L_{1}$ order 
                & $L_{\infty}$ error & $L_{\infty}$ order \\
\Xhline{0.65pt}
200                   & 2.08722e-01        & -
                      & 4.98364e-01        & -
                      & 1.18950e-01        & -
                      & 4.90750e-01        & - \\
400                   & 1.25884e-01        & 0.7295
                      & 5.94701e-01        & -0.2550
                      & 6.56664e-02        & 0.8571
                      & 4.81270e-01        & 0.0281\\
800                   & 8.52679e-02        & 0.5620
                      & 5.98136e-01        & -0.0083
                      & 3.80765e-02        & 0.7863
                      & 4.84771e-01        & -0.0105\\
\hline
\space    &\multicolumn{4}{l}{\cellcolor{gray!35}{WENO-A}}
          &\multicolumn{4}{l}{\cellcolor{gray!35}{MOP-GMWENO-A}}\\
\cline{2-5}  \cline{6-9}
$N$                   & $L_{1}$ error      & $L_{1}$ order 
                      & $L_{\infty}$ error & $L_{\infty}$ order
              & $L_{1}$ error      & $L_{1}$ order 
                & $L_{\infty}$ error & $L_{\infty}$ order \\
\Xhline{0.65pt}
200                   & 3.29771e-01        & -
                      & 5.51867e-01        & -
                      & 1.16473e-01        & -
                      & 4.80577e-01        & - \\
400                   & 2.02749e-01        & 0.7018
                      & 5.35326e-01        & 0.0439
                      & 6.65634e-02        & 0.8072
                      & 4.97352e-01        & -0.0495\\
800                   & 1.05379e-01        & 0.9441
                      & 5.54605e-01        & -0.0510
                      & 3.62328e-02        & 0.8774
                      & 4.87916e-01        & 0.0276\\
\hline
\space    &\multicolumn{4}{l}{\cellcolor{gray!35}{WENO-NIP}}
          &\multicolumn{4}{l}{\cellcolor{gray!35}{MOP-GMWENO-NIP}}\\
\cline{2-5}  \cline{6-9}
$N$                   & $L_{1}$ error      & $L_{1}$ order 
                      & $L_{\infty}$ error & $L_{\infty}$ order
              & $L_{1}$ error      & $L_{1}$ order 
                & $L_{\infty}$ error & $L_{\infty}$ order \\
\Xhline{0.65pt}
200                   & 1.19293e-01        & -
                      & 5.83513e-01        & -
                      & 1.14420e-01        & -
                      & 4.85124e-01        & - \\
400                   & 7.46960e-02        & 0.6754
                      & 5.61235e-01        & 0.0562
                      & 6.10704e-02        & 0.9058
                      & 4.82764e-01        & 0.0070\\
800                   & 4.56612e-02        & 0.7101
                      & 5.18618e-01        & 0.1139
                      & 3.35626e-02        & 0.8636
                      & 4.75671e-01        & 0.0214\\
\hline
\end{tabular*}
\end{myFontSize}
\end{table}

\begin{table}[!ht]
\begin{myFontSize}
\centering
\caption{Numerical errors and convergence orders of accuracy on 
Example \ref{ex:AccuracyTest:SLP} at $t = 2000.0$.}
\label{table:AccuracyTest:SLP:t2000}
\begin{tabular*}{\hsize}
{@{}@{\extracolsep{\fill}}cllllllll@{}}
\hline
\space    &\multicolumn{4}{l}{\cellcolor{gray!35}{WENO5-ILW}}  
          &\multicolumn{4}{l}{\cellcolor{gray!35}{WENO-JS}}  \\
\cline{2-5}  \cline{6-9}
$N$                   & $L_{1}$ error      & $L_{1}$ order 
                      & $L_{\infty}$ error & $L_{\infty}$ order
		  			  & $L_{1}$ error      & $L_{1}$ order 
      				  & $L_{\infty}$ error & $L_{\infty}$ order \\
\Xhline{0.65pt}
200                   & 2.27171E-01        & -
                      & 5.14236E-01        & -
                      & 6.12899E-01        & -
                      & 7.99265E-01        & - \\
400                   & 1.15918E-01        & 0.9707
                      & 4.77803E-01        & 0.1060
                      & 5.99215E-01        & 0.0326
                      & 8.20493E-01        & -0.0378 \\
800                   & 5.35871E-02        & 1.1131
                      & 4.74317E-01        & 0.0106
                      & 5.50158E-01        & 0.1232
                      & 8.14650E-01        & 0.0103 \\
\hline
\space    &\multicolumn{4}{l}{\cellcolor{gray!35}{WENO-Z}}
          &\multicolumn{4}{l}{\cellcolor{gray!35}{MOP-GMWENO-Z}}\\
\cline{2-5}  \cline{6-9}
$N$                   & $L_{1}$ error      & $L_{1}$ order 
                      & $L_{\infty}$ error & $L_{\infty}$ order
		  			  & $L_{1}$ error      & $L_{1}$ order 
      				  & $L_{\infty}$ error & $L_{\infty}$ order \\
\Xhline{0.65pt}
200                   & 3.86995e-01        & -
                      & 6.85835e-01        & -
                      & 4.51514e-01        & -
                      & 7.68680e-01        & - \\
400                   & 2.02287e-01        & 0.9359
                      & 5.18993e-01        & 0.4021
                      & 1.76765e-01        & 1.3529
                      & 5.01169e-01        & 0.6171\\
800                   & 1.66703e-01        & 0.2791 
                      & 5.04564e-01        & 0.0407
                      & 6.44772e-02        & 1.4550
                      & 5.04415e-01        & -0.0093\\
\hline
\space    &\multicolumn{4}{l}{\cellcolor{gray!35}{WENO-Z$\eta(\tau_{5})$}}
          &\multicolumn{4}{l}{\cellcolor{gray!35}{MOP-GMWENO-Z$\eta(\tau_{5})$}}\\
\cline{2-5}  \cline{6-9}
$N$                   & $L_{1}$ error      & $L_{1}$ order 
                      & $L_{\infty}$ error & $L_{\infty}$ order
		  			  & $L_{1}$ error      & $L_{1}$ order 
      				  & $L_{\infty}$ error & $L_{\infty}$ order \\
\Xhline{0.65pt}
200                   & 3.24546e-01        & -
                      & 6.51925e-01        & -
                      & 4.48811e-01        & -
                      & 7.65748e-01        & - \\
400                   & 2.31374e-01        & 0.4882
                      & 5.42432e-01        & 0.2653
                      & 1.82206e-01        & 1.3005
                      & 5.47587e-01        & 0.4838\\
800                   & 1.55489e-01        & 0.5734
                      & 5.12270e-01        & 0.0825
                      & 6.41444e-02        & 1.5062
                      & 4.98985e-01        & 0.1341\\
\hline
\space    &\multicolumn{4}{l}{\cellcolor{gray!35}{WENO-Z$\eta(\tau_{81})$}}
          &\multicolumn{4}{l}{\cellcolor{gray!35}{MOP-GMWENO-Z$\eta(\tau_{81})$}}\\
\cline{2-5}  \cline{6-9}
$N$                   & $L_{1}$ error      & $L_{1}$ order 
                      & $L_{\infty}$ error & $L_{\infty}$ order
		  			  & $L_{1}$ error      & $L_{1}$ order 
      				  & $L_{\infty}$ error & $L_{\infty}$ order \\
\Xhline{0.65pt}
200                   & 2.42963e-01        & -
                      & 6.39818e-01        & -
                      & 3.77388e-01        & -
                      & 7.39311e-01        & - \\
400                   & 1.33752e-01        & 0.8612
                      & 6.01344e-01        & 0.0895
                      & 1.61629e-01        & 1.2234
                      & 4.91776e-01        & 0.5882\\
800                   & 5.89144e-02        & 1.1829
                      & 5.73819e-01        & 0.0676
                      & 6.41956e-02        & 1.3321
                      & 4.95672e-01        & -0.0114\\
\hline
\space    &\multicolumn{4}{l}{\cellcolor{gray!35}{WENO-Z+}}
          &\multicolumn{4}{l}{\cellcolor{gray!35}{MOP-GMWENO-Z+}}\\
\cline{2-5}  \cline{6-9}
$N$                   & $L_{1}$ error      & $L_{1}$ order 
                      & $L_{\infty}$ error & $L_{\infty}$ order
		  			  & $L_{1}$ error      & $L_{1}$ order 
      				  & $L_{\infty}$ error & $L_{\infty}$ order \\
\Xhline{0.65pt}
200                   & 2.99492e-01        & -
                      & 5.45598e-01        & -
                      & 2.66825e-01        & -
                      & 7.08732e-01        & - \\
400                   & 2.42059e-01        & 0.3072
                      & 4.96267e-01        & 0.1367
                      & 1.57975e-01        & 0.7562
                      & 5.38619e-01        & 0.3960\\
800                   & 1.48193e-01        & 0.7079
                      & 5.49331e-01        & -0.1466
                      & 5.99206e-02        & 1.3986
                      & 4.94470e-01        & 0.1234\\
\hline
\space    &\multicolumn{4}{l}{\cellcolor{gray!35}{WENO-ZA}}
          &\multicolumn{4}{l}{\cellcolor{gray!35}{MOP-GMWENO-ZA}}\\
\cline{2-5}  \cline{6-9}
$N$                   & $L_{1}$ error      & $L_{1}$ order 
                      & $L_{\infty}$ error & $L_{\infty}$ order
		  			  & $L_{1}$ error      & $L_{1}$ order 
      				  & $L_{\infty}$ error & $L_{\infty}$ order \\
\Xhline{0.65pt}
200                   & 2.27174e-01        & -
                      & 6.35754e-01        & -
                      & 3.82186e-01        & -
                      & 7.40853e-01        & - \\
400                   & 2.03221e-01        & 0.1607
                      & 5.99664e-01        & 0.0843
                      & 1.81782e-01        & 1.0721
                      & 4.90704e-01        & 0.5943\\
800                   & 1.63892e-01        & 0.3103
                      & 5.32228e-01        & 0.1721
                      & 6.35776e-02        & 1.5156
                      & 4.94734e-01        & -0.0118\\
\hline
\space    &\multicolumn{4}{l}{\cellcolor{gray!35}{WENO-D}}
          &\multicolumn{4}{l}{\cellcolor{gray!35}{MOP-GMWENO-D}}\\
\cline{2-5}  \cline{6-9}
$N$                   & $L_{1}$ error      & $L_{1}$ order 
                      & $L_{\infty}$ error & $L_{\infty}$ order
		  			  & $L_{1}$ error      & $L_{1}$ order 
      				  & $L_{\infty}$ error & $L_{\infty}$ order \\
\Xhline{0.65pt}
200                   & 3.86995e-01        & -
                      & 6.85835e-01        & -
                      & 4.52027e-01        & -
                      & 7.69911e-01        & - \\
400                   & 2.02287e-01        & 0.9359
                      & 5.18995e-01        & 0.4021
                      & 1.79447e-01        & 1.3329
                      & 5.13353e-01        & 0.5947\\
800                   & 1.66552e-01        & 0.2804
                      & 5.04564e-01        & 0.0407
                      & 6.28595e-02        & 1.5134
                      & 4.94367e-01        & 0.0544\\
\hline
\space    &\multicolumn{4}{l}{\cellcolor{gray!35}{WENO-A}}
          &\multicolumn{4}{l}{\cellcolor{gray!35}{MOP-GMWENO-A}}\\
\cline{2-5}  \cline{6-9}
$N$                   & $L_{1}$ error      & $L_{1}$ order 
                      & $L_{\infty}$ error & $L_{\infty}$ order
		  			  & $L_{1}$ error      & $L_{1}$ order 
      				  & $L_{\infty}$ error & $L_{\infty}$ order \\
\Xhline{0.65pt}
200                   & 5.31200e-01        & -
                      & 7.70910e-01        & -
                      & 4.53035e-01        & -
                      & 7.77065e-01        & - \\
400                   & 4.08352e-01        & 0.3794
                      & 6.93282e-01        & 0.1531
                      & 1.73916e-01        & 1.3812
                      & 5.29427e-01        & 0.5536\\
800                   & 2.95123e-01        & 0.4685
                      & 5.90637e-01        & 0.2312
                      & 6.35906e-02        & 1.4515
                      & 4.87979e-01        & 0.1176\\
\hline
\space    &\multicolumn{4}{l}{\cellcolor{gray!35}{WENO-NIP}}
          &\multicolumn{4}{l}{\cellcolor{gray!35}{MOP-GMWENO-NIP}}\\
\cline{2-5}  \cline{6-9}
$N$                   & $L_{1}$ error      & $L_{1}$ order 
                      & $L_{\infty}$ error & $L_{\infty}$ order
		  			  & $L_{1}$ error      & $L_{1}$ order 
      				  & $L_{\infty}$ error & $L_{\infty}$ order \\
\Xhline{0.65pt}
200                   & 2.40800e-01        & -
                      & 5.64249e-01        & -
                      & 4.37570e-01        & -
                      & 7.64352e-01        & - \\
400                   & 1.34966e-01        & 0.8352
                      & 5.68093e-01        & -0.0098
                      & 1.58619e-01        & 1.4639
                      & 5.19438e-01        & 0.5573\\
800                   & 6.62377e-02        & 1.0269
                      & 5.18541e-01        & 0.1317
                      & 5.80606e-02        & 1.4499
                      & 4.84905e-01        & 0.0992\\
\hline
\end{tabular*}
\end{myFontSize}
\end{table}

\section{Numerical results}
\label{NumericalExperiments}
\subsection{Accuracy test}
To test the accuracy of different WENO schemes, we consider the 
smooth density perturbation advection problem \cite{AdaWENO-Z}. It 
is governed by the one-dimensional Euler equation
\begin{equation}\left\{
\begin{array}{l}
\dfrac{\partial\mathbf{U}}{\partial t} + 
\dfrac{\partial\mathbf{F}\big(\mathbf{U} \big)}{\partial x} = 0,\\
\mathbf{U} = \left( 
\begin{array}{l}
\rho\\
\rho u\\
E
\end{array}
\right),
\mathbf
{F}\big( \mathbf{U} \big) = \left(
\begin{array}{c}
\rho u\\
\rho u^{2} + p\\
u(E + p)
\end{array}
\right),
\end{array}\right.
\label{eq:Euler1D}
\end{equation}
where $\rho$ is the density, $u$ the velocity in the $x$ direction, $p$ the pressure, $E$ the total energy, and $p = (\gamma - 1)\Big( E - \rho u^{2}/2 \Big), \gamma = 1.4$. In all calculations of this subsection, the CFL 
number is set to be $\Delta x^{2/3}$ to rule out the effect of the 
time advancement on the convergence order of accuracy. The 
$L_{1}$, $L_{\infty}$ error of the density are computed by
\begin{equation*}
\displaystyle
\begin{array}{l}
L_{1} = h \cdot \displaystyle\sum\limits_{j} \big\lvert \rho_{j}^{
\mathrm{exact}} - (\rho_{h})_{j} \big\rvert, \quad
L_{\infty} = \displaystyle\max_{j} \big\lvert \rho_{j}^{
\mathrm{exact}} - (\rho_{h})_{j} \big\rvert,
\end{array}
\end{equation*}
where $h$ is the mesh size and $h = \Delta x$. $(\rho_{h})_{j}$ is 
the computing results and $\rho_{j}^{\mathrm{exact}}$ is the exact 
solution.

The initial conditions of this problem in our tests are given in 
Example \ref{ex:AccuracyTest:01} and Example \ref{ex:AccuracyTest:02}
below.

\begin{example}
\rm{The first accuracy test is defined by} 
\label{ex:AccuracyTest:01}
\end{example} 
\begin{equation}
\rho(x, 0) = 1 + 0.2\sin ( \pi x ), 
\quad u(x, 0) = 1, 
\quad p(x, 0) = 1.
\label{eq:Euler:IC1}
\end{equation}
The exact solution is computed by
\begin{equation*}
\rho(x, t) = 1 + 0.2\sin \big( \pi (x-ut) \big), 
\quad u(x, t) = 1, 
\quad p(x, t) = 1.
\end{equation*}
The computational domain is $[0, 2]$ and the output time is taken to 
be $t = 2.0$.

The numerical errors of the density and the convergence rorders of different WENO 
schemes are presented in Table \ref{tab:AccuracyTest:01}. For 
comparison purpose, we also present the results computed by the 
WENO5-ILW scheme. As expected, all the considered WENO 
schemes can achieve the designed convergence orders. In addition, 
the $L_{1}$ and $L_{\infty}$ errors of all considered schemes are 
very close to those of the WENO-ILW scheme for this smooth problem. 
Here, we point out that the MOP-GMWENO-X scheme and the 
WENO-X scheme are essentially the same method for 
smooth solution in the cases whose grid numbers are not too small (
e.g., $N \geq 20$ here). 
Actually, the minor accuracy loss of the MOP-GMWENO-X schemes that 
only occurs with very small grid numbers and will disappear 
immediately when the grid number increases very slightly is 
negligible and this has been discussed carefully and detailly in 
\cite{MOP-WENO-ACMk}. Furthermore, the MOP-GMWENO-X schemes 
have solutions with better accuracy than WENO-JS in general.

\begin{table}[ht]
\begin{myFontSize}
\centering
\caption{The $L_{1}$, $L_{\infty}$ errors of the density and 
the convergence rates of accuracy of different WENO 
schemes for Example \ref{ex:AccuracyTest:01}.}
\label{tab:AccuracyTest:01}
\begin{tabular*}{\hsize}
{@{}@{\extracolsep{\fill}}cllllllll@{}}
\hline
\space    &\multicolumn{4}{l}{\cellcolor{gray!35}{WENO5-ILW}}  
          &\multicolumn{4}{l}{\cellcolor{gray!35}{WENO-JS}}  \\
\cline{2-5}  \cline{6-9}
$N$                   & $L_{1}$ error      & $L_{1}$ order 
                      & $L_{\infty}$ error & $L_{\infty}$ order
		  			  & $L_{1}$ error      & $L_{1}$ order 
      				  & $L_{\infty}$ error & $L_{\infty}$ order \\
\Xhline{0.65pt}
10                    & 5.57925e-03        & -
                      & 4.31021e-03        & -
                      & 2.46709e-02        & -
                      & 1.76565e-02        & -\\
20                    & 1.87285e-04        & 4.8968
                      & 1.46004e-04        & 4.8837
                      & 1.30211e-03        & 4.2439
                      & 1.00423e-03        & 4.1360\\
40                    & 5.94483e-06        & 4.9775
                      & 4.66039e-06        & 4.9694
                      & 4.15463e-05        & 4.9700
                      & 3.74857e-05        & 4.7436\\
80                    & 1.86489e-07        & 4.9945
                      & 1.46400e-07        & 4.9925
                      & 1.29861e-06        & 4.9997
                      & 1.20927e-06        & 4.9541\\
160                   & 5.83322e-09        & 4.9987
                      & 4.58129e-09        & 4.9980
                      & 4.05543e-08        & 5.0010
                      & 3.81502e-08        & 4.9863\\
320                   & 1.81943e-10        & 5.0027
                      & 1.44265e-10        & 4.9890
                      & 1.26670e-09        & 5.0007
                      & 1.17287e-09        & 5.0236\\
\hline
\space    &\multicolumn{4}{l}{\cellcolor{gray!35}{WENO-Z}}
          &\multicolumn{4}{l}{\cellcolor{gray!35}{MOP-GMWENO-Z}}\\
\cline{2-5}  \cline{6-9}
$N$                   & $L_{1}$ error      & $L_{1}$ order 
                      & $L_{\infty}$ error & $L_{\infty}$ order
		  			  & $L_{1}$ error      & $L_{1}$ order 
      				  & $L_{\infty}$ error & $L_{\infty}$ order \\
\Xhline{0.65pt}
10                    & 6.08762e-03        & -
                      & 4.50504e-03        & -
                      & 1.35332e-02        & -
                      & 1.00260e-02        & -\\
20                    & 1.88487e-04        & 5.0133
                      & 1.47175e-04        & 4.9359
                      & 1.88487e-04        & 6.1659
                      & 1.47175e-04        & 6.0901\\
40                    & 5.94480e-06        & 4.9867
                      & 4.66387e-06        & 4.9799
                      & 5.94480e-06        & 4.9867
                      & 4.66387e-06        & 4.9799\\
80                    & 1.86489e-07        & 4.9945
                      & 1.46410e-07        & 4.9934
                      & 1.86489e-07        & 4.9945
                      & 1.46410e-07        & 4.9934\\
160                   & 5.83323e-09        & 4.9987
                      & 4.58132e-09        & 4.9981
                      & 5.83323e-09        & 4.9987
                      & 4.58133e-09        & 4.9981\\
320                   & 1.81941e-10        & 5.0028
                      & 1.44251e-10        & 4.9891
                      & 1.81941e-10        & 5.0028
                      & 1.44251e-10        & 4.9891\\
\hline
\space    &\multicolumn{4}{l}{\cellcolor{gray!35}{WENO-Z$\eta(\tau_{5})$}} 
          &\multicolumn{4}{l}{\cellcolor{gray!35}{MOP-GMWENO-Z$\eta(\tau_{5})$}}\\
\cline{2-5}  \cline{6-9}
$N$                   & $L_{1}$ error      & $L_{1}$ order 
                      & $L_{\infty}$ error & $L_{\infty}$ order
		  			  & $L_{1}$ error      & $L_{1}$ order 
      				  & $L_{\infty}$ error & $L_{\infty}$ order \\
\Xhline{0.65pt}
10                    & 6.02580e-03        & -
                      & 4.51422e-03        & -
                      & 1.34684e-02        & -
                      & 9.98726e-03        & -\\
20                    & 1.88232e-04        & 5.0006
                      & 1.47368e-04        & 4.9370
                      & 1.88232e-04        & 6.1609
                      & 1.47368e-04        & 6.0826\\
40                    & 5.94481e-06        & 4.9847
                      & 4.66425e-06        & 4.9816
                      & 5.94481e-06        & 4.9847
                      & 4.66425e-06        & 4.9816\\
80                    & 1.86489e-07        & 4.9945
                      & 1.46410e-07        & 4.9936
                      & 1.86489e-07        & 4.9945
                      & 1.46410e-07        & 4.9936\\
160                   & 5.83323e-09        & 4.9987
                      & 4.58133e-09        & 4.9981
                      & 5.83323e-09        & 4.9987
                      & 4.58133e-09        & 4.9981\\
320                   & 1.81938e-10        & 5.0028
                      & 1.44269e-10        & 4.9889
                      & 1.81946e-10        & 5.0027
                      & 1.44260e-10        & 4.9890\\
\hline
\space    &\multicolumn{4}{l}{\cellcolor{gray!35}{WENO-Z$\eta(\tau_{81})$}} 
          &\multicolumn{4}{l}{\cellcolor{gray!35}{MOP-GMWENO-Z$\eta(\tau_{81})$}}\\
\cline{2-5}  \cline{6-9}
$N$                   & $L_{1}$ error      & $L_{1}$ order 
                      & $L_{\infty}$ error & $L_{\infty}$ order
		  			  & $L_{1}$ error      & $L_{1}$ order 
      				  & $L_{\infty}$ error & $L_{\infty}$ order \\
\Xhline{0.65pt}
10                    & 5.57331e-03        & -
                      & 4.29858e-03        & -
                      & 1.23782e-02        & -
                      & 9.20236e-03        & -\\
20                    & 1.87285e-04        & 4.8952
                      & 1.46003e-04        & 4.8798
                      & 1.87285e-04        & 6.0464
                      & 1.46003e-04        & 5.9779\\
40                    & 5.94483e-06        & 4.9775
                      & 4.66039e-06        & 4.9694
                      & 5.94483e-06        & 4.9775
                      & 4.66039e-06        & 4.9694\\
80                    & 1.86489e-07        & 4.9945
                      & 1.46400e-07        & 4.9925
                      & 1.86489e-07        & 4.9945
                      & 1.46400e-07        & 4.9925\\
160                   & 5.83323e-09        & 4.9987
                      & 4.58130e-09        & 4.9980
                      & 5.83323e-09        & 4.9987
                      & 4.58130e-09        & 4.9980\\
320                   & 1.81942e-10        & 5.0027
                      & 1.44275e-10        & 4.9889
                      & 1.81944e-10        & 5.0027
                      & 1.44261e-10        & 4.9890\\
\hline
\space    &\multicolumn{4}{l}{\cellcolor{gray!35}{WENO-Z+}} 
          &\multicolumn{4}{l}{\cellcolor{gray!35}{MOP-GMWENO-Z+}}\\
\cline{2-5}  \cline{6-9}
$N$                   & $L_{1}$ error      & $L_{1}$ order 
                      & $L_{\infty}$ error & $L_{\infty}$ order
		  			  & $L_{1}$ error      & $L_{1}$ order 
      				  & $L_{\infty}$ error & $L_{\infty}$ order \\
\Xhline{0.65pt}
10                    & 5.91748e-03        & -
                      & 3.89100e-03        & -
                      & 4.61518e-03        & -
                      & 3.04605e-03        & -\\
20                    & 4.31540e-04        & 3.7774
                      & 3.83185e-04        & 3.3440
                      & 4.31540e-04        & 3.4188
                      & 3.83185e-04        & 2.9908\\
40                    & 1.28245e-05        & 5.0725
                      & 1.21272e-05        & 4.9817
                      & 1.28245e-05        & 5.0725
                      & 1.21272e-05        & 4.9817\\
80                    & 3.78662e-07        & 5.0818
                      & 3.90521e-07        & 4.9567
                      & 3.78662e-07        & 5.0818
                      & 3.90521e-07        & 4.9567\\
160                   & 1.15394e-08        & 5.0363
                      & 1.21818e-08        & 5.0026
                      & 1.15394e-08        & 5.0363
                      & 1.21818e-08        & 5.0026\\
320                   & 3.60494e-10        & 5.0004
                      & 3.72739e-10        & 5.0304
                      & 3.60491e-10        & 5.0005
                      & 3.72757e-10        & 5.0303\\
\hline
\space    &\multicolumn{4}{l}{\cellcolor{gray!35}{WENO-ZA}} 
          &\multicolumn{4}{l}{\cellcolor{gray!35}{MOP-GMWENO-ZA}}\\
\cline{2-5}  \cline{6-9}
$N$                   & $L_{1}$ error      & $L_{1}$ order 
                      & $L_{\infty}$ error & $L_{\infty}$ order
		  			  & $L_{1}$ error      & $L_{1}$ order 
      				  & $L_{\infty}$ error & $L_{\infty}$ order \\
\Xhline{0.65pt}
10                    & 5.78359e-03        & -
                      & 4.46507e-03        & -
                      & 1.26046e-02        & -
                      & 9.33088e-03        & -\\
20                    & 1.87474e-04        & 4.9472
                      & 1.46141e-04        & 4.9332
                      & 1.87474e-04        & 6.0711
                      & 1.46141e-04        & 5.9966\\
40                    & 5.94486e-06        & 4.9789
                      & 4.66039e-06        & 4.9708
                      & 5.94486e-06        & 4.9789
                      & 4.66039e-06        & 4.9708\\
80                    & 1.86489e-07        & 4.9945
                      & 1.46400e-07        & 4.9925
                      & 1.86489e-07        & 4.9945
                      & 1.46400e-07        & 4.9925\\
160                   & 5.83323e-09        & 4.9987
                      & 4.58130e-09        & 4.9980
                      & 5.83322e-09        & 4.9987
                      & 4.58128e-09        & 4.9980\\
320                   & 1.81941e-10        & 5.0028
                      & 1.44262e-10        & 4.9890
                      & 1.81942e-10        & 5.0027
                      & 1.44257e-10        & 4.9890\\
\hline
\space    &\multicolumn{4}{l}{\cellcolor{gray!35}{WENO-D}} 
          &\multicolumn{4}{l}{\cellcolor{gray!35}{MOP-GMWENO-D}}\\
\cline{2-5}  \cline{6-9}
$N$                   & $L_{1}$ error      & $L_{1}$ order 
                      & $L_{\infty}$ error & $L_{\infty}$ order
		  			  & $L_{1}$ error      & $L_{1}$ order 
      				  & $L_{\infty}$ error & $L_{\infty}$ order \\
\Xhline{0.65pt}
10                    & 6.08762e-03        & -
                      & 4.50504e-03        & -
                      & 1.35332e-02        & -
                      & 1.00260e-02        & -\\
20                    & 1.88487e-04        & 5.0133
                      & 1.47175e-04        & 4.9359
                      & 1.88487e-04        & 6.1659
                      & 1.47175e-04        & 6.0901\\
40                    & 5.94480e-06        & 4.9867
                      & 4.66387e-06        & 4.9799
                      & 5.94480e-06        & 4.9867
                      & 4.66387e-06        & 4.9799\\
80                    & 1.86489e-07        & 4.9945
                      & 1.46410e-07        & 4.9934
                      & 1.86489e-07        & 4.9945
                      & 1.46410e-07        & 4.9934\\
160                   & 5.83323e-09        & 4.9987
                      & 4.58132e-09        & 4.9981
                      & 5.83323e-09        & 4.9987
                      & 4.58133e-09        & 4.9981\\
320                   & 1.81941e-10        & 5.0028
                      & 1.44251e-10        & 4.9891
                      & 1.81941e-10        & 5.0028
                      & 1.44251e-10        & 4.9891\\
\hline
\space    &\multicolumn{4}{l}{\cellcolor{gray!35}{WENO-A}} 
          &\multicolumn{4}{l}{\cellcolor{gray!35}{MOP-GMWENO-A}}\\
\cline{2-5}  \cline{6-9}
$N$                   & $L_{1}$ error      & $L_{1}$ order 
                      & $L_{\infty}$ error & $L_{\infty}$ order
		  			  & $L_{1}$ error      & $L_{1}$ order 
      				  & $L_{\infty}$ error & $L_{\infty}$ order \\
\Xhline{0.65pt}
10                    & 5.57925e-03        & -
                      & 4.31021e-03        & -
                      & 1.23660e-02        & -
                      & 9.21706e-03        & -\\
20                    & 1.87285e-04        & 4.8968
                      & 1.46004e-04        & 4.8837
                      & 1.87285e-04        & 6.0450
                      & 1.46004e-04        & 5.9802\\
40                    & 5.94483e-06        & 4.9775
                      & 4.66039e-06        & 4.9694
                      & 5.94483e-06        & 4.9775
                      & 4.66039e-06        & 4.9694\\
80                    & 1.86489e-07        & 4.9945
                      & 1.46400e-07        & 4.9925
                      & 1.86489e-07        & 4.9945
                      & 1.46400e-07        & 4.9925\\
160                   & 5.83323e-09        & 4.9987
                      & 4.58129e-09        & 4.9980
                      & 5.83323e-09        & 4.9987
                      & 4.58130e-09        & 4.9980\\
320                   & 1.81943e-10        & 5.0027
                      & 1.44279e-10        & 4.9888
                      & 1.81941e-10        & 5.0028
                      & 1.44249e-10        & 4.9891\\
\hline
\space    &\multicolumn{4}{l}{\cellcolor{gray!35}{WENO-NIP}} 
          &\multicolumn{4}{l}{\cellcolor{gray!35}{MOP-GMWENO-NIP}}\\
\cline{2-5}  \cline{6-9}
$N$                   & $L_{1}$ error      & $L_{1}$ order 
                      & $L_{\infty}$ error & $L_{\infty}$ order
		  			  & $L_{1}$ error      & $L_{1}$ order 
      				  & $L_{\infty}$ error & $L_{\infty}$ order \\
\Xhline{0.65pt}
10                    & 1.59200e-03        & -
                      & 1.88378e-03        & -
                      & 1.04077e-02        & -
                      & 7.92518e-03        & -\\
20                    & 1.81848e-04        & 3.1300
                      & 1.34479e-04        & 3.8082
                      & 1.81848e-04        & 5.8388
                      & 1.34479e-04        & 5.8810\\
40                    & 5.94395e-06        & 4.9352
                      & 4.64675e-06        & 4.8550
                      & 5.94395e-06        & 4.9352
                      & 4.64675e-06        & 4.8550\\
80                    & 1.86489e-07        & 4.9943
                      & 1.46376e-07        & 4.9885
                      & 1.86489e-07        & 4.9943
                      & 1.46376e-07        & 4.9885\\
160                   & 5.83322e-09        & 4.9987
                      & 4.58125e-09        & 4.9978
                      & 5.83323e-09        & 4.9987
                      & 4.58125e-09        & 4.9978\\
320                   & 1.81943e-10        & 5.0027
                      & 1.44258e-10        & 4.9890
                      & 1.81943e-10        & 5.0027
                      & 1.44241e-10        & 4.9892\\
\hline
\end{tabular*}
\end{myFontSize}
\end{table}

\begin{example}
\rm{The second accuracy test is defined by} 
\label{ex:AccuracyTest:02}
\end{example} 
\begin{equation}
\rho(x,0)=1+0.2\sin \bigg( \pi x - \dfrac{\sin (\pi x)}{\pi}\bigg), 
\quad u(x, 0) = 1, 
\quad p(x, 0) = 1.
\label{eq:Euler:IC2}
\end{equation}
The exact solution is computed by
\begin{equation*}
\rho(x, t) = 1+0.2\sin \bigg( \pi (x-ut) - 
\dfrac{\sin \big(\pi (x-ut)\big)}{\pi}\bigg), 
\quad u(x, t) = 1, 
\quad p(x, t) = 1.
\end{equation*}
Also, the computational domain is $[0, 2]$ and the output time is 
taken to be $t = 2.0$.

Table \ref{tab:AccuracyTest:02} shows the numerical errors of the 
density and the convergence orders of different WENO schemes. Again, 
we give the results computed by the WENO5-ILW scheme. Firstly, 
it can be observed that WENO-JS provides the lowest accurate 
results among all considered schemes. Its $L_{\infty}$ convergence 
rate of accuracy drops by nearly $2$ orders that leads to the 
noticeable accuracy loss shown with the $L_{1}$ convergence orders. However, it can be seen that the MOP-GMWENO-X and WENO-X schemes can recover the designed convergence orders even in the presence of critical points.

\begin{table}[ht]
\begin{myFontSize}
\centering
\caption{The $L_{1}$, $L_{\infty}$ errors of the density and 
the convergence rates of accuracy of different WENO 
schemes for Example \ref{ex:AccuracyTest:02}.}
\label{tab:AccuracyTest:02}
\begin{tabular*}{\hsize}
{@{}@{\extracolsep{\fill}}cllllllll@{}}
\hline
\space    &\multicolumn{4}{l}{\cellcolor{gray!35}{WENO5-ILW}}  
          &\multicolumn{4}{l}{\cellcolor{gray!35}{WENO-JS}}  \\
\cline{2-5}  \cline{6-9}
$N$                   & $L_{1}$ error      & $L_{1}$ order 
                      & $L_{\infty}$ error & $L_{\infty}$ order
		  			  & $L_{1}$ error      & $L_{1}$ order 
      				  & $L_{\infty}$ error & $L_{\infty}$ order \\
\Xhline{0.65pt}
10                    & 2.65885e-02        & -
                      & 2.42757e-02        & -
                      & 3.96420e-02        & -
                      & 4.33175e-02        & -\\
20                    & 1.85214e-03        & 3.8435
                      & 2.25478e-03        & 3.4285
                      & 4.12024e-03        & 3.2662
                      & 4.83249e-03        & 3.1641\\
40                    & 6.72815e-05        & 3.7828
                      & 9.26702e-05        & 4.6047
                      & 2.73458e-04        & 3.9133
                      & 3.72118e-04        & 3.6989\\
80                    & 2.15464e-06        & 4.9647
                      & 3.01580e-06        & 4.9415
                      & 1.42262e-05        & 4.2647
                      & 3.21426e-05        & 3.5332\\
160                   & 6.77495e-08        & 4.9911
                      & 9.51825e-08        & 4.9857
                      & 6.77367e-07        & 4.3925
                      & 3.02134e-06        & 3.4112\\
320                   & 2.12023e-09        & 4.9979
                      & 2.98236e-09        & 4.9962
                      & 3.38590e-08        & 4.3223
                      & 3.05139e-07        & 3.3077\\
\hline
\space    &\multicolumn{4}{l}{\cellcolor{gray!35}{WENO-Z}}
          &\multicolumn{4}{l}{\cellcolor{gray!35}{MOP-GMWENO-Z}}\\
\cline{2-5}  \cline{6-9}
$N$                   & $L_{1}$ error      & $L_{1}$ order 
                      & $L_{\infty}$ error & $L_{\infty}$ order
		  			  & $L_{1}$ error      & $L_{1}$ order 
      				  & $L_{\infty}$ error & $L_{\infty}$ order \\
\Xhline{0.65pt}
10                    & 2.04632e-02        & -
                      & 2.06144e-02        & -
                      & 3.04124e-02        & -
                      & 3.38563e-02        & -\\
20                    & 1.18306e-03        & 4.1124
                      & 1.13020e-03        & 4.1890
                      & 2.04082e-03        & 3.8974
                      & 3.07733e-03        & 3.4597\\
40                    & 5.62173e-05        & 4.3954
                      & 9.47028e-05        & 3.5770
                      & 7.54893e-05        & 4.7567
                      & 1.06191e-04        & 4.8569\\
80                    & 2.04617e-06        & 4.7800
                      & 3.01102e-06        & 4.9751
                      & 2.04617e-06        & 5.2053
                      & 3.01102e-06        & 5.1403\\
160                   & 6.68022e-08        & 4.9369
                      & 9.51811e-08        & 4.9834
                      & 6.68022e-08        & 4.9369
                      & 9.51810e-08        & 4.9834\\
320                   & 2.11297e-09        & 4.9826
                      & 2.98203e-09        & 4.9963
                      & 2.11297e-09        & 4.9826
                      & 2.98203e-09        & 4.9963\\
\hline
\space    &\multicolumn{4}{l}{\cellcolor{gray!35}{WENO-Z$\eta(\tau_{5})$}} 
          &\multicolumn{4}{l}{\cellcolor{gray!35}{MOP-GMWENO-Z$\eta(\tau_{5})$}}\\
\cline{2-5}  \cline{6-9}
$N$                   & $L_{1}$ error      & $L_{1}$ order 
                      & $L_{\infty}$ error & $L_{\infty}$ order
		  			  & $L_{1}$ error      & $L_{1}$ order 
      				  & $L_{\infty}$ error & $L_{\infty}$ order \\
\Xhline{0.65pt}
10                    & 2.09002e-02        & -
                      & 2.09558e-02        & -
                      & 3.04489e-02        & -
                      & 3.39484e-02        & -\\
20                    & 1.19359e-03        & 4.1301
                      & 1.17806e-03        & 4.1529
                      & 2.01481e-03        & 3.9177
                      & 3.04597e-03        & 3.4784\\
40                    & 5.64972e-05        & 4.4010
                      & 9.45530e-05        & 3.6391
                      & 7.53576e-05        & 4.7407
                      & 1.05810e-04        & 4.8474\\
80                    & 2.04941e-06        & 4.7849
                      & 3.01125e-06        & 4.9727
                      & 2.04941e-06        & 5.2005
                      & 3.01125e-06        & 5.1350\\
160                   & 6.68307e-08        & 4.9386
                      & 9.51809e-08        & 4.9835
                      & 6.68307e-08        & 4.9386
                      & 9.51809e-08        & 4.9835\\
320                   & 2.11313e-09        & 4.9831
                      & 2.98254e-09        & 4.9961
                      & 2.11313e-09        & 4.9831
                      & 2.98253e-09        & 4.9961\\
\hline
\space    &\multicolumn{4}{l}{\cellcolor{gray!35}{WENO-Z$\eta(\tau_{81})$}} 
          &\multicolumn{4}{l}{\cellcolor{gray!35}{MOP-GMWENO-Z$\eta(\tau_{81})$}}\\
\cline{2-5}  \cline{6-9}
$N$                   & $L_{1}$ error      & $L_{1}$ order 
                      & $L_{\infty}$ error & $L_{\infty}$ order
		  			  & $L_{1}$ error      & $L_{1}$ order 
      				  & $L_{\infty}$ error & $L_{\infty}$ order \\
\Xhline{0.65pt}
10                    & 2.62185e-02        & -
                      & 2.41549e-02        & -
                      & 2.86051e-02        & -
                      & 3.25261e-02        & -\\
20                    & 1.85177e-03        & 3.8236
                      & 2.25439e-03        & 3.4215
                      & 2.27118e-03        & 3.6548
                      & 3.22579e-03        & 3.3339\\
40                    & 6.72815e-05        & 4.7826
                      & 9.26702e-05        & 4.6045
                      & 7.70601e-05        & 4.8813
                      & 1.04861e-04        & 4.9431\\
80                    & 2.15464e-06        & 4.9647
                      & 3.01580e-06        & 4.9415
                      & 2.15464e-06        & 5.1605
                      & 3.01580e-06        & 5.1198\\
160                   & 6.77495e-08        & 4.9911
                      & 9.51825e-08        & 4.9857
                      & 6.77495e-08        & 4.9911
                      & 9.51825e-08        & 4.9857\\
320                   & 2.12024e-09        & 4.9979
                      & 2.98237e-09        & 4.9962
                      & 2.12023e-09        & 4.9979
                      & 2.98236e-09        & 4.9962\\
\hline
\space    &\multicolumn{4}{l}{\cellcolor{gray!35}{WENO-Z+}} 
          &\multicolumn{4}{l}{\cellcolor{gray!35}{MOP-GMWENO-Z+}}\\
\cline{2-5}  \cline{6-9}
$N$                   & $L_{1}$ error      & $L_{1}$ order 
                      & $L_{\infty}$ error & $L_{\infty}$ order
		  			  & $L_{1}$ error      & $L_{1}$ order 
      				  & $L_{\infty}$ error & $L_{\infty}$ order \\
\Xhline{0.65pt}
10                    & 7.73904e-03        & -
                      & 7.48121e-03        & -
                      & 1.27380e-02        & -
                      & 1.14846e-02        & -\\
20                    & 4.95302e-04        & 3.9658
                      & 5.45534e-04        & 3.7775
                      & 1.92014e-03        & 2.7299
                      & 2.50856e-03        & 2.1948\\
40                    & 5.00701e-05        & 3.3063
                      & 8.02140e-05        & 2.7657
                      & 5.50151e-05        & 5.1252
                      & 1.16909e-04        & 4.4234\\
80                    & 3.04854e-06        & 4.0378
                      & 8.56431e-06        & 3.2274
                      & 3.04854e-06        & 4.1736
                      & 8.56431e-06        & 3.7709\\
160                   & 1.77633e-07        & 4.1011
                      & 9.41566e-07        & 3.1852
                      & 1.77633e-07        & 4.1011
                      & 9.41566e-07        & 3.1852\\
320                   & 1.10103e-08        & 4.0120
                      & 1.06615e-07        & 3.1427
                      & 1.10103e-08        & 4.0120
                      & 1.06615e-07        & 3.1427\\
\hline
\space    &\multicolumn{4}{l}{\cellcolor{gray!35}{WENO-ZA}} 
          &\multicolumn{4}{l}{\cellcolor{gray!35}{MOP-GMWENO-ZA}}\\
\cline{2-5}  \cline{6-9}
$N$                   & $L_{1}$ error      & $L_{1}$ order 
                      & $L_{\infty}$ error & $L_{\infty}$ order
		  			  & $L_{1}$ error      & $L_{1}$ order 
      				  & $L_{\infty}$ error & $L_{\infty}$ order \\
\Xhline{0.65pt}
10                    & 2.62478e-02        & -
                      & 2.43179e-02        & -
                      & 2.89498e-02        & -
                      & 3.24003e-02        & -\\
20                    & 1.84035e-03        & 3.8341
                      & 2.24260e-03        & 3.4388
                      & 2.10030e-03        & 3.7849
                      & 2.99185e-03        & 3.4369\\
40                    & 6.72478e-05        & 4.7743
                      & 9.26639e-05        & 4.5970
                      & 7.70265e-05        & 4.7691
                      & 1.04867e-04        & 4.8344\\
80                    & 2.15457e-06        & 4.9640
                      & 3.01580e-06        & 4.9414
                      & 2.15457e-06        & 5.1599
                      & 3.01580e-06        & 5.1199\\
160                   & 6.77493e-08        & 4.9911
                      & 9.51825e-08        & 4.9857
                      & 6.77493e-08        & 4.9911
                      & 9.51825e-08        & 4.9857\\
320                   & 2.12024e-09        & 4.9979
                      & 2.98237e-09        & 4.9962
                      & 2.12024e-09        & 4.9979
                      & 2.98237e-09        & 4.9962\\
\hline
\space    &\multicolumn{4}{l}{\cellcolor{gray!35}{WENO-D}} 
          &\multicolumn{4}{l}{\cellcolor{gray!35}{MOP-GMWENO-D}}\\
\cline{2-5}  \cline{6-9}
$N$                   & $L_{1}$ error      & $L_{1}$ order 
                      & $L_{\infty}$ error & $L_{\infty}$ order
		  			  & $L_{1}$ error      & $L_{1}$ order 
      				  & $L_{\infty}$ error & $L_{\infty}$ order \\
\Xhline{0.65pt}
10                    & 2.04632e-02        & -
                      & 2.06144e-02        & -
                      & 3.04124e-02        & -
                      & 3.38563e-02        & -\\
20                    & 1.18306e-03        & 4.1124
                      & 1.13020e-03        & 4.1890
                      & 2.04082e-03        & 3.8974
                      & 3.07733e-03        & 3.4597\\
40                    & 5.62173e-05        & 4.3954
                      & 9.47028e-05        & 3.5770
                      & 7.54893e-05        & 4.7567
                      & 1.06191e-04        & 4.8567\\
80                    & 2.04617e-06        & 4.7800
                      & 3.01102e-06        & 4.9751
                      & 2.04617e-06        & 5.2053
                      & 3.01102e-06        & 5.1403\\
160                   & 6.68022e-08        & 4.9369
                      & 9.51811e-08        & 4.9834
                      & 6.68022e-08        & 4.9369
                      & 9.51810e-08        & 4.9834\\
320                   & 2.11297e-09        & 4.9826
                      & 2.98203e-09        & 4.9963
                      & 2.11297e-09        & 4.9826
                      & 2.98203e-09        & 4.9963\\
\hline
\space    &\multicolumn{4}{l}{\cellcolor{gray!35}{WENO-A}} 
          &\multicolumn{4}{l}{\cellcolor{gray!35}{MOP-GMWENO-A}}\\
\cline{2-5}  \cline{6-9}
$N$                   & $L_{1}$ error      & $L_{1}$ order 
                      & $L_{\infty}$ error & $L_{\infty}$ order
		  			  & $L_{1}$ error      & $L_{1}$ order 
      				  & $L_{\infty}$ error & $L_{\infty}$ order \\
\Xhline{0.65pt}
10                    & 1.83037e-02        & -
                      & 1.90477e-02        & -
                      & 2.99262e-02        & -
                      & 3.35183e-02        & -\\
20                    & 1.25386e-03        & 3.8677
                      & 1.77100e-03        & 3.4270
                      & 1.92453e-03        & 3.9588
                      & 2.85230e-03        & 3.5548\\
40                    & 6.72815e-05        & 4.2200
                      & 9.26702e-05        & 4.2563
                      & 7.70601e-05        & 4.6424
                      & 1.04861e-04        & 4.7656\\
80                    & 2.15464e-06        & 4.9647
                      & 3.01580e-06        & 4.9415
                      & 2.15464e-06        & 5.1605
                      & 3.01580e-06        & 5.1198\\
160                   & 6.77495e-08        & 4.9911
                      & 9.51825e-08        & 4.9857
                      & 6.77495e-08        & 4.9911
                      & 9.51825e-08        & 4.9857\\
320                   & 2.12024e-09        & 4.9979
                      & 2.98238e-09        & 4.9962
                      & 2.12023e-09        & 4.9979
                      & 2.98237e-09        & 4.9962\\
\hline
\space    &\multicolumn{4}{l}{\cellcolor{gray!35}{WENO-NIP}} 
          &\multicolumn{4}{l}{\cellcolor{gray!35}{MOP-GMWENO-NIP}}\\
\cline{2-5}  \cline{6-9}
$N$                   & $L_{1}$ error      & $L_{1}$ order 
                      & $L_{\infty}$ error & $L_{\infty}$ order
		  			  & $L_{1}$ error      & $L_{1}$ order 
      				  & $L_{\infty}$ error & $L_{\infty}$ order \\
\Xhline{0.65pt}
10                    & 2.33776e-02        & -
                      & 2.39846e-02        & -
                      & 2.01935e-02        & -
                      & 2.86901e-02        & -\\
20                    & 1.64174e-03        & 3.8318
                      & 2.00222e-03        & 3.5824
                      & 2.19785e-03        & 3.1997
                      & 3.03570e-03        & 3.2405\\
40                    & 6.70209e-05        & 4.6145
                      & 9.18945e-05        & 4.4455
                      & 7.68682e-05        & 4.8376
                      & 1.04304e-04        & 4.8632\\
80                    & 2.15445e-06        & 4.9592
                      & 3.01546e-06        & 4.9295
                      & 2.15445e-06        & 5.1570
                      & 3.01546e-06        & 5.1123\\
160                   & 6.77491e-08        & 4.9910
                      & 9.51826e-08        & 4.9855
                      & 6.77491e-08        & 4.9910
                      & 9.51826e-08        & 4.9855\\
320                   & 2.12024e-09        & 4.9979
                      & 2.98236e-09        & 4.9962
                      & 2.12023e-09        & 4.9979
                      & 2.98236e-09        & 4.9962\\
\hline
\end{tabular*}
\end{myFontSize}
\end{table}

\subsection{Two-dimensional Euler equations}
To examine the enhancement of the MOP-GMWENO-X schemes, we solve the 
two-dimensional Euler equations in this subsection. Its strong 
conservative form is given as
\begin{equation}\left\{
\begin{array}{l}
\dfrac{\partial\mathbf{U}}{\partial t} + 
\dfrac{\partial\mathbf{F}\big(\mathbf{U} \big)}{\partial x} +
\dfrac{\partial\mathbf{G}\big(\mathbf{U} \big)}{\partial y}= 0,\\
\mathbf{U} = \left( 
\begin{array}{l}
\rho\\
\rho u\\
\rho v\\
E
\end{array}
\right),
\mathbf
{F}\big( \mathbf{U} \big) = \left(
\begin{array}{c}
\rho u\\
\rho u^{2} + p\\
\rho uv\\
u(E + p)
\end{array}
\right),
\mathbf
{G}\big( \mathbf{U} \big) = \left(
\begin{array}{c}
\rho u\\
\rho vu\\
\rho v^{2} + p\\
v(E + p)
\end{array}
\right),
\end{array}\right.
\label{eq:Euler2D}
\end{equation}
where $\rho, u, p$ and $E$ are the same as in Eq. \eqref{eq:Euler1D} and $v$ is velocity in the $y$ direction. 
Similarly, we use the following equations of state for ideal gases to close Eq. \eqref{eq:Euler2D}
\begin{equation*}
p = (\gamma - 1)\Big( E - \dfrac{1}{2}\rho (u^{2} + v^{2}) \Big), 
\quad \gamma = 1.4.
\end{equation*}
For brevity in the discussion, we only focus on two standard 
examples. The first one is the 2D Riemann problem 
\cite{Riemann-2D-01,Riemann2D-02,Riemann2D-03}, and the second one is
the shock-vortex interaction problem 
\cite{Shock-vortex_interaction-1,Shock-vortex_interaction-2,
Shock-vortex_interaction-3}. We choose the CFL number to be $0.5$ in 
all calculations of this subsection.

\subsubsection{2D Riemann problem}
\begin{example}
\rm{We consider the 2D Riemann problem. Since the first test 2D 
Riemann problem was introduced by \cite{Riemann-2D-01}, it has 
become a widely-used two-dimensional test case for the 
high-resolution numerical methods \cite{Riemann2D-03,WENO-PPM5,
Riemann2D-04}. In our present example, the Configuration 9 of 
\cite{Riemann2D-03} is taken. It is defined on the rectangular 
computational domain $[0,1]\times[0,1]$ that is divided into four 
quadrants by lines $x = 0.5$ and $y = 0.5$, and the following 
initial constant states in each quadrant are specified}
\label{ex:Riemann2D}
\end{example}
\begin{equation*}
\mathbf{U}(x, y, 0) = \left\{
\begin{aligned}
\begin{array}{ll}
(1.0, 0.0, 0.3, 1.0)^{\mathrm{T}},         & 0.5 \leq x \leq 1.0, 
                                0.5 \leq y \leq 1.0, \\
(2.0, 0.0, -0.3, 1.0)^{\mathrm{T}},        & 0.0 \leq x \leq 0.5, 
                                0.5 \leq y \leq 1.0, \\
(1.039, 0.0, -0.8133, 0.4)^{\mathrm{T}},   & 0.0 \leq x \leq 0.5, 
                                0.0 \leq y \leq 0.5, \\
(0.5197, 0.0, -0.4259, 0.4)^{\mathrm{T}},  & 0.5 \leq x \leq 1.0, 
                                0.0 \leq y \leq 0.5. \\
\end{array}
\end{aligned}
\right.
\label{eq:initial_Euer2D:Riemann2D}
\end{equation*}
The outflow condition is employed on all edges. We 
discretize the computational domain into $800 \times 800$ cells and 
set $t = 0.3$.

\begin{figure}[ht]
\centering
  \includegraphics[height=0.432\textwidth]
  {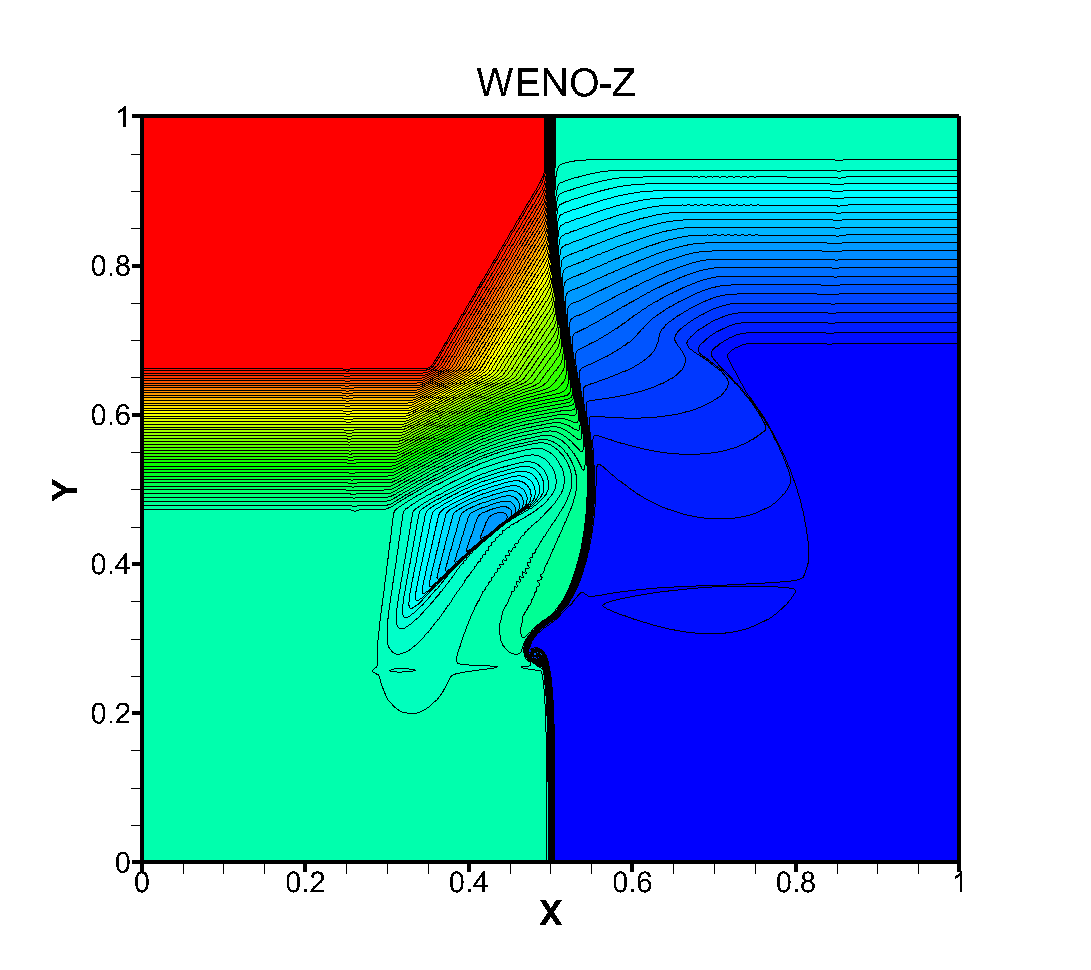}
  \includegraphics[height=0.432\textwidth]
  {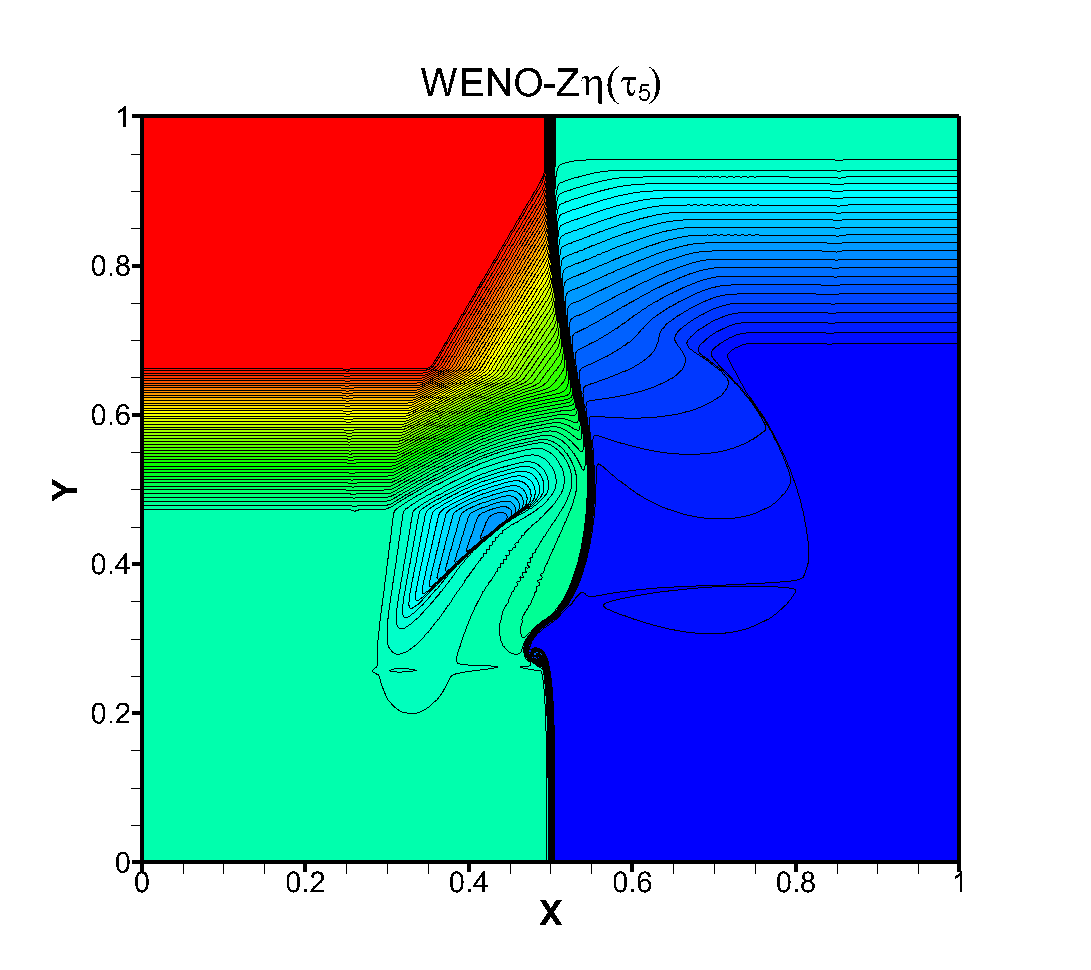}\\
  \includegraphics[height=0.432\textwidth]
  {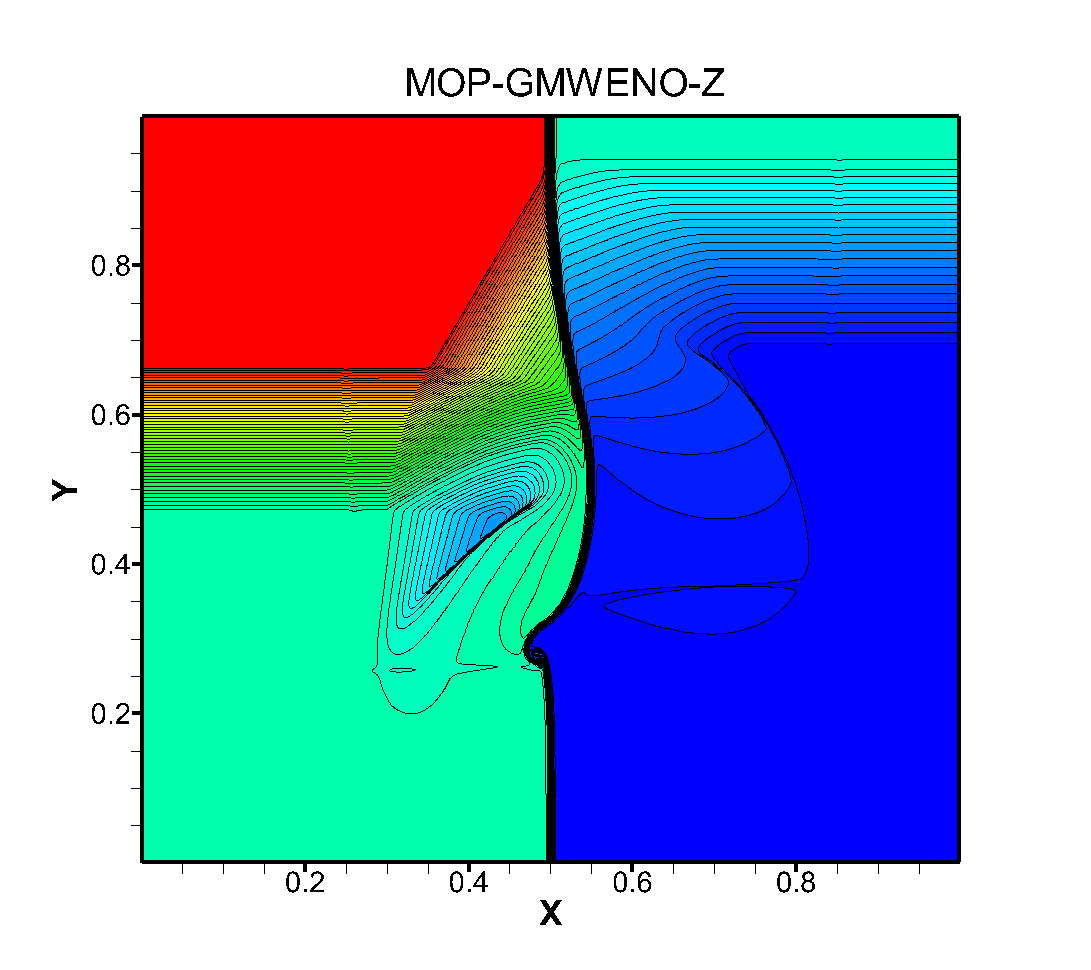}
  \includegraphics[height=0.432\textwidth]
  {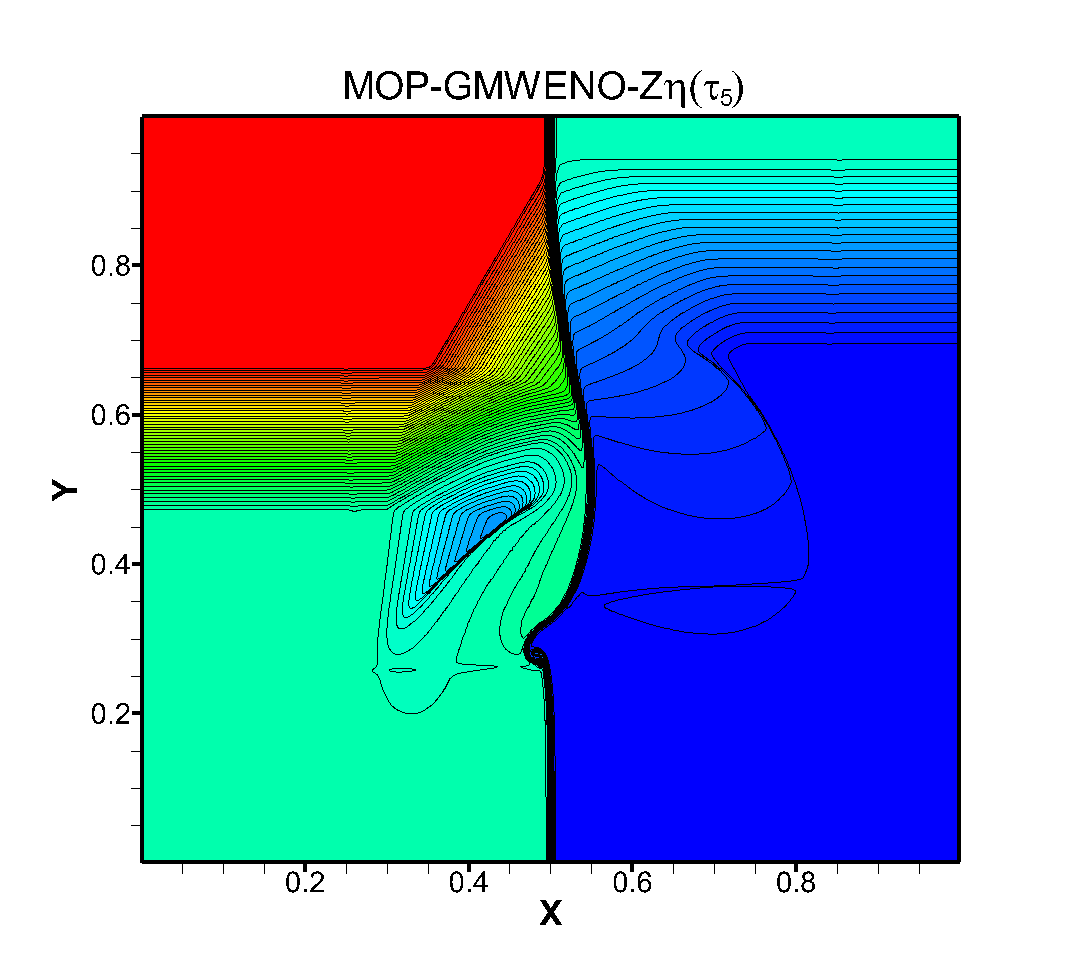}\\
  \includegraphics[height=0.375\textwidth]
  {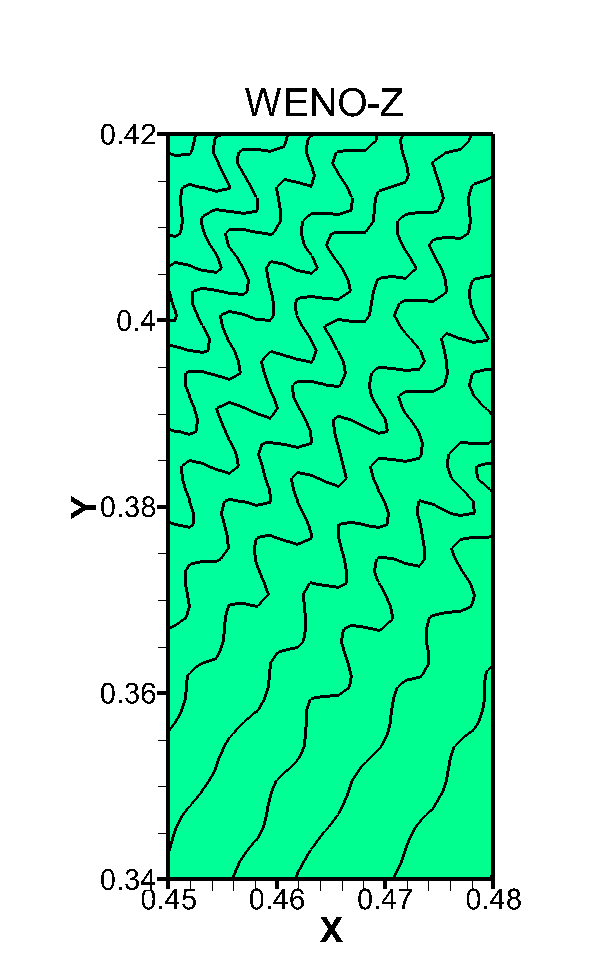} 
  \includegraphics[height=0.375\textwidth]
  {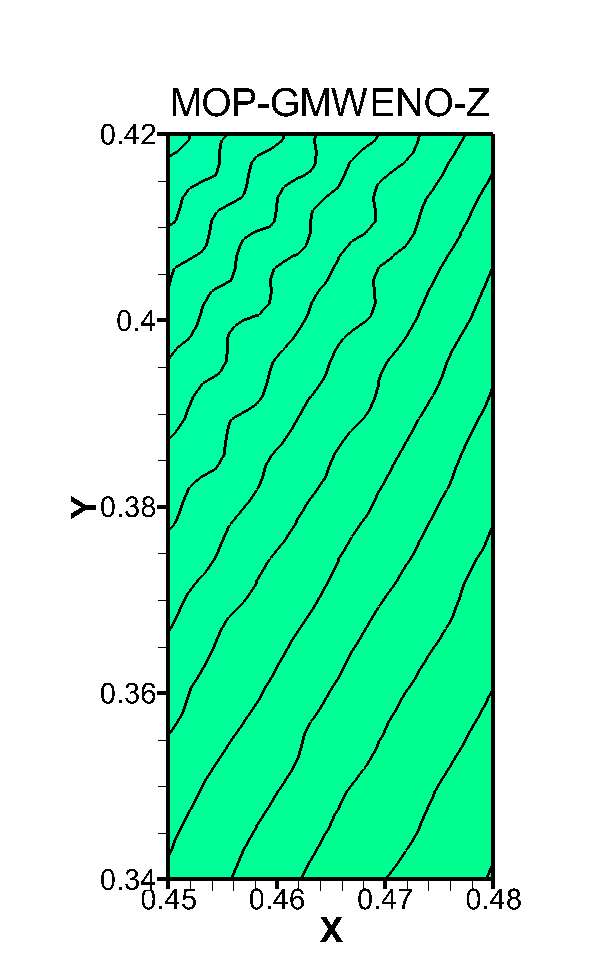} 
  \includegraphics[height=0.375\textwidth]
  {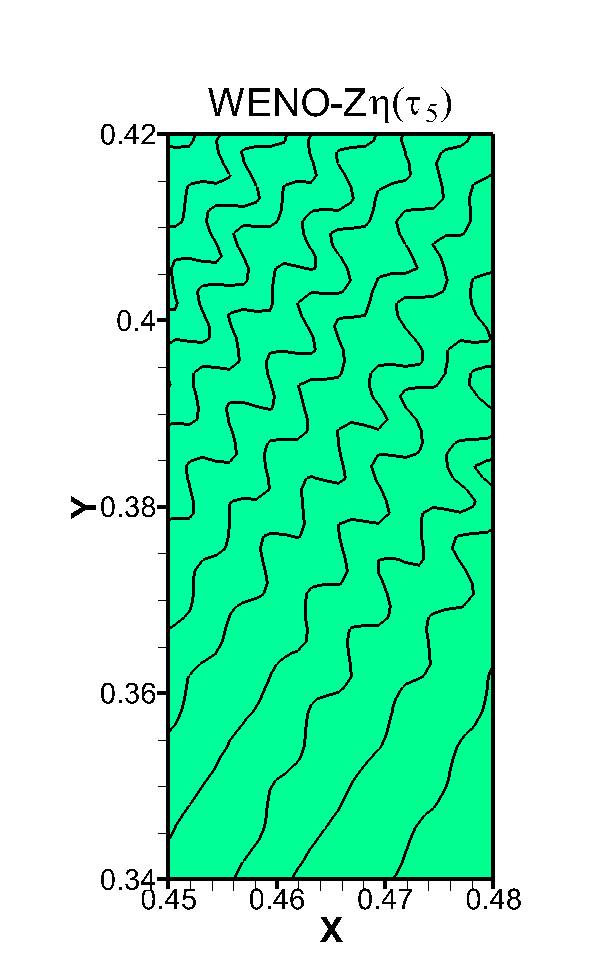} 
  \includegraphics[height=0.375\textwidth]
  {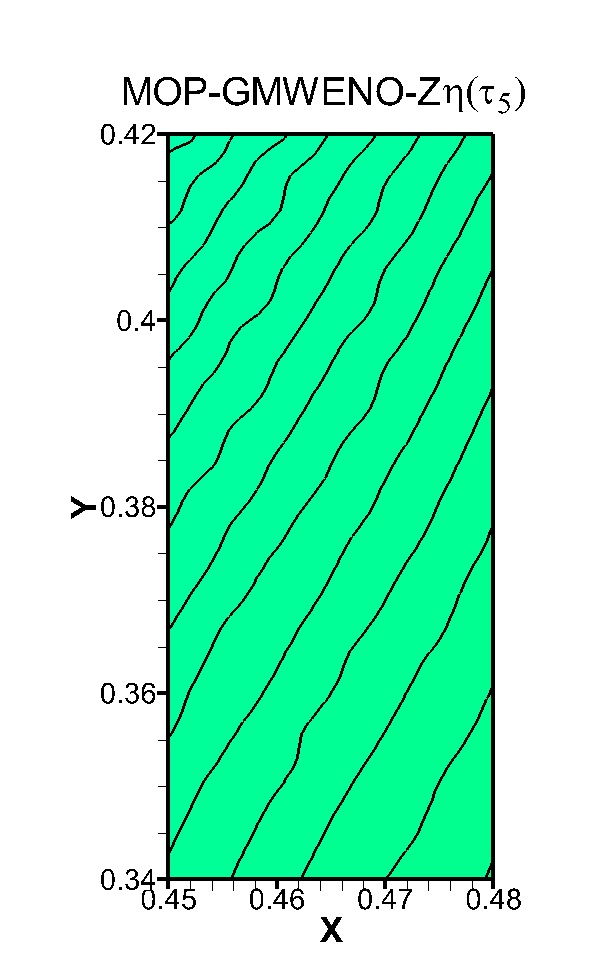}     
\caption{Density contours of Example \ref{ex:Riemann2D}. Left: WENO-Z and MOP-GMWENO-Z; Right:  
WENO-Z$\eta(\tau_{5})$ and MOP-GMWENO-Z$\eta(\tau_{5})$.}
\label{fig:ex:Riemann2D:1}
\end{figure}

\begin{figure}[ht]
\centering
  \includegraphics[height=0.432\textwidth]
  {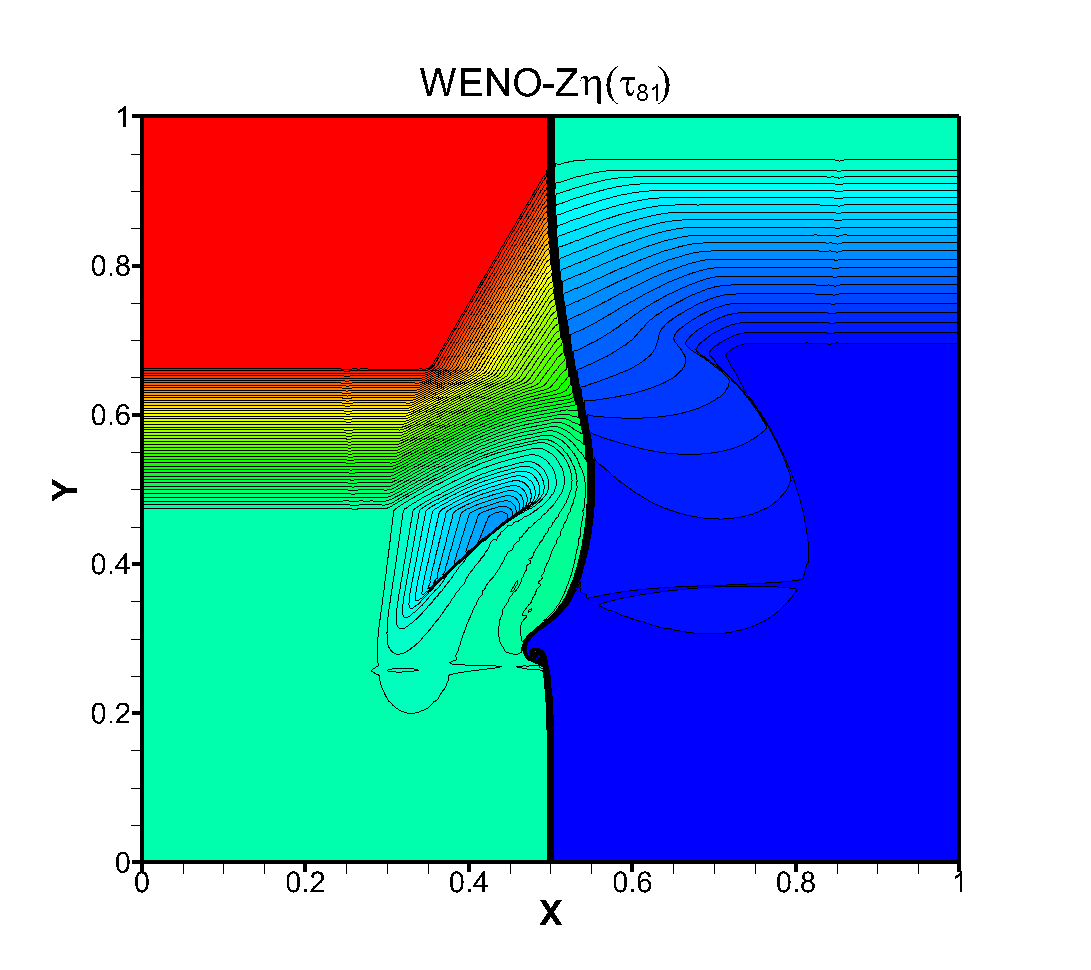}
  \includegraphics[height=0.432\textwidth]
  {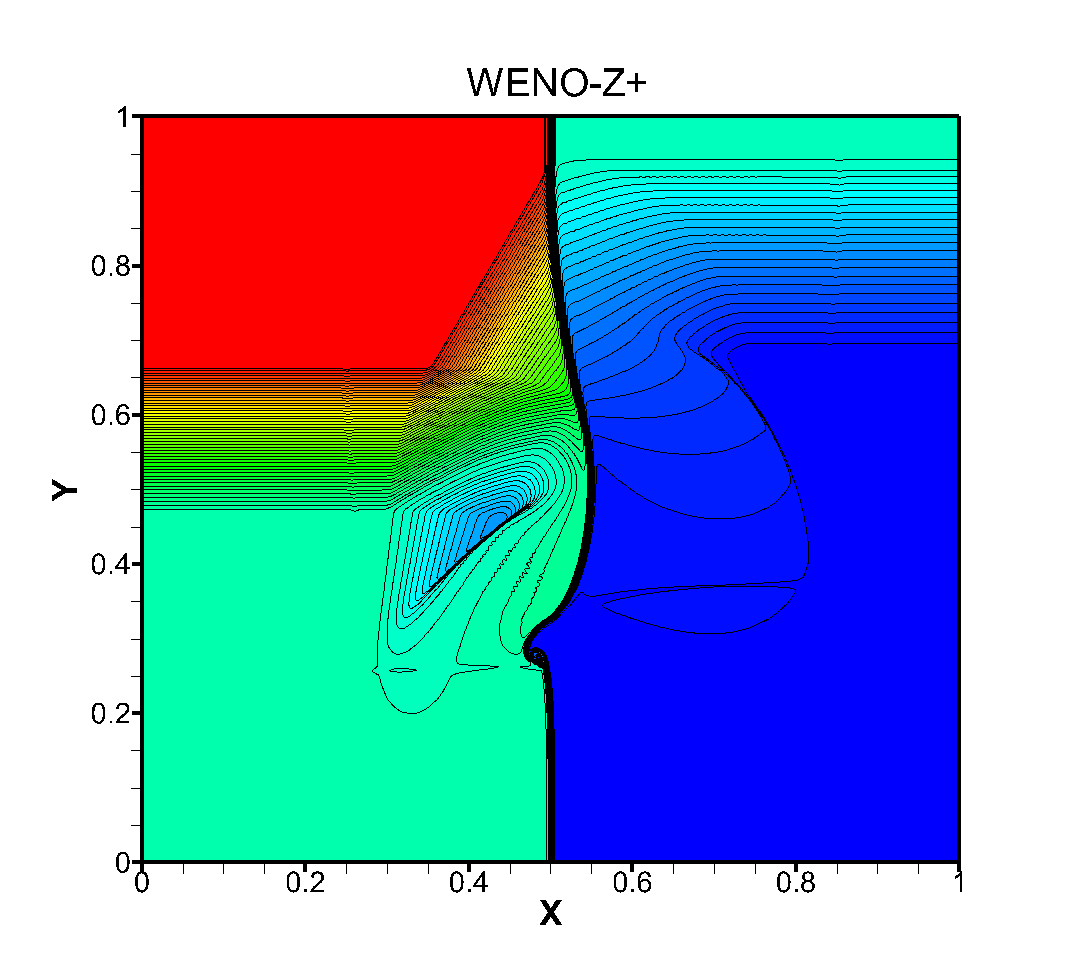}\\
  \includegraphics[height=0.432\textwidth]
  {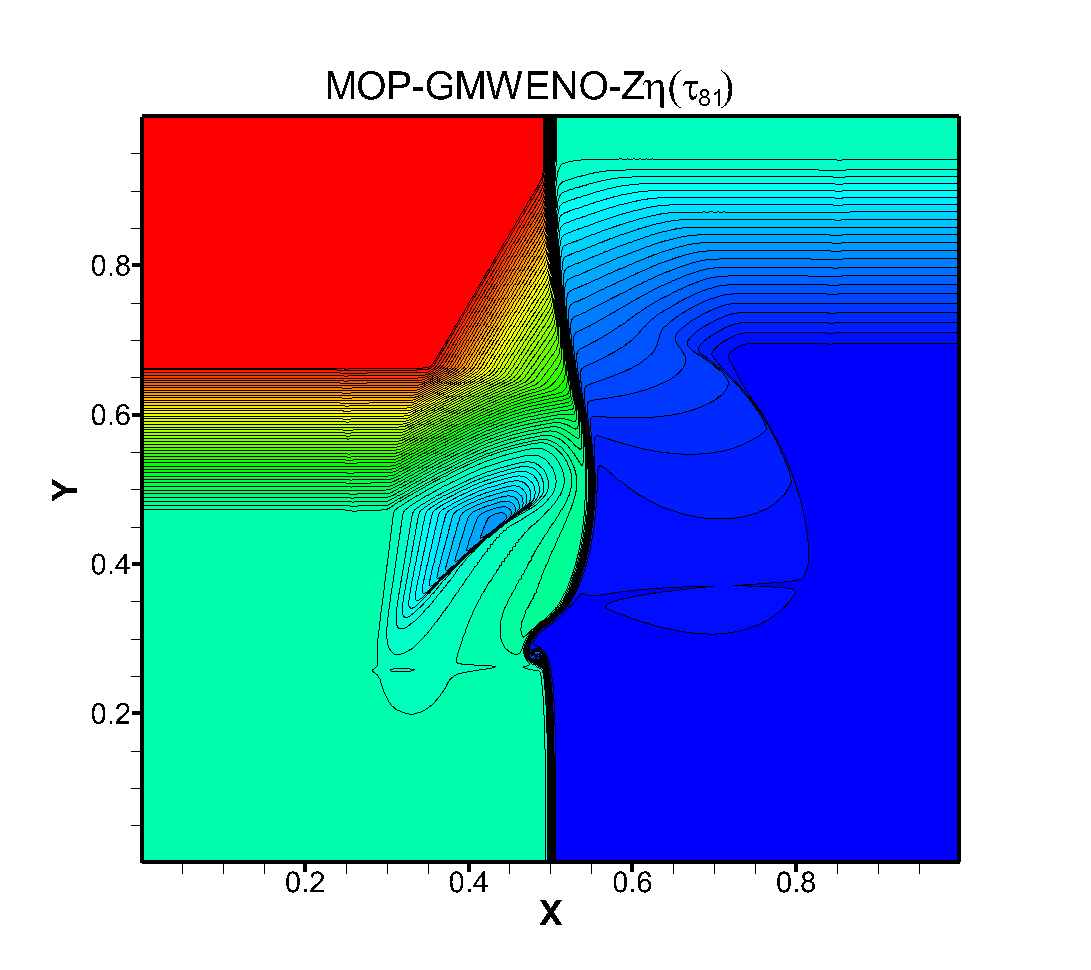}
  \includegraphics[height=0.432\textwidth]
  {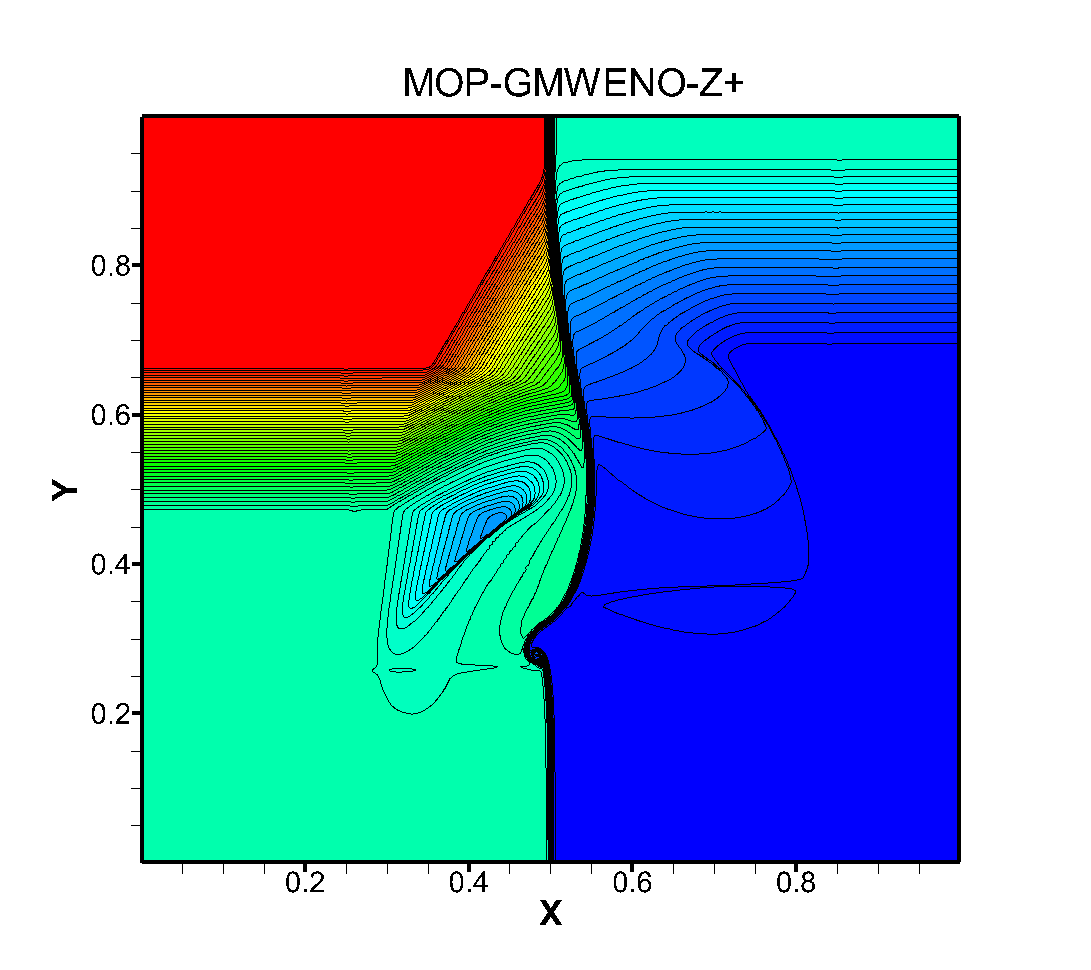}\\
  \includegraphics[height=0.375\textwidth]
  {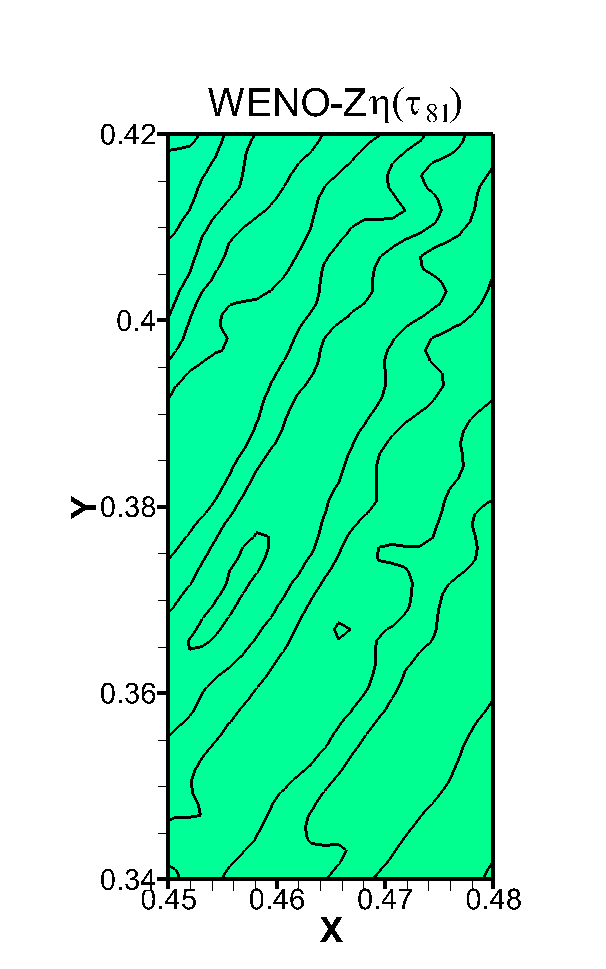} 
  \includegraphics[height=0.375\textwidth]
  {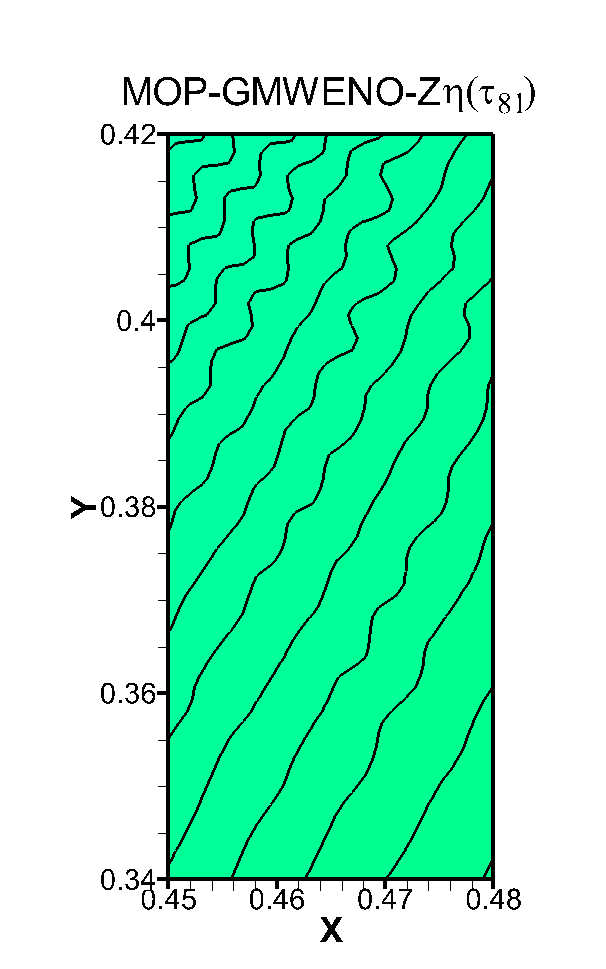}     
  \includegraphics[height=0.375\textwidth]
  {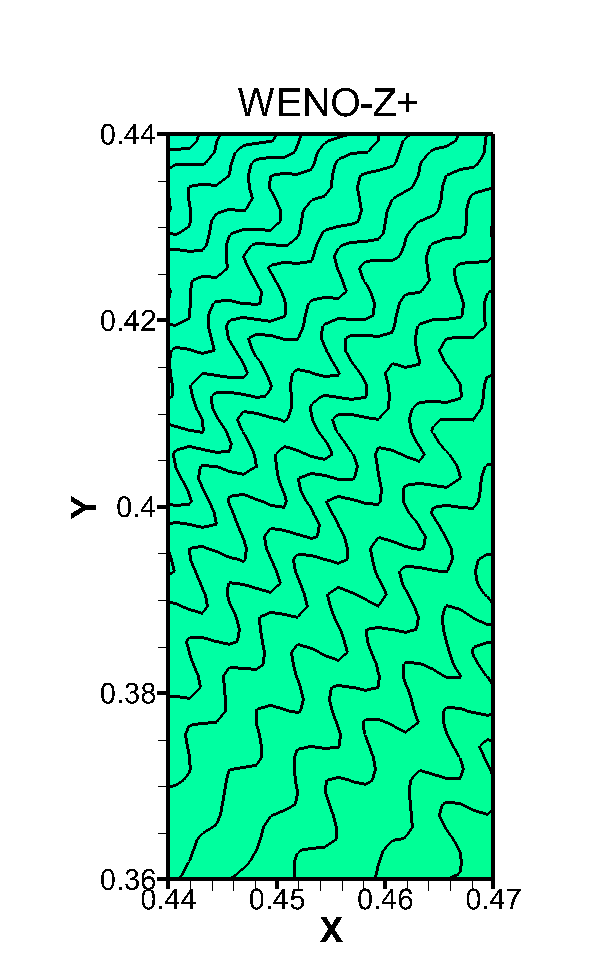} 
  \includegraphics[height=0.375\textwidth]
  {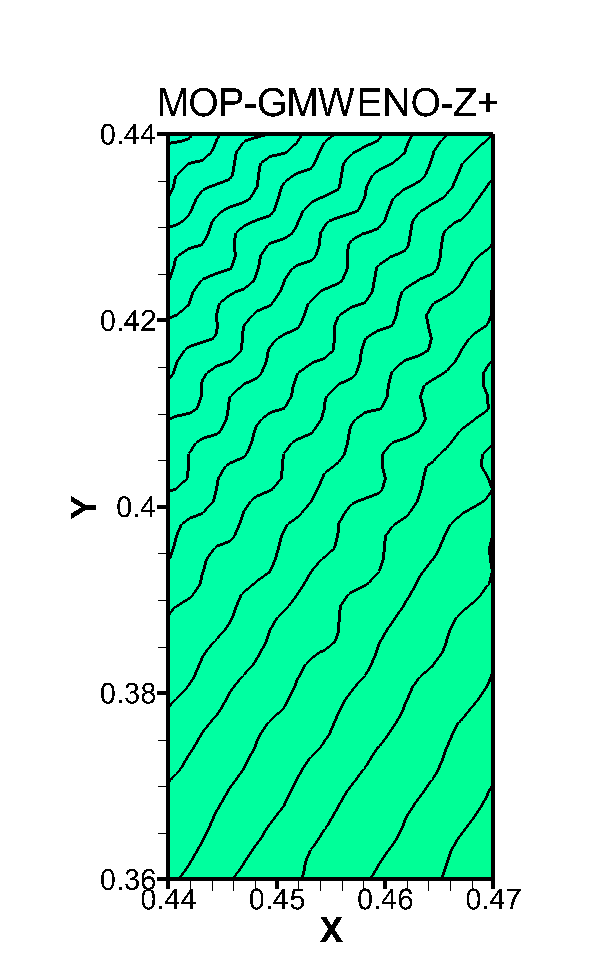} 
\caption{Density contours of Example \ref{ex:Riemann2D}. Left 
column: WENO-Z$\eta(\tau_{81})$ and MOP-GMWENO-Z$\eta(\tau_{81})$; 
Right: WENO-Z+ and MOP-GMWENO-Z+.}
\label{fig:ex:Riemann2D:2}
\end{figure}

\begin{figure}[ht]
\centering
  \includegraphics[height=0.432\textwidth]
  {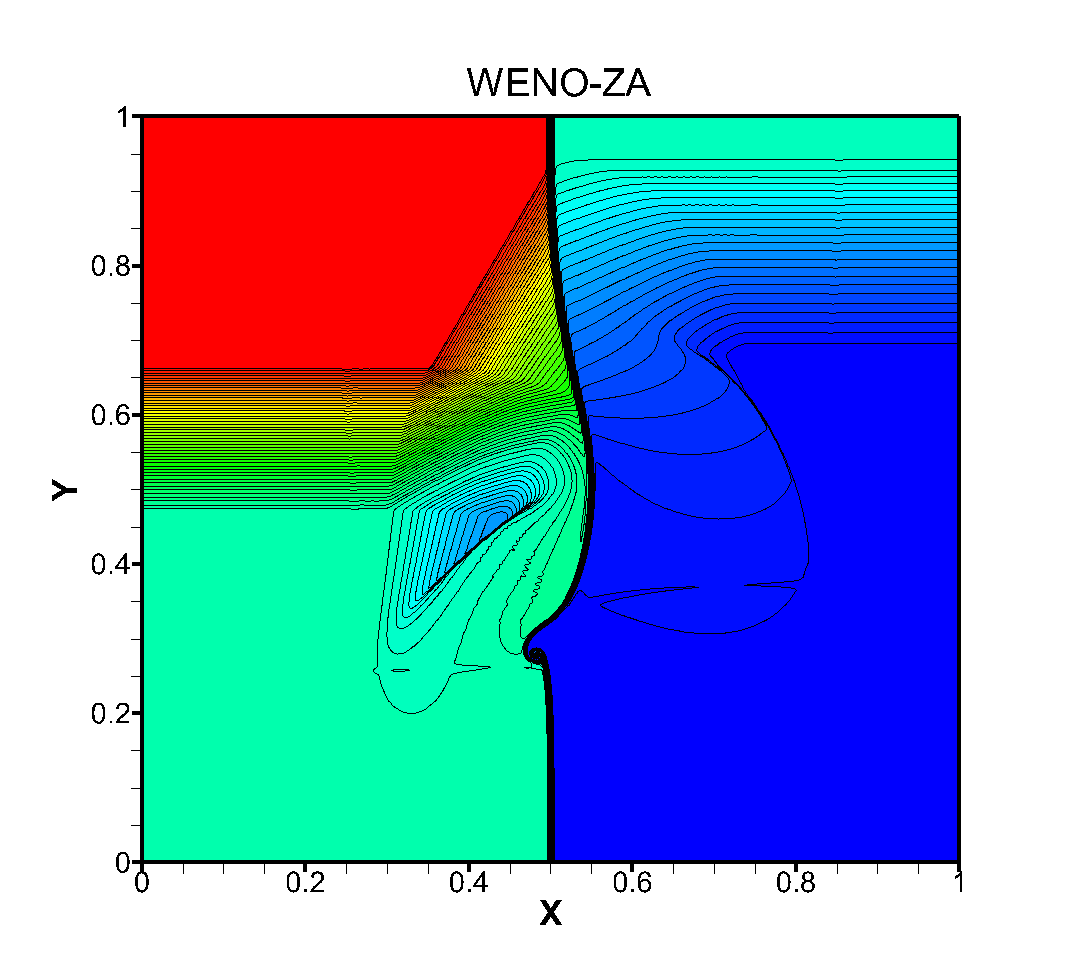}
  \includegraphics[height=0.432\textwidth]
  {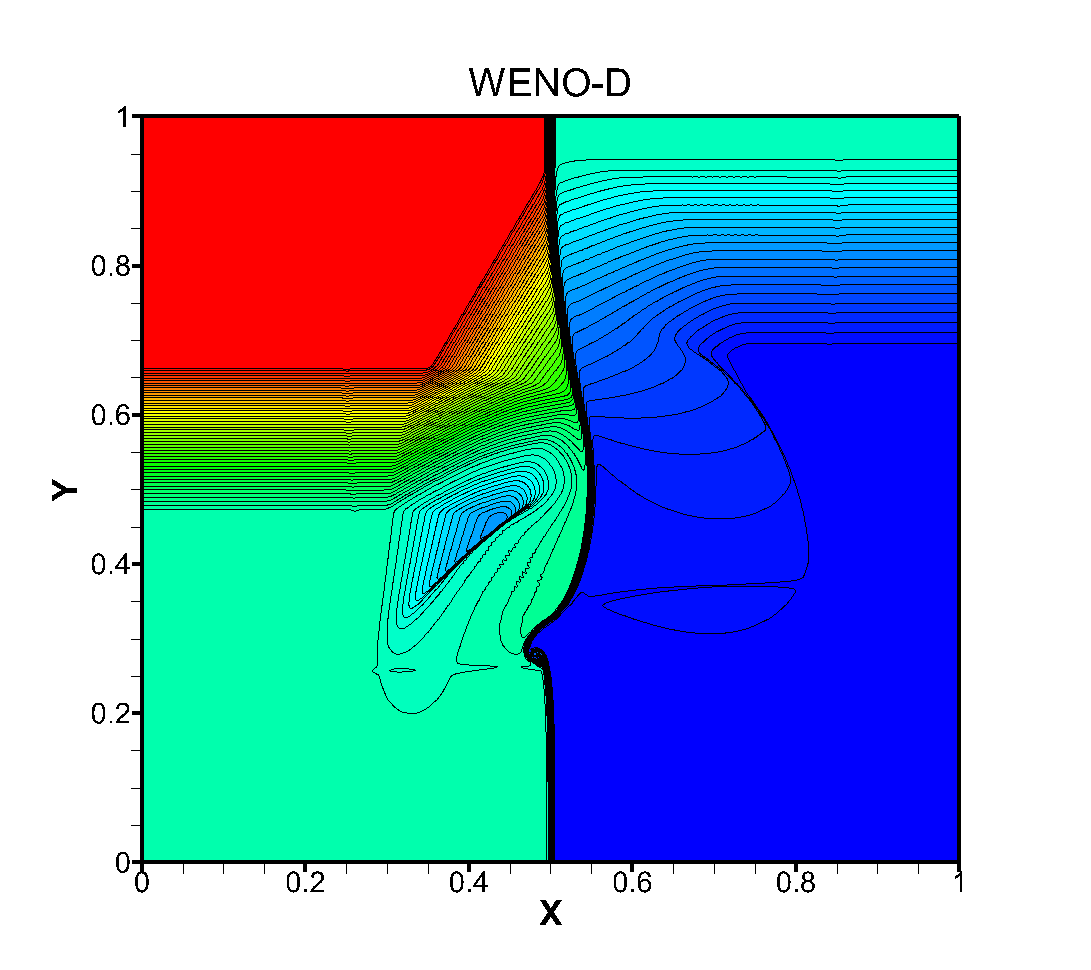}\\
  \includegraphics[height=0.432\textwidth]
  {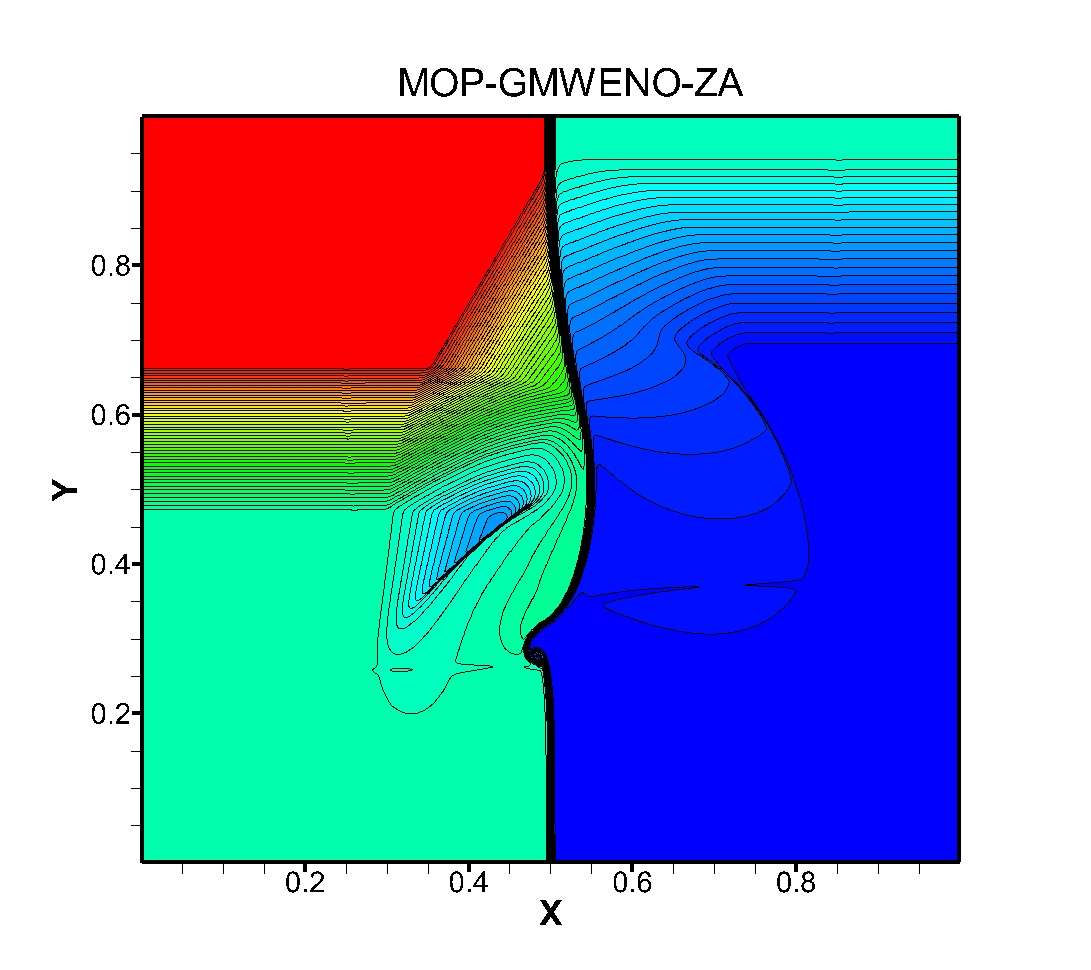}
  \includegraphics[height=0.432\textwidth]
  {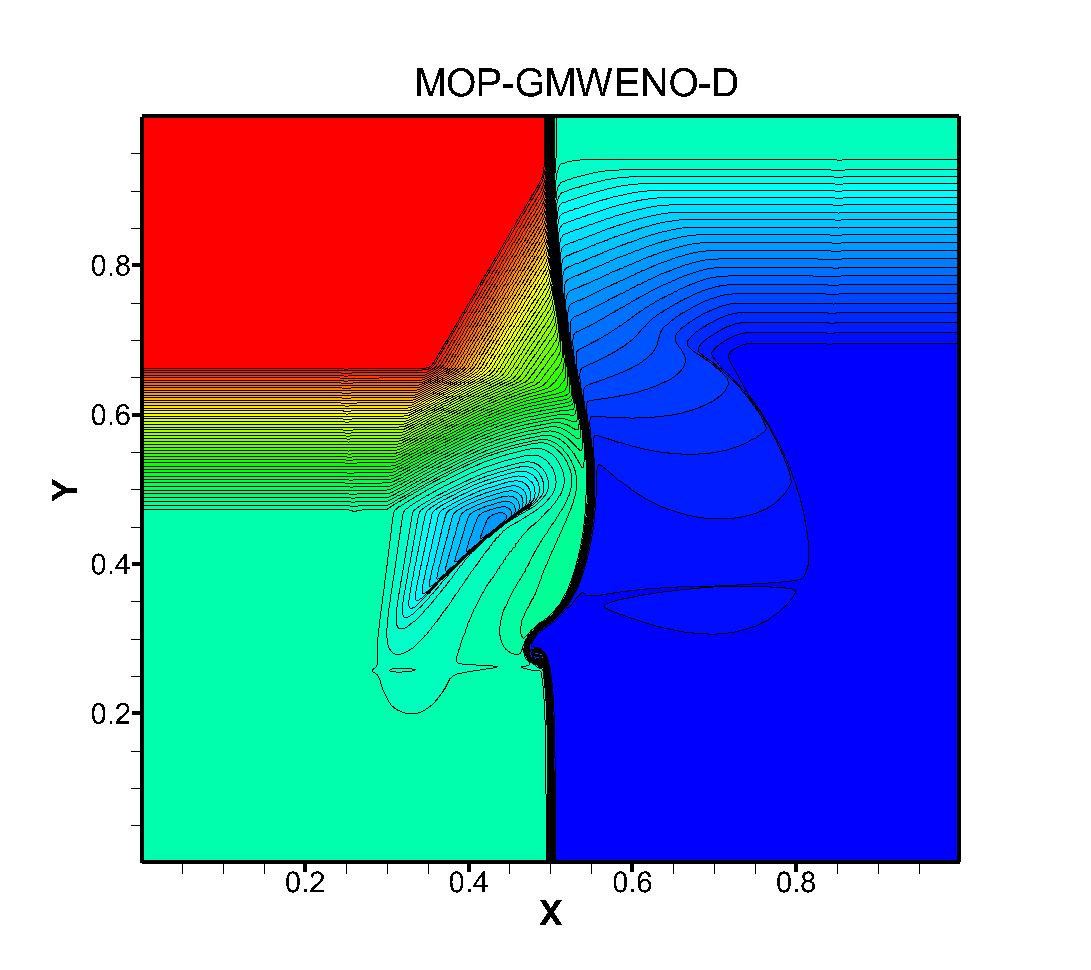}\\
  \includegraphics[height=0.375\textwidth]
  {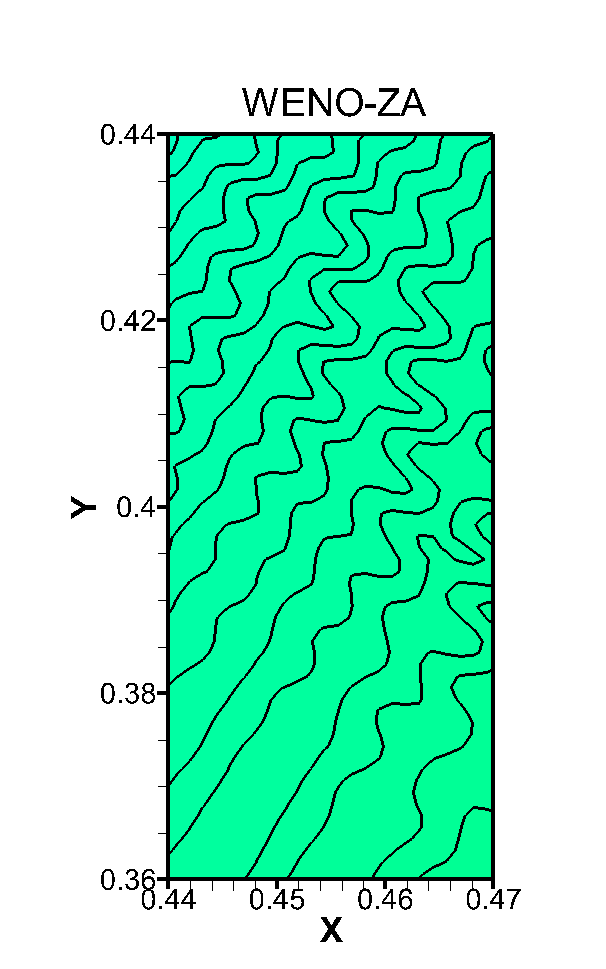} 
  \includegraphics[height=0.375\textwidth]
  {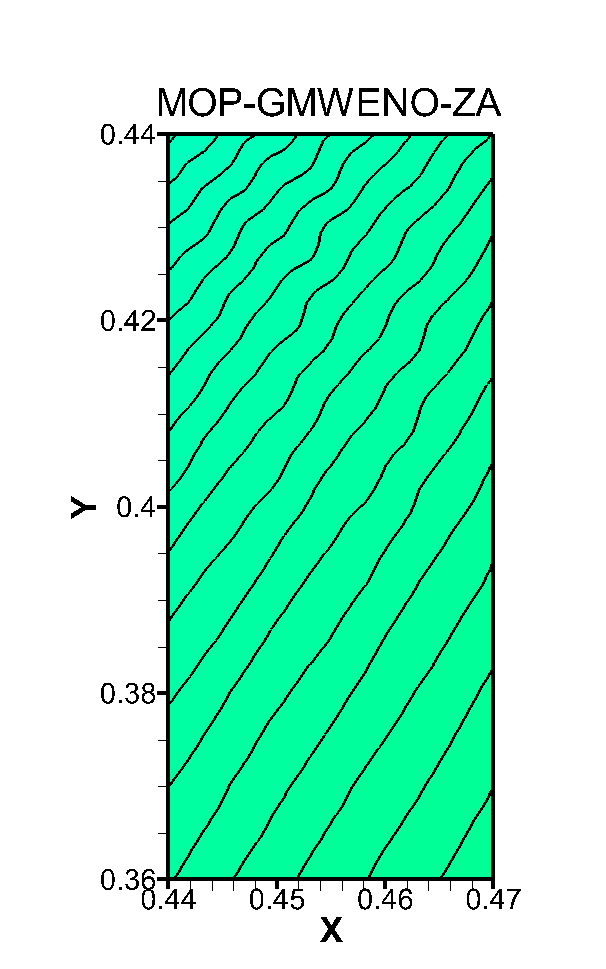} 
  \includegraphics[height=0.375\textwidth]
  {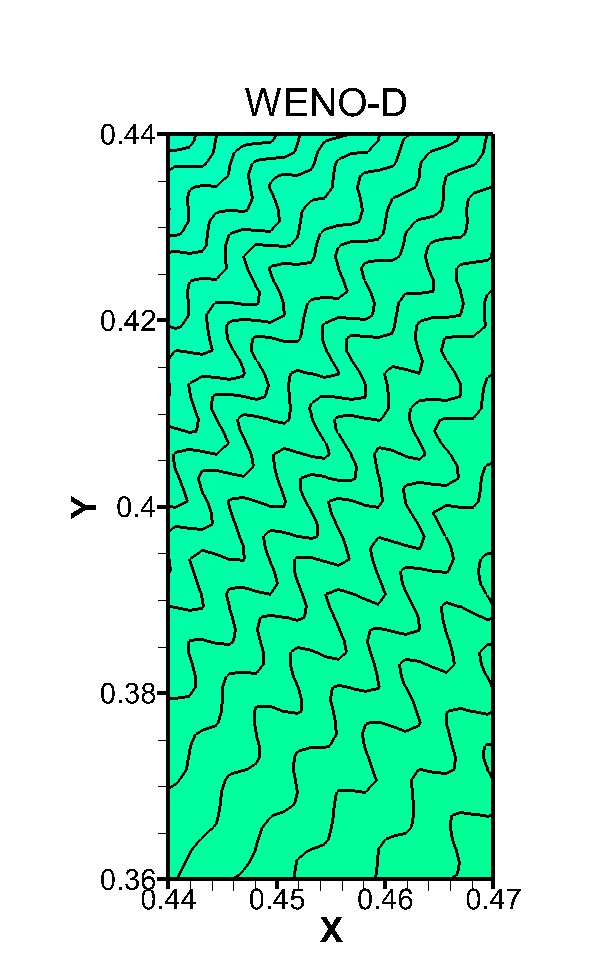} 
  \includegraphics[height=0.375\textwidth]
  {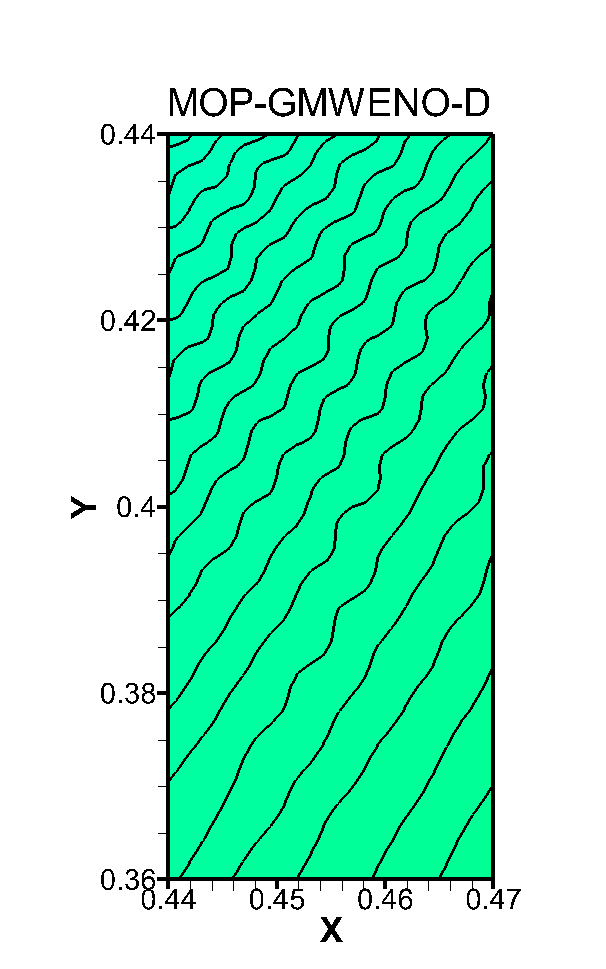}     
\caption{Density contours of Example \ref{ex:Riemann2D}. Left 
column: WENO-ZA and MOP-GMWENO-ZA; Right: WENO-D and MOP-GMWENO-D.}
\label{fig:ex:Riemann2D:3}
\end{figure}

\begin{figure}[ht]
\centering
  \includegraphics[height=0.432\textwidth]
  {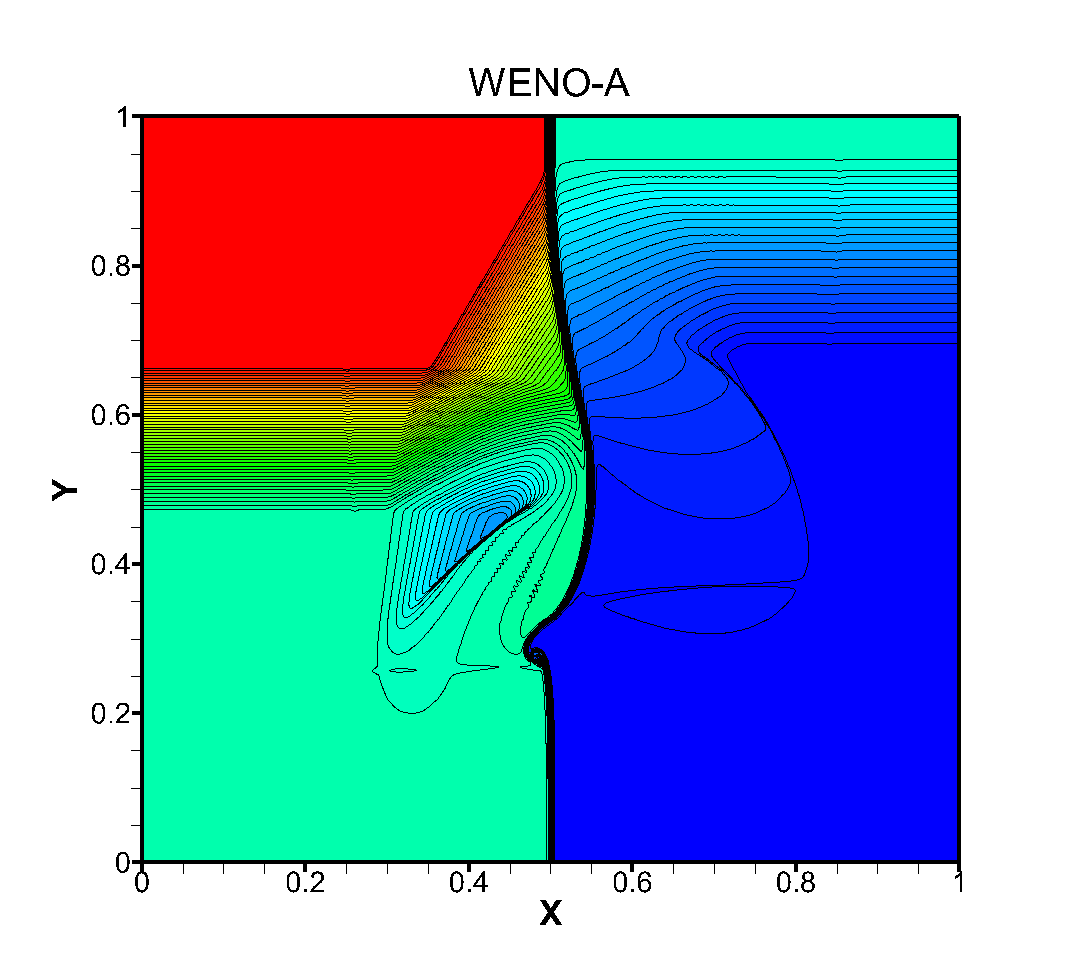}
  \includegraphics[height=0.432\textwidth]
  {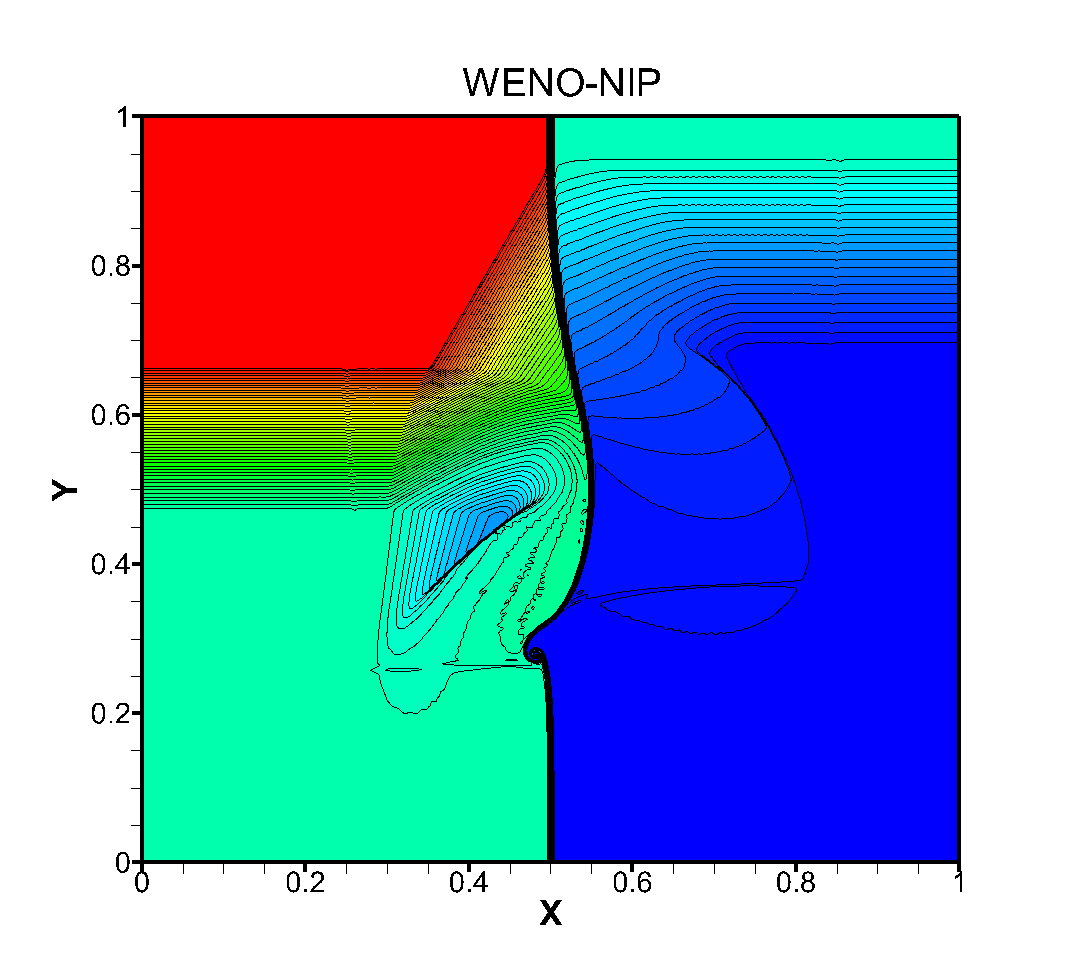}\\
  \includegraphics[height=0.432\textwidth]
  {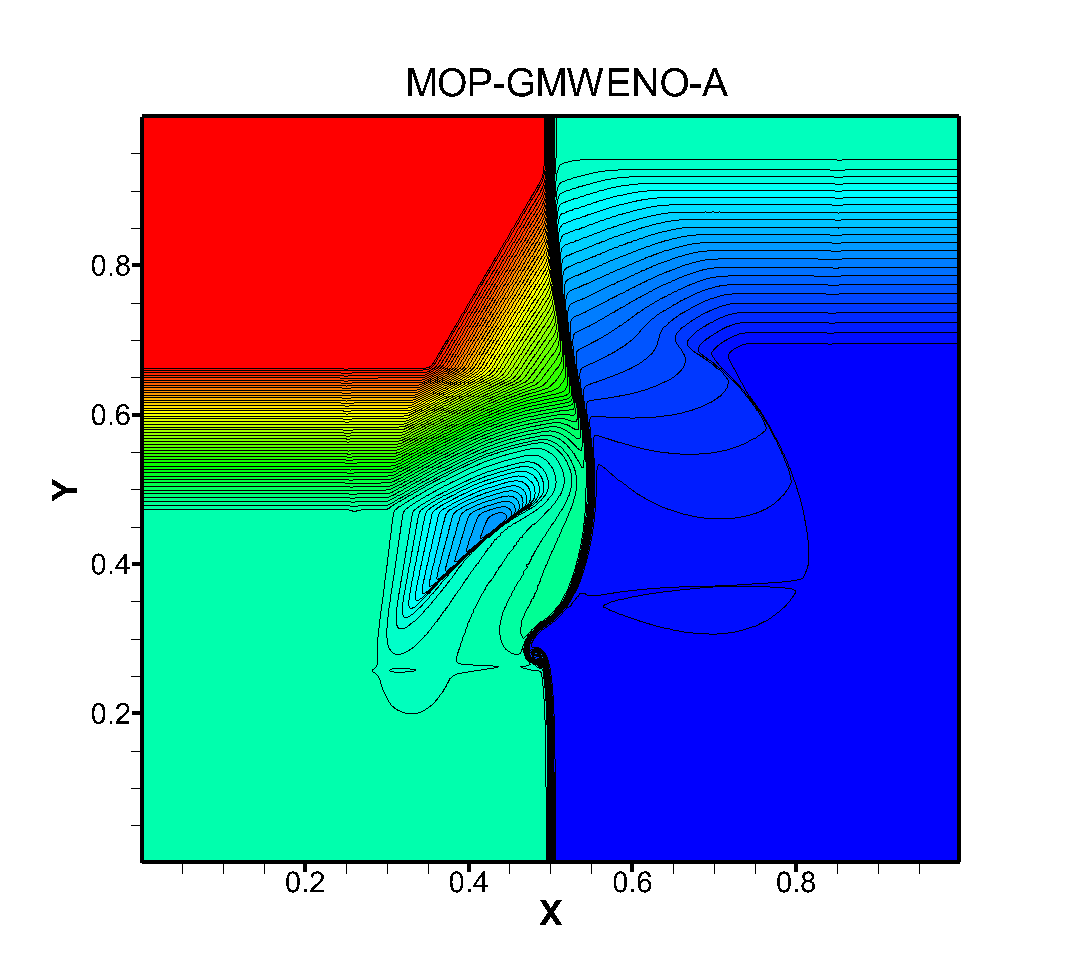}
  \includegraphics[height=0.432\textwidth]
  {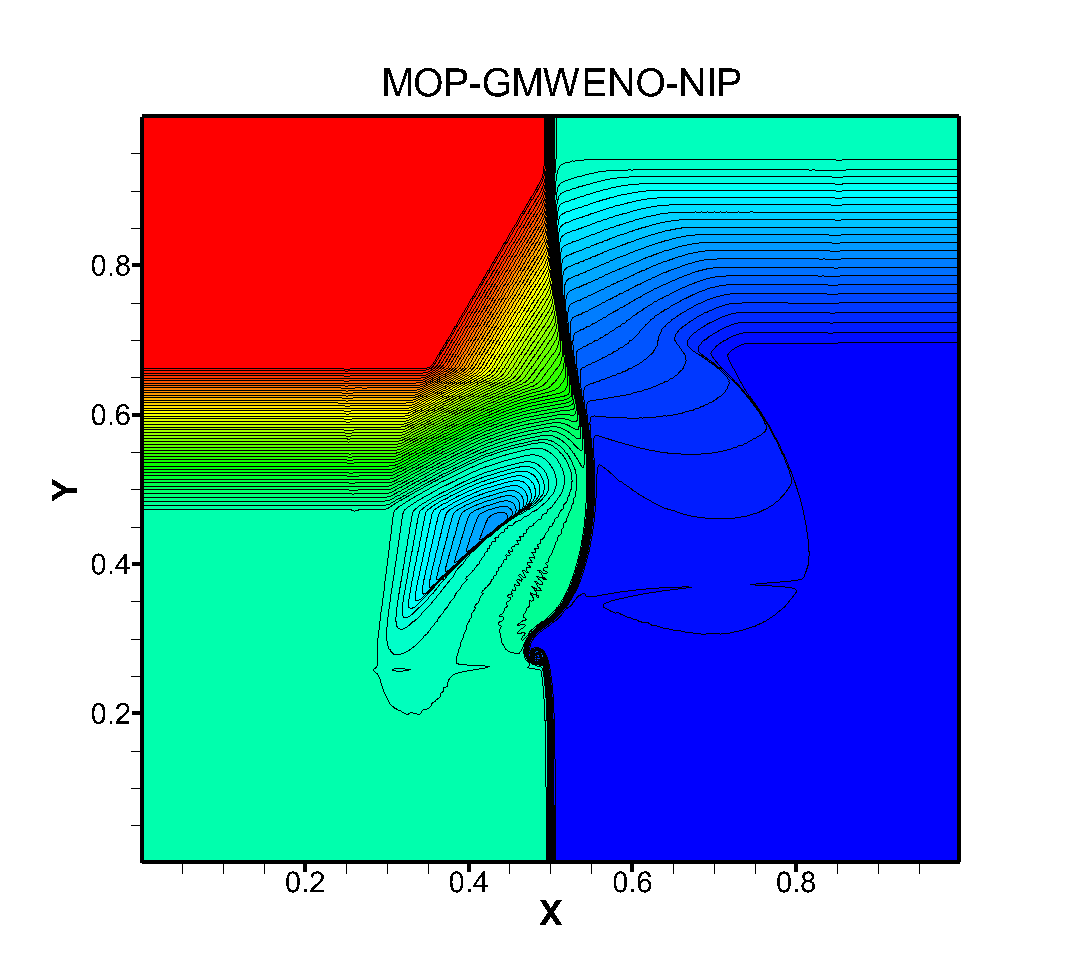}\\
  \includegraphics[height=0.375\textwidth]
  {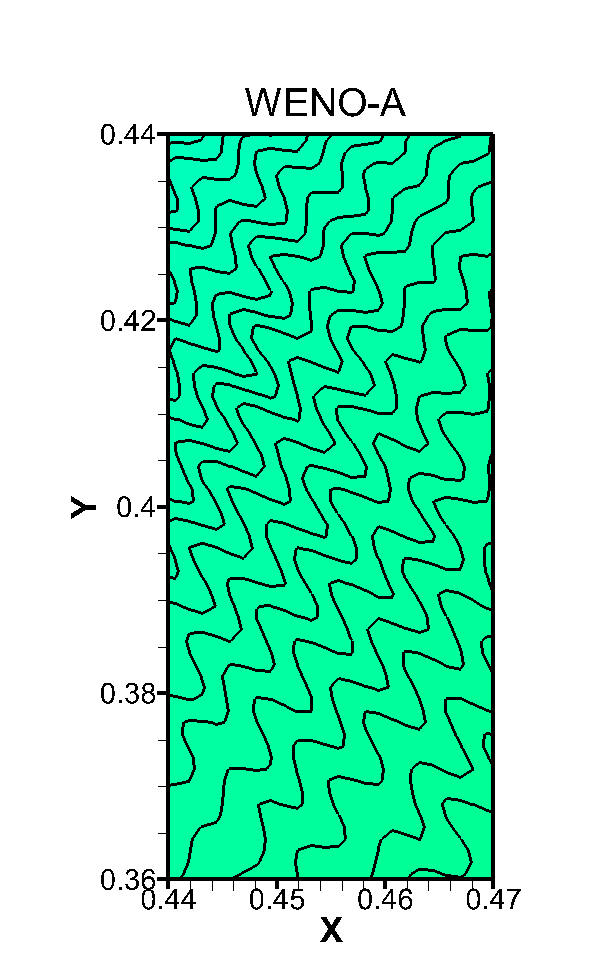} 
  \includegraphics[height=0.375\textwidth]
  {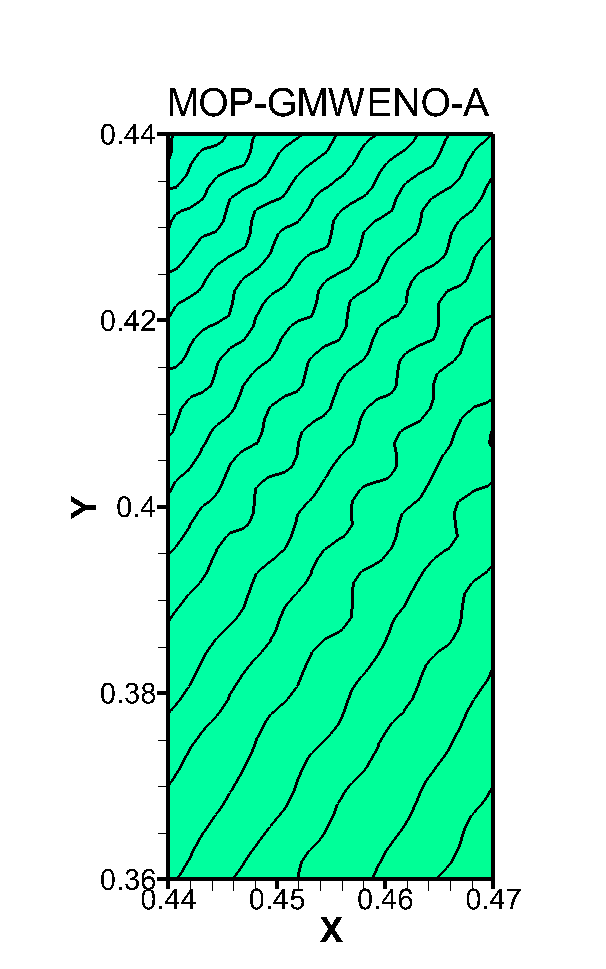} 
  \includegraphics[height=0.375\textwidth]
  {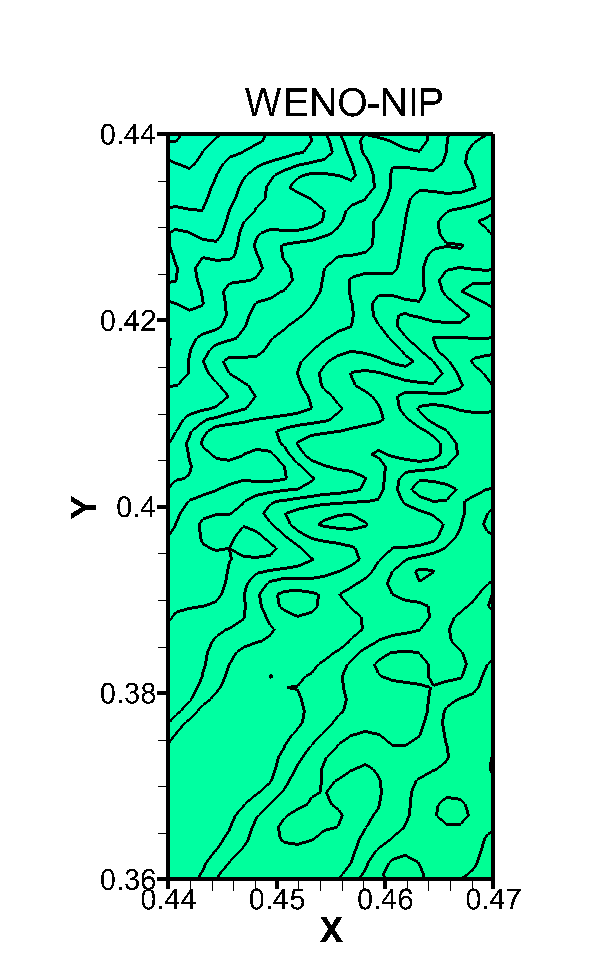} 
  \includegraphics[height=0.375\textwidth]
  {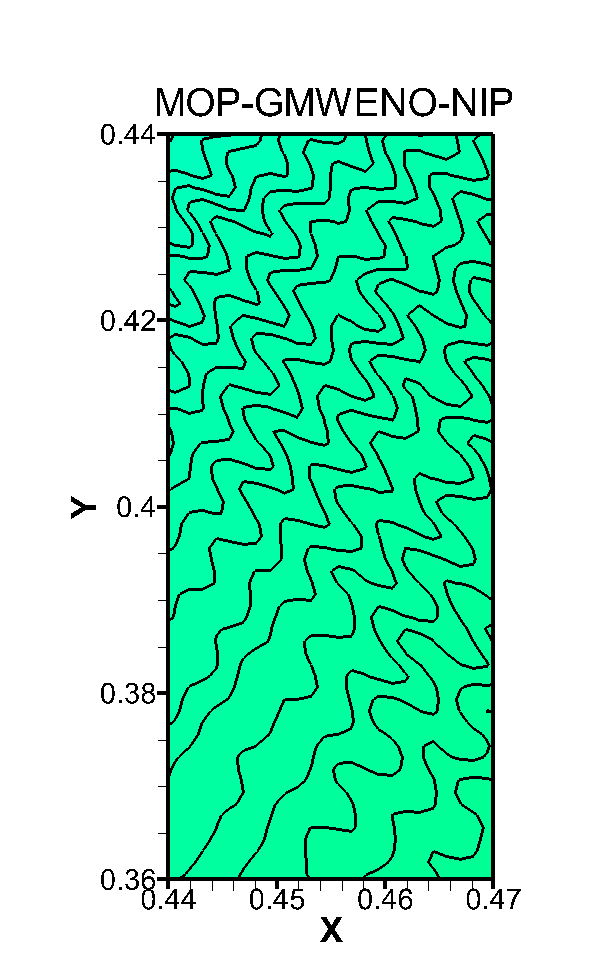}     
\caption{Density contours of Example \ref{ex:Riemann2D}. Left 
column: WENO-A and MOP-GMWENO-A; Right: WENO-NIP and MOP-GMWENO-NIP.}
\label{fig:ex:Riemann2D:4}
\end{figure}

The solutions are shown in Fig. \ref{fig:ex:Riemann2D:1} to Fig. 
\ref{fig:ex:Riemann2D:4}. In the first two rows, we give the density profiles of the 2D Riemann problem simulated by the WENO-X schemes and MOP-GMWENO-X schemes respectively. In the last rows, we display the close-up view of the zone where the post-shock oscillations are 
generated. We can observe that: (1) the main structure of the 2D 
Riemann problem was captured properly by all the considered schemes; 
(2) in the solutions of the WENO-X schemes, evident post-shock 
oscillations can be observed, however, in almost all the solutions 
of the MOP-GMWENO-X schemes except MOP-GMWENO-NIP, the post-shock 
oscillations are dramatically reduced; (3) moreover, from the 
close-up views in the third rows, it is easy to see that the 
amplitudes of the post-shock oscillations produced by the WENO-X 
schemes (except WENO-NIP) are much larger than the MOP-GMWENO-X schemes; (4) although both the WENO-NIP and MOP-GMWENO-NIP schemes produce evident 
post-shock oscillations, it appears that the MOP-GMWENO-NIP scheme 
performs slightly better than the WENO-NIP scheme, and we will 
discuss this further in the next example. In summary, we claim that 
this should be an additional benefit of the new proposed WENO 
schemes that satisfy the \textit{OP} property. 

\subsubsection{Shock-vortex interaction}
\begin{example}
\rm{This is a very favorable two-dimensional test case for 
high-resolution methods \cite{Shock-vortex_interaction-1,
Shock-vortex_interaction-2,Shock-vortex_interaction-3}. The 
computational domain is $[0, 1]\times[0, 1]$, and the initial condition is given by}
\label{ex:shock-vortex}
\end{example}
\begin{equation*}
\mathbf{U}(x, y, 0) = \left\{
\begin{aligned}
\begin{array}{ll}
\mathbf{U}_{\mathrm{L}} = \Big(\rho_{\mathrm{L}} + \delta \rho, u_{\mathrm{L}} + \delta u, 
v_{\mathrm{L}} + \delta v, p_{\mathrm{L}} + \delta p\Big)^{\mathrm{T}},  & x < 0.5, \\
\mathbf{U}_{\mathrm{R}} = \Big(\rho_{\mathrm{R}}, u_{\mathrm{R}}, 
v_{\mathrm{R}}, p_{\mathrm{R}}\Big)^{\mathrm{T}}, & x \geq 0.5, \\
\end{array}
\end{aligned}
\right.
\label{eq:initial_Euer2D:shock-vortex-interaction}
\end{equation*}
with
\begin{equation*}
\begin{array}{l}
\begin{aligned}
&\rho_{\mathrm{L}} = 1, u_{\mathrm{L}} = \sqrt{\gamma}, 
v_{\mathrm{L}} = 0, p_{\mathrm{L}} = 1, \\
&\rho_{\mathrm{R}} = \rho_{\mathrm{L}}\bigg(\dfrac{\gamma - 1 + 
(\gamma + 1)p_{\mathrm{R}}}{\gamma + 1 + (\gamma - 1)p_{\mathrm{R}}} 
\bigg), u_{\mathrm{R}} = \dfrac{u_{\mathrm{L}}(1 - 
p_{\mathrm{R}})}{\Big(\gamma-1 + p_{\mathrm{R}}(\gamma + 1)\Big)^{1/2}}, 
v_{\mathrm{R}} = 0, p_{\mathrm{R}} = 1.3,\\
&\delta \rho = \dfrac{\rho_{\mathrm{L}}^{2}}{
p_{\mathrm{L}}}\dfrac{\delta T}{\gamma - 1}, 
\delta u = \dfrac{\epsilon}{r_\mathrm{c}}(y - y_{\mathrm{c}})
\mathrm{e}^{\alpha(1-r^{2})}, 
\delta v =  \dfrac{- \epsilon}{r_\mathrm{c}}(x - x_{\mathrm{c}})
\mathrm{e}^{\alpha(1-r^{2})}, 
\delta p = \dfrac{\gamma \rho_{\mathrm{L}}^{2}}{
\rho_{\mathrm{L}}}\dfrac{\delta T}{\gamma - 1},
\end{aligned}
\end{array}
\end{equation*}
where $r = \sqrt{((x - x_{\mathrm{c}})^{2} + (y - y_{\mathrm{c}})^{2})/r_{\mathrm{c}}^{2}}, \delta T = - (\gamma - 1)\epsilon^{2}\mathrm{e}^{2\alpha (1 - r^{2})}/(4\alpha \gamma)$, and $\epsilon$ is 0.3, $r_{\mathrm{c}}$ is 0.05, $\alpha$ is 0.204, $x_{\mathrm{c}}$ is 0.25, $y_{\mathrm{c}}$ is 0.5. The outflow condition is used 
on all edges. We discretize the computational domain into
$800 \times 800$ cells and set $t = 0.35$.

\begin{figure}[ht]
\centering
  \includegraphics[height=0.44\textwidth]
  {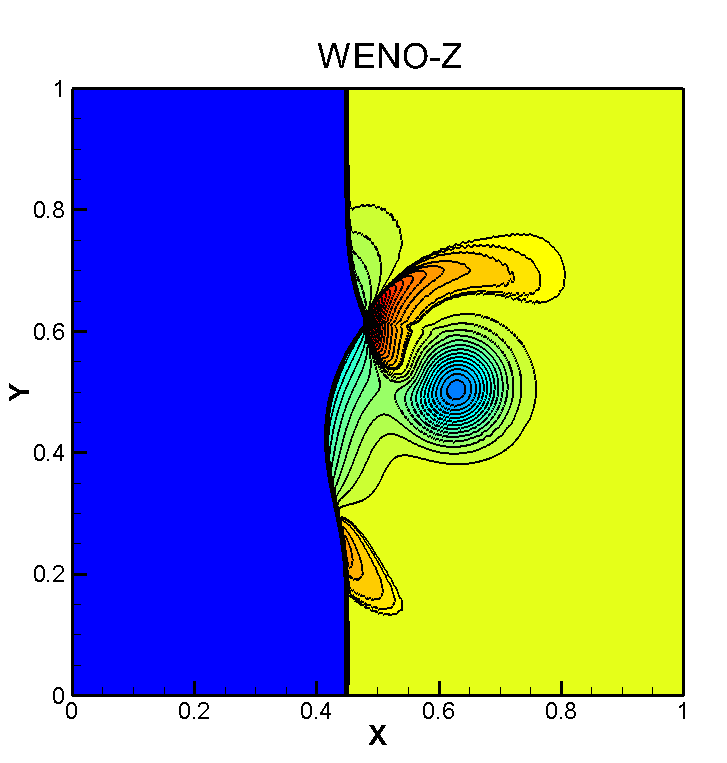}          \hspace{1.2ex}
  \includegraphics[height=0.44\textwidth]
  {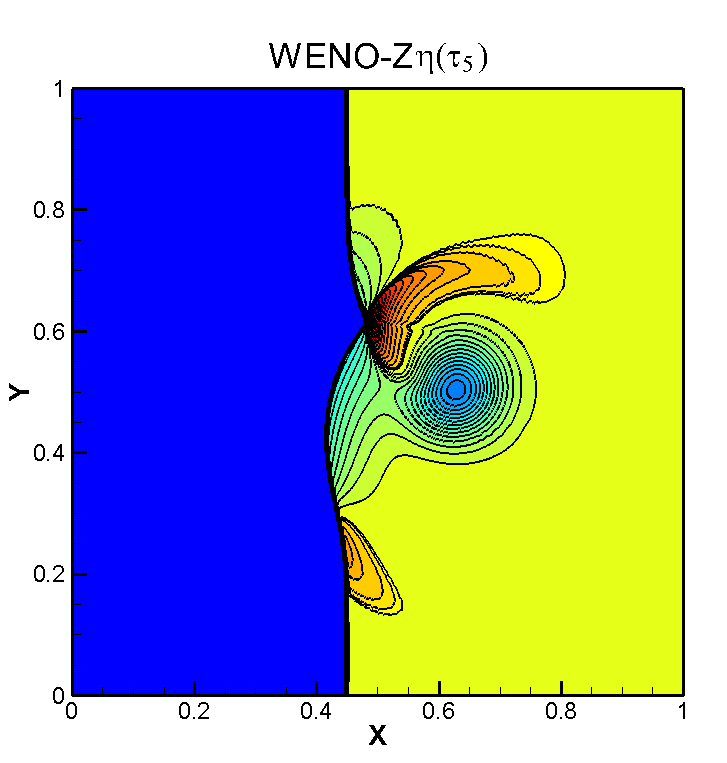}\\
  \includegraphics[height=0.44\textwidth]
  {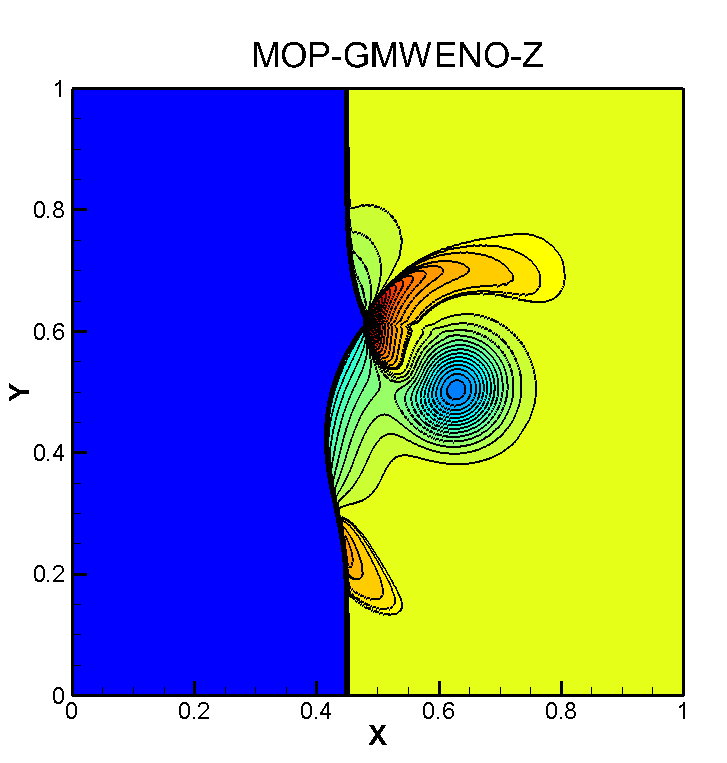}      \hspace{1.2ex}
  \includegraphics[height=0.44\textwidth]
  {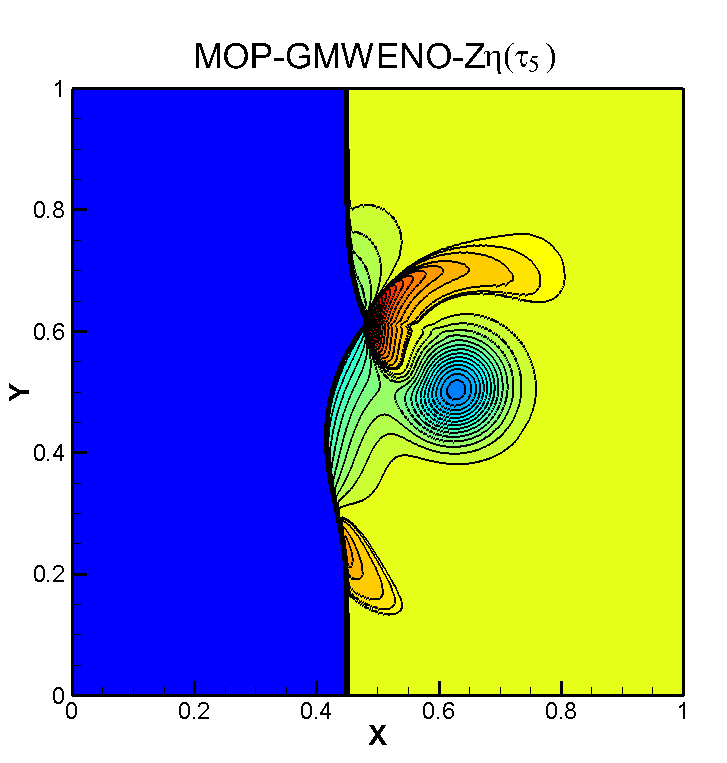}\\
\hspace{-3.5ex}  
  \includegraphics[height=0.345\textwidth]
  {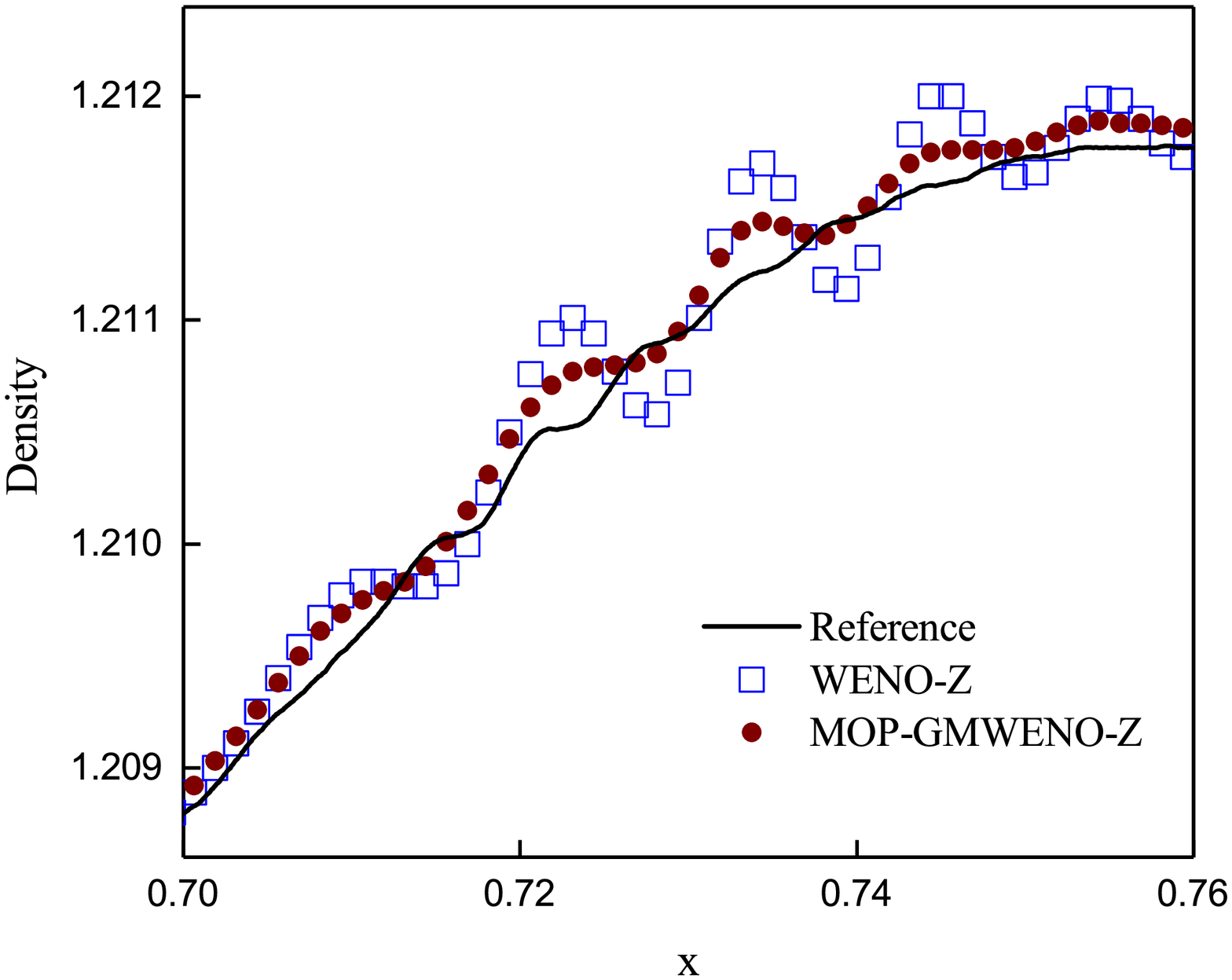}
  \includegraphics[height=0.345\textwidth]
  {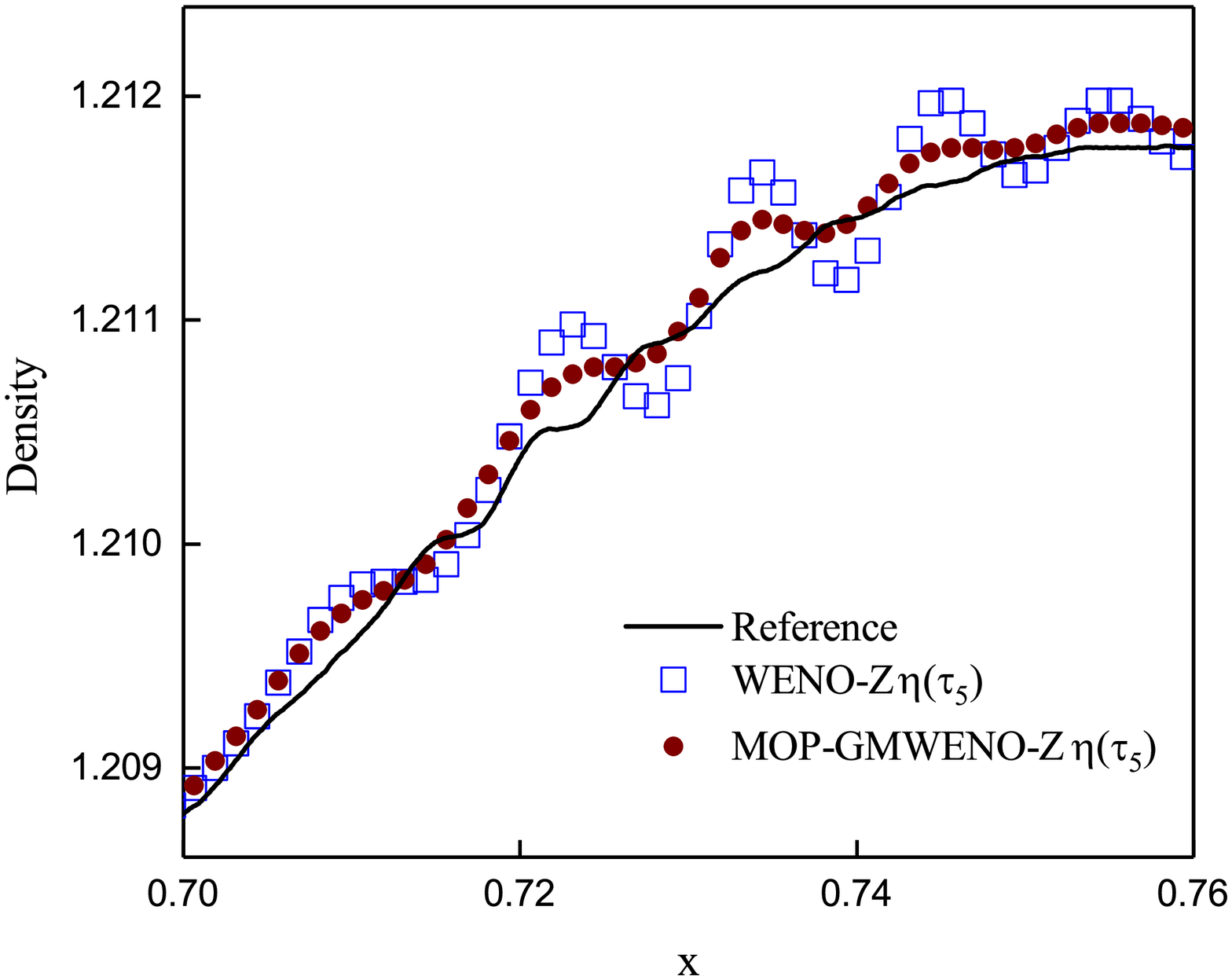}    
\caption{Density contours of Example \ref{ex:shock-vortex} (Row 1 and Row 2) and the density cross-sectional slices at $y = 0.65$ (Row 3w). Left: WENO-Z and MOP-GMWENO-Z; Right: WENO-Z$\eta(\tau_{5})$ and MOP-GMWENO-Z$\eta(\tau_{5})$.}
\label{fig:ex:SVI:1}
\end{figure}

\begin{figure}[ht]
\centering
  \includegraphics[height=0.44\textwidth]
  {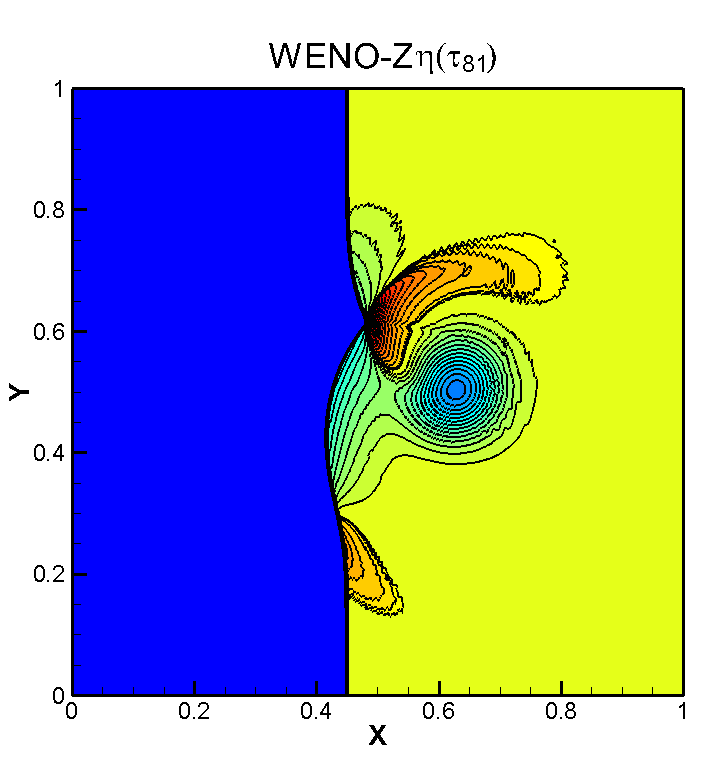}          \hspace{1.2ex}
  \includegraphics[height=0.44\textwidth]
  {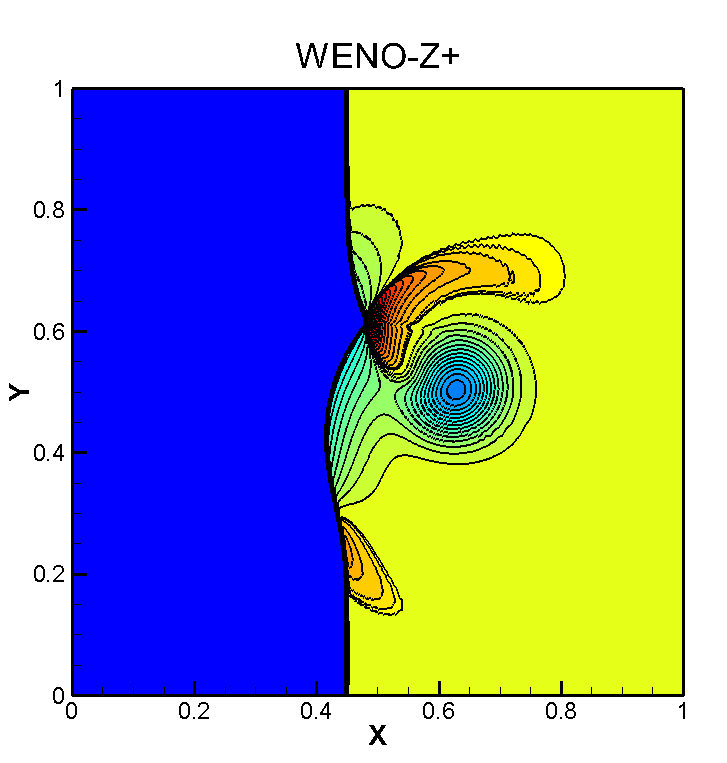}\\
  \includegraphics[height=0.44\textwidth]
  {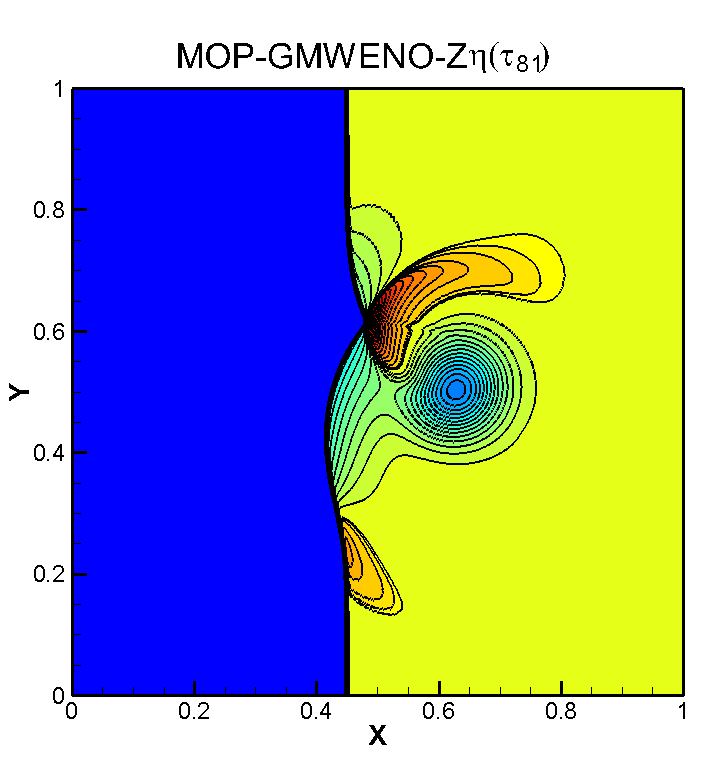}      \hspace{1.2ex}
  \includegraphics[height=0.44\textwidth]
  {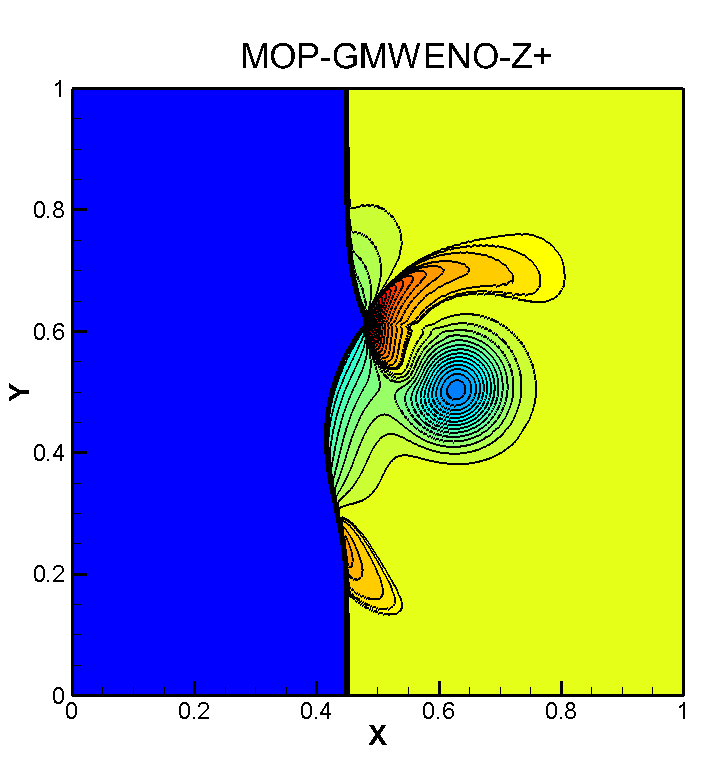}\\
\hspace{-3.5ex}  
  \includegraphics[height=0.345\textwidth]
  {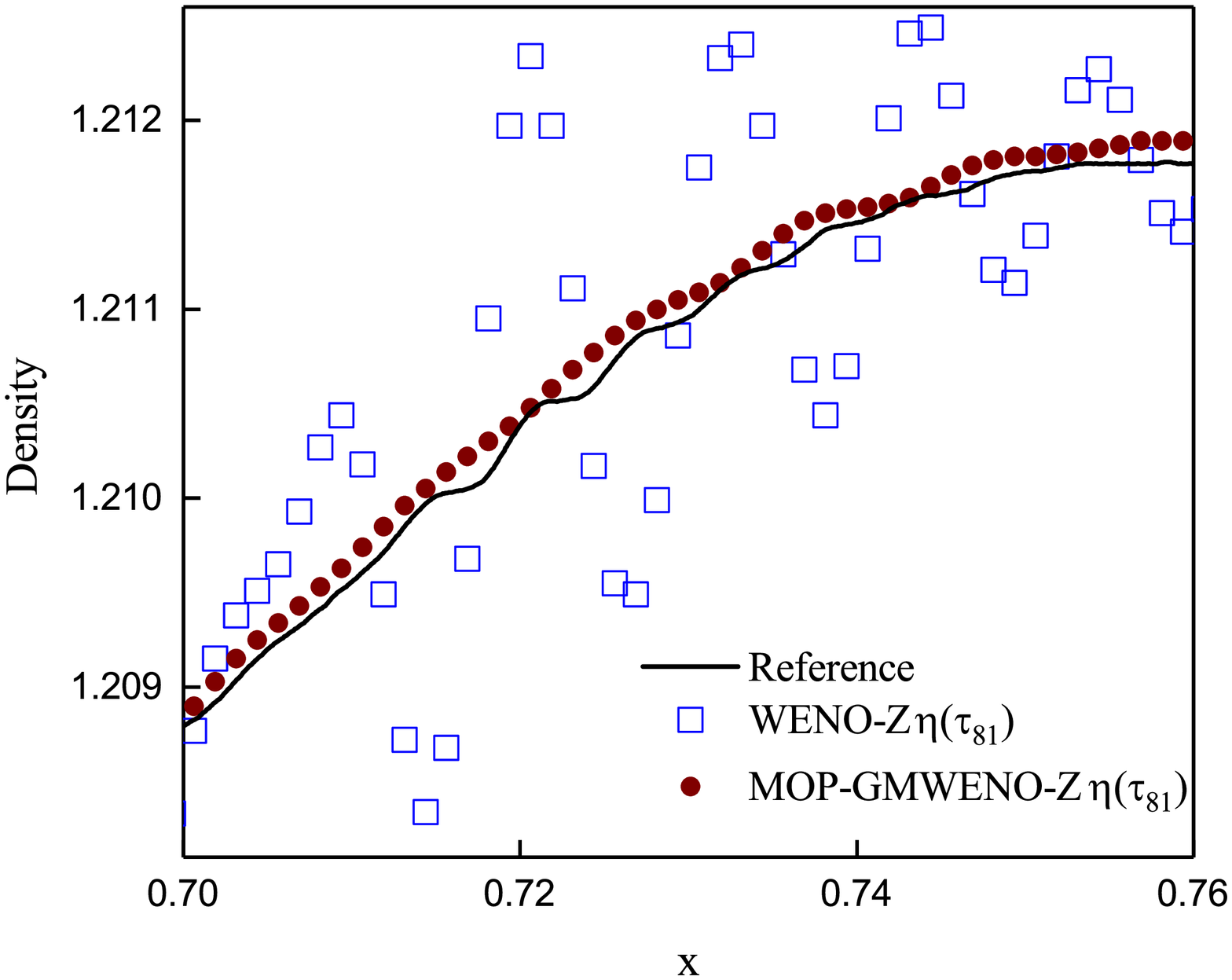}
  \includegraphics[height=0.345\textwidth]
  {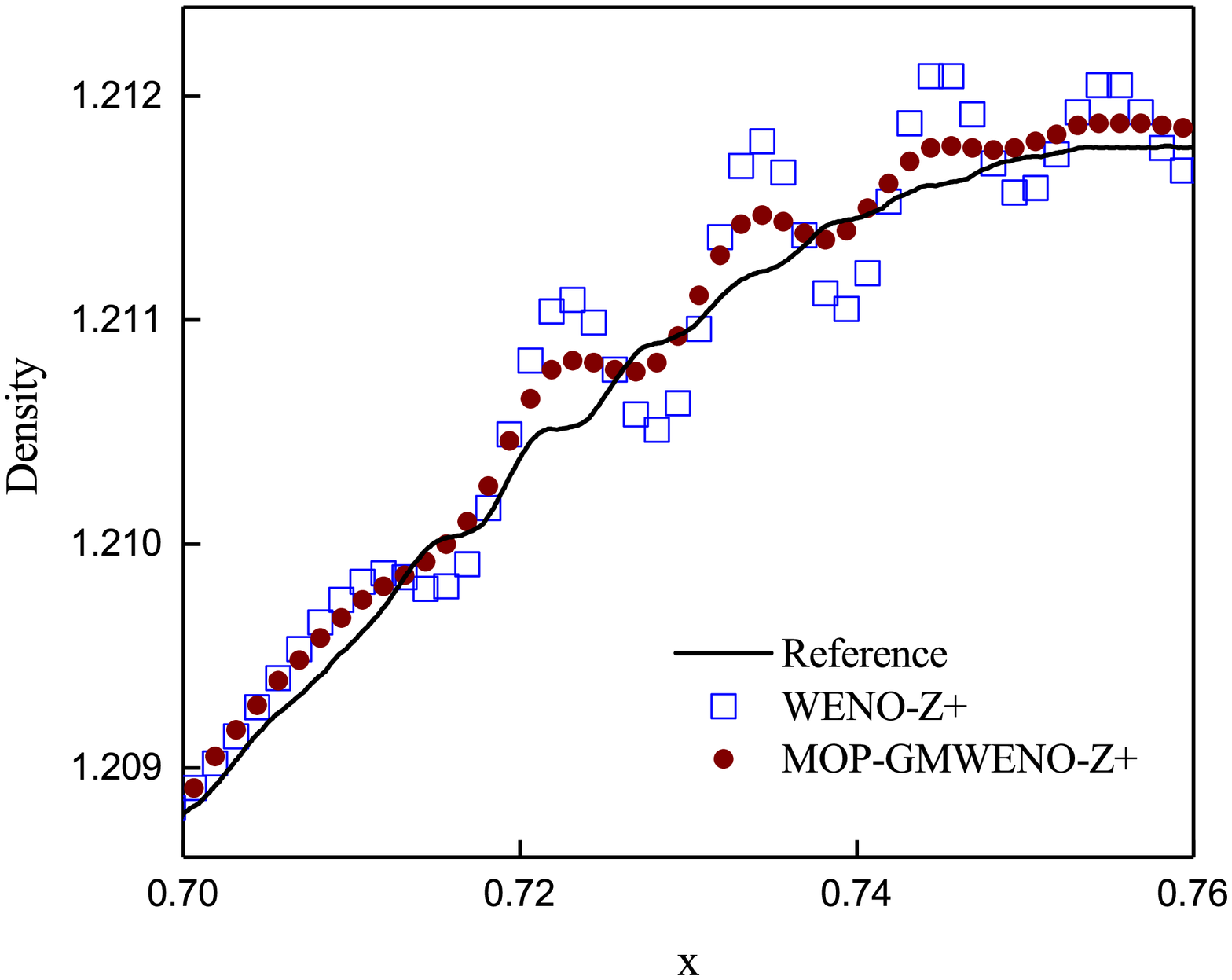}    
\caption{Density contours of Example \ref{ex:shock-vortex} (Row 1 and Row 2) and the density cross-sectional slices at $y = 0.65$ (the last row). Left: WENO-Z$\eta(\tau_{81})$ and MOP-GMWENO-Z$\eta(\tau_{81})$; Right: WENO-Z+ and MOP-GMWENO-Z+.}
\label{fig:ex:SVI:2}
\end{figure}

\begin{figure}[ht]
\centering
  \includegraphics[height=0.44\textwidth]
  {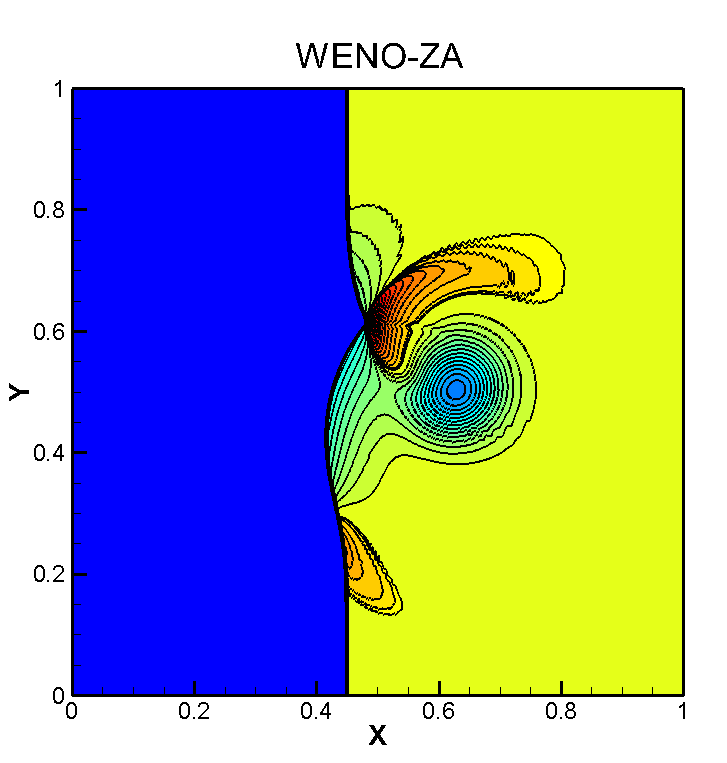}          \hspace{1.2ex}
  \includegraphics[height=0.44\textwidth]
  {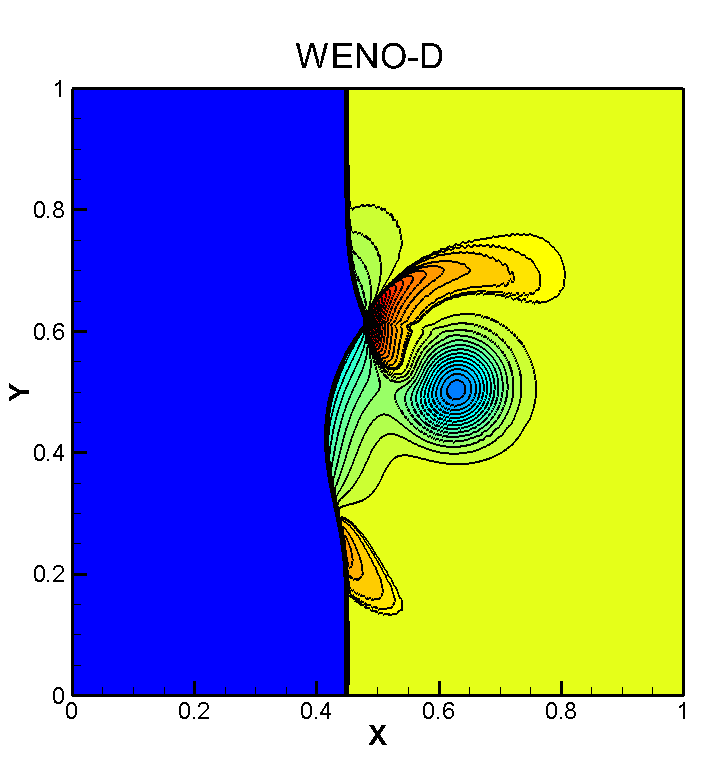}\\
  \includegraphics[height=0.44\textwidth]
  {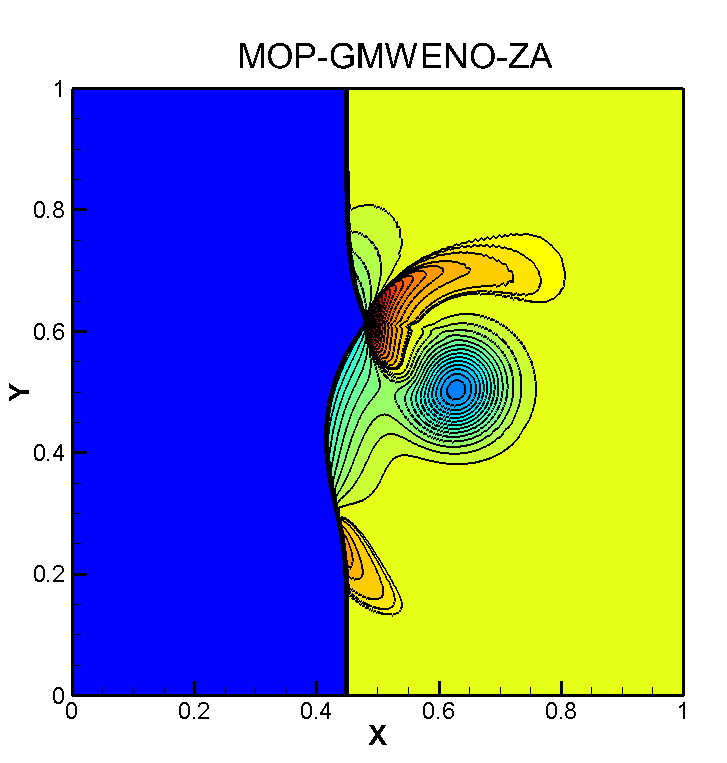}      \hspace{1.2ex}
  \includegraphics[height=0.44\textwidth]
  {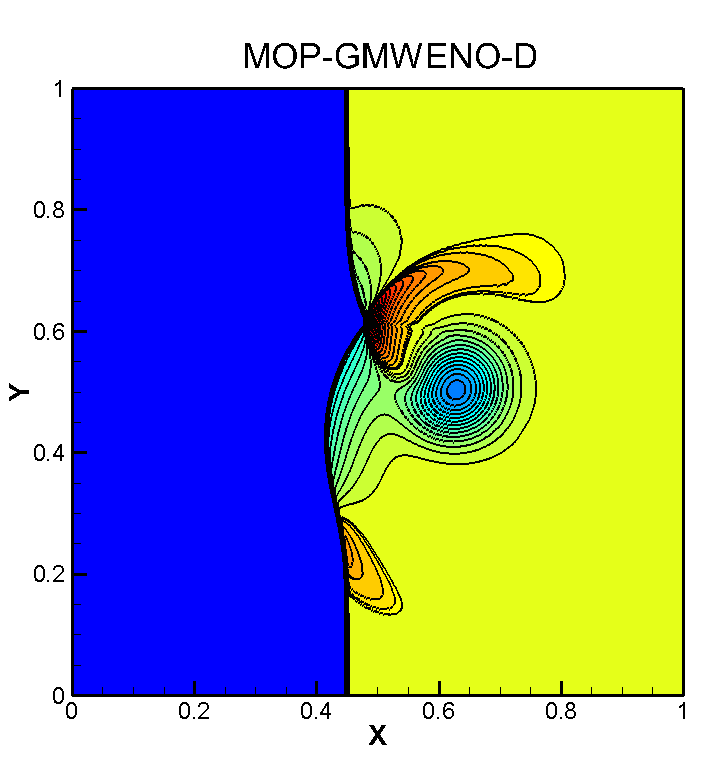}\\
\hspace{-3.5ex}  
  \includegraphics[height=0.345\textwidth]
  {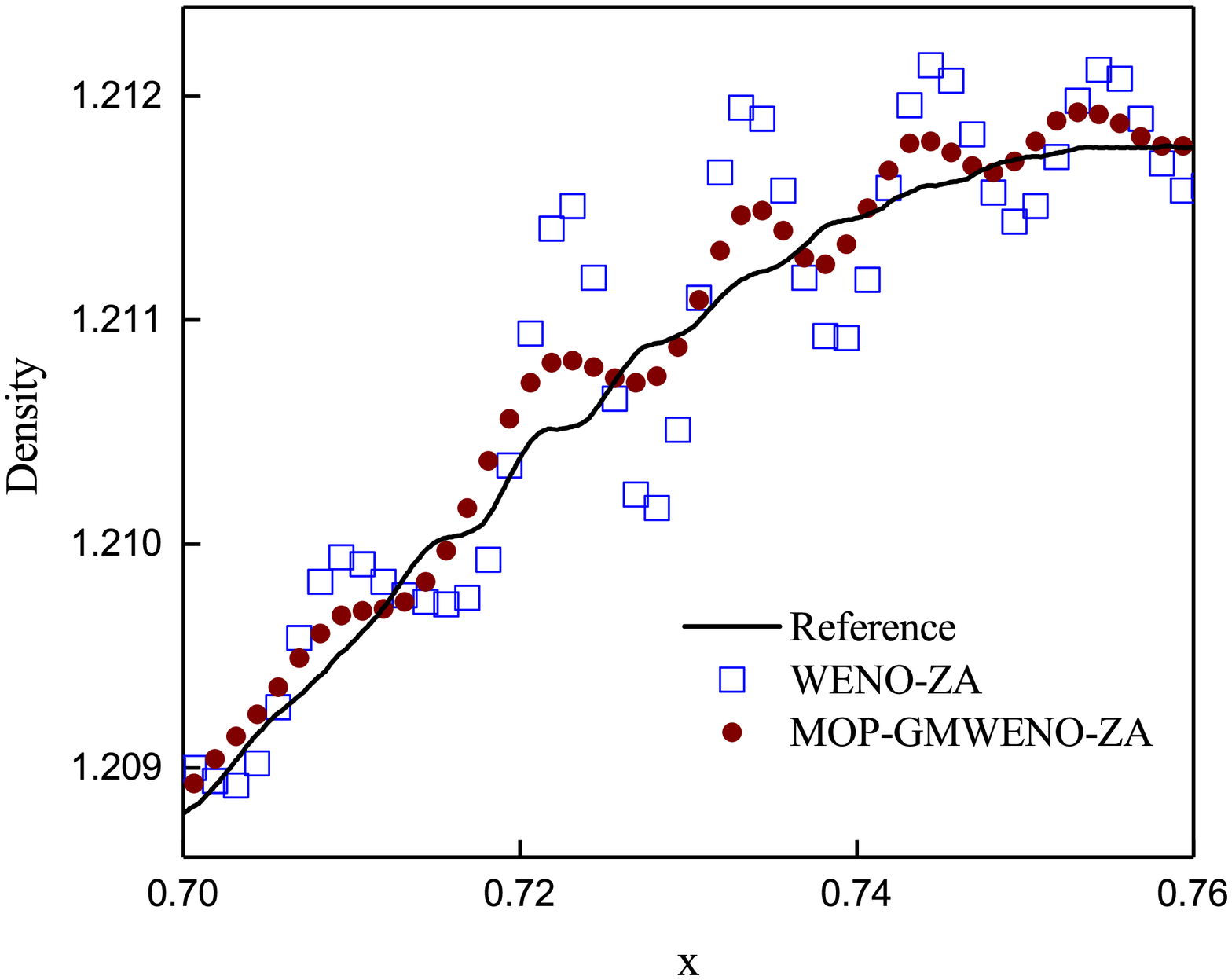}
  \includegraphics[height=0.345\textwidth]
  {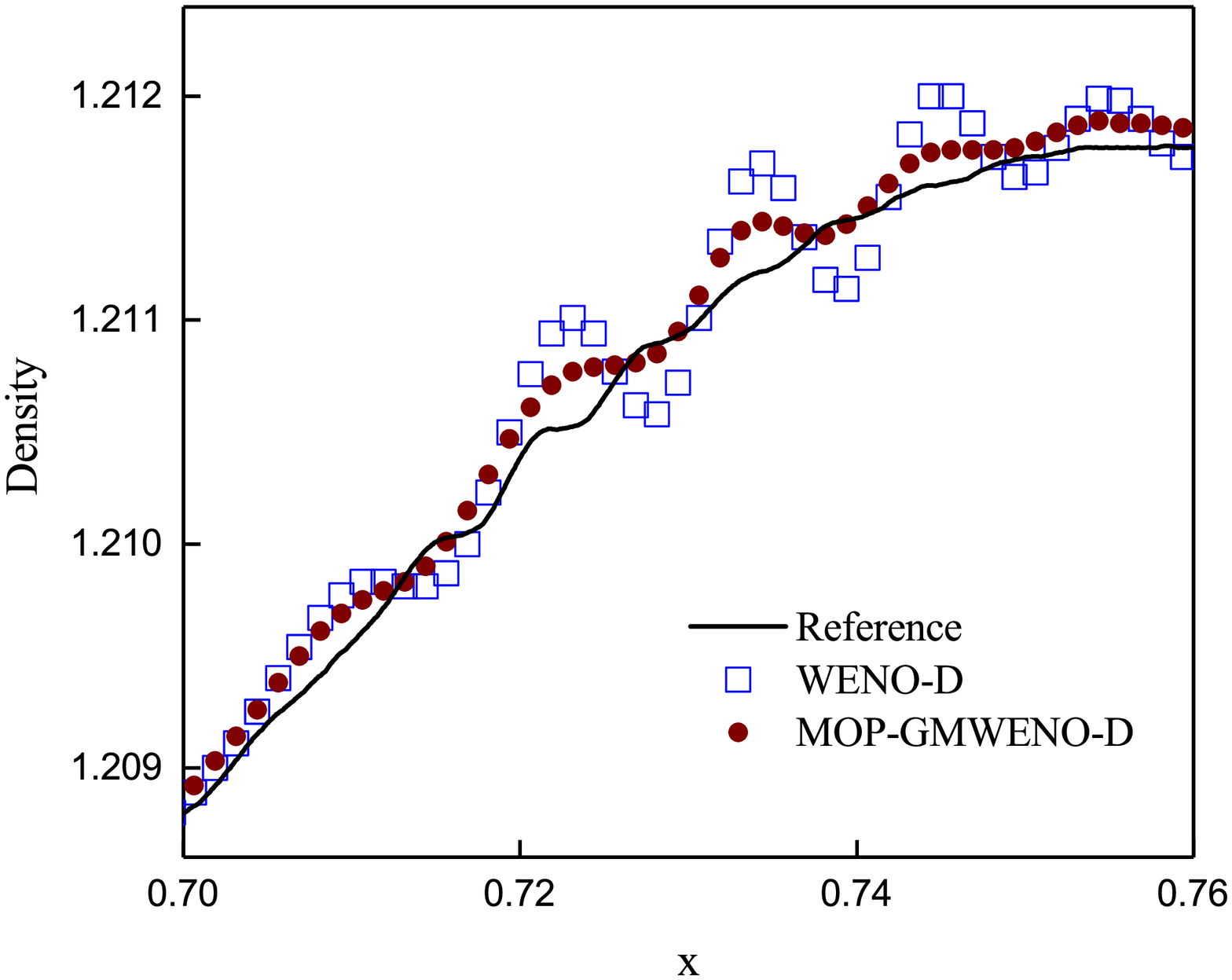}    
\caption{Density contours of Example \ref{ex:shock-vortex} (Row 1 and Row 2) and the density cross-sectional slices at $y = 0.65$ (the last row). Left:  WENO-ZA and MOP-GMWENO-ZA; Right: WENO-D and MOP-GMWENO-D.}
\label{fig:ex:SVI:3}
\end{figure}

\begin{figure}[ht]
\centering
  \includegraphics[height=0.44\textwidth]
  {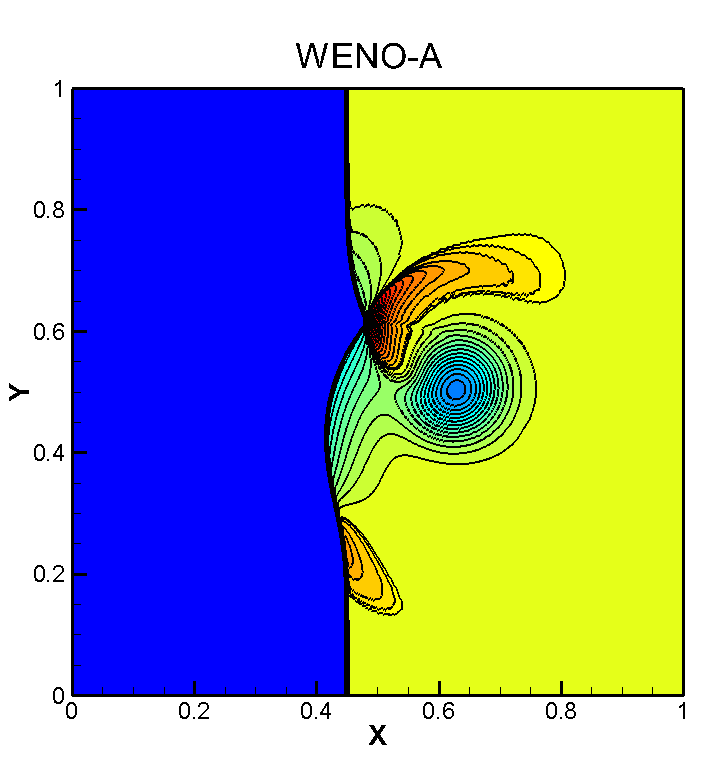}          \hspace{1.2ex}
  \includegraphics[height=0.44\textwidth]
  {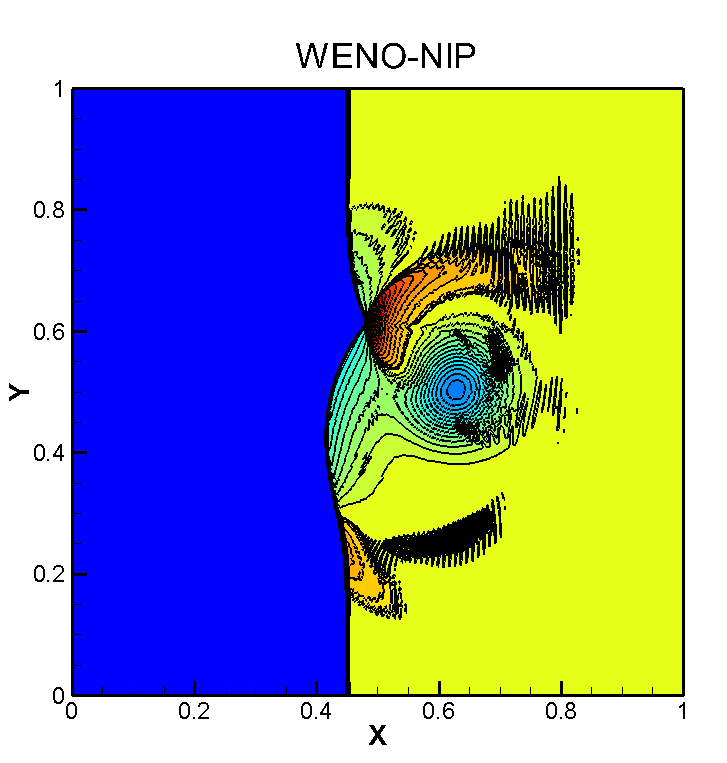}\\
  \includegraphics[height=0.44\textwidth]
  {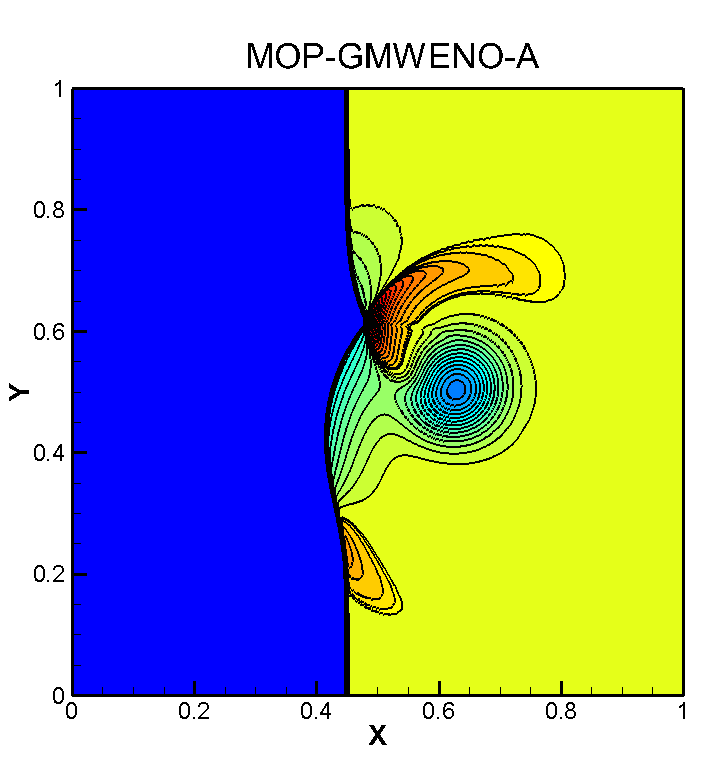}      \hspace{1.2ex}
  \includegraphics[height=0.44\textwidth]
  {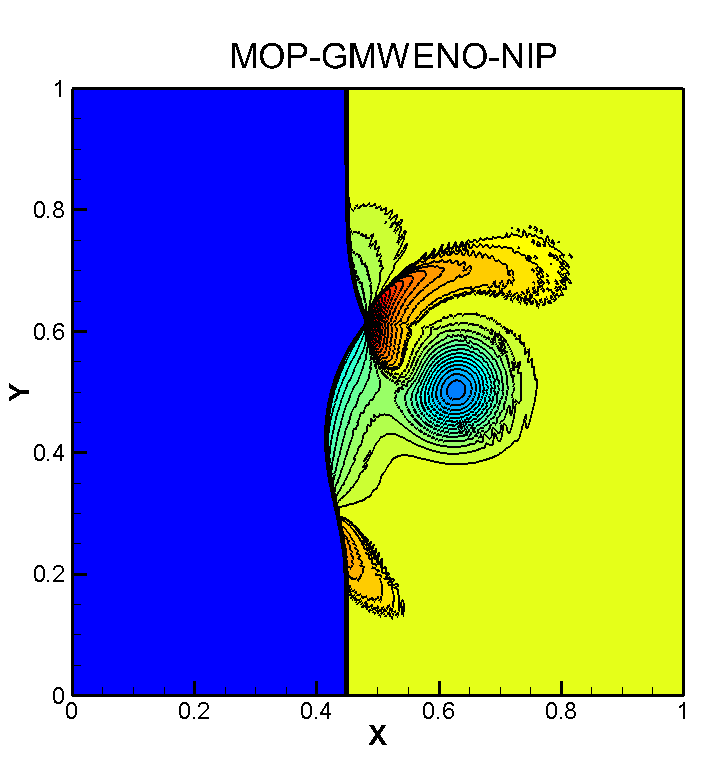}\\
\hspace{-3.5ex}  
  \includegraphics[height=0.345\textwidth]
  {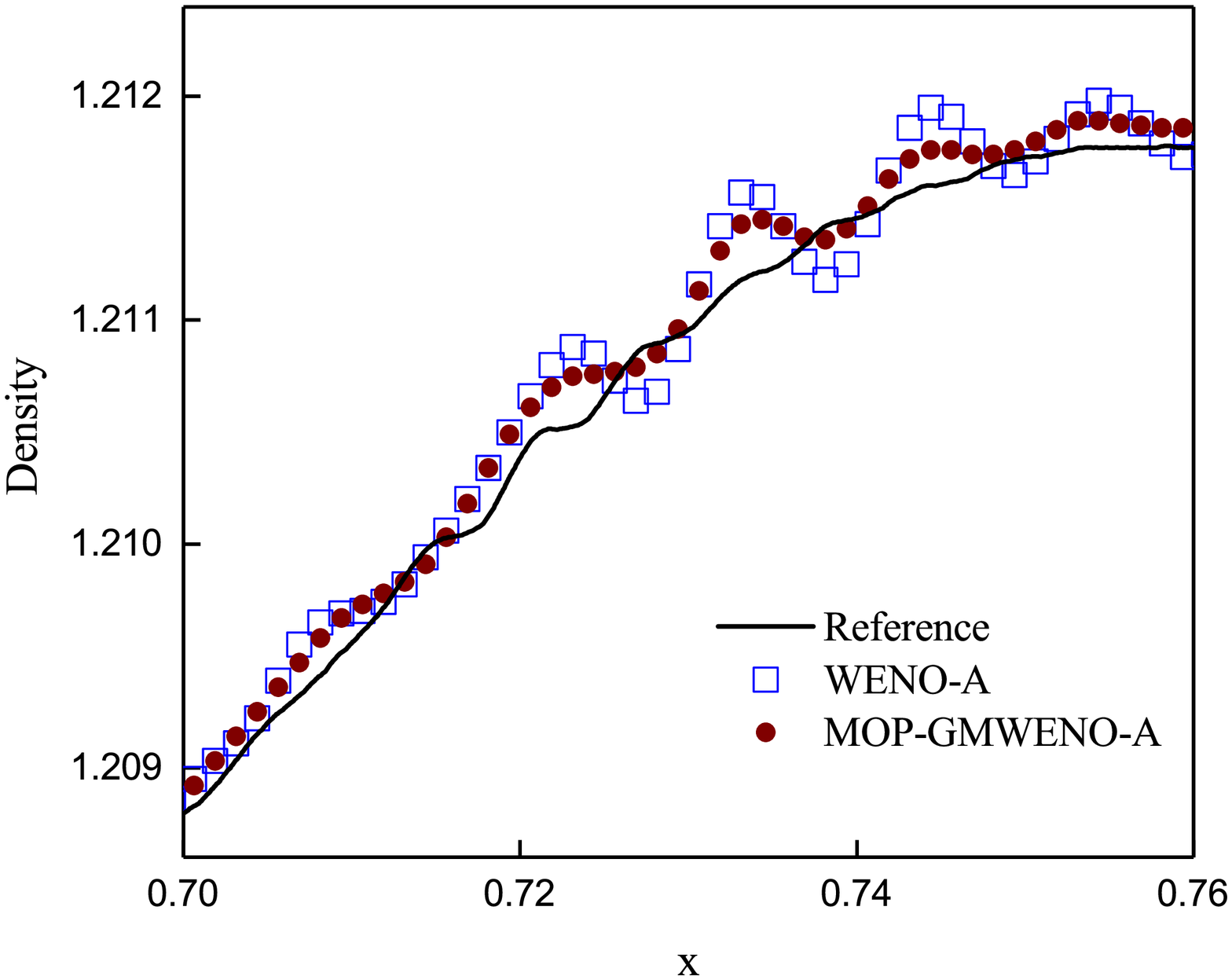}
  \includegraphics[height=0.345\textwidth]
  {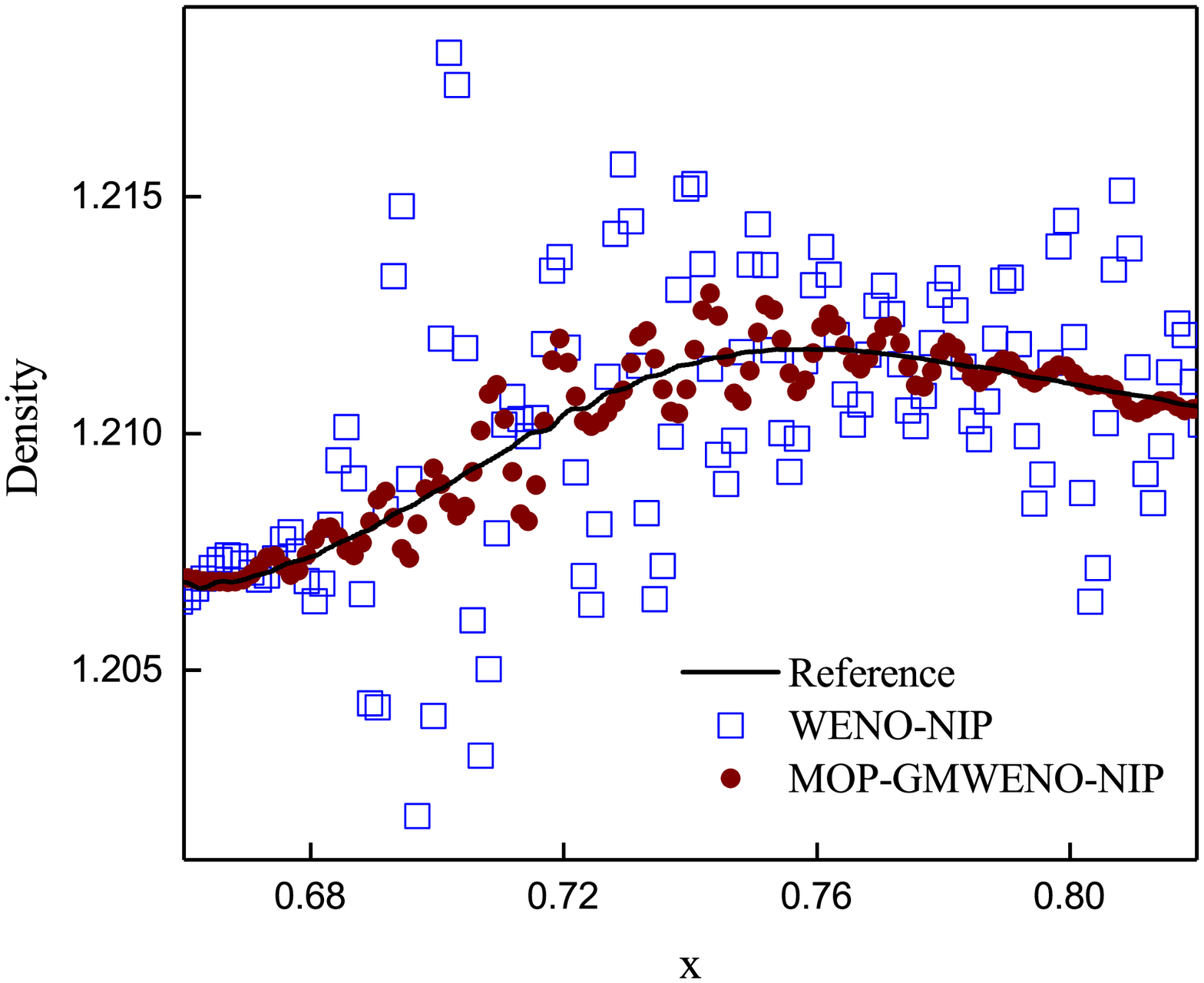}    
\caption{Density contours of Example \ref{ex:shock-vortex} (Row 1 and Row 2) and the density cross-sectional slices at $y = 0.65$ (the last row). Left:  WENO-A and MOP-GMWENO-A; Right: WENO-NIP and MOP-GMWENO-NIP.}
\label{fig:ex:SVI:4}
\end{figure}

We show the results computed by the MOP-GMWENO-X schemes and the WENO-X schemes in Fig. \ref{fig:ex:SVI:1} to 
Fig. \ref{fig:ex:SVI:4}. In the first two rows, 
the solutions in the density profile 
of the WENO-X schemes and the MOP-GMWENO-X schemes are 
given, respectively. In order to show the enhancement of the 
MOP-GMWENO-X schemes more clearly, in the last rows, we show the density cross-sectional slices at $y = 0.65$. 
We can observe that: (1) the main structure of this complicated flow
are captured properly by all considered schemes; 
(2) the WENO-X schemes generate evident numerical oscillations, 
however, the MOP-GMWENO-X schemes, except MOP-GMWENO-NIP, can 
considerably decrease these oscillations; (3) it is interesting that the oscillation produced by WENO-NIP is extremely violent, leading to the fact that the MOP-GMWENO-NIP scheme produces oscillations with larger amplitudes than that of the other 
MOP-GMWENO-X schemes, even larger than the other WENO-X schemes; (4) 
in spite of this, we can easily find that the MOP-GMWENO-NIP scheme
significantly reduces the oscillations comparing with the WENO-NIP 
scheme; (5) also, from the last rows, we can easily 
find that the amplitudes of the oscillations produced by 
the WENO-X schemes are much greater than the
MOP-GMWENO-X schemes. As mentioned before, this should be a benefit 
of the WENO-Z-type schemes with \textit{OP} generalized mappings.

\section{Conclusions}
\label{secConclusions} 

In this paper we improved the family of the WENO-Z-type schemes. We 
have extended the order-preserving (OP) criterion to the WENO-Z-type 
schemes. As reported in the literature, the OP property plays a critical part in preserving high resolutions and meanwhile avoiding spurious 
oscillations for the traditional mapped WENO schemes on long 
simulations. As the fact that the real-time one-to-one relationships 
between the nonlinear weights of the WENO-JS scheme and those of the 
WENO-Z-type schemes are very similar to the designed mapping 
relationships of the traditional mapped WENO schemes, we generated 
the idea of establishing the concept of generalized mapped WENO 
scheme. Accordingly, we devised a uniform formula for the Z-type 
weights. And then, we designed a general algorithm to implement our 
idea. Therefore, the improved WENO-Z-type schemes, dubbed 
MOP-GMWENO-X, have been proposed by naturally introducing the OP 
criterion. Compared to the traditional WENO-Z-type schemes, the 
improved WENO-Z-type schemes have the following enhancements: (1) 
they can amend the drawback of the traditional WENO-Z-type schemes 
of suffering from either generating spurious oscillations or losing 
high resolutions in long simulations of hyperbolic problems; (2) 
they can significantly decreasing the post-shock oscillations in the 
simulations of the 2D Euler problems with strong shock waves.



\bibliographystyle{model1b-shortjournal-num-names}
\bibliography{refs}

\begin{thebibliography}{45}
\expandafter\ifx\csname natexlab\endcsname\relax\def\natexlab#1{#1}\fi
\providecommand{\bibinfo}[2]{#2}
\ifx\xfnm\relax \def\xfnm[#1]{\unskip,\space#1}\fi
\bibitem[{Acker et~al.(2016)Acker, de~R.~Borges and Costa}]{WENO-Zplus}
\bibinfo{author}{F.~Acker}, \bibinfo{author}{R.B. de~R.~Borges},
  \bibinfo{author}{B.~Costa}, \bibinfo{title}{An improved {WENO-Z} scheme},
  \bibinfo{shortjournal}{J. Comput. Phys.} \bibinfo{volume}{313}
  (\bibinfo{year}{2016}) \bibinfo{pages}{726--753}.
\bibitem[{Borges et~al.(2008)Borges, Carmona, Costa and Don}]{WENO-Z}
\bibinfo{author}{R.~Borges}, \bibinfo{author}{M.~Carmona},
  \bibinfo{author}{B.~Costa}, \bibinfo{author}{W.S. Don}, \bibinfo{title}{An
  improved weighted essentially non-oscillatory scheme for hyperbolic
  conservation laws}, \bibinfo{shortjournal}{J. Comput. Phys.}
  \bibinfo{volume}{227} (\bibinfo{year}{2008}) \bibinfo{pages}{3191--3211}.
\bibitem[{Castro et~al.(2011)Castro, Costa and Don}]{WENO-Z01}
\bibinfo{author}{M.~Castro}, \bibinfo{author}{B.~Costa}, \bibinfo{author}{W.S.
  Don}, \bibinfo{title}{High order weighted essentially non-oscillatory
  {WENO-Z} schemes for hyperbolic conservation laws}, \bibinfo{shortjournal}{J.
  Comput. Phys.} \bibinfo{volume}{230} (\bibinfo{year}{2011})
  \bibinfo{pages}{1766--1792}.
\bibitem[{Chatterjee(1999)}]{Shock-vortex_interaction-1}
\bibinfo{author}{A.~Chatterjee}, \bibinfo{title}{Shock wave deformation in
  shock-vortex interactions}, \bibinfo{shortjournal}{Shock Waves}
  \bibinfo{volume}{9} (\bibinfo{year}{1999}) \bibinfo{pages}{95--105}.
\bibitem[{Fan(2014)}]{WENO-eta-02}
\bibinfo{author}{P.~Fan}, \bibinfo{title}{High order weighted essentially
  non-oscillatory {WENO}-$\eta$ schemes for hyperbolic conservation laws},
  \bibinfo{shortjournal}{J. Comput. Phys.} \bibinfo{volume}{269}
  (\bibinfo{year}{2014}) \bibinfo{pages}{355--385}.
\bibitem[{Fan et~al.(2014)Fan, Shen, Tian and Yang}]{WENO-eta}
\bibinfo{author}{P.~Fan}, \bibinfo{author}{Y.~Shen}, \bibinfo{author}{B.~Tian},
  \bibinfo{author}{C.~Yang}, \bibinfo{title}{A new smoothness indicator for
  improving the weighted essentially non-oscillatory scheme},
  \bibinfo{shortjournal}{J. Comput. Phys.} \bibinfo{volume}{269}
  (\bibinfo{year}{2014}) \bibinfo{pages}{329--354}.
\bibitem[{Feng et~al.(2012)Feng, Hu and Wang}]{WENO-PM}
\bibinfo{author}{H.~Feng}, \bibinfo{author}{F.~Hu}, \bibinfo{author}{R.~Wang},
  \bibinfo{title}{A new mapped weighted essentially non-oscillatory scheme},
  \bibinfo{shortjournal}{J. Sci. Comput.} \bibinfo{volume}{51}
  (\bibinfo{year}{2012}) \bibinfo{pages}{449--473}.
\bibitem[{Feng et~al.(2014)Feng, Huang and Wang}]{WENO-IM}
\bibinfo{author}{H.~Feng}, \bibinfo{author}{C.~Huang},
  \bibinfo{author}{R.~Wang}, \bibinfo{title}{An improved mapped weighted
  essentially non-oscillatory scheme}, \bibinfo{shortjournal}{Appl. Math.
  Comput.} \bibinfo{volume}{232} (\bibinfo{year}{2014})
  \bibinfo{pages}{453--468}.
\bibitem[{Ha et~al.(2013)Ha, Kim, Lee and Yoon}]{WENO-NS}
\bibinfo{author}{Y.~Ha}, \bibinfo{author}{C.H. Kim}, \bibinfo{author}{Y.J.
  Lee}, \bibinfo{author}{J.~Yoon}, \bibinfo{title}{An improved weighted
  essentially non-oscillatory scheme with a new smoothness indicator},
  \bibinfo{shortjournal}{J. Comput. Phys.} \bibinfo{volume}{232}
  (\bibinfo{year}{2013}) \bibinfo{pages}{68--86}.
\bibitem[{Harten(1989)}]{ENO1987JCP83}
\bibinfo{author}{A.~Harten}, \bibinfo{title}{{ENO} schemes with subcell
  resolution}, \bibinfo{shortjournal}{J. Comput. Phys.} \bibinfo{volume}{83}
  (\bibinfo{year}{1989}) \bibinfo{pages}{148--184}.
\bibitem[{Harten et~al.(1987)Harten, Engquist, Osher and
  Chakravarthy}]{ENO1987JCP71}
\bibinfo{author}{A.~Harten}, \bibinfo{author}{B.~Engquist},
  \bibinfo{author}{S.~Osher}, \bibinfo{author}{S.R. Chakravarthy},
  \bibinfo{title}{Uniformly high order accurate essentially non-oscillatory
  schemes {III}}, \bibinfo{shortjournal}{J. Comput. Phys.} \bibinfo{volume}{71}
  (\bibinfo{year}{1987}) \bibinfo{pages}{231--303}.
\bibitem[{Harten and Osher(1987)}]{ENO1987V24}
\bibinfo{author}{A.~Harten}, \bibinfo{author}{S.~Osher},
  \bibinfo{title}{Uniformly high order accurate essentially non-oscillatory
  schemes {I}}, \bibinfo{shortjournal}{SIAM J. Numer. Anal.}
  \bibinfo{volume}{24} (\bibinfo{year}{1987}) \bibinfo{pages}{279--309}.
\bibitem[{Harten et~al.(1986)Harten, Osher, Engquist and
  Chakravarthy}]{ENO1986}
\bibinfo{author}{A.~Harten}, \bibinfo{author}{S.~Osher},
  \bibinfo{author}{B.~Engquist}, \bibinfo{author}{S.R. Chakravarthy},
  \bibinfo{title}{Some results on uniformly high order accurate essentially
  non-oscillatory schemes}, \bibinfo{shortjournal}{Appl. Numer. Math.}
  \bibinfo{volume}{2} (\bibinfo{year}{1986}) \bibinfo{pages}{347--377}.
\bibitem[{Henrick et~al.(2005)Henrick, Aslam and Powers}]{WENO-M}
\bibinfo{author}{A.K. Henrick}, \bibinfo{author}{T.D. Aslam},
  \bibinfo{author}{J.M. Powers}, \bibinfo{title}{Mapped weighted essentially
  non-oscillatory schemes: Achieving optimal order near critical points},
  \bibinfo{shortjournal}{J. Comput. Phys.} \bibinfo{volume}{207}
  (\bibinfo{year}{2005}) \bibinfo{pages}{542--567}.
\bibitem[{Hu(2017)}]{WENO-NIS}
\bibinfo{author}{F.~Hu}, \bibinfo{title}{The weighted {ENO} scheme based on the
  modified smoothness indicator}, \bibinfo{shortjournal}{Comput. Fluids}
  \bibinfo{volume}{150} (\bibinfo{year}{2017}) \bibinfo{pages}{1--7}.
\bibitem[{Jiang and Shu(1996)}]{WENO-JS}
\bibinfo{author}{G.S. Jiang}, \bibinfo{author}{C.W. Shu},
  \bibinfo{title}{Efficient implementation of weighted {ENO} schemes},
  \bibinfo{shortjournal}{J. Comput. Phys.} \bibinfo{volume}{126}
  (\bibinfo{year}{1996}) \bibinfo{pages}{202--228}.
\bibitem[{Kim et~al.(2016)Kim, Ha and Yoon}]{WENO-P}
\bibinfo{author}{C.H. Kim}, \bibinfo{author}{Y.~Ha}, \bibinfo{author}{J.~Yoon},
  \bibinfo{title}{Modified {Non-linear} {Weights} for {Fifth-Order} {Weighted}
  {Essentially} {Non-oscillatory} {Schemes}}, \bibinfo{shortjournal}{J. Sci.
  Comput.} \bibinfo{volume}{67} (\bibinfo{year}{2016})
  \bibinfo{pages}{299--323}.
\bibitem[{Lax and Liu(1998)}]{Riemann2D-03}
\bibinfo{author}{P.D. Lax}, \bibinfo{author}{X.D. Liu},
  \bibinfo{title}{Solution of two-dimensional {R}iemann problems of gas
  dynamics by positive schemes}, \bibinfo{shortjournal}{SIAM J. Sci. Comput.}
  \bibinfo{volume}{19} (\bibinfo{year}{1998}) \bibinfo{pages}{319--340}.
\bibitem[{Li et~al.(2015)Li, Liu and Zhang}]{WENO-PPM5}
\bibinfo{author}{Q.~Li}, \bibinfo{author}{P.~Liu}, \bibinfo{author}{H.~Zhang},
  \bibinfo{title}{Piecewise {Polynomial} {Mapping} {Method} and {Corresponding}
  {WENO} {Scheme} with {Improved} {Resolution}}, \bibinfo{shortjournal}{Commun.
  Comput. Phys.} \bibinfo{volume}{18} (\bibinfo{year}{2015})
  \bibinfo{pages}{1417--1444}.
\bibitem[{Li and Zhong(2021{\natexlab{a}})}]{WENO-ACM}
\bibinfo{author}{R.~Li}, \bibinfo{author}{W.~Zhong}, \bibinfo{title}{An
  efficient mapped {WENO} scheme using approximate constant mapping},
  \bibinfo{shortjournal}{Numer. Math. Theor. Meth. Appl.}
  (\bibinfo{year}{2021}{\natexlab{a}}) \bibinfo{pages}{Published online,
  doi:10.4208/nmtma.OA--2021--0074}.
\bibitem[{Li and Zhong(2021{\natexlab{b}})}]{PoAOP-WENO-X}
\bibinfo{author}{R.~Li}, \bibinfo{author}{W.~Zhong}, \bibinfo{title}{Locally
  {Order-Preserving} {Mapping} for {WENO} {Methods}},
  \bibinfo{shortjournal}{East Asian Appl. Math.}
  (\bibinfo{year}{2021}{\natexlab{b}}) \bibinfo{pages}{Under review}.
\bibitem[{Li and Zhong(2021{\natexlab{c}})}]{WENO-MAIMi}
\bibinfo{author}{R.~Li}, \bibinfo{author}{W.~Zhong}, \bibinfo{title}{A modified
  adaptive improved mapped {WENO} method}, \bibinfo{shortjournal}{Commun.
  Comput. Phys.}  (\bibinfo{year}{2021}{\natexlab{c}}) \bibinfo{pages}{Accepted
  for publication}.
\bibitem[{Li and Zhong(2021{\natexlab{d}})}]{MOP-WENO-ACMk}
\bibinfo{author}{R.~Li}, \bibinfo{author}{W.~Zhong}, \bibinfo{title}{A new
  mapped {WENO} scheme using order-preserving mapping},
  \bibinfo{shortjournal}{Commun. Comput. Phs.}
  (\bibinfo{year}{2021}{\natexlab{d}}) \bibinfo{pages}{Accepted for
  publication}.
\bibitem[{Li and Zhong(2021{\natexlab{e}})}]{MOP-WENO-X}
\bibinfo{author}{R.~Li}, \bibinfo{author}{W.~Zhong}, \bibinfo{title}{Towards
  building the {OP-Mapped} {WENO} schemes: {A} general methodology},
  \bibinfo{shortjournal}{Math. Comput. Appl.} \bibinfo{volume}{26}
  (\bibinfo{year}{2021}{\natexlab{e}}) \bibinfo{pages}{67}.
\bibitem[{Liu et~al.(2018)Liu, Shen, Zeng and Yu}]{WENO-ZA}
\bibinfo{author}{S.~Liu}, \bibinfo{author}{Y.~Shen}, \bibinfo{author}{F.~Zeng},
  \bibinfo{author}{M.~Yu}, \bibinfo{title}{A new weighting method for improving
  the {WENO-Z} scheme}, \bibinfo{shortjournal}{Int. J. Numer. Meth. Fluids}
  \bibinfo{volume}{87} (\bibinfo{year}{2018}) \bibinfo{pages}{271--291}.
\bibitem[{Liu et~al.(1994)Liu, Osher and Chan}]{WENO-LiuXD}
\bibinfo{author}{X.D. Liu}, \bibinfo{author}{S.~Osher},
  \bibinfo{author}{T.~Chan}, \bibinfo{title}{Weighted essentially
  non-oscillatory schemes}, \bibinfo{shortjournal}{J. Comput. Phys.}
  \bibinfo{volume}{115} (\bibinfo{year}{1994}) \bibinfo{pages}{200--212}.
\bibitem[{Luo and Wu(2021)}]{WENO-ZplusI}
\bibinfo{author}{X.~Luo}, \bibinfo{author}{S.~Wu}, \bibinfo{title}{Improvement
  of the weno-z+ scheme}, \bibinfo{shortjournal}{Comput. Fluids}
  \bibinfo{volume}{218} (\bibinfo{year}{2021}) \bibinfo{pages}{104855}.
\bibitem[{Pao and Salas(1981)}]{Shock-vortex_interaction-2}
\bibinfo{author}{S.P. Pao}, \bibinfo{author}{M.D. Salas}, \bibinfo{title}{A
  numerical study of two-dimensional shock-vortex interaction}, in:
  \bibinfo{booktitle}{AIAA 14th Fluid and Plasma Dynamics Conference},
  \bibinfo{address}{California, Palo Alto, 1981}.
\bibitem[{Peng et~al.(2019)Peng, Zhai, Ni, Yong and Shen}]{AdaWENO-Z}
\bibinfo{author}{J.~Peng}, \bibinfo{author}{C.~Zhai}, \bibinfo{author}{G.~Ni},
  \bibinfo{author}{H.~Yong}, \bibinfo{author}{Y.~Shen}, \bibinfo{title}{An
  adaptive characteristic-wise reconstruction {WENO-Z} scheme for gas dynamic
  euler equations}, \bibinfo{shortjournal}{Comput. Fluids}
  \bibinfo{volume}{179} (\bibinfo{year}{2019}) \bibinfo{pages}{34--51}.
\bibitem[{Pirozzoli(2011)}]{Riemann2D-04}
\bibinfo{author}{S.~Pirozzoli}, \bibinfo{title}{Numerical methods for
  high-speed flows}, \bibinfo{shortjournal}{Annu. Rev. Fluid Mech.}
  \bibinfo{volume}{43} (\bibinfo{year}{2011}) \bibinfo{pages}{163--194}.
\bibitem[{Ren et~al.(2003)Ren, Liu and Zhang}]{Shock-vortex_interaction-3}
\bibinfo{author}{Y.X. Ren}, \bibinfo{author}{M.~Liu},
  \bibinfo{author}{H.~Zhang}, \bibinfo{title}{A characteristic-wise hybrid
  compact-{WENO} scheme for solving hyperbolic conservation laws},
  \bibinfo{shortjournal}{J. Comput. Phys.} \bibinfo{volume}{192}
  (\bibinfo{year}{2003}) \bibinfo{pages}{365--386}.
\bibitem[{Samala and Raju(2018)}]{MWENO-P}
\bibinfo{author}{R.~Samala}, \bibinfo{author}{G.N. Raju}, \bibinfo{title}{A
  modified fifth-order {WENO} scheme for hyperbolic conservation laws},
  \bibinfo{shortjournal}{Comput. Math. Appl.} \bibinfo{volume}{75}
  (\bibinfo{year}{2018}) \bibinfo{pages}{1531--1549}.
\bibitem[{Schulz-Rinne(1993)}]{Riemann2D-02}
\bibinfo{author}{C.W. Schulz-Rinne}, \bibinfo{title}{Classification of the
  {R}iemann problem for two-dimensional gas dynamics},
  \bibinfo{shortjournal}{SIAM J. Math. Anal.} \bibinfo{volume}{24}
  (\bibinfo{year}{1993}) \bibinfo{pages}{76--88}.
\bibitem[{Schulz-Rinne et~al.(1993)Schulz-Rinne, Collins and
  Glaz}]{Riemann-2D-01}
\bibinfo{author}{C.W. Schulz-Rinne}, \bibinfo{author}{J.P. Collins},
  \bibinfo{author}{H.M. Glaz}, \bibinfo{title}{Numerical solution of the
  {R}iemann problem for two-dimensional gas dynamics},
  \bibinfo{shortjournal}{SIAM J. Sci. Comput.} \bibinfo{volume}{14}
  (\bibinfo{year}{1993}) \bibinfo{pages}{1394--1414}.
\bibitem[{Shen and Zha(2008)}]{WENO-eta-LSI}
\bibinfo{author}{Y.~Shen}, \bibinfo{author}{G.~Zha}, \bibinfo{title}{A robust
  seventh-order {WENO} scheme and its application}, \bibinfo{shortjournal}{AIAA
  2008-757}, in: \bibinfo{booktitle}{46th AIAA Aerospace Sciences Meeting and
  Exhibit}, \bibinfo{address}{Reno, Nevada, 2008}.
\bibitem[{Shu(1998)}]{WENOoverview}
\bibinfo{author}{C.W. Shu}, \bibinfo{title}{Essentially non-oscillatory and
  weighted essentially non-oscillatory schemes for hyperbolic conservation
  laws}, in: \bibinfo{booktitle}{{Advanced Numerical Approximation of Nonlinear
  Hyperbolic Equations. Lecture Notes in Mathematics}}, volume
  \bibinfo{volume}{1697}, \bibinfo{publisher}{Springer},
  \bibinfo{address}{Berlin}, \bibinfo{year}{1998}, pp.
  \bibinfo{pages}{325--432}.
\bibitem[{Shu and Osher(1988)}]{ENO-Shu1988}
\bibinfo{author}{C.W. Shu}, \bibinfo{author}{S.~Osher},
  \bibinfo{title}{Efficient implementation of essentially non-oscillatory
  shock-capturing schemes}, \bibinfo{shortjournal}{J. Comput. Phys.}
  \bibinfo{volume}{77} (\bibinfo{year}{1988}) \bibinfo{pages}{439--471}.
\bibitem[{Shu and Osher(1989)}]{ENO-Shu1989}
\bibinfo{author}{C.W. Shu}, \bibinfo{author}{S.~Osher},
  \bibinfo{title}{Efficient implementation of essentially non-oscillatory
  shock-capturing schemes {II}}, \bibinfo{shortjournal}{J. Comput. Phys.}
  \bibinfo{volume}{83} (\bibinfo{year}{1989}) \bibinfo{pages}{32--78}.
\bibitem[{Vevek et~al.(2018)Vevek, Zang and New}]{WENO-RM-Vevek2018}
\bibinfo{author}{U.S. Vevek}, \bibinfo{author}{B.~Zang}, \bibinfo{author}{T.H.
  New}, \bibinfo{title}{A {New} {Mapped} {WENO} {Method} for {Hyperbolic}
  {Problems}}, \bibinfo{shortjournal}{ICCFD10}, in: \bibinfo{booktitle}{Tenth
  International Conference on Computational Fluid Dynamics},
  \bibinfo{address}{Barcelona, Spain, 2018}.
\bibitem[{Vevek et~al.(2019)Vevek, Zang and New}]{WENO-AIM}
\bibinfo{author}{U.S. Vevek}, \bibinfo{author}{B.~Zang}, \bibinfo{author}{T.H.
  New}, \bibinfo{title}{Adaptive mapping for high order {WENO} methods},
  \bibinfo{shortjournal}{J. Comput. Phys.} \bibinfo{volume}{381}
  (\bibinfo{year}{2019}) \bibinfo{pages}{162--188}.
\bibitem[{Wang et~al.(2016)Wang, Feng and Huang}]{WENO-RM260}
\bibinfo{author}{R.~Wang}, \bibinfo{author}{H.~Feng},
  \bibinfo{author}{C.~Huang}, \bibinfo{title}{A {New} {Mapped} {Weighted}
  {Essentially} {Non-oscillatory} {Method} {Using} {Rational} {Function}},
  \bibinfo{shortjournal}{J. Sci. Comput.} \bibinfo{volume}{67}
  (\bibinfo{year}{2016}) \bibinfo{pages}{540--580}.
\bibitem[{Wang et~al.(2019)Wang, Wang and Don}]{WENO-D_WENO-A}
\bibinfo{author}{Y.~Wang}, \bibinfo{author}{B.S. Wang}, \bibinfo{author}{W.S.
  Don}, \bibinfo{title}{Generalized {Sensitivity} {Parameter} {Free} {Fifth}
  {Order} {WENO} {Finite} {Difference} {Scheme} with {Z-Type} {Weights}},
  \bibinfo{shortjournal}{J. Sci. Comput.} \bibinfo{volume}{81}
  (\bibinfo{year}{2019}) \bibinfo{pages}{1329--1358}.
\bibitem[{Yuan(2020)}]{WENO-NIP}
\bibinfo{author}{M.~Yuan}, \bibinfo{title}{A new weighted essentially
  non-oscillatory {WENO-NIP} scheme for hyperbolic conservation laws},
  \bibinfo{shortjournal}{Comput. Fluids} \bibinfo{volume}{197}
  (\bibinfo{year}{2020}) \bibinfo{pages}{104168}.
\bibitem[{Zeng et~al.(2018)Zeng, Shen and Liu}]{P-WENO}
\bibinfo{author}{F.~Zeng}, \bibinfo{author}{Y.~Shen}, \bibinfo{author}{S.~Liu},
  \bibinfo{title}{A perturbational weighted essentially non-oscillatory
  scheme}, \bibinfo{shortjournal}{Comput. Fluids} \bibinfo{volume}{172}
  (\bibinfo{year}{2018}) \bibinfo{pages}{196--208}.
\bibitem[{Zhang and Shu(2007)}]{WENO-ZS}
\bibinfo{author}{S.~Zhang}, \bibinfo{author}{C.W. Shu}, \bibinfo{title}{A new
  smoothness indicator for the weno schemes and its effect on the convergence
  to steady state solutions}, \bibinfo{shortjournal}{J. Sci. Comput.}
  \bibinfo{volume}{31} (\bibinfo{year}{2007}) \bibinfo{pages}{273--305}.

\end{thebibliography}

\end{document}